\documentclass[a4paper,11pt]{book}

\usepackage[a4paper]{geometry}
\usepackage[latin1]{inputenc}
\usepackage[T1]{fontenc}

\usepackage[english,francais]{babel}
\usepackage{lmodern}
\rmfamily
\DeclareFontShape{T1}{lmr}{b}{sc}{<->ssub*cmr/bx/sc}{}
\DeclareFontShape{T1}{lmr}{bx}{sc}{<->ssub*cmr/bx/sc}{}

\usepackage{amssymb,amsmath,amsthm,amscd}
\usepackage{mathrsfs}

\usepackage{subfig}

\usepackage{epsfig}
\usepackage{graphicx} 
\usepackage{pstricks}

\usepackage{hyperref}

\usepackage{url}

\usepackage[nottoc]{tocbibind}

\usepackage{fancyhdr}
\usepackage{enumitem}

\makeatletter
\let\ps@plain=\ps@empty
\makeatother

\let\origdoublepage\cleardoublepage
\newcommand{\clearemptydoublepage}{%
  \clearpage
  {\pagestyle{empty}\origdoublepage}%
}
\let\cleardoublepage\clearemptydoublepage

\FrenchFootnotes 
\AddThinSpaceBeforeFootnotes 

\newcommand*\chapterstar[1]{%
  \chapter*{#1}%
  \addcontentsline{toc}{chapter}{#1}%
  \markboth{#1}{#1}}

\theoremstyle{plain}
\newtheorem{theo}{Théorème}[chapter]

\newtheorem{prop}[theo]{Proposition}
\newtheorem{lemm}[theo]{Lemme}

\newtheorem*{theo*}{Théorème} 
\newtheorem*{conj*}{Conjecture} 
\newtheorem*{prop*}{Proposition}

\theoremstyle{definition}
\newtheorem{defi}[theo]{Définition}
\newtheorem{rema}[theo]{Remarque}
\newtheorem*{exem}{Exemple}

\def\ie{\emph{i.e. }}

\def\R{\mathbb{R}}
\def\C{\mathbb{C}}
\def\N{\mathbb{N}}
\def\Z{\mathbb{Z}}
\def\P{\mathbb{P}^1}

\def\Rinfty{\R^3\cup\{\infty\}}

\def\l{\lambda}
\def\a{\alpha}
\def\b{\beta}
\def\t{\theta}
\def\e{\varepsilon}
\def\d{\delta}
\def\ga{\gamma}
\def\s{\sigma}

\def\X{\mathcal{X}}
\def\B{\mathcal{B}}
\def\D{\mathcal{D}}
\def\p{\mathcal{P}}
\def\E{\mathcal{E}}
\def\A{\mathcal{A}}

\def\w{\widetilde}

\def\g{\left(G,H \right) }

\def\gxt{\left(G(x,t),H(x,t) \right) }

\def\Y{\mathbf{Y}}
\def\wY{\mathbf{\widetilde Y}}
\def\S{SU(2)}
\def\ssm{\smallsetminus}
\def\eq{E_\t(\l,\mu,t)}

\DeclareMathOperator{\tr}{Tr}
\DeclareMathOperator{\res}{Res}
\DeclareMathOperator{\id}{id}
\DeclareMathOperator{\PI}{P_I}
\DeclareMathOperator{\PVI}{P_{VI}}
\DeclareMathOperator{\I}{I}

\title{Probl\`eme de Plateau, \'equations fuchsiennes et probl\`eme de Riemann--Hilbert}
\author{Laura \sc{Desideri}}
\date{10 mars 2011}

\begin{document}

\maketitle

\newpage
\thispagestyle{empty}

\noindent
Laura \textsc{Desideri}\vspace{0.1cm}\noindent\\
Universit\"at T\"ubingen\\ 
Mathematisches Institut\\
Auf der Morgenstelle 10\\
72 076 T\"ubingen, Germany\vspace{0.1cm}\noindent\\
E-mail : \texttt{desideri@mathematik.uni-tuebingen.de}\\
Url : \url{http://www.mathematik.uni-tuebingen.de/ab/Differentialgeometrie/desideri.html}

\pdfbookmark[0]{Résumé}{resume}
    
\chapter*{R\'esum\'e}

Ce m\'emoire est consacr\'e \`a la r\'esolution du probl\`eme de Plateau \`a bord polygonal dans l'espace euclidien de dimension trois. Il s'appuie sur la m\'ethode de r\'esolution propos\'ee par Ren\'e Garnier dans un article publi\'e en 1928 et qui a \'et\'e oubli\'e depuis, voire ignor\'e \`a l'\'epoque. L'approche de Garnier est tr\`es diff\'erente de la m\'ethode variationnelle, elle est plus g\'eom\'etrique et constructive, et permet d'obtenir des disques minimaux sans point de branchement. Cependant, elle est parfois tr\`es compliqu\'ee, voire obscure et incompl\`ete. En s'inspirant des id\'ees de Garnier, on propose une nouvelle d\'emonstration, qui est non seulement compl\`ete, mais \'egalement plus simple et plus moderne que la sienne. Ce travail repose principalement sur l'utilisation plus syst\'ematique des syst\`emes fuchsiens et la mise en \'evidence du lien entre la r\'ealit\'e d'un syst\`eme et sa monodromie.

La m\'ethode de Garnier repose sur le fait que, par la repr\'esentation de Weierstrass spinorielle des surfaces minimales, on peut associer une \'equation fuchsienne r\'eelle du second ordre, d\'efinie sur la sph\`ere de Riemann, \`a tout disque minimal \`a bord polygonal. La monodromie de cette \'equation est d\'etermin\'ee par les directions orient\'ees des c\^ot\'es du bord. Le bon point de vue consiste \`a consid\'erer des polygones pouvant avoir un sommet en l'infini. Pour r\'esoudre le probl\`eme de Plateau, on est donc amen\'e \`a r\'esoudre un probl\`eme de Riemann--Hilbert. On proc\`ede ensuite en deux \'etapes : tout d'abord, on d\'ecrit explicitement, par d\'eformations isomonodromiques, la famille de tous les disques minimaux dont le bord est un polygone de directions orient\'ees donn\'ees. Puis on utilise cette description pour \'etudier les longueurs des c\^ot\'es des bords polygonaux, et on montre ainsi que tout polygone est le bord d'un disque minimal.

\subsubsection*{Mots-clefs} Surfaces minimales, syst\`emes compl\`etement int\'egrables, \'equations fuchsiennes et syst\`emes fuchsiens, probl\`eme de Riemann--Hilbert, d\'eformations isomonodromiques, syst\`eme de Schlesinger.

\paragraph*{Classification math\'ematique par sujets (2010)} 53A10, 34M03, 34M35, 34M50, 34M55, 34M56.

\newpage

\begin{otherlanguage}{english}

  \begin{center} \Large \bf The Plateau problem, Fuchsian equations and the Riemann--Hilbert problem \end{center}

\bigskip

\section*{Abstract}

This dissertation is devoted to the resolution of the Plateau problem in the case of polygonal boundary curves in the three-dimensional Euclidean space. It relies on the method developed by Ren\'e Garnier and published in 1928 in a paper which seems today to be totally forgotten. Garnier's approach is more geometrical and constructive than the variational one, and it provides minimal disks without branch point. However, it is sometimes really complicated, and even obscure or incomplete. Following Garnier's initial ideas, we propose a new proof, which intends not only to be complete, but also simpler and moderner than his one. This work mainly relies on a systematic use of Fuchsian systems and on the relation that we establish between the reality of such systems and their monodromy.

Garnier's method is based on the following fact: using the spinor Weierstrass representation for minimal surfaces, we can associate a real Fuchsian second-order equation, defined on the Riemann sphere, with each minimal disk with a polygonal boundary curve. The monodromy of the equation is determined by the oriented directions of the edges of the boundary. To solve the Plateau problem, we are thus led to solve a Riemann--Hilbert problem. We then proceed in two steps: first, by means of isomonodromic deformations, we construct and describe the family of all minimal disks with a polygonal boundary curve of given oriented directions. Then we use this description to study the edges's lengths of their boundary curves, and we show that every polygon is the boundary of a minimal disk.

\subsubsection*{Keywords} 
Minimal surfaces, integrable systems, Fuchsian equations and Fuchsian systems, the Riemann--Hilbert problem, isomonodromic deformations, Schlesinger system.

\paragraph*{Mathematics Subject Classification (2010)} 53A10, 34M03, 34M35, 34M50, 34M55, 34M56.

\end{otherlanguage}

\setcounter{tocdepth}{1} 
\pdfbookmark[0]{Table des mati\`eres}{tablematieres}

\tableofcontents

\chapterstar{Introduction}

Ce m\'emoire a pour but de pr\'esenter une r\'esolution du probl\`eme de Plateau \`a bord polygonal, qui est tr\`es diff\'erente de la m\'ethode variationnelle, et qui repose sur une m\'ethode \'elabor\'ee par Ren\'e Garnier. Garnier a expos\'e cette m\'ethode dans l'article \emph{Le Probl\`eme de Plateau} \cite{Garnier28}. %
Publi\'e en 1928, c'est-\`a-dire environ deux ans avant les d\'emonstrations du probl\`eme de Plateau obtenues ind\'ependamment par T. Rad\'o~\cite{Rado} et J. Douglas~\cite{Douglas}, cet article semble avoir \'et\'e compl\`etement oubli\'e, voire ignor\'e \`a l'\'epoque. M\^eme si l'existence de cette r\'esolution est aujourd'hui connue de certains sp\'ecialistes, lorsque j'ai commenc\'e ma th\`ese (dont ce m\'emoire est un des r\'esultats), personne ne semblait \^etre en mesure de dire comment elle fonctionnait, ni m\^eme si elle \'etait correcte ou non. Sa d\'emonstration est en effet tr\`es compliqu\'ee, parfois elliptique et obscure, et certains passages en sont m\^eme peu convaincants. En s'inspirant des id\'ees de Garnier, on propose ici une nouvelle preuve de ce r\'esultat, qui soit non seulement compl\`ete et compr\'ehensible, mais aussi plus simple, et  qui apporte un point de vue nouveau sur la m\'ethode de Garnier. Ce travail repose principalement sur l'utilisation plus syst\'ematique des syst\`emes fuchsiens et la mise en \'evidence du lien entre la r\'ealit\'e d'un tel syst\`eme et sa monodromie. Cette clarification des fondements de la m\'ethode de Garnier m'a permis de l'\'etendre au cas o\`u l'espace ambiant est l'espace de Minkowski de dimension trois~\cite{Desideri-minkowski}.

\bigskip

Les surfaces minimales sont les surfaces dont la courbure moyenne est partout nulle. Elles constituent les points critiques de la fonctionnelle d'aire pour les variations fixant le bord. La th\'eorie des surfaces minimales a commenc\'e au \textsc{xviii}$^\text e$ si\`ecle, avec les d\'ebuts du calcul des variations, 
et conna\^it d'importantes avanc\'ees dans la seconde moiti\'e du \textsc{xix}$^\text e$ si\`ecle, avec notamment la repr\'esentation due \`a Weierstrass de toute immersion conforme minimale \`a partir de deux fonctions holomorphes. \`A la fin du \textsc{xix}$^\text e$ si\`ecle et au d\'ebut du \textsc{xx}$^\text e$ si\`ecle, les math\'ematiciens s'int\'eressent au <<~probl\`eme de Plateau~>>, du nom du physicien belge Joseph Plateau qui en 1873, a \'etabli exp\'erimentalement, par de tr\`es nombreuses exp\'eriences sur les films de savon, que toute courbe ferm\'ee de l'espace est le bord d'une surface minimale. L'\'enonc\'e math\'ematique du probl\`eme de Plateau est le suivant : \emph{\'etant donn\'e une courbe ferm\'ee connexe de Jordan de l'espace euclidien de dimension trois, montrer qu'il existe une surface minimale r\'eguli\`ere et ayant la topologie d'un disque dont le bord soit la courbe ferm\'ee}. Au d\'ebut des ann\'ees 1930, Tibor Rad\'o~\cite{Rado} et Jesse Douglas~\cite{Douglas} obtiennent ind\'ependamment par la m\'ethode variationnelle les premiers r\'esultats g\'en\'eraux (reconnus !) du probl\`eme de Plateau. Cependant, ils ne parviennent pas \`a exclure l'existence de points de branchement isol\'es \`a l'int\'erieur ou au bord du disque minimal. Il faut attendre les ann\'ees 1970, et les travaux de R. Osserman~\cite{Osserman}, R. Gulliver~\cite{Gulliver} et R. Osserman, R. Gulliver et H. L. Royden~\cite{GulliverOssermanRoyden} pour obtenir une d\'emonstration du probl\`eme de Plateau qui soit absolument compl\`ete.

La m\'ethode de Garnier pour r\'esoudre le probl\`eme de Plateau est tr\`es diff\'erente de la m\'ethode variationnelle. M\^eme si elle para\^it moins puissante, elle permet d'obtenir des surfaces, qui, contrairement aux solutions de Douglas--Rad\'o, sont r\'eguli\`eres partout. De plus, l'approche de Garnier est plus g\'eom\'etrique, s'inscrivant dans la continuation des travaux de K. Weierstrass, B. Riemann, H.-A. Schwarz et G. Darboux. Elle est \'egalement plus constructive que la m\'ethode variationnelle.

La m\'ethode de Garnier repose sur la correspondance de tout disque
minimal \`a bord polygonal avec une \'equation fuchsienne r\'eelle du
second ordre d\'efinie sur la sph\`ere de Riemann. Cette
correspondance est ant\'erieure aux travaux de Garnier. Elle est
donn\'ee par la repr\'esentation de Weierstrass, aujourd'hui dite
spinorielle, des immersions conformes minimales. Cette \'equation
fuchsienne semble \^etre mentionn\'ee pour la premi\`ere fois, de
mani\`ere ind\'ependante et presque simultan\'ee, dans un bref article
de Karl Weierstrass~\cite{Weierstrass1} publi\'e au mois de d\'ecembre
1866, et lors d'une pr\'esentation posthume des travaux de Bernhard
Riemann~\cite{Riemann} par Hattendorf le 6 janvier 1867 \`a la
Soci\'et\'e Royale de G\"ottingen. Riemann n'utilise pas la
repr\'esentation de Weierstrass, mais deux repr\'esentations conformes
(sph\'erique et plane) du m\^eme disque minimal. Gaston Darboux \'etudie en d\'etail cette
\'equation associ\'ee \`a un disque minimal \`a bord polygonal (\cite{Darboux}, chapitre \textsc{xiii}), et expose
les difficult\'es \`a surmonter pour \^etre en mesure de r\'esoudre le
probl\`eme de Plateau. Au premier rang de celles-ci figure la
d\'etermination d'une \'equation fuchsienne \`a partir de sa monodromie
: c'est le <<~probl\`eme de Riemann--Hilbert~>>, qui deviendra
bient\^ot le vingt-et-uni\`eme des vingt-trois probl\`emes propos\'es par
David Hilbert au Congr\`es International de Paris en 1900. C'est
seulement une vingtaine d'ann\'ees apr\`es ces observations de Darboux
que seront obtenues les premi\`eres solutions du probl\`eme de
Riemann--Hilbert, par J. Plemelj~\cite{Plemelj} et G. Birkhoff~\cite{Birkhoff} -- solutions dont A. A. Bolibruch a montr\'e des d\'ecennies plus tard par une s\'erie de contre-exemples~\cite{Bolibruch1}, \cite{Bolibruch2} qu'elles contiennent une erreur.

Garnier est un \'etudiant de Paul Painlev\'e. En 1912, il publie un article~\cite{Garnier12} qui rassemble les r\'esultats de sa th\`ese et dans lequel il \'etudie en particulier les d\'eformations isomonodromiques d'\'equations fuchsiennes ayant un nombre arbitraire de singularit\'es et aucune singularit\'e logarithmique. Le syst\`eme diff\'erentiel qui gouverne ces d\'eformations, connu aujourd'hui sous sa forme hamiltonienne sous le nom de \emph{syst\`eme de Garnier}, est en un sens une g\'en\'eralisation de la sixi\`eme \'equation de Painlev\'e $\PVI$. En 1926, il propose une r\'esolution du probl\`eme de Riemann--Hilbert~\cite{Garnier26} bas\'ee sur l'\'etude du syst\`eme Schlesinger au voisinage de ses singularit\'es non mobiles, et de ses liens avec le syst\`eme de Garnier. Les r\'esultats obtenus dans ces deux articles lui permettent d'esp\'erer \^etre en mesure de lever les difficult\'es mises en \'evidence par Darboux pour la r\'esolution du probl\`eme de Plateau. Il lui reste n\'eanmoins encore beaucoup de travail \`a accomplir pour obtenir cette r\'esolution~\cite{Garnier28}.

Depuis les ann\'ees 1970, leurs liens avec des probl\`emes issus de la
physique sont \`a l'origine de l'int\'er\^et nouveau que suscitent les
\'equations de Painlev\'e, et cons\'ecutivement, le syst\`eme de Garnier.
C'est \`a Kazuo Okamoto et \`a Hironobu Kimura que l'on doit la
<<~red\'ecouverte~>> du syst\`eme de Garnier au d\'ebut des ann\'ees 1980
et, en particulier, la mise en \'evidence de sa structure
hamiltonienne~\cite{Okamoto}. Dans ce contexte, et gr\^ace notamment
aux travaux de Mikio Sato, Tetsuji Miwa et Michio Jimbo~\cite{SMJ}
sur le probl\`eme de Riemann--Hilbert et le syst\`eme de Schlesinger,
la r\'esolution du probl\`eme de Plateau par Garnier rev\^et elle aussi
un int\'er\^et nouveau, avec entre autre la possibilit\'e d'une
simplification.

\subsection*{R\'esum\'e des chapitres}

L'objet de ce m\'emoire est la d\'emonstration du th\'eor\`eme suivant. 

\begin{theo}[Probl\`eme de Plateau \`a bord polygonal]
Tout polygone $P\subset\R^3$ en position g\'en\'erique, ayant \'eventuellement un sommet en l'infini, est le bord d'au moins un disque minimal immerg\'e. De plus, si $P$ a un sommet en l'infini, alors le disque minimal a un bout h\'elico\"idal en ce sommet.
\label{thm-Plateau}
\end{theo}

On dit ici qu'un polygone $P$ \`a $n+3$ c\^ot\'es est \emph{en position g\'en\'erique} si le $(n+3)$-uplet des directions orient\'ees de ses c\^ot\'es $D=(D_1,\ldots,D_{n+3})$ est dans l'ensemble $\D^n$ (d\'efinition~\ref{def-Dn}), \ie si deux directions quelconques de $P$ ne sont pas colin\'eaires et trois directions quelconques ne sont pas coplanaires.

Pour toute direction orient\'ee $D\in\D^n$, on introduit l'ensemble $\p^n_D$ des polygones \`a $n+3$ c\^ot\'es de direction $D$ ayant \'eventuellement un sommet en l'infini (\ie des lignes bris\'ees \'eventuellement infinies), d\'efinis \`a translation et homoth\'etie de rapport positif pr\`es (d\'efinition~\ref{def-PnD}) : ces polygones sont caract\'eris\'es par $n$ rapports de longueurs de c\^ot\'es, entre leurs $n+1$ longueurs finies, et l'ensemble $\p^n_D$ est ainsi isomorphe \`a $]0,+\infty[\,^n$. On d\'efinit \'egalement l'ensemble $\X^n_D$ des immersions conformes minimales $X$ qui repr\'esentent des disques minimaux ayant un bord polygonal $P\in\p^n_D$, et un bout h\'elico\"idal si $P$ a un sommet en l'infini, \'egalement \`a translation et homoth\'etie de rapport positif pr\`es (d\'efinition~\ref{def-XnD}). On peut toujours supposer qu'une telle immersion est d\'efinie sur le demi-plan sup\'erieur
\[
\C_+ = \left\{x\in \C \ | \ \Im (x)>0\right\}.
\]
On peut alors paraphraser ainsi le th\'eor\`eme~\ref{thm-Plateau} : il revient \`a montrer que pour toute direction $D\in\D^n$, l'application suivante est surjective
\begin{align*}
\X^n_D &\longrightarrow \p^n_D\\
X & \longmapsto   \partial X\left(\C_+\right).
\end{align*}
Pour cela, la m\'ethode que propose Garnier repose sur une correspondance bijective explicite entre une classe ad\'equate d'\'equations fuchsiennes, not\'ee $\E^n_D$, et l'ensemble $\X^n_D$. On cherchera donc plut\^ot \`a montrer que la composition suivante est surjective
\[
\E^n_D \stackrel{1:1}{\longrightarrow}\X^n_D \stackrel{\partial}{\longrightarrow} \p^n_D \stackrel{\sim}{\longrightarrow}  (0,+\infty)^n.
\]

Apr\`es deux premiers chapitres introductifs, on d\'efinit et on caract\'erise au chapitre~\ref{chapitre-equ-fu} l'ensemble d'\'equations $\E^n_D$, en constituant une sorte de dictionnaire entre les ensembles $\X^n_D$ et $\E^n_D$. Au chapitre~\ref{chapitre-isomono}, on consid\`ere l'ensemble analogue $\A^n_D$ de syst\`emes fuchsiens, et on d\'ecrit au moyen de d\'eformations isomonodromiques l'ensemble $\X^n_D$. Le chapitre~\ref{chapitre-longueur} est consacr\'e \`a la r\'esolution du probl\`eme de Plateau proprement dite : on utilise la description pr\'ec\'edente pour \'etudier les rapports de longueurs des bords polygonaux des immersions de $\X^n_D$, et on montre ainsi que tout polygone de directions orient\'ees $D$  est le bord d'au moins un disque minimal.

\paragraph*{Chapitre~\ref{chapitre-prelim1}. Surfaces minimales}

On expose des aspects g\'en\'eraux sur les surfaces minimales de l'espace euclidien de dimension trois. Le point essentiel est la repr\'esentation de Weierstrass que l'on appelle aujourd'hui spinorielle : tout couple de fonctions $\g$ holomorphes sur une le demi-plan sup\'erieur et sans z\'ero commun d\'efinit une immersion conforme minimale de $\C_+$ dans $\R^3$, et r\'eciproquement, toute immersion de ce type est obtenue par un couple de fonctions holomorphes sans z\'ero commun.

\paragraph*{Chapitre~\ref{chapitre-prelim2}.  \'Equations fuchsiennes et syst\`emes fuchsiens}

On donne une introduction assez d\'etaill\'ee des notions de base telles que le comportement local au voisinage des singularit\'es, la monodromie, le probl\`eme de Riemann--Hilbert, les d\'eformations isomonodromiques et, en particulier, le syst\`eme de Schlesinger. On explicite aussi les liens entre \'equations et syst\`emes fuchsiens.

\paragraph*{Chapitre~\ref{chapitre-equ-fu}. L'\'equation associ\'ee \`a un disque minimal \`a bord polygonal}

Ce chapitre n'est pas consacr\'e \`a la r\'esolution du probl\`eme de Plateau proprement dite, mais plut\^ot \`a l'\'etude de la correspondance entre disques minimaux \`a bord polygonal et \'equations fuchsiennes. Cette correspondance est ant\'erieure aux travaux de Garnier sur le probl\`eme de Plateau, elle est d\'ej\`a \'etudi\'ee par Darboux (\cite{Darboux}, chapitre \textsc{xiii}).

On consid\`ere une immersion conforme minimale $X : \C_+\to\R^3$ qui repr\'esente un disque minimal \`a bord polygonal de direction $D$, c'est-\`a-dire un \'el\'ement de $\X^n_D$. Cette immersion est caract\'eris\'ee par ses donn\'ees de Weierstrass $G$ et $H$, qui sont des fonctions holomorphes dans $\C_+$, et qui sont lin\'eairement ind\'ependantes d\`es que l'image de $X$ n'est pas plane. Elles sont donc solutions d'une unique \'equation diff\'erentielle ordinaire lin\'eaire du second ordre
\begin{equation}
y''+p(x)y'+q(x)y=0.
\tag{$E$}\label{E-intro}
\end{equation}
L'\'equation~\eqref{E-intro} est l'\'equation associ\'ee \`a l'immersion $X$. On note $\E^n_D$ l'ensemble des \'equations qui sont associ\'ees en ce sens \`a une immersion appartenant \`a $\X^n_D$. Le but de ce chapitre est d'obtenir une caract\'erisation de l'ensemble $\E^n_D$, en traduisant des propri\'et\'es g\'eom\'etriques des immersions $X$ en terme de propri\'et\'es analytiques des \'equations~\eqref{E-intro}. Une \'equation~\eqref{E-intro} de $\E^n_D$ a deux types de singularit\'es : les ant\'ec\'edents par l'immersion $X$ des sommets du bord polygonal $P$, qui sont r\'eels
\[
 t_1<\cdots<t_n< t_{n+1}=0, \ t_{n+2}=1,\ t_{n+3}=\infty,
\]
et les ombilics de $X$, qui sont des singularit\'es apparentes. En appliquant le principe de r\'eflexion de Schwarz, on montre que l'\'equation~\eqref{E-intro} s'\'etend \`a la sph\`ere de Riemann, sur laquelle c'est une \'equation fuchsienne r\'eelle, et on d\'etermine comment les donn\'ees de Weierstrass sont transform\'ees autour des singularit\'es $t_i$. On montre ainsi que la monodromie de l'\'equation~\eqref{E-intro} est enti\`erement d\'etermin\'ee par la direction orient\'ee $D$ du bord polygonal de $X$ : l'ensemble $\E_D^n$ est isomonodromique. Il n'y a par contre aucune traduction naturelle des longueurs des c\^ot\'es de $P$ en terme de propri\'et\'es de l'\'equation~\eqref{E-intro}.

On obtient ainsi que les \'equations de $\E^n_D$ sont caract\'eris\'ees par trois conditions : une condition~\ref{cond-SdR} qui est d'ordre local (nature et position des singularit\'es, valeurs des exposants), une condition~\ref{cond-mono} qui impose la monodromie \`a partir de la direction $D$, et une condition de r\'ealit\'e~\ref{cond-realite}. Finalement, on montre que l'ensemble $\E_D^n$ est en bijection avec l'ensemble $\X^n_D$.

\paragraph*{Chapitre~\ref{chapitre-isomono}. D\'eformations isomonodromiques}

\'Etant donn\'e un $(n+3)$-uplet de directions orient\'ees $D\in\D^n$, le but de ce chapitre est d'utiliser l'ensemble $\E^n_D$ pour d\'ecrire explicitement l'ensemble $\X_D^n$. Contrairement \`a Garnier, pour obtenir cette description, on va plut\^ot utiliser des syst\`emes fuchsiens, \`a la place des \'equations fuchsiennes de $\E_D^n$. Cette approche apporte un point de vue nouveau \`a la m\'ethode de Garnier et la simplifie notablement.

On commence donc par introduire l'ensemble analogue  $\A_D^n$ des syst\`emes fuchsiens du premier ordre de taille $2\times2$, qui sont associ\'es, dans un sens que l'on pr\'ecisera, aux immersions de l'ensemble $\X^n_D$. On \'etablit une caract\'erisation de ces syst\`emes par des conditions~\ref{cond-sys-SdR}, \ref{cond-sys-mono}, et \ref{cond-sys-realite}, qui sont les analogues des conditions~\ref{cond-SdR}, \ref{cond-mono}, et \ref{cond-realite} pr\'ec\'edentes. Notamment, les conditions~\ref{cond-mono} et \ref{cond-sys-mono}, qui portent sur la monodromie, sont identiques. L'ensemble $\A^n_D$ n'est pas en bijection avec l'ensemble $\X^n_D$, puisque des syst\`emes fuchsiens diff\'erents peuvent d\'efinir la m\^eme \'equation.

Pour d\'ecrire l'ensemble $\A^n_D$, on l\`eve ensuite une difficult\'e ignor\'ee par Garnier, qu'est la condition de r\'eali\'e~\ref{cond-sys-realite}. On montre que la <<~r\'ealit\'e~>> d'un syst\`eme fuchsien peut \^etre caract\'eris\'ee par sa monodromie : on \'etablit une condition n\'ecessaire et suffisante portant sur la monodromie d'un syst\`eme pour que celui-ci satisfasse la condition~\ref{cond-sys-realite}. En particulier, cette condition est v\'erifi\'ee par une monodromie satisfaisant la condition~\ref{cond-sys-mono} : l'ensemble $\A^n_D$ est donc simplement l'ensemble des syst\`emes v\'erifiant les conditions~\ref{cond-sys-SdR} et~\ref{cond-sys-mono}.

Enfin, on utilise des d\'eformations isomonodromiques pour d\'ecrire les syst\`emes de $\A^n_D$. On obtient que l'ensemble $\A^n_D$ contient une famille isomonodromique de syst\`emes fuchsiens $\left(A_D(t), t\in\pi^n\right)$ param\'etr\'ee par la position des singularit\'es $t=(t_1,\ldots,t_n)$ variant dans le simplexe
\[
 \pi^n=\left\{ (t_1,\ldots,t_n)\in\R^n\ \big|\ t_1<\cdots<t_n<0\right\} ,
\]
d\'ecrite par le syst\`eme de Schlesinger et qui est en bijection avec l'ensemble $\X^n_D$. On obtient \'egalement un r\'esultat de r\'egularit\'e en $t$ pour cette famille. On en d\'eduit une description explicite de l'ensemble $\X^n_D=\left(X_D(t), t\in\pi^n\right)$, et de la famille $\left(P_D(t), t\in\pi^n\right)\subset\p^n_D$ des polygones de direction $D$ qui sont le bord d'au moins un disque minimal.

\paragraph*{Chapitre~\ref{chapitre-longueur}. Rapports de longueurs des c\^ot\'es}

Le but de chapitre est de montrer que la famille de polygones $\left(P_D(t), t\in\pi^n\right)$ d\'ecrit enti\`erement l'ensemble $\p^n_D$. Un syst\`eme de coordonn\'ees sur $\p_D^n$ est donn\'e par $n$ rapports de longueurs de c\^ot\'es. Pour chaque valeur de $t\in\pi^n$, les donn\'ees de Weierstrass $\gxt$ de l'immersion $X_D(t)$ sont obtenues \`a partir d'une solution fondamentale du syst\`eme fuchsien $\left(A_D(t)\right)$. Les rapports de longueurs des c\^ot\'es du polygone $P_D(t)$ s'\'ecrivent donc
\[
 r_i(t)  = \frac{\displaystyle\int_{t_i}^{t_{i+1}} \left(|G(x,t)|^2 + |H(x,t)|^2\right) dx}{\displaystyle\int_0^1 \left(|G(x,t)|^2 + |H(x,t)|^2 \right) dx}
\]
($i=1,\ldots,n$), et on obtient ainsi la fonction <<~rapports des longueurs~>> $F_D(t)$ associ\'ee \`a la direction $D$
\[
 F_D:\pi^n \to \,]0,+\infty[\,^n, \qquad F_D(t)=(r_1(t),\ldots,r_n(t)).
\]
Le but de ce chapitre est donc d'\'etablir le th\'eor\`eme suivant, qui conclut la d\'emonstration du th\'eor\`eme~\ref{thm-Plateau}, et qui en est la partie la plus difficile.

\begin{theo*}
 \'Etant donn\'e un $(n+3)$-uplet de directions orient\'ees $D\in\D^n$, la fonction <<~rapports des longueurs~>> $F_D:\pi^n \to \,]0,+\infty[\,^n$ est surjective.
\end{theo*}

On propose une d\'emonstration de ce th\'eor\`eme tr\`es diff\'erente de celle Garnier, bas\'ee sur l'\'etude de la famille $\left(A_D(t), t\in\pi^n\right)$ au bord du simplexe $\pi^n$ et une r\'ecurrence portant sur le nombre $n+3$ de c\^ot\'es des polygones. Par identification naturelle des simplexes $\pi^n$ et $]0,+\infty[\,^n$, on obtient une fonction $\w F_D:\,]0,+\infty[\,^n \to\,]0,+\infty[\,^n$. Pour montrer que la fonction $F_D$ est surjective, on montre que la fonction $\w F_D$ est de degr\'e $1$, c'est-\`a-dire homotope \`a l'identit\'e. On \'etablit un r\'esultat de topologie qui nous permet de nous ramener \`a montrer que la fonction $\w F_D$ est continue et de degr\'e $1$ au bord de $]0,+\infty[\,^n$. Pour obtenir cela, il faut interpr\'eter la fonction $F_D\big|_{\partial \pi^n}$ en terme de nouvelles fonctions <<~rapports des longueurs~>> de dimension inf\'erieure : c'est l'objet de la proposition~\ref{prop-F-cont-au-bord} dont l'\'enonc\'e para\^it naturel et qui est l'\'etape la plus importante de la d\'emonstration : la fonction $F_D(t)$ s'\'etend contin\^ument \`a chacune des faces du bord du simplexe $\pi^n$ (qui sont des simplexes de dimension inf\'erieure). Chaque face est caract\'eris\'ee par la <<~disparition~>> de certains $t_i$, qui ont fusionn\'e avec la singularit\'e suivante $t_{i+1}$. On affirme qu'alors la fonction $F_D(t)$ restreinte \`a chaque face est, \`a hom\'eomorphisme pr\`es, la fonction <<~rapports des longueurs~>> $F_{D'}:\pi^k \to \,]0,+\infty[\,^k$ ($1\leq k\leq n-1$) d\'efinie par les directions orient\'ees $D'\in\D^k$ obtenues \`a partir de $D$ en <<~enlevant~>> les composantes $D_i$ correspondant aux $t_i$ qui ont disparu. Une fois que l'on a obtenu la proposition~\ref{prop-F-cont-au-bord}, il suffit pour conclure de faire une r\'ecurrence sur le nombre $n+3$ de c\^ot\'es, dont l'h\'er\'edit\'e est assur\'ee par le r\'esultat de topologie mentionn\'e plus haut, et dont l'initialisation au rang $n=1$ (cas d'un bord quadrilat\'eral) est imm\'ediate une fois que l'on a obtenu la proposition~\ref{prop-F-cont-au-bord}.

La majeure partie de ce chapitre est donc consacr\'ee \`a la d\'emonstration de la proposition~\ref{prop-F-cont-au-bord}. La partie la plus difficile est d'obtenir la continuit\'e de la fonction $F_D(t)$ au bord, et non pas son interpr\'etation g\'eom\'etrique. On s'appuie sur des r\'esultats g\'en\'eraux sur les singularit\'es fixes du syst\`eme de Schlesinger, que Garnier appelle les pseudo-chocs, c'est-\`a-dire en les points tels que $t_i=t_j$, $i\neq j$. Ces r\'esultats sont une partie plus connue du travail de Garnier~\cite{Garnier26}, et ont \'et\'e d\'evelopp\'es et g\'en\'eralis\'es par Sato, Miwa et Jimbo~\cite{SMJ}. On reprend ces r\'esultats en en approfondissant des aspects qui nous seront utiles pour \'etudier l'holomorphie de la fonction $F_D(t)$ en les pseudo-chocs. On applique ensuite cette \'etude g\'en\'erale aux solutions particuli\`eres du syst\`eme de Schlesinger qui nous int\'eresse, c'est-\`a-dire au cas r\'eel, et on \'etablit la proposition~\ref{prop-F-cont-au-bord}.

\vspace{.5cm}

\noindent \textsc{Remerciements.} Je souhaite remercier mon directeur de th\`ese Fr\'ed\'eric H\'elein de m'avoir sugg\'er\'e de travailler sur la r\'esolution du probl\`eme de Plateau par R. Garnier, et pour son aide tout au long de ce travail. 

\chapter{Surfaces minimales}
\label{chapitre-prelim1}

On expose dans ce chapitre des aspects g\'en\'eraux sur les surfaces minimales de l'espace euclidien de dimension trois $\left(\R^3,\left\langle \, ,\right\rangle \right) $. On note $\left(O,e_1,e_2,e_3 \right) $ un rep\`ere orthonormal de $\R^3$. Une immersion conforme $X:\Sigma\to\R^3$ d'une surface de Riemann $\Sigma$ dans $\R^3$ est dite minimale si sa courbure moyenne est partout nulle. Rappelons que la courbure moyenne d'une immersion est la moiti\'e de la trace de sa deuxi\`eme forme fondamentale.

\section{Repr\'esentation de Weierstrass}

La repr\'esentation de Weierstrass est un outil fondamental dans l'\'etude des surfaces minimales. Elle permet \`a la fois de caract\'eriser et de construire des surfaces minimales. Donnons tout d'abord une premi\`ere forme, classique, de cette repr\'esentation.

\begin{theo}
Soient $\Sigma$ une surface de Riemann et $x_0$ un point de $\Sigma$.

Soient une fonction $g$ m\'eromorphe dans $\Sigma$ et une $1$-forme diff\'erentielle $\omega$ holomorphe dans $\Sigma$ telles que
\begin{itemize}
 \item  les z\'eros de $\omega$ sont d'ordre pair,
 \item  $g$ a un p\^ole d'ordre $m$ en un point $a\in\Sigma$ si et seulement si $\omega$ a un z\'ero d'ordre $2m$ en $a$.
 \end{itemize}
Alors l'application $X$ d\'efinie sur le rev\^etement universel $\w\Sigma$ de $\Sigma$ par
\[
 X(x) = \Re \int_{x_0}^x \left(1-g^2,i(1+g^2),2g \right) \omega
\]
est une immersion conforme minimale de $\widetilde\Sigma$ dans $\R^3$.

R\'eciproquement, si $X:\Sigma\to\R^3$ est une immersion conforme minimale, alors il existe un point $X_0\in\R^3$, une fonction $g$ m\'eromorphe dans $\Sigma$ et une $1$-forme diff\'erentielle $\omega$ holomorphe dans $\Sigma$ v\'erifiant les deux conditions ci-dessus tels que
\[
 X(x) = X_0 + \Re \int_{x_0}^x \left(1-g^2,i(1+g^2),2g \right) \omega.
\]
\end{theo}

La diff\'erentielle de Hopf de l'immersion $X$ est, par d\'efinition, la $2$-forme diff\'erentielle
\[
    Q=\left\langle \frac{d^2X}{dx^2},N\right\rangle dx^2,
\]
et elle s'exprime en fonction des donn\'ees $(g,\omega)$ par $Q=-\omega dg$. On peut voir facilement que la fonction $g$ est le projet\'e st\'er\'eographique par rapport au p\^ole nord du vecteur de Gauss $N:\Sigma\to\mathbb S^2$ de l'immersion $X$. Les donn\'ees g\'eom\'etriques de l'immersion $X$ sont caract\'eris\'ees par les donn\'ees $(g,\omega)$~: sa m\'etrique induite et sa seconde forme fondamentale sont
\[
 ds^2=\left(1+|g|^2 \right)^2|\omega|^2 , \qquad \text{II}=Q+\bar Q.
\]

Cependant, la repr\'esentation qu'utilise Garnier, et que l'on va utiliser exclusivement dans ce m\'emoire, est la repr\'esentation aujourd'hui dite spinorielle des surfaces minimales. Bien que soit probablement sous cette forme que la repr\'esentation de Weierstrass ait \'et\'e donn\'ee pour la premi\`ere fois --- par K. Weierstrass lui-m\^eme~\cite{Weierstrass1} ---, elle n'est pas consid\'er\'ee aujourd'hui comme la repr\'esentation \emph{classique}. Par souci de simplicit\'e, comme on ne s'int\'eresse dans ce m\'emoire qu'aux disques minimaux, on n'\'enonce cette repr\'esentation que dans le cas des immersions $X$ d\'efinie dans le demi-plan sup\'erieur ou demi-plan de Poincar\'e
\begin{equation}
 \C_+ = \left\{x\in \C \ | \ \Im (x)>0\right\},
\label{def-C+}
\end{equation}
o\`u $\Im (x)$ d\'esigne la partie imaginaire du nombre complexe $x$. Il n'y a pas alors de probl\`eme de p\'eriode, et de passage au rev\^etement universel. On pourra se reporter \`a~\cite{KuSch} pour un \'enonc\'e plus g\'en\'eral et pour plus de d\'etails sur la repr\'esentation spinorielle.


\begin{theo}
Soit $x_0$ un point du demi-plan sup\'erieur $\C_+$.

Pour tout couple $\g :\C_+\to\C^2\ssm\left\{(0,0)\right\}$ de fonctions holomorphes dans $\C_+$ sans z\'ero commun, l'application $X:\C_+\to\R^3$ d\'efinie par
\begin{equation}
 X(x) = \Re \int_{x_0}^x
\begin{pmatrix}
 i\left(G(\xi)^2 - H(\xi)^2\right) \\
 G(\xi)^2 + H(\xi)^2\\
 2iG(\xi)H(\xi)
\end{pmatrix}
d\xi
\label{rep-W}
\end{equation}
est une immersion conforme minimale.

R\'eciproquement, si $X:\C_+\to\R^3$ est une immersion conforme minimale, alors il existe un point $X_0\in\R^3$, et un couple $\g :\Sigma\to\C^2\ssm\left\{(0,0)\right\}$ de fonctions holomorphes tels que
\[
 X(x) = X_0 + \Re \int_{x_0}^x
\begin{pmatrix}
 i\left(G(\xi)^2 - H(\xi)^2\right) \\
 G(\xi)^2 + H(\xi)^2\\
 2iG(\xi)H(\xi)
\end{pmatrix}
d\xi.
\]
\end{theo}

Comme on utilisera exclusivement cette repr\'esentation, on l'appellera, contrairement \`a l'usage actuel, la repr\'esentation de Weierstrass, et le couple de fonctions holomorphes $\left(G,H \right)$ les donn\'ees de Weierstrass de l'immersion $X$. La correspondance entre les deux repr\'esentations pr\'ec\'edentes est donn\'ee par
\[
 g=-\frac{G}{H}, \qquad \omega=-iH^2dx.
\]
Le projet\'e st\'er\'eographique nord du vecteur de Gauss $N$ est $-G/H$. La diff\'erentielle de Hopf est donn\'ee par le Wronskien des fonctions $G$ et $H$
\begin{equation}
 Q=i\left( GH'-HG'\right) dx^2,
\label{Q}
\end{equation}
et la m\'etrique induite et la seconde forme fondamentale par
\begin{equation}
 ds^2=\left(|G|^2+|H|^2 \right)^2|dx|^2 , \qquad \text{II}=Q+\bar Q.
\label{metrique-ff2}
\end{equation}

\begin{figure}
 \begin{minipage}[b]{.4\linewidth}
  \centering\epsfig{figure=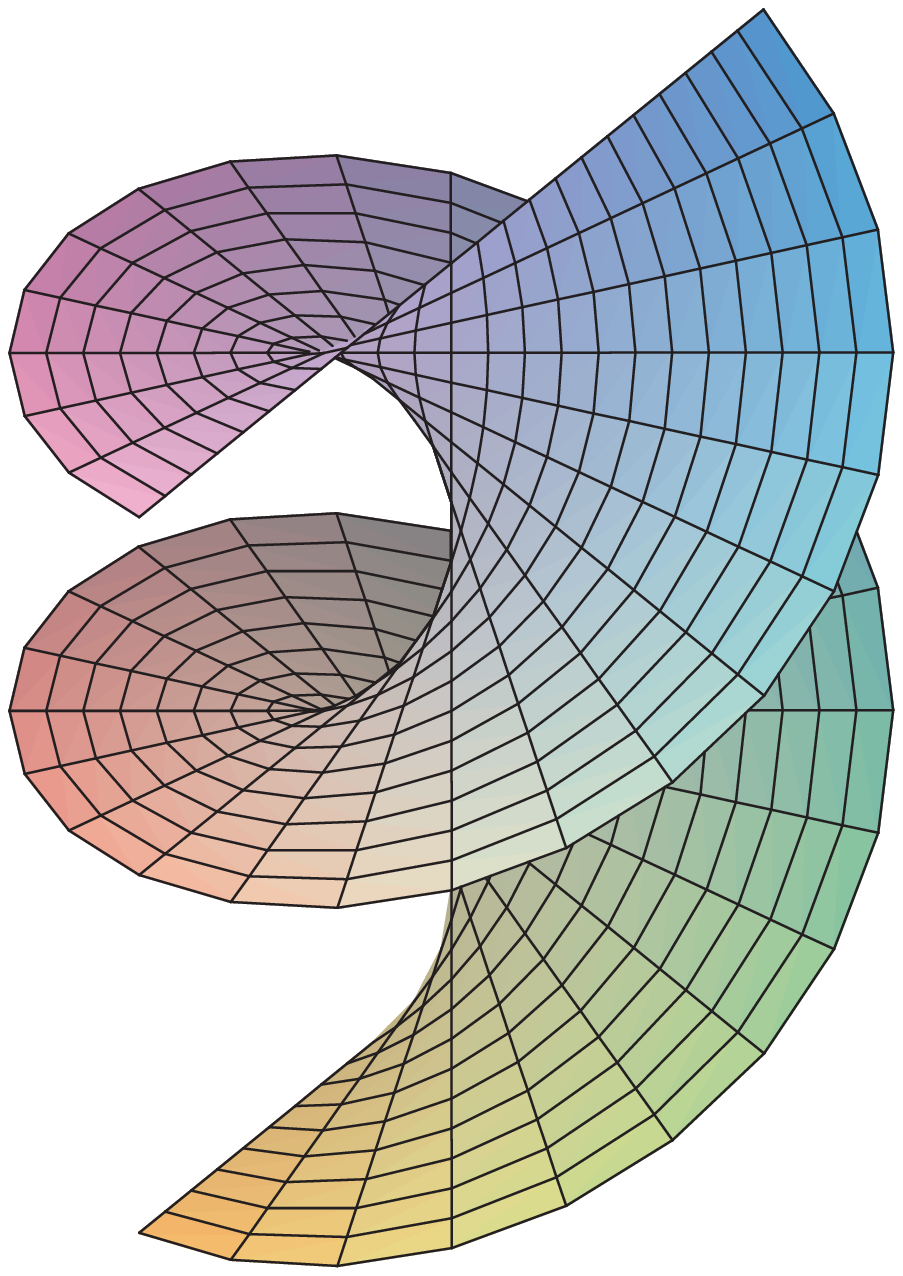,width=\linewidth}
  \caption{Une h\'elico\"ide}
  \label{fig-helicoide}
 \end{minipage} \hfill
 \begin{minipage}[b]{.4\linewidth}
  \centering\epsfig{figure=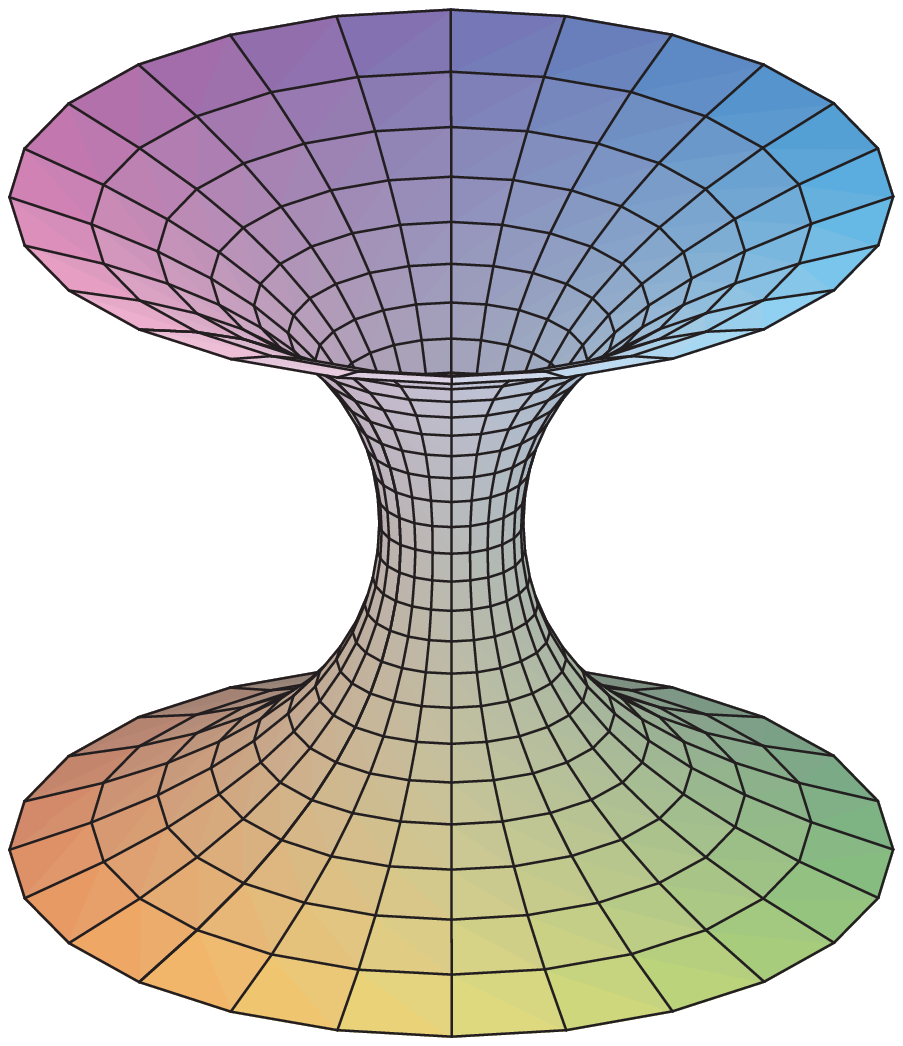,width=\linewidth}
  \caption{Une cat\'eno\"ide}
  \label{fig-catenoide}
 \end{minipage}
\end{figure}

\begin{exem} Voici les exemples les plus classiques de surfaces minimales.
\begin{enumerate}
 \item Si les fonctions $G$ et $H$ sont proportionnelles, alors l'immersion associ\'ee d\'efinit une surface minimale contenue dans un plan (c'est m\^eme une \'equivalence). Si $\Sigma=\C$ et si les fonctions $G$ et $H$ sont constantes, on obtient un plan entier.
 \item Si on choisit $\Sigma=\C^*$, $G(x)=1$, $H(x)=1/x$, on obtient une \emph{h\'elico\"ide}. L'immersion $X$ est d\'efinie dans le rev\^etement universel de $\C^*$. Les h\'elico\"ides sont des surfaces r\'egl\'ees (figure~\ref{fig-helicoide}).
 \item Si on choisit $\Sigma=\C^*$, $G(x)=e^{i\frac{\pi}{4}}$, $H(x)=e^{i\frac{\pi}{4}}/x$, on obtient une \emph{cat\'eno\"ide}. On peut montrer qu'alors l'immersion $X$ est bien d\'efinie dans $\C^*$. Les cat\'eno\"ides sont les seules surfaces minimales de r\'evolution (figure~\ref{fig-catenoide}).
\end{enumerate}
\end{exem}

Une application diff\'erentiable $X:\C_+\to\R^3$ donn\'ee par~\eqref{rep-W} o\`u les fonctions $G$ et $H$ sont seulement suppos\'ees holomorphes, repr\'esente une surface minimale \emph{g\'en\'eralis\'ee}, c'est-\`a-dire 
qui peut avoir des points de branchement. Ces points de branchement sont les points o\`u la d\'eriv\'ee $\partial X/\partial x$ s'annule, et o\`u donc la surface minimale n'est plus immerg\'ee. Ce sont exactement les z\'eros communs des fonctions $G$ et $H$.

On voit que l'immersion $X$ ne change pas si on change le signe du couple $\left( G,H\right) $. En fait, les donn\'ees de Weierstrass $\left(G,H \right)$ associ\'ees \`a une immersion conforme minimale $X$ sont uniques au signe pr\`es. Par ailleurs, si on consid\`ere deux repr\'esentations conformes sur $\C_+$ du m\^eme disque minimal, elles se d\'eduisent l'une de l'autre par composition \`a droite par une repr\'esentation conforme du demi-plan $\C_+$ dans lui-m\^eme, \ie par une application de M\"obius
\[
 x\mapsto \frac{ax+b}{cx+d} \qquad \quad \text{ o\`u }
 \begin{pmatrix}
  a & b\\
  c & d
 \end{pmatrix}
 \in PSL(2,\R).
\]
Il suffit donc de fixer l'image de trois points par une immersion $X:\C_+\to\R^3$ pour la d\'eterminer enti\`erement \`a partir de son image.

Remarquons que si la repr\'esentation de Weierstrass donne une description locale tr\`es simple des immersions conformes minimales, elle paraît \emph{a priori} peu utile \`a la r\'esolution du probl\`eme de Plateau. Il semble en effet difficile de d\'eduire d'une courbe que l'on s'est fix\'ee \`a l'avance des conditions sur les donn\'ees de Weierstrass $\g$ qui assurent que l'immersion conforme minimale associ\'ee passe par cette courbe. On verra au chapitre~\ref{chapitre-equ-fu} comment l'\'equation associ\'ee \`a un disque minimal \`a bord polygonal permet de d\'eduire de cette description locale des contraintes globales sur les donn\'ees de Weierstrass.

\section{Surface minimale conjugu\'ee et famille associ\'ee}

Les coordonn\'ees d'une immersion conforme minimale sont les parties r\'eelles de fonctions holomorphes : elles sont donc harmoniques. Rappelons qu'\`a toute application harmonique $f$ d\'efinie sur une surface de Riemann $\Sigma$, on peut associer une autre application harmonique $f^*$, qui est \emph{a priori} d\'efinie dans le rev\^etement universel $\w\Sigma$ de $\Sigma$, telle que la fonction $f+if^*$ soit holomorphe dans $\w\Sigma$ ($f^*$ est d\'efinie \`a une constante additive pr\`es). L'application $f^*$ est appel\'ee \emph{l'application harmonique conjugu\'ee} de $f$. On peut ainsi introduire la d\'efinition suivante.

\begin{defi}
 Soit $X:\Sigma\to\R^3$ une immersion conforme minimale. Alors l'immersion conforme minimale $X^*:\widetilde\Sigma\to\R^3$ dont les coordonn\'ees sont les applications harmoniques conjugu\'ees de celles de $X$ est appel\'ee \emph{l'immersion conjugu\'ee} de $X$. Elle est d\'efinie \`a une translation pr\`es.
\end{defi}

Si l'immersion $X:\C_+\to\R^3$ a pour donn\'ees de Weierstrass $(G,H)$, alors l'immersion conjugu\'ee $X^*$ s'\'ecrit
\[
 X^*(x) =  \Im \int_{x_0}^x
\begin{pmatrix}
 i\left(G(\xi)^2 - H(\xi)^2\right) \\
 G(\xi)^2 + H(\xi)^2\\
 2iG(\xi)H(\xi)
\end{pmatrix}
d\xi,
\]
et ses donn\'ees de Weierstrass sont
\[
 e^{i\frac{\pi}{4}} G, \quad e^{i\frac{\pi}{4}} H.
\]
Les immersions $X$ et $X^*$ ont la m\^eme application de Gauss, et elles sont localement isom\'etriques. Par exemple, la surface conjugu\'ee d'une cat\'eno\"ide est une h\'elico\"ide, bien qu'elles ne soient pas globalement isom\'etriques. L'\'equation diff\'erentielle des lignes de courbure de $X$ est donn\'ee par
\[
 \Re\left(GH'-HG' \right) dx^2=0,
\]
et celle des lignes asymptotiques par
\[
 \Im\left(GH'-HG' \right) dx^2=0.
\]
Les lignes de courbure et les lignes asymptotiques sont donc \'echang\'ees entre une surface minimale et sa conjugu\'ee. Comme une surface minimale et sa conjugu\'ee ont les m\^emes g\'eod\'esiques et la m\^eme application de Gauss, on en d\'eduit donc le lemme suivant.

\begin{lemm}
Si une surface minimale de $\R^3$ contient un segment de droite de vecteur directeur $v$, alors ce segment correspond sur la surface minimale conjugu\'ee \`a une courbe plane contenue dans un plan normal \`a $v$ et que la surface coupe perpendiculairement.
\label{lemme-SM-conj}
\end{lemm}

En effet, si $\left( \mathcal S\right)$ est une surface immerg\'ee dans $\R^3$, alors les droites contenues dans $\left( \mathcal S\right)$ sont exactement les courbes qui sont \`a la fois des lignes asymptotiques et des g\'eod\'esiques de $\left( \mathcal S\right)$. De m\^eme, les courbes trac\'ees sur $\left( \mathcal S\right)$ et contenues dans un plan que la surface $\left( \mathcal S\right)$ coupe perpendiculairement sont exactement les courbes qui sont \`a la fois des lignes de courbure et des g\'eod\'esiques de $\left( \mathcal S\right)$.

Par exemple, les m\'eridiens d'une cat\'eno\"ide correspondent sur une h\'elico\"ide conjugu\'ee aux droites qui engendrent l'h\'elico\"ide. Le cercle m\'edian de la cat\'eno\"ide correspond \`a la droite centrale de l'h\'elico\"ide.

Plus g\'en\'eralement, pour tout $\l\in\C^*$, on peut d\'efinir l'immersion conforme minimale $X_\l:\C_+\to\R^3$ de donn\'ees de Weierstrass $\l(G,H)$. On a
\[
 X_\l(x) = \Re(\l^2) X(x) + \Im(\l^2) X^*(x).
\]
Si le scalaire $\l$ est r\'eel ou purement imaginaire, alors les immersions $X_\l$ sont homoth\'etiques \`a l'immersion $X$. Lorsque le scalaire $\l$ appartient au cercle unit\'e $\mathbb S^1$, les immersions  $X_\l$ sont localement isom\'etriques \`a l'immersion $X$. La famille d'immersions conformes minimales $\left(X_\l
\right)_{\l\in\mathbb S^1}$ est appel\'ee \emph{famille associ\'ee} \`a l'immersion $X$.

\section{Principes de r\'eflexion de Schwarz}

Les deux propositions suivantes mettent en \'evidence certaines sym\'etries apparaissant sur les surfaces minimales. Elles permettent \'egalement d'\'etendre les surfaces minimales ayant un bord au del\`a de celui-ci, lorsque ce bord contient un segment de droite ou une courbe contenue dans un plan que la surface coupe perpendiculairement. Ces r\'esultats nous seront tr\`es utiles par la suite. On note $\mathbb D$ le disque unit\'e ouvert de $\C$, $\mathbb D^+ = \left\{x\in\mathbb D\ |\ \Im(x)>0 \right\} $ et $\mathbb D^- = \left\{x\in\mathbb D\ |\ \Im(x)<0 \right\} $.

\begin{prop}
 Soit une immersion conforme minimale $X:\mathbb D^+\to\R^3$. Si $X$ s'\'etend contin\^ument \`a l'intervalle $]-1,1[=\mathbb D\cap\R$, et si l'image par $X$ de l'intervalle $]-1,1[$ est un segment de droite, alors l'immersion $X$ se prolonge \`a $\mathbb D^-$ par r\'eflexion par rapport \`a cette droite et $X:\mathbb D\to\R^3$ est une immersion conforme minimale. De plus, deux points sym\'etriques sur l'image $X(\mathbb D)$ ont des ant\'ec\'edents conjugu\'es.
\end{prop}

\begin{prop}
 Soit une immersion conforme minimale $X:\mathbb D^+\to\R^3$. Si $X$ s'\'etend contin\^ument \`a l'intervalle $]-1,1[=\mathbb D\cap\R$, et si l'image par $X$ de l'intervalle $]-1,1[$ est une courbe contenue dans un plan que la surface $X(\mathbb D^+)$ coupe perpendiculairement, alors l'immersion $X$ se prolonge \`a $\mathbb D^-$ par r\'eflexion par rapport \`a ce plan et $X:\mathbb D\to\R^3$ est une immersion conforme minimale. De plus, deux points sym\'etriques sur l'image $X(\mathbb D)$ ont des ant\'ec\'edents conjugu\'es.
\end{prop}

On donnera une d\'emonstration de ces propositions au chapitre~\ref{chapitre-equ-fu}.

Par le lemme~\ref{lemme-SM-conj}, une r\'eflexion axiale sur une surface minimale correspond sur la surface minimale conjugu\'ee \`a une r\'eflexion par rapport \`a un plan orthogonal \`a cet axe, et r\'eciproquement.

\section{Description quaternionique}

Consid\'erons l'isomorphisme de $\R^3$ dans l'ensemble $E^3$ des matrices de $M(2,\C)$ hermitiennes \`a trace nulle, qui identifie un vecteur $X=(X_1,X_2,X_3)^t\in\R^3$ avec la matrice $\w X$ d\'efinie par
\[
 \w X =
 \begin{pmatrix}
  -X_3     & X_1-iX_2\\
  X_1+iX_2 & X_3
 \end{pmatrix}.
\]
Le produit scalaire de $\R^3$ induit sur $E^3$ le produit scalaire suivant
\[
 \left\langle X,Y\right\rangle = \frac12 \tr \left(\widetilde X\widetilde Y \right),
\]
et la norme euclidienne d'un vecteur $X$ est donn\'ee par l'oppos\'e du d\'eterminant de la matrice $\widetilde X$
\[
 X_1^2+X_2^2+X_3^2 =- \det \widetilde X.
\]
Pour toute matrice $A\in\S$, l'application
\[
 R_A : M\mapsto \bar A^t M A
\]
est une isom\'etrie directe de $E^3$ pour ce produit scalaire. On identifie $SO(E^3)$ avec le groupe $SO(3)$ des rotations de $\R^3$ : pour toute matrice $A\in\S$, on appelle aussi $R_A$ la rotation correspondante dans $SO(3)$ et pour tout vecteur $X\in \R^3$, on a
\[
 \w{(R_AX)} = \bar A^t \w X A.
\]
On obtient le morphisme de groupe
\begin{align*}
 R :\ \S &\to SO(3)\\
     A  &\mapsto R_A
\end{align*}
qui est le rev\^etement \`a deux feuillets de $SO(3)$ par le groupe $Spin(3)\simeq\S$. On peut expliciter ce morphisme : si $R\in SO(3)$ est une rotation d'angle $\varphi$ et d'axe unitaire $\delta=(\delta_1,\delta_2,\delta_3)$,
alors les deux relev\'es de $R$ sont $A$ et $-A$ avec
\begin{equation}
    A = \cos\left(\frac{\varphi}{2}\right)\I_2 -
    i\sin\left(\frac{\varphi}{2}\right)
 \begin{pmatrix}
  -\delta_3          & \delta_1-i\delta_2\\
  \delta_1+i\delta_2 & \delta_3
 \end{pmatrix}.
\label{def-R_A}
\end{equation}
Rappelons que si on pose
\[
 J =
\begin{pmatrix}
 0 & -1\\
 1 & 0
\end{pmatrix},
\]
alors pour toute matrice $M\in SU(2)$, on a
\begin{equation}
 MJ=J\bar{M}.
\label{SU2}
\end{equation}

La proposition suivante explicite le caract\`ere spinoriel de la repr\'esentation de Weierstrass~\eqref{rep-W}.

\begin{prop}
 Soit $X:\C_+\to\R^3$ une immersion conforme minimale de donn\'ees de Weierstrass $Y=(G,H)$. Soit une matrice $A$ dans $SU(2)$. Alors le vecteur $ Y A$ constitue les donn\'ees de Weierstrass de l'immersion conforme minimale $ R_A\left( X\right) $ image de l'immersion $X$ par la rotation $R_A$.
\label{prop-spin}
\end{prop}

\begin{proof}
Supposons que l'immersion $X=\left(X_1,X_2,X_3 \right) :\C_+\to\R^3$ soit donn\'ee par le vecteur $Y$ par la formule de Weierstrass~\eqref{rep-W} (\ie $X_0=O$). Il suffit d'\'ecrire l'immersion $X$ en terme de matrices $2\times2$~:
\[
 \widetilde X(x) =
 \begin{pmatrix}
  -X_3(x)         & X_1(x)-iX_2(x)\\
  X_1(x)+iX_2(x) & X_3(x)
 \end{pmatrix}.
\]
Calculons $X_1+iX_2$~:
\begin{align*}
 X_1(x)+iX_2(x) &= \frac i2 \int_{x_0}^x\left(G(\xi)^2 - H(\xi)^2\right)d\xi
    -\frac i2 \int_{x_0}^x\left(\overline G(\xi)^2 - \overline H(\xi)^2\right)d\bar \xi\\
          &\quad +\frac i2 \int_{x_0}^x\left(G(\xi)^2 + H(\xi)^2\right)d\xi
    +\frac i2 \int_{x_0}^x\left(\overline G(\xi)^2 + \overline H(\xi)^2\right)d\bar \xi\\
          &= i\int_{x_0}^x G(\xi)^2d\xi + i\int_{x_0}^x \overline H(\xi)^2d\bar \xi.
\end{align*}
On obtient donc
\[
 \widetilde X(x) = i\int_{x_0}^x
 \begin{pmatrix}
  -GH & -H^2\\
  G^2 & GH
 \end{pmatrix}
 d\xi
 +i\int_{x_0}^x
 \begin{pmatrix}
  \overline G\overline H  & -\overline G^2\\
   \overline H^2          & -\overline G\overline H
 \end{pmatrix}
 d\bar \xi,
\]
ce que l'on peut \'ecrire sous la forme
\[
 \widetilde X(x) = i\int_{x_0}^x
 J\cdot Y(\xi)^t\cdot Y(\xi)
 d\xi
 +i\int_{x_0}^x
 \overline Y(\xi)^t\cdot \overline Y(\xi)\cdot J
 d\bar \xi.
\]
Par l'identit\'e~\eqref{SU2}, on trouve
\[
 \bar A^t \widetilde X(x) A = i\int_{x_0}^x
 J\cdot \left(Y(\xi)A\right)^t \cdot \left(Y(\xi)A\right)
 d\xi
 +i\int_{x_0}^x
 \overline{\left(Y(\xi)A\right)}^t\cdot \overline{\left(Y(\xi)A\right)}\cdot J
 d\bar \xi.
\]
Les donn\'ees de Weierstrass $YA$ d\'efinissent donc l'immersion conforme minimale $R_A\left( X\right)$.
\end{proof}

On reprend les notations de la section pr\'ec\'edente.

\begin{lemm}
Soit $X:\mathbb D^+\to\R^3$ une immersion conforme minimale de donn\'ees de Weierstrass $Y:\mathbb D^+\to\C^2$. On suppose que $Y$ s'\'etend contin\^ument \`a $]-1,1[$. Alors
\begin{itemize}
 \item l'image par $X$ de l'intervalle $]-1,1[$ est un segment de droite si et seulement s'il existe une matrice $A\in SU(2)$ telle que le vecteur $YA$ soit \`a valeurs r\'eelles ou purement imaginaires sur $]-1,1[$ ;
 \item l'image par $X$ de l'intervalle $]-1,1[$ est une courbe contenue dans un plan que la surface coupe perpendiculairement si et seulement s'il existe une matrice $A\in SU(2)$ telle que le vecteur $e^{i\frac{\pi}{4}}YA$ soit \`a valeurs r\'eelles ou purement imaginaires sur $]-1,1[$.
\end{itemize}
\label{lemme-realite}
\end{lemm}

\begin{proof}
Soit $Y=(G,H)$ les donn\'ees de Weierstrass de l'immersion $X$. Pour la premi\`ere assertion, on va montrer que l'image de $]-1,1[$ par l'immersion $X$ est un segment de droite dirig\'e par le vecteur de base $e_2=(0,1,0)$ si et seulement si les fonctions $G^2(x)$, $H^2(x)$ et $G(x)H(x)$ sont r\'eelles sur $]-1,1[$, c'est-\`a-dire si et seulement si les fonctions $G(x)$ et $H(x)$ sont toutes les deux r\'eelles ou purement imaginaires. On en d\'eduit alors la premi\`ere assertion par la proposition~\ref{prop-spin}.

La condition suffisante est imm\'ediate. Pour la n\'ecessit\'e, il faut exprimer par exemple que sur $]-1,1[$, la troisi\`eme composante $X_3(x)$ de l'immersion est constante et que son application de Gauss $N(x)$ est orthogonale au vecteur $e_2$. Comme la projection st\'er\'eographique nord de $N(x)$ est $-G(x)/H(x)$, on obtient que sur $]-1,1[$
\[
    \left\{
        \begin{array}{ll}
            -G/H \in \R \\
            GH\in \R
        \end{array}
    \right., \qquad
    \text{i.e.:} \quad
    \left\{
        \begin{array}{ll}
            G\overline{H}=\overline{G}H \\
            GH=\overline{GH}
        \end{array}.
    \right.
\]
Ceci donne le r\'esultat annonc\'e, puisque les donn\'ees de Weierstrass $G(x)$ et $H(x)$ ne peuvent pas \^etre simultan\'ement nulles.

Pour la deuxi\`eme assertion, il suffit de consid\'erer l'immersion conjugu\'ee $X^*$, qui a pour donn\'ees de Weierstrass $e^{i\frac{\pi}{4}}Y$. Alors le lemme~\ref{lemme-SM-conj} nous permet de nous ramener au cas pr\'ec\'edent.
\end{proof}

Comme on va le voir \`a la section~\ref{section-mono}, le lemme~\ref{lemme-realite} permet de retrouver les principes de r\'eflexion de Schwarz.


\chapter{\'Equations fuchsiennes et syst\`emes fuchsiens}
\label{chapitre-prelim2}

On pr\'esente dans ce chapitre les notions de base de la th\'eorie des \'equations et syst\`emes fuchsiens sur la sph\`ere de Riemann. On commence par \'etudier les \'equations fuchsiennes, on donne ensuite les r\'esultats analogues pour les syst\`emes d'\'equations, et enfin, on pr\'ecise les liens entre syst\`emes fuchsiens et \'equations fuchsiens (dans le cas non r\'esonnant), dont on aura besoin au chapitre~\ref{chapitre-isomono}. Pour une approche plus compl\`ete, ainsi que pour conna\^itre les d\'emonstrations des r\'esultats \'enonc\'es, on pourra se reporter \`a~\cite{IKSY} --- particuli\`erement pour ce qui concerne les transformations isomonodromiques, que ce soit le syst\`eme de Garnier ou le syst\`eme de Schlesinger. Pour le probl\`eme de Riemann--Hilbert pour les syst\`emes fuchsiens, on pourra se r\'ef\'erer \`a Anosov et Bolibruch~\cite{AB}, ou plus simplement \`a~\cite{Beauville} pour une pr\'esentation g\'en\'erale du probl\`eme et des r\'esultats de Bolibruch.

\section{\'Equations fuchsiennes}

On consid\`ere une \'equation diff\'erentielle lin\'eaire du second ordre d\'efinie sur la sph\`ere de Riemann $\P=\C\cup\{\infty\}$
\begin{equation}
 D^2y + p(x)Dy+ q(x)y=0
\tag{$E$}
\label{E-chap1}
\end{equation}
o\`u $D=\frac{d}{dx}$ d\'esigne la d\'erivation par rapport \`a la variable complexe $x\in\C$. On suppose que les c\oe fficients $p(x)$ et $q(x)$ sont des fonctions m\'eromorphes sur $\P$. On note $S$ l'ensemble des singularit\'es de l'\'equation~\eqref{E-chap1}, \ie  des points en lesquels $p(x)$ ou $q(x)$ a un p\^ole
\[
 S=\{x_1,\ldots,x_n\}.
\]
Les solutions de l'\'equation~\eqref{E-chap1} sont des fonctions \emph{multi-formes} dans $\P\ssm S$, c'est-\`a-dire des fonctions holomorphes dans le rev\^etement universel de $\P\ssm S$. Par abus de langage, on notera encore $y(x)$ une telle fonction. Ces solutions forment un espace vectoriel de dimension $2$. On appelle \emph{syst\`eme fondamental de solutions} un vecteur $Y(x)=\left( y_1(x),y_2(x)\right) $ dont les composantes forment une base de cet espace.

\subsection{\'Etude locale}

On commence par \'etudier le comportement des solutions de l'\'equation~\eqref{E-chap1} au voisinage de ses singularit\'es. On en d\'eduira ensuite une caract\'erisation globale des \'equations fuchsiennes.

\subsubsection*{Singularit\'es r\'eguli\`eres et singularit\'es fuchsiennes}

En g\'en\'eral, les solutions de l'\'equation~\eqref{E-chap1} ne sont pas uniformes au voisinage d'une singularit\'e. On distingue certains types de singularit\'es.

\begin{defi}
 On dit qu'une singularit\'e $x=x_0$ de l'\'equation~\eqref{E-chap1} est \emph{fuchsienne} si la fonction $p(x)$ a en $x=x_0$ un p\^ole d'ordre au plus $1$ et la fonction $q(x)$ un p\^ole d'ordre au plus $2$.
\label{def-sg-fu}
\end{defi}

On distingue une autre cat\'egorie de singularit\'es : on consid\`ere les singularit\'es $x=x_0$ au voisinage desquelles toute solution a une croissance au plus polynomiale en $1/|x-x_0|$ quand $x\to x_0$. Comme \emph{a priori} une solution de l'\'equation~\eqref{E-chap1} a un point de branchement logarithmique en une singularit\'e, il faut \^etre plus pr\'ecis dans cette d\'efinition.

\begin{defi}
 On dit qu'une singularit\'e $x=x_0$ de l'\'equation~\eqref{E-chap1} est \emph{r\'eguli\`ere} si pour tout secteur $S$ centr\'e en $x_0$, pour tout rev\^etement $\widetilde S$ de ce secteur dans le rev\^etement de $\P\ssm S$ et pour toute solution $y$ de l'\'equation~\eqref{E-chap1}, la restriction $y_{|\widetilde S}$ a une croissance polynomiale en $1/|x-x_0|$ quand $x\to x_0, \ x\in S$.
\label{def-sg-reg}
\end{defi}

Comme on va le voir, une singularit\'e fuchsienne est toujours r\'eguli\`ere. Pour les \'equations, la r\'eciproque est \'egalement vraie (\cite{Hartman}), mais elle est fausse en g\'en\'eral pour les syst\`emes d'\'equations.

\subsubsection*{M\'ethode de Fr\"obenius}

La m\'ethode de Fr\"obenius permet de d\'ecrire le comportement local des solutions de l'\'equation~\eqref{E-chap1} au voisinage d'une singularit\'e fuchsienne. On se place au point $x=0$ en supposant qu'il est une telle singularit\'e.

Si on cherche les solutions formelles de l'\'equation~\eqref{E-chap1} de la forme
\[
 y(x) = x^s \sum_{n=0}^\infty b_nx^n,
\]
on se rend compte que le nombre complexe $s$ ne peut prendre au plus que deux valeurs, qui sont les racines de l'\'equation quadratique
\begin{equation}
 s^2 + (a-1)s + b,
\label{eq-cara}
\end{equation}
avec
\[
 a =\lim_{x\to0}xp(x), \quad
 b =\lim_{x\to0}x^2q(x).
\]
L'\'equation~\eqref{eq-cara} s'appelle l'\emph{\'equation caract\'eristique} de l'\'equation~\eqref{E-chap1} en la singularit\'e fuchsienne $x=0$. Ses racines s'appellent les \emph{exposants} en $x=0$. Si on les note $s_1$ et $s_2$ avec
\[
 \Re s_2 \leq \Re s_1,
\]
alors on peut v\'erifier qu'il existe toujours une solution convergente (multi-valu\'ee) $y_1(x)$ de l'\'equation~\eqref{E-chap1} de la forme
\[
 y_1(x) = x^{s_1} \sum_{n=0}^\infty b_nx^n, \qquad b_0=1.
\]
Pour expliciter une autre solution lin\'eairement ind\'ependante de $y_1(x)$, il faut distinguer deux cas~:
\begin{itemize}
 \item s'il existe \'egalement une solution convergente $y_2(x)$ de la forme
\[
 y_2(x) = x^{s_2} \sum_{n=0}^\infty b_nx^n, \qquad b_0=1,
\]
alors la singularit\'e fuchsienne $x=0$ est dite \emph{non logarithmique}. En particulier, c'est toujours le cas si $s_1-s_2$ n'est pas un entier naturel ;
 \item sinon, la singularit\'e fuchsienne $x=0$ est dite \emph{logarithmique}, et la deuxi\`eme solution canonique en $x=0$ est de la forme
\[
 y_1(x)\log x + x^{s_1} \sum_{n=0}^\infty c_nx^n + x^{s_2} \sum_{n=0}^\infty d_nx^n.
\]
\end{itemize}

On peut observer que la singularit\'e fuchsienne $x=0$ est non logarithmique si et seulement s'il existe un syst\`eme fondamental de solutions $Y(x)$ dont la matrice de monodromie en $x=0$ soit diagonale.

Les expressions que l'on vient de donner pour les solutions de l'\'equation~\eqref{E-chap1} au voisinage d'une singularit\'e fuchsienne montrent qu'une singularit\'e fuchsienne est r\'eguli\`ere.

\subsubsection*{\'Equations fuchsiennes}

Il nous reste \`a \'etudier le point $x=\infty$. Pour cela, on fait le changement de variable $w=1/x$ dans l'\'equation~\eqref{E-chap1}, et la nature du point $x=\infty$ est celle du point $w=0$ dans la nouvelle \'equation. On montre ainsi facilement que le point $x=\infty$ est une singularit\'e fuchsienne de l'\'equation~\eqref{E-chap1} si et seulement si les fonctions
\[
    w^{-1}p\left(w^{-1}\right), \quad w^{-2}q\left(w^{-1}\right)
\]
sont holomorphes au point $w=0$. On note alors $a_\infty$ et $b_\infty$ leurs valeurs respectives en $w=0$, et l'\'equation caract\'eristique au point $x=\infty$ est
\[
    s^2+(1-a_\infty)s+b_\infty=0.
\]

\begin{defi}
 On dit que l'\'equation~\eqref{E-chap1} est une \emph{\'equation fuchsienne} sur la sph\`ere de Riemann $\P$ si toutes ses singularit\'es, y compris \'eventuellement le point en l'infini, sont fuchsiennes.
\end{defi}

On obtient alors la caract\'erisation suivante des \'equations fuchsiennes.

\begin{prop}
 L'\'equation~\eqref{E-chap1} est fuchsienne sur la sph\`ere de Riemann $\P$, de singularit\'es $x_1,\ldots,x_{n-1},x_n=\infty$, si et seulement si ses c\oe fficients sont de la forme
\[
 p(x) = \sum_{i=1}^{n-1} \frac{a_i}{x-x_i}, \qquad
 q(x) = \sum_{i=1}^{n-1} \frac{b_i}{(x-x_i)^2} + \sum_{i=1}^{n-1} \frac{c_i}{x-x_i},
\]
avec
\[
 \sum_{i=1}^{n-1} c_i =0.
\]
\label{prop-cara-eq-fu}
\end{prop}

On range dans un tableau appel\'e \emph{sch\'ema de Riemann} les singularit\'es fuchsiennes de l'\'equation~\eqref{E-chap1}, et les exposants $\t^+_i$ et $\t^-_i$ en chaque singularit\'e $x=x_i$~:
\begin{equation}
  \begin{pmatrix}
  x=x_1 & \cdots & x=x_n\\
  \t^+_1  & \cdots & \t^+_n\\
  \t^-_1  & \cdots & \t^-_n\\
 \end{pmatrix}.
\label{SdR-chap1}
\end{equation}

\begin{prop}[Relation de Fuchs]
Supposons que l'\'equation~\eqref{E-chap1} soit fuchsienne et que son sch\'ema de Riemann soit donn\'e par~\eqref{SdR-chap1}. Alors la somme de tous les exposants de \eqref{E-chap1} ne d\'epend que du nombre de singularit\'es, et plus pr\'ecis\'ement
\begin{equation}
 \sum_{i=1}^n (\t^+_i+\t^-_i) = n-2.
\label{rel-Fu}
\end{equation}
\label{prop-rel-Fu}
\end{prop}

\begin{proof}
Il suffit d'\'ecrire que la somme des r\'esidus du c\oe fficient $p(x)$ est nulle. Par la proposition~\ref{prop-cara-eq-fu}, on a
\[
 p(x) = \sum_{i=1}^{n-1} \frac{a_i}{x-x_i}
\]
et par d\'efinition du r\'esidu $a_\infty$, on a $ a_\infty = \sum_{i=1}^{n-1} a_i$. D'apr\`es les \'equations caract\'eristiques en chacune des singularit\'es, on d\'eduit
\[
 a_i = 1-\t^+_i-\t^-_i \ (i=1,\ldots,n-1),\quad a_\infty = 1+\t^+_n+\t^-_n,
\]
ce qui permet de conclure.
\end{proof}

\subsection{\'Equations projectivement \'equivalentes et schwarzien}

\'Etant donn\'e une fonction $u$ non constante et m\'eromorphe dans un ouvert $U$ d'une surface de Riemann, le \emph{schwarzien} de $u$ par rapport \`a une coordonn\'ee conforme $x$ est donn\'e par
\[
 S_x(u) = \left( \frac{u''}{u'}\right) ' - \frac12 \left( \frac{u''}{u'}\right)^2
\]
o\`u $u'=\frac{du}{dx}$. Si $z$ est une autre coordonn\'ee conforme, alors $S_z(u)=S_z(x)+S_x(u)\left(\frac{dx}{dz} \right)^2 $. De plus, le schwarzien est invariant sous l'action de $PGL(2,\C)$~:
\[
 S_x\left( \frac{au+b}{cu+d}\right) = S_x(u) \qquad \text{pour tout }
\begin{pmatrix}
 a & b\\
 c & d
\end{pmatrix}
\in GL(2,\C).
\]
Ces deux propri\'et\'es assurent en particulier que le schwarzien $S_x(u)$ est identiquement nul si et seulement si la fonction $u$ est une homographie $u(x)=\frac{ax+b}{cx+d}$. Une fonction $u$ est dite $PGL(2,\C)$-multi-forme si deux branches arbitraires de $u(x)$ sont reli\'ees par une homographie. Si une fonction est $PGL(2,\C)$-multi-forme, alors son schwarzien est uniforme.

\bigskip

Pour tout syst\`eme fondamental de solutions $Y(x)=\left( y_1(x),y_2(x)\right) $ de l'\'equation~\eqref{E-chap1}, le schwarzien du rapport $u=\frac{y_1}{y_2}$ est ind\'ependant du choix de $Y(x)$ et vaut
\begin{equation}
  S_x\left(\frac{y_1}{y_2} \right) = 2q(x) - \frac12 p(x)^2 - Dp(x).
\label{schwarzien}
\end{equation}
Le rapport $\frac{y_1}{y_2}$ est d\'efini \`a partir de l'\'equation~\eqref{E-chap1} \`a une homographie pr\`es.

\begin{defi}
 La classe d'\'equivalence du rapport de deux solutions lin\'eairement ind\'ependantes de l'\'equation~\eqref{E-chap1} est appel\'ee \emph{la solution projective} de l'\'equation~\eqref{E-chap1}. Deux \'equations diff\'erentielles lin\'eaires du second ordre \`a c\oe fficients m\'eromorphes dans la sph\`ere de Riemann sont dites \emph{projectivement \'equivalentes} si elles ont la m\^eme solution projective.
\end{defi}

Soient deux \'equations $(E_1)$ et $(E_2)$ ayant le m\^eme ensemble de singularit\'es $S$. Alors elles sont projectivement \'equivalentes si et seulement s'il existe une fonction $\Phi(x)$ holomorphe et jamais nulle dans le rev\^etement universel de l'ensemble $\P\ssm S$ telle que toute solution $y_2(x)$ de l'\'equation $(E_2)$ soit obtenue par la multiplication d'une solution $y_1(x)$ de l'\'equation $(E_1)$ par la fonction $\Phi(x)$. La fonction $\Phi(x)$ est alors de la forme
\[
 \Phi(x) = \prod_{a\in S\ssm\{\infty\}} (x-a)^{\t_a}.
\]

\subsection{Monodromie}

On ne suppose pas que l'\'equation~\eqref{E-chap1} est fuchsienne. On a vu qu'en g\'en\'eral, les solutions de l'\'equation~\eqref{E-chap1} sont des fonctions multi-formes dans $\P\ssm S$. Pour mesurer ce d\'efaut d'uniformit\'e de ses solutions, on introduit la \emph{monodromie} de l'\'equation~\eqref{E-chap1}, qui est une classe d'\'equivalence de repr\'esentations du groupe fondamental de l'ensemble $\P\ssm S$.

Soient un point $x_0\in\P\ssm S$ et un ouvert simplement connexe $U$ de $\P\ssm S$ contenant $x_0$. On consid\`ere un syst\`eme fondamental de solutions $Y(x)$ de l'\'equation~\eqref{E-chap1} d\'efini dans $U$. On peut prolonger analytiquement le syst\`eme $Y(x)$ le long de tout lacet de point de base $x_0$ et contenu dans $\P\ssm S$, et ce prolongement ne d\'epend que de la classe d'homotopie du lacet. Pour toute classe $\a$ dans le groupe fondamental $\pi_1(\P\ssm S,x_0)$, on peut donc noter $\a\ast Y(x)$ le prolongement du syst\`eme $Y(x)$ le long de tout repr\'esentant de $\a$. Alors le syst\`eme $\a\ast Y(x)$ est d\'efini dans $U$ et il est aussi un syst\`eme fondamental de solutions de l'\'equation~\eqref{E-chap1}. Il existe donc une unique matrice $M_\a(Y)\in GL(2,\C)$ qui v\'erifie
\[
 \a\ast Y(x) = Y(x) M_\a(Y).
\]
On appelle la matrice $M_\a(Y)$ la matrice de monodromie du syst\`eme $Y(x)$ le long de $\a$. On d\'efinit ainsi une application
\[
 \rho_Y : \pi_1(\P\ssm S,x_0)\to GL(2,\C), \quad \a\mapsto M_\a(Y).
\]
On choisit un ordre dans le groupe fondamental $\pi_1(\P\ssm S,x_0)$ de la fa\c con suivante : on d\'efinit le produit $\b\a$ de deux \'el\'ements $\a, \b\in\pi_1(\P\ssm S,x_0)$ comme \'etant la classe du lacet qui suit d'abord $\a$ puis $\b$ (dans le sens naturel). On a alors $(\b\a)\ast Y(x) = \b\ast\left( \a\ast Y\right) (x)$, donc
\[
 M_{\b\a}(Y) = M_\b(Y)M_\a(Y),
\]
et l'application $\rho_Y$ est un hom\'eomorphisme du groupe $\pi_1(\P\ssm S,x_0)$ dans $GL(2,\C)$ : c'est une \emph{repr\'esentation lin\'eaire} de rang $2$ (si on inverse l'ordre dans $\pi_1(\P\ssm S,x_0)$, on obtient une anti-repr\'esentation). On appelle l'application $\rho_Y$ la \emph{repr\'esentation de monodromie} de l'\'equation~\eqref{E-chap1} par rapport au syst\`eme fondamental $Y(x)$.

Consid\'erons \`a pr\'esent un autre syst\`eme fondamental de solutions $Z(x)$ d\'efini dans l'ouvert $U$. Il existe une unique matrice $C\in GL(2,\C)$, appel\'ee \emph{matrice de connexion} entre les syst\`emes $Y(x)$ et $Z(x)$, telle que
\[
 Z(x)=Y(x) C.
\]
Alors pour tout $\a\in\pi_1(\P\ssm S,x_0)$, on a
\[
 \a\ast Z(x) = \a\ast Y(x)\cdot C = Y(x) M_\a(Y)C = Z(x) C^{-1}M_\a(Y)C,
\]
c'est-\`a-dire
\begin{equation}
 M_\a(Z) = C^{-1}M_\a(Y)C.
\label{mono-conjuguees}
\end{equation}
Les deux repr\'esentations $\rho_Y$ et $\rho_Z$ sont donc conjugu\'ees. La relation de conjugaison entre repr\'esentations est une relation d'\'equivalence. On voit donc que l'ensemble de toutes les repr\'esentations de monodromie de l'\'equation~\eqref{E-chap1} (par rapport \`a chaque syst\`eme fondamental) constitue une classe de conjugaison. Cette classe est canoniquement associ\'ee \`a l'\'equation~\eqref{E-chap1} : on l'appelle la \emph{monodromie} de l'\'equation~\eqref{E-chap1}.

Le groupe fondamental $\pi_1(\P\ssm S,x_0)$ est engendr\'e par les classes
de lacets $\ga_1,\ldots,\ga_n$ tournant respectivement une fois dans
le sens direct autour de la singularit\'e $x=x_i$, en laissant les
autres singularit\'es \`a l'ext\'erieur, soumises \`a la relation
$\ga_n\cdots\ga_1=1$. La repr\'esentation de monodromie $\rho_Y$ par
rapport \`a un syst\`eme $Y(x)$ est donc d\'etermin\'ee par la famille
$(M_1,\ldots,M_n)$, o\`u
\[
 M_i = M_{\ga_i}(Y).
\]
Les matrices $M_i$ v\'erifient aussi
\[
 M_n\cdots M_1 = \I_2.
\]
On appelle la famille $(M_1,\ldots,M_n)$ un \emph{syst\`eme de g\'en\'erateurs} de la monodromie de l'\'equation~\eqref{E-chap1}.

\begin{defi}
 Une repr\'esentation $\rho : G\to GL(m,\C)$ d'un groupe $G$ est dite \emph{irr\'eductible} si les sous-espaces vectoriels de $\C^m$ invariants par $\rho$ sont exactement $\{0\}$ et $\C^m$.
\end{defi}

La monodromie de l'\'equation~\eqref{E-chap1} est dite \emph{irr\'eductible} si elle admet un repr\'esentant irr\'eductible, c'est-\`a-dire si elle admet un syst\`eme de g\'en\'erateurs $(M_1,\ldots,M_n)$ constitu\'e de matrices qui ne soient pas simultan\'ement trigonalisables. Si l'\'equation~\eqref{E-chap1} est fuchsienne, alors le fait qu'elle ait une monodromie irr\'eductible est \'equivalent \`a ce qu'elle soit elle-m\^eme irr\'eductible, \ie que l'op\'erateur diff\'erentiel
\[
 L=D^2 + p(x)D + q(x)
\]
n'admette que des factorisations triviales.

\subsection{Le probl\`eme de Riemann--Hilbert pour les \'equations lin\'eaires du second ordre}

On ne consid\`ere pour l'instant le probl\`eme de Riemann--Hilbert que dans le cas des \'equations du second ordre. Il n'y a pas de diff\'erence fondamentale avec les \'equations d'ordre sup\'erieur. Par contre, la discussion est diff\'erente dans le cas des syst\`emes fuchsiens. Le probl\`eme de Riemann--Hilbert pour les \'equations fuchsiennes est exactement le vingt-et-uni\`eme des vingt-trois probl\`emes propos\'es par Hilbert au Congr\`es International de Paris en 1900~:

\bigskip

\noindent
\textbf{Le probl\`eme de Riemann--Hilbert.} \emph{Trouver une \'equation fuchsienne ayant des singularit\'es donn\'ees et une monodromie donn\'ee.}

\bigskip
Formul\'e ainsi, on peut facilement voir que le probl\`eme de Riemann--Hilbert n'a en g\'en\'eral pas de solution. En effet, soit $S=\{x_1,\ldots,x_{n-1},x_n=\infty\}\subset\P$ un ensemble de singularit\'es. D'apr\`es la proposition~\ref{prop-cara-eq-fu}, une \'equation fuchsienne du second ordre dont l'ensemble des singularit\'es soit $S$ d\'epend de $e(S)$ param\`etres, avec
\[
 e(S) = 3n-4.
\]
Par ailleurs, on peut montrer que l'ensemble des classes de conjugaison de repr\'esentations $\rho : \pi_1(\P\ssm S,x_0)\to GL(2,\C)$ d\'epend de $m(S)$ param\`etres, avec
\[
 m(S) = 4(n-2)+1.
\]
D\`es que $n>3$, on a donc
\[
 m(S)-e(S)>0.
\]
\`A singularit\'es fix\'ees, l'application qui \`a une \'equation
fuchsienne du second ordre associe sa monodromie n'est donc pas
surjective d\`es que $n>3$. Ce calcul remonte \`a
Poincar\'e~\cite{Poincare}. Si on veut pouvoir construire une
\'equation fuchsienne dont la monodromie est donn\'ee, il faut donc
s'autoriser \`a ajouter des param\`etres suppl\'ementaires : les
singularit\'es apparentes sont les seuls param\`etres possibles.

\begin{defi}
 Une singularit\'e fuchsienne de l'\'equation~\eqref{E-chap1} est dite \emph{apparente} si elle n'est pas logarithmique et si ses exposants sont des entiers relatifs.
\label{def-sg-app}
\end{defi}

Une singularit\'e fuchsienne $x=a$ est apparente si et seulement si toutes les solutions de l'\'equation~\eqref{E-chap1} sont m\'eromorphes en $x=a$. Il n'y a donc pas de monodromie en ces singularit\'es.

On v\'erifie alors qu'une \'equation fuchsienne du second ordre ayant ses singularit\'es dans $S$, et ayant au plus $N$ singularit\'es apparentes \`a l'ext\'erieur de $S$ d\'epend de $e(S)+N$ param\`etres. Pourtant, il n'est pas \'evident qu'autoriser $N=m(S)-e(S)=n-3$ singularit\'es apparentes soit suffisant pour obtenir une r\'eponse positive au probl\`eme de Riemann--Hilbert. Lorsque la monodromie est irr\'eductible, Ohtsuki~\cite{Oht} a obtenu la bonne majoration du nombre de singularit\'es apparentes, \`a la condition qu'un des g\'en\'erateurs de la monodromie soit diagonalisable. Mais le r\'esultat le plus g\'en\'eral est d\^u \`a Bolibruch.

\begin{theo}[\cite{Bolibruch}]
 \'Etant donn\'e un ensemble fini $S\subset\P$ \`a $n$ \'el\'ements et une repr\'esentation irr\'eductible $\rho : \pi_1(\P\ssm S)\to GL(2,\C)$, il existe une \'equation fuchsienne du second ordre dont l'ensemble des singularit\'es soit $S$, dont la monodromie soit la classe de $\rho$ et ayant au plus $n-3$ singularit\'es apparentes.
\label{thm-boli}
\end{theo}

\section{Syst\`emes fuchsiens}
\label{section-sys-fu}

\subsection{D\'efinitions}

Consid\'erons un syst\`eme diff\'erentiel lin\'eaire du premier ordre
\begin{equation}
 DY=A(x)Y
 \tag{$A_0$}\label{A0}
\end{equation}
o\`u $D=\frac{d}{dx}$ et la fonction $A(x)$ est m\'eromorphe sur la sph\`ere de Riemann $\P$, \`a valeur dans $M(2,\C)$. On suppose que le syst\`eme~\eqref{A0} est \emph{fuchsien}, c'est-\`a-dire que tous les p\^oles de $A(x)$ sont simples%
\footnote{
 contrairement \`a ce qui se passe pour les \'equations, les notions de singularit\'es r\'eguli\`eres et fuchsiennes ne co\"incident pas pour les syst\`emes d'\'equations. Une singularit\'e fuchsienne, c'est-\`a-dire un p\^ole simple, est r\'eguli\`ere (\emph{cf} d\'efinition~\ref{def-sg-reg}), mais la r\'eciproque est fausse.
}
. Comme l'ensemble des syst\`eme fuchsiens sur la sph\`ere de Riemann est stable par transformation de M\"obius, on peut choisir comme pr\'ec\'edemment
\[
 t_1,\ldots,t_n,\ t_{n+1}=0,\ t_{n+2}=1,\ t_{n+3}=\infty
\]
les singularit\'es du syst\`eme~\eqref{A0}, et on a donc
\begin{equation}
 A(x)=\sum_{i=1}^{n+2}\frac{A_i}{x-t_i}.
\end{equation}
Comme on suppose que l'infini est un point singulier, le r\'esidu
\[
 A_\infty:=-\sum_{i=1}^{n+2}A_i
\]
n'est pas la matrice nulle (on note parfois $A_{n+3}$ pour $A_\infty$). On note $S(t)$ l'ensemble des singularit\'es :
\[
 S(t):=\{t_1,\ldots, t_{n+3}\}.
\]
Le syst\`eme~\eqref{A0} est donc d\'efini dans l'ensemble $ \P\ssm S(t)$. Ses solutions, qui sont des couples de fonctions d\'efinies sur le rev\^etement universel de $ \P\ssm S(t)$, forment un espace vectoriel de dimension $2$. On appelle \emph{matrice fondamentale de solutions} une matrice $\Y(x)$ dont les colonnes $Y_1(x),Y_2(x)$ forment une base de cet espace. Une telle matrice v\'erifie l'\'equation $D\Y=A(x)\Y $. On d\'efinit la monodromie du syst\`eme~\eqref{A0} comme on l'a fait pour les \'equations du second ordre.

\bigskip

On suppose de plus que le syst\`eme~\eqref{A0} v\'erifie les deux hypoth\`eses suivantes~:
\begin{itemize}
 \item le syst\`eme~\eqref{A0} est \emph{non r\'esonnant} : les valeurs propres $\t_i^+$ et $\t_i^-$ de la matrice $A_i$ satisfont $ \t_i^+-\t_i^- \notin \Z$ ($i=1,\ldots,n+3$);
 \item le syst\`eme~\eqref{A0} est \emph{normalis\'e en l'infini} :
\[
 A_\infty=-\sum_{i=1}^{n+2}A_i=
    \begin{pmatrix}
     \t_\infty^+ & 0\\
     0 & \t_\infty^-\\
    \end{pmatrix}.
\]
\end{itemize}

Comme le syst\`eme~\eqref{A0} est non r\'esonnant, les singularit\'es $x=t_i$ ne sont pas logarithmiques. Ceci assure l'existence au voisinage de chaque singularit\'e d'une matrice fondamentale de la forme suivante.

\begin{prop}
On suppose le syst\`eme~\eqref{A0} non r\'esonnant. Alors, pour tout $i=1,\ldots,n+2$, il existe une unique matrice $P_i(x)$ holomorphe au point $x=t_i$ v\'erifiant $P_i(t_i)=\I_2$ et telle que
\[
 P_i(x) (x-t_i)^{A_i}
\]
soit une matrice fondamentale de solutions du syst\`eme~\eqref{A0}, o\`u
\[
 (x-t_i)^{A_i} = \exp\left( A_i \log(x-t_i)\right) .
\]
\label{prop-sol-cano}
\end{prop}

On ne donne pas la d\'emonstration de cette proposition, mais remarquons simplement que la matrice $P_i(x)$ est solution de l'\'equation
\[
 DP_i = A(x)P_i - P_i \frac{A_i}{x-t_i}.
\]
Soit $L_i$ la matrice diagonalis\'ee de $A_i$
\[
 L_i=
\begin{pmatrix}
 \t_i^+&0\\
 0&\t_i^-
\end{pmatrix}.
\]
Alors, il existe des matrices fondamentales de solutions de la forme
\[
 R_i(x) (x-t_i)^{L_i}
\]
o\`u la matrice $R_i(x)$ est holomorphe et inversible au point $x=t_i$ et $R_i(t_i) \in GL(2,\C)$ diagonalise $A_i$
\[
 A_i = R_i(t_i)L_iR_i(t_i)^{-1}.
\]
Ces solutions sont dites \emph{canoniques} au point $x=t_i$, parce que leur matrice de monodromie en ce point est diagonale~:
\[
 \begin{pmatrix}
  e^{2i\pi\t_i^+} & 0\\
  0 & e^{2i\pi\t_i^-}
 \end{pmatrix}.
\]
En l'infini, comme le syst\`eme~\eqref{A0} est normalis\'e en l'infini, il existe une unique solution canonique de la forme
\[
 \Y_\infty(x)= R_\infty\left(\frac1x\right) x^{-L_\infty}
\]
o\`u la matrice $R_\infty(w)$ est holomorphe en $w=0$ et $R_\infty(0)=\I_2$.

\subsection{D\'eformations isomonodromiques}

On s'int\'eresse \`a pr\'esent au probl\`eme suivant : si on consid\`ere que le syst\`eme~\eqref{A0} d\'epend d'un param\`etre variable, comment d\'ecrire l'ensemble des syst\`emes fuchsiens (ou des \'equations fuchsiennes) ayant une monodromie donn\'ee ? On pr\'esente d'abord la th\'eorie g\'en\'erale des d\'eformations isomonodromiques, et on en d\'eduira le syst\`eme de Schlesinger \`a la section suivante (le syst\`eme de Garnier est quant \`a lui introduit \`a l'appendice~\ref{annexe-garnier}).

On consid\`ere une famille de syst\`emes diff\'erentiels lin\'eaires $2\times 2$ d\'ependant d'un param\`etre $t$ variant dans un ouvert simplement connexe $U$ de $\C^n$~:
\begin{equation}
 DY = A(x,t) Y
\label{sys-dep-t}
\end{equation}
o\`u la fonction $A(x,t)$ est d\'efinie dans $\P\times U$, \`a valeurs dans $M(2,\C)$. On suppose que pour tout $t\in U$ fix\'e, la fonction $x\mapsto A(x,t)$ est holomorphe en dehors d'un ensemble fini $S(t)\subset \P$ de points singuliers, et que les points de $S(t)$ sont des fonctions holomorphes de $t$. On d\'efinit le sous-ensemble $S$ de $\P\times U$ des singularit\'es du syst\`eme
\[
 S := \bigcup_{t\in U} S(t)\times\{t\},
\]
qui est donc une hypersurface. Localement, l'ensemble $S$ a autant
de composantes connexes qu'il y a de points dans les ensembles
$S(t)$ et chacune de ces composantes connexes est un graphe de
$\P\times U$ au dessus de l'ouvert $U$. Quitte \`a restreindre
l'ouvert simplement connexe $U$, on suppose que ceci est vrai dans
$U$ entier. Sans entrer dans des d\'etails techniques de topologie,
on voit que les classes d'homotopie des lacets de $\P\ssm S(t)$
bas\'es en un point $x_0(t)$ de $\P\ssm S(t)$ sont alors
ind\'ependantes de $t$. Quitte \`a restreindre de nouveau l'ouvert
$U$, on peut choisir un point de base $x_0\in\P$ ind\'ependant de
$t$. Il suffit pour cela que $x_0$ et $U$ v\'erifient
\[
 \left( \{x_0\} \times U\right) \cap S = \emptyset.
\]
On note $\pi_1\left(\P\ssm S(t), x_0 \right) $ le groupe d'homotopie correspondant.

On peut ainsi d\'efinir la monodromie de la famille de
syst\`emes~\eqref{sys-dep-t}. Soit une solution fondamentale
$\Y(x,t)$, \ie une matrice solution de~\eqref{sys-dep-t},
holomorphe et inversible en tout point $(x_0,t)$ ($t\in U$). Pour
toute classe d'homotopie $\a\in\pi_1\left(\P\ssm S(t), x_0
\right)$, le prolongement analytique $\a\ast\Y(x,t)$ de $\Y(x,t)$
le long de $\a$ est encore une solution fondamentale en $(x_0,t)$
: il existe une unique matrice $\rho_\Y(t,\a)\in GL(2,\C)$ telle
que
\[
 \a\ast\Y(x,t) = \Y(x,t) \rho_\Y(t,\a).
\]
On obtient donc une famille analytique de repr\'esentations de monodromie
\[
 \rho_\Y(t,\cdot) : \pi_1\left(\P\ssm S(t), x_0 \right) \to GL(2,\C).
\]

\begin{defi}
 Une solution fondamentale $\Y(x,t)$ est dite \emph{$M$-invariante} si sa repr\'esentation de monodromie $\rho_\Y(t,\cdot)$ est ind\'ependante de $t$.
\label{def-M-inv}
\end{defi}

\begin{defi}
 La famille~\eqref{sys-dep-t} de syst\`emes diff\'erentiels est dite \emph{isomonodromique} si elle admet une solution fondamentale qui soit $M$-invariante.
\label{def-fam-isomono}
\end{defi}

On note $d$ la diff\'erentiation par rapport \`a la variable $t=(t_1,\ldots,t_n)$
\[
 d=\sum_{i=1}^n \frac{\partial}{\partial t_i}dt_i.
\]
On a les r\'esultats suivants.

\begin{lemm}
 Une solution fondamentale $\Y(x,t)$ est $M$-invariante si et seulement si la $1$-forme \`a valeurs matricielles
\[
 \Omega(x,t) := d\Y(x,t)\Y(x,t)^{-1}
\]
est uniforme dans $\left( \P\times U\right) \ssm S$.
\label{lemme-Omega}
\end{lemm}

\begin{prop}
 Le syst\`eme de Pfaff
\begin{equation}
 \begin{split}
  DY &= A(x,t) Y \\
  dY &= \Omega(x,t) Y
 \end{split}
\label{(sys-pfaff)}
\end{equation}
est compl\`etement int\'egrable si et seulement si le syst\`eme suivant est v\'erifi\'e
\begin{equation}
 \begin{split}
  dA(x,t) &= D\Omega(x,t) + \left[\Omega(x,t), A(x,t) \right]\\
  d\Omega(x,t) &= \Omega(x,t) \wedge \Omega(x,t).
 \end{split}
\label{defo-eq}
\end{equation}
\label{prop-defo-isomono}
\end{prop}

Le syst\`eme~\eqref{(sys-pfaff)} s'\'ecrit
\[
 d_{(x,t)} Y = \omega Y
\]
o\`u la $1$-forme $\omega$ est d\'efinie par $\omega = Adx+\Omega$, et $d_{(x,t)}$ est la diff\'erentiation par rapport \`a la variable $(x,t)$. S'il existe une matrice inversible $\Y(x,t)$ telle que  $\omega = d_{(x,t)}\Y\cdot \Y^{-1}$, alors on a de mani\`ere imm\'ediate
\[
 d_{(x,t)}\omega = \omega \wedge \omega,
\]
o\`u le produit ext\'erieur $\a \wedge \a$ d'une $1$-forme
$\a=(\a_{ij})$ \`a valeurs dans $M(2,\C)$ est la matrice dont
l'\'el\'ement $(i,j)$ est $ \a_{i1} \wedge\a_{1j} + \a_{i2}
\wedge\a_{2j} $. La r\'eciproque constitue le th\'eor\`eme de Fr\"obenius.
La condition n\'ecessaire et suffisante d'int\'egrabilit\'e
$d_{(x,t)}\omega = \omega \wedge \omega$ est exactement le
syst\`eme~\eqref{defo-eq}. Le syst\`eme~\eqref{defo-eq} s'appelle
l'\emph{\'equation de d\'eformation} de~\eqref{sys-dep-t}.

La proposition~\ref{prop-defo-isomono} nous dit donc que le syst\`eme~\eqref{sys-dep-t} admet une solution fondamentale $M$-invariante $\Y(x,t)$ si et seulement si le syst\`eme~\eqref{defo-eq} admet une solution $\Omega(x,t)$ uniforme dans $\left( \P\times U\right) \ssm S$. La solution fondamentale $\Y(x,t)$ v\'erifie alors
\[
   D\Y = A(x,t)\Y, \qquad  d\Y = \Omega(x,t) \Y.
\]

La $1$-forme $\Omega$ d\'epend du choix d'une solution fondamentale $M$-invariante. La proposition suivante permet de comparer entre elles les solutions fondamentales $M$-invariantes.

\begin{prop}
On suppose que la famille de syst\`emes~\eqref{sys-dep-t} est
isomonodromique, de monodromie irr\'eductible. Soit une
solution fondamentale $\Y_1(x,t)$ $M$-invariante. Alors une
solution fondamentale $\Y_2(x,t)$ est aussi $M$-invariante si et
seulement s'il existe une fonction holomorphe $\mu :U\to \C^*$ et
une matrice $C\in GL(2,\C)$ ind\'ependante de $t$ telles que
\[
 \Y_2(x,t) = \mu(t)\Y_1(x,t)\cdot C.
\]
\label{prop-sol-M-inv}
\end{prop}

\subsection{Le syst\`eme de Schlesinger}

On applique les r\'esultats pr\'ec\'edents \`a la d\'eformation d'un syst\`eme fuchsien non r\'esonnant. On pose
\begin{equation}
 \B^n= \left\{ (t_1,\ldots,t_n) \in (\C^* \ssm \{1\})^n \quad |\quad \forall i\neq j \quad t_i \neq t_j \right\},
\label{def-B}
\end{equation}
et on consid\`ere \`a pr\'esent la position des singularit\'es $t=(t_1,\ldots,t_n)\in\B^n$ comme un param\`etre du syst\`eme~\eqref{A0}, dont d\'ependent les matrices $A_i=A_i(t)$. On suppose que les valeurs propres $\t_i^+$ et $\t_i^-$ ($i=1,\ldots,n+3$) sont ind\'ependantes de $t$. Soit $U$ un ouvert simplement connexe de l'ensemble $\B^n$. Les d\'eformations de param\`etre $t \in U$ du syst\`eme~\eqref{A0} qui pr\'eservent la monodromie sont gouvern\'ees par le syst\`eme de Schlesinger :

\begin{theo}
 On suppose le syst\`eme fuchsien~\eqref{A0} non r\'esonnant et normalis\'e en l'infini. Alors la matrice fondamentale de solutions $\mathbf Y_\infty(x,t)$ est $M$-invariante si et seulement si les matrices $A_i(t)$, $i=1,\ldots,n+2$, satisfont le syst\`eme aux d\'eriv\'ees partielles (\emph{syst\`eme de Schlesinger})
\begin{equation}
 dA_i = \sum_{\substack{j=1\\ j\neq i}}^{n+2} [A_j,A_i] d\log(t_i-t_j),  \qquad i=1,\ldots,n+2.
\label{schlesinger}
\end{equation}

De plus, le syst\`eme de Schlesinger~\eqref{schlesinger} est compl\`etement int\'egrable.
\label{thm-schlesinger}
\end{theo}

De mani\`ere plus d\'etaill\'ee, le syst\`eme de Schlesinger s'\'ecrit
\begin{eqnarray*}
  &\dfrac{\partial A_i}{\partial t_j} = \dfrac{[A_i,A_j]}{t_i-t_j} \qquad & i=1,\ldots,n+2,\; j=1,\ldots,n, \; i\neq j\\
  &\displaystyle\sum_{i=1}^{n+2}\dfrac{\partial A_i}{\partial t_j}=0 \qquad & j=1,\ldots,n.
\end{eqnarray*}

La premi\`ere partie du th\'eor\`eme~\ref{thm-schlesinger} est obtenue en appliquant la proposition~\ref{prop-defo-isomono}. La premi\`ere \'etape consiste \`a calculer la $1$-forme $\Omega(x,t)$ associ\'ee \`a la matrice fondamentale $\Y_\infty(x,t)$ et d\'efinie au lemme~\ref{lemme-Omega}. Elle est obtenue par une \'etude locale au voisinage de chaque singularit\'e $x=t_i$ gr\^ace aux matrices fondamentales canoniques $R_i(x) (x-t_i)^{L_i}$.

\begin{lemm}
Si la matrice fondamentale de solutions $\Y_\infty(x,t)$ est $M$-invariante, alors la $1$-forme $\Omega(x,t)=d\Y_\infty(x,t) \Y_\infty(x,t)^{-1}$ s'\'ecrit
\[
 \Omega(x,t) = -\sum_{i=1}^{n+2} \frac{A_i(t)}{x-t_i} dt_i.
\]
\label{lemme-omega-sch}
\end{lemm}

On montre ensuite facilement que l'\'equation de d\'eformation~\eqref{defo-eq}
\[
   dA= D\Omega + \left[\Omega, A \right], \qquad  d\Omega = \Omega \wedge \Omega,
\]
avec
\[
 A(x,t)=\sum_{i=1}^{n+2}\frac{A_i(t)}{x-t_i}, \qquad \Omega(x,t) = -\sum_{i=1}^{n+2} \frac{A_i(t)}{x-t_i} dt_i
\]
est \'equivalente au syst\`eme de Schlesinger~\eqref{schlesinger}.

\subsection{La propri\'et\'e de Painlev\'e}

Soit une \'equation diff\'erentielle
\begin{equation}
 F\left(t,y,\frac{dy}{dt}, \ldots,\frac{d^ny}{dt^n}\right) = 0
\label{(eq-prop-P)}
\end{equation}
o\`u la fonction $F\left(t,y_0,y_1, \ldots,y_n\right)$ est polynomiale en $(y_0,y_1, \ldots,y_n)$ \`a c\oe fficients m\'eromorphes en $t$.

\begin{defi}
 On dit que l'\'equation~\eqref{(eq-prop-P)} a des points de branchement (respectivement des singularit\'es essentielles) \emph{mobiles} si ses solutions ont des points de branchement (respectivement des singularit\'es essentielles) dont la position d\'epend des constantes d'int\'egration.

On dit que l'\'equation~\eqref{(eq-prop-P)} a la \emph{propri\'et\'e de Painlev\'e} si elle n'a ni point de branchement mobile, ni singularit\'e essentielle mobile.
\label{def-prop-P}
\end{defi}

Quand $n=2$, les six \'equations de Painlev\'e $\PI,\ldots,\PVI$ constituent, \`a changement de variables pr\`es, l'ensemble des \'equations~\eqref{(eq-prop-P)} rationnelles ayant la propri\'et\'e de Painlev\'e qui ne sont ni lin\'eaires, ni int\'egrables par une quadrature.

\begin{theo}[\cite{Malgrange}, \cite{Miwa}]
 Le syst\`eme de Schlesinger~\eqref{schlesinger} a la propri\'et\'e de Painlev\'e. De plus, toute solution du syst\`eme de Schlesinger~\eqref{schlesinger} s'\'etend au rev\^etement universel de l'ensemble $\B^n$ de mani\`ere m\'eromorphe.
 \label{thm-sch-prop-P}
\end{theo}

Par contre, le syst\`eme de Garnier~\eqref{Gn}, qui d\'ecrit les d\'eformations isomonodromiques d'\'equations fuchsiennes sans singularit\'e logarithmique (voir l'appendice~\ref{annexe-garnier}), n'a pas la propri\'et\'e de Painlev\'e.

\section{Passage d'une \'equation \`a un syst\`eme d'\'equations}
\label{section-eq->sys}

Comme on va le voir au chapitre suivant, les \'equations fuchsiennes sont les objets naturellement associ\'es aux disques minimaux \`a bord polygonal. Cependant, le syst\`eme de Garnier~\eqref{Gn}, qui d\'ecrit les d\'eformations isomonodromiques de ces \'equations, n'a pas la propri\'et\'e de Painlev\'e. On va donc choisir, contrairement \`a l'approche suivie par Garnier, de transformer les \'equations fuchsiennes du second ordre que l'on va obtenir au chapitre suivant en syst\`emes fuchsiens du premier ordre de taille $2\times2$. On donne ici une description des relations entre \'equations et syst\`emes fuchsiens, dans le cas non r\'esonnant (\ie sans singularit\'e logarithmique), qui est celui qui nous int\'eresse.

\subsection{D'un syst\`eme du premier ordre \`a une \'equation du second ordre}

C'est le sens imm\'ediat. On consid\`ere un syst\`eme diff\'erentiel $2\times 2$ du premier ordre
\begin{equation}
  DY=A(x)Y, \qquad A(x)=
    \begin{pmatrix}
     A_{11}(x) & A_{12}(x)\\
     A_{21}(x) & A_{22}(x)
    \end{pmatrix},
\label{sysqcq}
\end{equation}
o\`u les fonctions $A_{ij}(x)$ sont m\'eromorphes sur la sph\`ere de Riemann.

\begin{lemm}
 Si la fonction $A_{12}(x)$ n'est pas identiquement nulle, alors la premi\`ere composante $y_1$ de toute solution $Y=(y_1,y_2)^t$ du syst\`eme~\eqref{sysqcq} v\'erifie l'\'equation du second ordre
\begin{equation}
  D^2y+p(x)Dy+q(x)y=0,
\label{eqcq}
\end{equation}
avec
\begin{align}
 p(x)&=-\frac{DA_{12}(x)}{A_{12}(x)}-\tr A(x)\label{p(A)}\\
 q(x)&=-DA_{11}(x)+A_{11}(x)\frac{DA_{12}(x)}{A_{12}(x)}+\det A(x).
\label{q(A)}
\end{align}

De plus, si $\Y(x)=(Y(x),Z(x))$ est une matrice fondamentale de solutions du syst\`eme~\eqref{sysqcq}, alors sa premi\`ere ligne $(y_1(x),z_1(x))$ est un syst\`eme fondamental de solutions de l'\'equation~\eqref{eqcq}.
\label{lemme-sys->eq}
\end{lemm}

Il est donc imm\'ediat que si le syst\`eme~\eqref{sysqcq} est fuchsien, alors l'\'equation qui lui est associ\'ee est fuchsienne. De plus, on a~:

\begin{lemm}
 Si $x=\l$ est un z\'ero de $A_{12}(x)$ d'ordre $m$, mais n'est pas une singularit\'e du syst\`eme~\eqref{sysqcq}, alors $x=\l$ est une singularit\'e apparente de l'\'equation~\eqref{eqcq}, d'exposants $0$ et $m+1$.
\label{lemme-sg-app-sys}
\end{lemm}

Consid\'erons l'\'equation associ\'ee au syst\`eme fuchsien~\eqref{A0}, toujours suppos\'e non r\'esonnant et normalis\'e en l'infini. Comme la fraction rationnelle $A_{12}(x)$ a exactement $n+2$ p\^oles simples et, par la normalisation en l'infini, un z\'ero d'ordre deux en l'infini, alors elle a exactement $n$ z\'eros dans $\C$ compt\'es avec multiplicit\'e. Supposons \`a pr\'esent que les z\'eros de la
fonction $A_{12}(x)$ sont simples. On les note $\l_1,\ldots,\l_n$, et on a donc
\[
 A_{12}(x) = \xi \frac{\Lambda(x)}{T(x)},
\]
o\`u
\begin{equation}
 \xi =\sum_{i=1}^{n+2}t_iA^i_{12}, \quad \Lambda(x)=\prod_{k=1}^n(x-\l_k), \quad
        T(x)=\prod_{i=1}^{n+2}(x-t_i).
\label{def-X-T-L}
\end{equation}
\'Etant donn\'ee la derni\`ere partie du lemme~\ref{lemme-sys->eq}, le sch\'ema de Riemann de l'\'equation~\eqref{eqcq} est donn\'e par
\begin{eqnarray}
            &\left(%
                \begin{array}{ccc}
                    x=t_i   & x=\infty       & x=\l_k\\
                    \t_i^+  & \t_\infty^+    & 0\\
                    \t_i^-  & \t_\infty^-+1  & 2\\
                \end{array}
             \right) \label{SdR-eq-sys}\\
            &\begin{array}{ccc}
                    i=1,\dotsc,n+2, & \,  & k=1,\dotsc,n \nonumber\\
             \end{array}
\end{eqnarray}
et les singularit\'es $x=\l_k$ sont apparentes. La diff\'erence entre les exposants en l'infini du syst\`eme~\eqref{A0} et de l'\'equation~\eqref{eqcq} provient de la normalisation en l'infini : puisque la matrice $A_\infty$ est diagonale, la solution canonique en l'infini $\Y_\infty(x)$ du syst\`eme~\eqref{A0} s'\'ecrit
\[
\Y_\infty(x) = \left( \I_2 + \mathcal O \left(\tfrac1x\right)\right) x^{-A_\infty},
\]
et donc la fonction $(\Y_\infty(x))_{1,1}$ est d'exposant $\t_\infty^+$, mais la fonction$(\Y_\infty(x))_{1,2}$ est d'exposant $\t_\infty^-+1$.

\subsection{Les syst\`emes fuchsiens associ\'es \`a une \'equation fuchsienne}

On consid\`ere une \'equation fuchsienne du second ordre d\'efinie sur la sph\`ere de Riemann
\begin{equation}
 D^2y+p(x)Dy+q(x)y=0,
\label{F}
\end{equation}
de sch\'ema de Riemann~\eqref{SdR-eq-sys}, dont les singularit\'es sont distinctes, telle que ses exposants v\'erifient $ \t_i^+ - \t_i^- \notin \Z$ ($i=1,\ldots,n+3$), et que les singularit\'es $x=\l_k$ sont apparentes. On peut caract\'eriser l'ensemble des syst\`emes fuchsiens~\eqref{A0} normalis\'es en l'infini qui d\'efinissent au sens du lemme~\ref{lemme-sys->eq} l'\'equation~\eqref{F}. On vient de voir que si un tel syst\`eme existe, alors son c\oe fficient $ A_{12}(x)$ est enti\`erement d\'etermin\'e par les param\`etres $t_i$ et $\l_k$ de l'\'equation~\eqref{F} et le param\`etre suppl\'ementaire $\xi $, qui est ind\'ependant de l'\'equation. Il en est en fait de m\^eme pour les autres c\oe fficients de $A(x)$. Dans~\cite{IKSY} est donn\'ee l'expression explicite des matrices $A_i$ en fonction de ces param\`etres (proposition 6.3.1. p. 208). Comme un r\'esultat aussi pr\'ecis ne nous sera pas utile par la suite, on se contente ici de donner l'existence de ces syst\`emes et de pr\'eciser leur d\'ependance en $\xi $. Comme on n'impose \`a l'avance aucune normalisation en l'infini, on obtient <<~deux fois plus~>> de syst\`emes que dans~\cite{IKSY}, \ie on obtient deux familles \`a un param\`etre de syst\`emes, \`a la place d'une seule. Dans la r\'esolution du probl\`eme de Plateau, on aura en effet besoin de pouvoir choisir la normalisation en l'infini. On ne donne pas la d\'emonstration de la proposition suivante (on pourra se reporter \`a~\cite{IKSY} ou \`a \cite{Desideri-these}, proposition~3.8).

\begin{prop}
L'ensemble des syst\`emes fuchsiens~\eqref{A0} normalis\'es en l'infini associ\'es au sens du lemme~\ref{lemme-sys->eq} \`a l'\'equation~\eqref{F} est constitu\'e de deux familles \`a un param\`etre~:
\[
DY = A^+_\xi (x) Y \quad \left(\xi \in\C^*\right)
\]
et
\[
DY = A^-_\xi (x) Y \quad \left(\xi \in\C^*\right).
\]
Ces deux familles se caract\'erisent par leur normalisation en l'infini~: pour tout $\xi \in\C^*$
\[
 \left(A_\xi ^+\right)_\infty =
    \begin{pmatrix}
     \t_\infty^+ & 0\\
     0      & \t_\infty^-\\
    \end{pmatrix}
\]
et
\[
 \left(A_\xi ^-\right)_\infty =
    \begin{pmatrix}
     \t_\infty^--1 & 0\\
     0      & \t_\infty^+-1\\
    \end{pmatrix}.
\]
De plus,
\[
 A_\xi ^\pm(x) =
    \begin{pmatrix}
    1 & 0\\
    0 & \xi 
    \end{pmatrix}
A_0^\pm(x)
    \begin{pmatrix}
    1 & 0\\
    0 & \frac1\xi 
    \end{pmatrix}
\]
o\`u les matrices $A_0^+(x)$ et $A_0^-(x)$ sont explicitement d\'etermin\'ees par l'\'equation~\eqref{F}.
\label{prop-eq->sys}
\end{prop}

Remarquons que la proposition~\ref{prop-eq->sys} permet d'\'etudier les liens entre le syst\`eme de Garnier et le syst\`eme de Schlesinger, qui sont \'etudi\'es en d\'etail dans~\cite{IKSY}, mais qui ne nous seront finalement pas utiles dans la suite.

\chapter{L'\'equation associ\'ee \`a un disque minimal \`a bord polygonal}
\label{chapitre-equ-fu}

Dans ce chapitre, on se donne une immersion conforme minimale $X:\C_+\to \R^3$ du demi-plan sup\'erieur $\C_+$ dont l'image est limit\'ee par un polygone $P$ \`a $n+3$ sommets. On note $Y_0=(G,H):\C_+\to \C^2$ ses donn\'ees de Weierstrass spinorielles. On suppose que $X$ n'est pas contenue dans un plan, et on voit alors facilement que les fonctions $G$ et $H$ sont lin\'eairement ind\'ependantes : la fonction $Y_0$ constitue un syst\`eme fondamental de solutions d'une unique \'equation diff\'erentielle lin\'eaire du second ordre
\begin{equation}
D^2y+p(x)Dy+q(x)y=0
\tag{$E$}\label{E*}
\end{equation}
o\`u $D=\frac{d}{dx}$ d\'esigne la d\'erivation par rapport \`a $x$. Les solutions de~\eqref{E*} sont les fonctions $y$ d\'efinies dans $\C_+$ telles que le d\'eterminant suivant soit identiquement nul
\[
 \left|
\begin{array}{ccc}
  G  & H   & y\\
 G'  & H'  & y'\\
 G'' & H'' & y''
\end{array}
 \right| = 0.
\]
En d\'eveloppant ce d\'eterminant par rapport \`a sa troisi\`eme colonne, on obtient que les c\oe fficients de l'\'equation~\eqref{E*}, qui sont d\'efinis dans le demi-plan sup\'erieur $\C_+$, s'expriment en fonction des donn\'ees $G$ et $H$ par
\[
    p(x) = -\frac{GH''-HG''}{GH'-HG'}, \qquad
    q(x) = \frac{G'H''-H'G''}{GH'-HG'}.
\]
Rappelons que le projet\'e st\'er\'eographique nord du vecteur de Gauss de l'immersion $X$ est donn\'e par $g = - G /H$, et, par~\eqref{schwarzien}, le schwarzien de $g$ est donc reli\'e aux c\oe fficients $p(x)$ et $q(x)$ par
\[
 S_x(g) = 2q(x) - \frac12 p(x)^2 - p'(x).
\]
La diff\'erentielle de Hopf~\eqref{Q} est donn\'ee par le Wronskien du syst\`eme fondamental $Y_0$
\[
 Q=i\left( GH'-HG'\right) dx^2 = i\exp\left(-\textstyle\int p\right)dx^2.
\]
On peut tout de suite observer que les fonctions $p(x)$ et $q(x)$, qui sont m\'eromorphes dans $\C_+$, ont deux types de singularit\'es :
\begin{itemize}
\item les ant\'ec\'edents $t_i$ des sommets du polygone, en lesquels $Y_0$ est singuli\`ere,
\item les ombilics de l'immersion $X$, \ie les z\'eros de sa diff\'erentielle de Hopf, en lesquels la fonction $Y_0$, et donc toute solution de l'\'equation~\eqref{E*}, est holomorphe. 
\end{itemize}
Les ombilics sont donc des singularit\'es fuchsiennes \emph{apparentes} (d\'efinition~\ref{def-sg-app}). On verra que les $t_i$ sont \'egalement des singularit\'es fuchsiennes.

On peut d\'efinir une \'equation~\eqref{E*} \`a partir de toute surface minimale qui n'est pas contenue dans un plan. Diff\'erentes immersions conformes minimales peuvent d\'efinir la m\^eme \'equation. \`A la proposition~\ref{prop-spin}, on a vu qu'une rotation de la surface repr\'esent\'ee par $X$ se traduit par une transformation lin\'eaire sur $Y_0$. Une telle transformation ne change donc pas l'\'equation~\eqref{E*}. De m\^eme, la famille associ\'ee d'immersions conformes minimales $X_\l$ ($\l\in\mathbb S^2$), qui ont pour donn\'ees de Weierstrass $\l\cdot Y_0$, et en particulier l'immersion conjugu\'ee \`a $X$, d\'efinissent la m\^eme \'equation que l'immersion $X$. Pour \'etudier l'\'equation~\eqref{E*}, on pourra donc transformer le syst\`eme $Y_0$ par toute application lin\'eaire inversible, et par exemple changer la position du rep\`ere orthonormal $(O,e_1,e_2,e_3)$ de $\R^3$.

\bigskip

Le but de ce chapitre est d'obtenir une caract\'erisation des \'equations diff\'erentielles lin\'eaires du second ordre qui proviennent, dans le sens que l'on vient de donner, d'une surface minimale \`a bord polygonal. On va voir que certaines propri\'et\'es g\'eom\'etriques de l'immersion $X$ se traduisent \'el\'egamment en terme de propri\'et\'es analytiques de l'\'equation~\eqref{E*}, comme la nature des singularit\'es et leurs exposants (proposition~\ref{prop-exp-ti} et lemme~\ref{lemme-sg-app}). On va montrer que l'\'equation~\eqref{E*} est fuchsienne r\'eelle et que sa monodromie est enti\`erement d\'etermin\'ee par la direction des c\^ot\'es du polygone $P$ (proposition~\ref{prop-mono}). Le contenu de ce chapitre \'etait connu avant que Garnier ne s'attaque au probl\`eme de Plateau. Les r\'esultats connus \`a la fin du \textsc{xix}$^{\textrm e}$ si\`ecle sont rassembl\'es par Darboux au chapitre \textsc{xiii} de~\cite{Darboux}. On y ajoute, et ceci ne figure pas non plus dans l'article de Garnier, des pr\'ecisions sur l'orientation du polygone, l'expression de la monodromie de l'\'equation et surtout la d\'emonstration de la proposition~\ref{prop-cara-eq},  qui assure la validit\'e de la m\'ethode de r\'esolution propos\'ee par Garnier. On d\'ecrit \'egalement plus pr\'ecis\'ement les ensembles de surfaces que l'on va construire, et les ensembles correspondants d'\'equations.

\section{Disques minimaux \`a bord polygonal}

On commence par introduire les espaces et les notations appropri\'es pour les disques minimaux que l'on souhaite construire, et pour leurs bords polygonaux. On va voir tout d'abord que l'on doit imposer certaines conditions naturelles sur ces polygones, ainsi que d'autres conditions qui sont peut-\^etre moins naturelles, mais dont on aura besoin dans la r\'esolution du probl\`eme de Plateau.

Soit un polygone $P$ \`a $n+3$ sommets distincts de $\R^3$ ($n\in\N^*$). On note $a_1,\ldots, a_{n+3}$ ses sommets, et pour tout $i=1,\ldots,n+3$
\[
 \ell_i = ||a_ia_{i+1}|| >0
\]
la longueur du $i$-i\`eme c\^ot\'e, et
\[
 u_i = \frac{\overrightarrow{a_ia_{i+1}}}{\ell_i}
\]
le vecteur unitaire qui dirige et oriente le $i$-i\`eme c\^ot\'e de $P$. On note \'egalement par $D_i$ la direction vectorielle orient\'ee du vecteur $u_i$. On a la condition de fermeture du polygone
\begin{equation}
 \sum_{i=1}^{n+3} \ell_i u_i =0.
\label{poly-ferme}
\end{equation}
L'ensemble des polygones non plans \`a $n+3$ sommets est param\'etr\'e par un point de $\R^3$, $n$ nombres r\'eels non nuls et $n+3$ vecteurs unitaires formant une famille g\'en\'eratrice de $\R^3$. Comme on peut extraire une base de cette famille g\'en\'eratrice, les trois longueurs manquantes seront d\'efinies de mani\`ere unique par la condition de fermeture~\eqref{poly-ferme}, mais les c\^ot\'es correspondant du polygone ne seront pas n\'ecessairement orient\'es par les vecteurs unitaires que l'on s'est donn\'es (les longueurs $\ell_i$ peuvent \^etre n\'egatives). Il ne para\^it donc pas tr\`es naturel de param\'etrer un polygone par ses directions \emph{orient\'ees}. Pourtant, la m\'ethode de Garnier permet de prescrire la direction et l'orientation des c\^ot\'es des bords polygonaux des surfaces minimales que l'on construit. En contrepartie, elle ne nous permettra pas de contr\^oler la fermeture de ces bords : on va obtenir des polygones qui ne sont pas n\'ecessairement des courbes ferm\'ees, ce sont des lignes bris\'ees \'eventuellement infinies.

\begin{defi}
 On appelle polygone \`a $n+3$ sommets de $\Rinfty$ la donn\'ee de $n+2$ points $a_1,\ldots, a_{n+2}$ de $\R^3$ et de deux directions orient\'ees $D_{n+2}$ et $D_{n+3}$.
\label{def-poly-ouvert}
\end{defi}

En quelque sorte, un polygone de $\Rinfty$ est un polygone dont le dernier sommet $a_{n+3}$ peut \^etre en l'infini. Les polygones de $\R^3$ sont les polygones de $\Rinfty$ dont le premier et le dernier c\^ot\'es sont s\'ecants, c'est-\`a-dire tels que les demi-droites affines $(a_{n+2},D_{n+2})$ et $(a_1,-D_{n+3})$ sont s\'ecantes ; le point d'intersection est le sommet suppl\'ementaire $a_{n+3}\in\R^3$. Par abus de langage, on appellera  simplement polygone tout polygone de $\Rinfty$.

On dit qu'un polygone $P$ est \emph{non d\'eg\'en\'er\'e} si aucun des produits vectoriels $u_{i-1}\times u_i$ n'est nul ($i=1,\ldots,n+3$, les indices se comprennent modulo $n+3$). On peut alors d\'efinir en chacun de ses sommets $a_i$~:
\begin{itemize}
 \item la mesure $\t_i\pi$ de l'angle ext\'erieur \`a $P$ (\ie l'angle entre les vecteurs $u_{i-1}$ et $u_i$) telle que $0<\t_i<1$ ;
 \item le vecteur unitaire normal au polygone $P$ au sommet $a_i$
    \[
     v_i = \frac{-u_{i-1}\times u_i}{||u_{i-1}\times u_i||}.
    \]
\end{itemize}
Tous les r\'esultats des chapitres~\ref{chapitre-equ-fu} et~\ref{chapitre-isomono} s'appliquent \`a l'ensemble des polygones non plans et non d\'eg\'en\'er\'es. Mais pour r\'esoudre le probl\`eme de Plateau, on sera amen\'e, au chapitre~\ref{chapitre-longueur}, \`a imposer des restrictions suppl\'ementaires sur les polygones que l'on consid\`ere. Comme on va proc\'eder par r\'ecurrence, il faut introduire une famille de polygones telle que les conditions sur les directions des c\^ot\'es se transmettent \`a des sous-ensembles de directions.

\begin{defi}
 On d\'efinit l'ensemble $\D^n$ des $(n+3)$-uplets $D=(D_1,\ldots,D_{n+3})$ de directions orient\'ees de $\R^3$ qui v\'erifient les deux conditions suivantes
\begin{itemize}
 \item deux directions quelconques $D_i$ et $D_j$ ($i\neq j$) ne sont pas colin\'eaires ;
 \item pour tout $i\neq n+1,\ n+2$, les directions $D_i$, $D_{n+1}$ et $D_{n+2}$ ne sont pas coplanaires.
\end{itemize}
On appellera un \'el\'ement de $\D^n$ un jeu de directions orient\'ees.
\label{def-Dn}
\end{defi}

Si les directions $D=(D_1,\ldots,D_{n+3})$ d'un polygone $P$ sont dans $\D^n$, alors tous ses <<~sous-polygones~>> --- obtenus en \'eliminant des c\^ot\'es de $P$ en faisant fusionner des sommets successifs --- seront non plans et non d\'eg\'en\'er\'es.

Dans la r\'esolution du probl\`eme de Plateau, on va construire des surfaces minimales, et donc des polygones, d\'efinies \`a translations et homoth\'eties de rapport positif pr\`es. Les directions orient\'ees sont invariantes par l'action de $\R^3\times\R^*_+$. On introduit donc~:

\begin{defi}
 Pour tout jeu $D\in\D^n$, on d\'efinit le quotient $\p^n_D$ de l'ensemble des polygones \`a $n+3$ sommets distincts de $\Rinfty$ dont le jeu de directions orient\'ees soit $D$ par le groupe $\R^3\times\R^*_+$ des translations et des homoth\'eties de rapport positif.
\label{def-PnD}
\end{defi}

Les ensembles $\p^n_D$ ne sont jamais vides, puisqu'il n'y a pas de condition de fermeture. Pour tout jeu $D\in\D^n$, l'ensemble $\p^n_D$ contient en particulier tous les polygones ferm\'es de $\R^3$ de directions orient\'ees $D$. Sur chaque ensemble $\p^n_D$, un syst\`eme de coordonn\'ees est donn\'e par le choix de $n$ rapports de longueur entre les $n+1$ longueurs qui sont toujours finies. On choisit le syst\`eme de coordonn\'ees d\'efini par
\begin{equation}
 (r_1,\ldots,r_n):\p^n_D\to \, ]0,+\infty[\,^n, \qquad r_i(P)=\frac{||a_ia_{i+1}||}{||a_{n+1}a_{n+2}||},
\label{coord-poly}
\end{equation}
o\`u $a_1,\ldots,a_{n+2}$ sont les sommets d'un repr\'esentant de $P\in\p^n_D$. On a alors l'isomorphisme
\[
 \p^n_D \simeq \, ]0,+\infty[\,^n.
\]

\bigskip

D\'ecrivons \`a pr\'esent l'ensemble des surfaces minimales que l'on va construire par la m\'ethode de Garnier, et dont les bords sont des \'el\'ements de $\p^n_D$. On souhaite construire des surfaces minimales qui ne se recouvrent pas elles-m\^emes aux sommets de leur bord polygonal, et qui seront donc localement plong\'ee au voisinage des sommets. En conservant les notations pr\'ec\'edentes, cela signifie qu'elles font au sommet $a_i$ ou bien un angle saillant (\ie compris entre $0$ et $\pi$) de $(1-\t_i)\pi$ ou bien un angle rentrant (\ie compris entre $\pi$ et $2\pi$) de $(1+\t_i)\pi$. Au sommet $a_{n+3}$, puisqu'on autorise un bout h\'elico\"idal, on suppose que les surfaces ont n\'ecessairement un angle saillant, de mani\`ere \`a ce qu'elles puissent <<~se refermer correctement~>> au cours de la d\'eformation isomonodromique.

Les surfaces que l'on va construire sont les \'el\'ements des ensembles suivants. Comme  on ne consid\`ere que des surfaces ayant la topologie du disque, on peut toujours supposer qu'elles sont repr\'esent\'ees sur le demi-plan sup\'erieur $\C_+$. 

\begin{defi}
Pour tout jeu $D\in\D^n$, on d\'efinit le quotient $\X^n_D$ par le groupe $\R^3\times\R^*_+$ des translations et des homoth\'eties de rapport positif de l'ensemble des immersions conformes minimales $X:\C_+\to\R^3$ telles que
\begin{itemize}
\item $X$ s'\'etend contin\^ument \`a $\overline\R = \R\cup\{\infty\}$, $X\big|_{\overline\R}$ repr\'esente un polygone $P\in\p_D^n$, et $X$ n'a pas de point de branchement au bord, except\'e peut-\^etre en les sommets de $P$,
\item $X$ a au sommet $a_i$ ($i=1,\dotsc,n+2$) un angle de $(1-\e_i\t_i)\pi$, o\`u $ \e_i=\pm 1$, et au sommet  $a_{n+3}$ un angle de $(1-\t_{n+3})\pi$,
\item si le dernier sommet, $a_{n+3}$, du polygone $P$ est en l'infini, alors $X$ a un bout h\'elico\"idal,
\item si $a_{n+3}\in\R^3$, \ie si les demi-droites issues de $a_1$ et de $a_{n+2}$ et dirig\'ees respectivement par $-D_{n+3}$ et $D_{n+2}$ sont s\'ecantes, alors la surface repr\'esent\'ee par $X$ est born\'ee dans $\R^3$. 
\end{itemize}
\label{def-XnD}
\end{defi}

On continue \`a appeler immersions les \'el\'ements des ensembles $\X^n_D$, m\^eme s'il s'agit de classes d'\'equivalence d'immersions. Soit $X$ une immersion de $\X^n_D$. On note $P\in\p^n_D$ son bord polygonal, et
\[
 Y_0=(G,H):\C_+\to \C^2\ssm\{(0,0)\}
\]
ses donn\'ees de Weierstrass. La fonction $Y_0$ est holomorphe dans le demi-plan sup\'erieur $\C_+$ et l'immersion $X$ est donn\'ee par
\[
 X(x) = \Re \int_{x_0}^x
\begin{pmatrix}
 i\left(G(\xi)^2 - H(\xi)^2\right) \\
 G(\xi)^2 + H(\xi)^2\\
 2iG(\xi)H(\xi)
\end{pmatrix}
d\xi
\]
o\`u $x_0$ est un point arbitraire de $\C_+$ (puisque $X$ est d\'efinie \`a translation pr\`es). On d\'efinit les points
\[
    t_1<\cdots<t_{n+3}
\]
de $\overline\R$ qui sont les ant\'ec\'edents des sommets de $P$ par l'immersion $X$. Quitte \`a composer $X$ par une homographie, on peut toujours supposer
\[
    t_{n+1}=0, \quad t_{n+2}=1, \quad t_{n+3}=\infty.
\]
D'apr\`es la premi\`ere des conditions de la d\'efinition pr\'ec\'edente, la fonction $Y_0$ est continue et non nulle sur chacun des intervalles $]t_i,t_{i+1}[$. Cette hypoth\`ese est naturelle si l'on veut pouvoir prolonger la surface \`a travers chacun des c\^ot\'es du polygone $P$, et appliquer le principe de r\'eflexion de Schwarz. Sous cette hypoth\`ese, l'application de Gauss $N(x)$ de l'immersion $X$ admet une limite en chaque sommet de $P$, qui est orthogonale aux c\^ot\'es adjacents au sommet. On note $N(t_i)$ le vecteur de Gauss limite en $x=t_i$, il v\'erifie $N(t_i)=\pm v_i$. On verra \`a la section~\ref{section-exp} que la deuxi\`eme des conditions implique que l'immersion $X$ a un point de branchement au bord en un sommet $a_i$ si et seulement si elle a un angle rentrant, \ie si $\e_i=-1$. L'ordre du point de branchement est alors $1$.

\section{Monodromie et propri\'et\'es de r\'ealit\'e}
\label{section-mono}

On note $S(t)$ l'ensemble des singularit\'es de l'immersion $X$
\[
 S(t):=\{t_1,\ldots, t_{n+3}\} \subset \overline\R
\]
o\`u $\overline\R=\R\cup\{\infty\}$. On va voir que l'\'equation~\eqref{E*} est bien d\'efinie dans la sph\`ere de Riemann, tandis que les donn\'ees de Weierstrass $G(x)$ et $H(x)$ ont des points de branchement en les points $x=t_i$, et sont donc holomorphes dans le rev\^etement universel de l'ensemble $ \P\ssm S(t)$. On va d\'eterminer, par des consid\'erations g\'eom\'etriques, le comportement et la monodromie des fonctions $G(x)$ et $H(x)$ en ces singularit\'es. On va voir que ceux-ci sont reli\'es aux propri\'et\'es de r\'ealit\'e de l'immersion $X$.

\subsection{Propri\'et\'es de r\'ealit\'e}

La proposition suivante est une cons\'equence directe du lemme~\ref{lemme-realite}. Elle assure en particulier que les points $x=t_i$ ne sont donc pas des points de branchement pour les fonctions $p(x)$ et $q(x)$.

\begin{prop}
 Les c\oe fficients $p(x)$ et $q(x)$ de l'\'equation~\eqref{E*} sont \`a valeurs r\'eelles dans $\overline\R\ssm S(t)$ et s'\'etendent en des fonctions m\'eromorphes et uniformes dans $\P\ssm S(t)$.
\label{prop-p-q-reels}
\end{prop}

\begin{proof}
Pour montrer que les c\oe fficients  $p(x)$ et $q(x)$ sont r\'eels sur l'axe r\'eel, il suffit de trouver pour tout $i=1,\dotsc,n+3$ un syst\`eme fondamental de solutions $\left(G_i(x),H_i(x)\right)$ dont les composantes soient toutes les deux r\'eelles ou toutes les deux purement imaginaires sur l'intervalle $]t_i,t_{i+1}[$. Par le
lemme~\ref{lemme-realite}, on sait qu'il existe une matrice $S_i\in SU(2)$ telle que le syst\`eme fondamental
$\left(G_i(x),H_i(x)\right) = Y_0(x)\cdot S_i$ convienne. On peut choisir la matrice $S_i$ telle que le syst\`eme $Y_0(x)\cdot S_i$ soit r\'eel sur $]t_i,t_{i+1}[$. La matrice $S_i$ est un relev\'e d'une rotation envoyant le vecteur $u_i$ sur le deuxi\`eme vecteur de base $e_2=(0,1,0)$, ou sur son oppos\'e $(0,-1,0)$.

On peut donc prolonger les fonctions $p(x)$ et $q(x)$ au demi-plan inf\'erieur $\C_-=\{x\in \C \ | \ \Im (x)<0\}$ en posant pour tout $x\in\C_-$
\[
 p(x) = \overline{p(\bar{x})},\qquad
 q(x) = \overline{q(\bar{x})},
\]
et on obtient ainsi qu'elles sont m\'eromorphes dans $\P\ssm S(t)$.
\end{proof}

Comme les propri\'et\'es de r\'ealit\'e jouent un r\^ole essentiel dans l'\'etude de l'\'equation~\eqref{E*}, on introduit l'application suivante $\tau$ d\'efinie sur le faisceau des fonctions m\'eromorphes $\mathcal M_{\P}$, qui \`a une fonction m\'eromorphe dans un ouvert $\Omega$ associe sa <<~conjugu\'ee~>> d\'efinie dans $\bar{\Omega}$, dans le sens suivant :
\begin{equation}
\begin{split}
 \tau :\
  \mathcal M_{\P}\left( \Omega\right)  & \to \mathcal M_{\P}\left( \bar{\Omega}\right) \\
  f                        & \mapsto \tau(f)= (x\mapsto \overline{f(\bar{x})}).
\end{split}
\label{def-tau}
\end{equation}
L'application $\tau$ est anti-lin\'eaire. Si $\Omega$ est un domaine de $\P$ stable par conjugaison (\ie sym\'etrique par rapport \`a l'axe r\'eel), alors pour toute fonction $f$ m\'eromorphe dans $\Omega$, on a
\begin{align*}
 \tau(f)=f  &\Leftrightarrow f\left(\Omega\cap\R\right) \subset \overline\R\\
 \tau(f)=-f &\Leftrightarrow f\left(\Omega\cap\R\right) \subset i\overline\R.
\end{align*}

La fonction holomorphe $\tau(Y_0)=\left(\tau(G),\tau(H)\right):\C_- \to\C^2$ constitue \'egalement les donn\'ees de Weierstrass d'une immersion conforme minimale $X^-:\C_-\to\R^3$. Un calcul rapide montre que cette immersion repr\'esente la surface minimale sym\'etrique de $X\left( \C_+\right)$ par rapport au second axe de coordonn\'ees $(O,e_2)$. Comme la matrice
\[
J=\begin{pmatrix}
   0 & -1\\
   1 & 0
  \end{pmatrix}
\]
est un relev\'e du demi-tour par rapport au second axe de coordonn\'ees, on obtient~:

\begin{lemm}
Soit une fonction holomorphe $Y:\C_+ \to\C^2\ssm\{(0,0)\}$. Alors, les deux fonctions
\[
 Y:\C_+ \to\C^2 \quad \text{ et } \quad \tau(Y)\cdot J:\C_- \to\C^2
\]
sont les donn\'ees de Weierstrass de la m\^eme surface minimale.
\label{lemme-sym-e2}
\end{lemm}

On obtient ainsi le principe de r\'eflexion de Schwarz. En effet, de m\^eme que les c\oe fficients $p(x)$ et $q(x)$, pour tout $i=1,\ldots,n+3$, le syst\`eme fondamental $\left(G_i(x),H_i(x)\right)$ introduit \`a la d\'emonstration de la proposition~\ref{prop-p-q-reels} se prolonge analytiquement au demi-plan inf\'erieur $\C_-$ \`a travers l'intervalle $]t_i,t_{i+1}[$ en posant pour tout $x\in\C_-$
\[
 \left(G_i,H_i\right)(x) = \tau\left((G_i,H_i)\right)(x).
\]
Le syst\`eme $\left(G_i,H_i\right)$ est alors holomorphe dans l'ouvert simplement connexe $U_i$
\[
 U_i = \C_+\cup\C_-\cup]t_i,t_{i+1}[.
\]
L'immersion de donn\'ees de Weierstrass $\left(G_i,H_i\right)$ se
prolonge donc \'egalement en une immersion d\'efinie dans l'ouvert
$U_i$, et le lemme~\ref{lemme-sym-e2} nous dit qu'elle d\'efinie une
surface minimale sym\'etrique par rapport au second axe de
coordonn\'ees $(O,e_2)$. Comme on a
\[
    \begin{pmatrix}
        G & H
    \end{pmatrix}
    =
    \begin{pmatrix}
        G_i & H_i
    \end{pmatrix}
    S_i^{-1},
\]
on obtient ainsi $n+3$ prolongements $Y_i(x)$ du syst\`eme $Y_0(x)$ \`a travers chacun des intervalles $]t_i,t_{i+1}[$~:
\begin{equation}
  Y_i:U_i\to\C^2, \qquad Y_i\big|_{\C_+} = Y_0.
\label{def-Yi}
\end{equation}
Chacun de ces prolongements induit un prolongement $X^i : U_i\to\R^3$ de l'immersion $X$, qui repr\'esente dans $\C_-$ la surface minimale sym\'etrique de $X\left( \C_+\right)$ par rapport au $i$-i\`eme c\^ot\'e du polygone $P$. De plus, les points sym\'etriques sur la surface minimale ont des ant\'ec\'edents par l'immersion $X^i$ qui sont conjugu\'es. Ceci nous permet de d\'eterminer la monodromie de l'\'equation~\eqref{E*}.

\subsection{Monodromie}

\begin{figure}
\centering
\begin{pspicture}(0,0)(12.5,8)
\rput(7,4){\includegraphics[height=7cm]{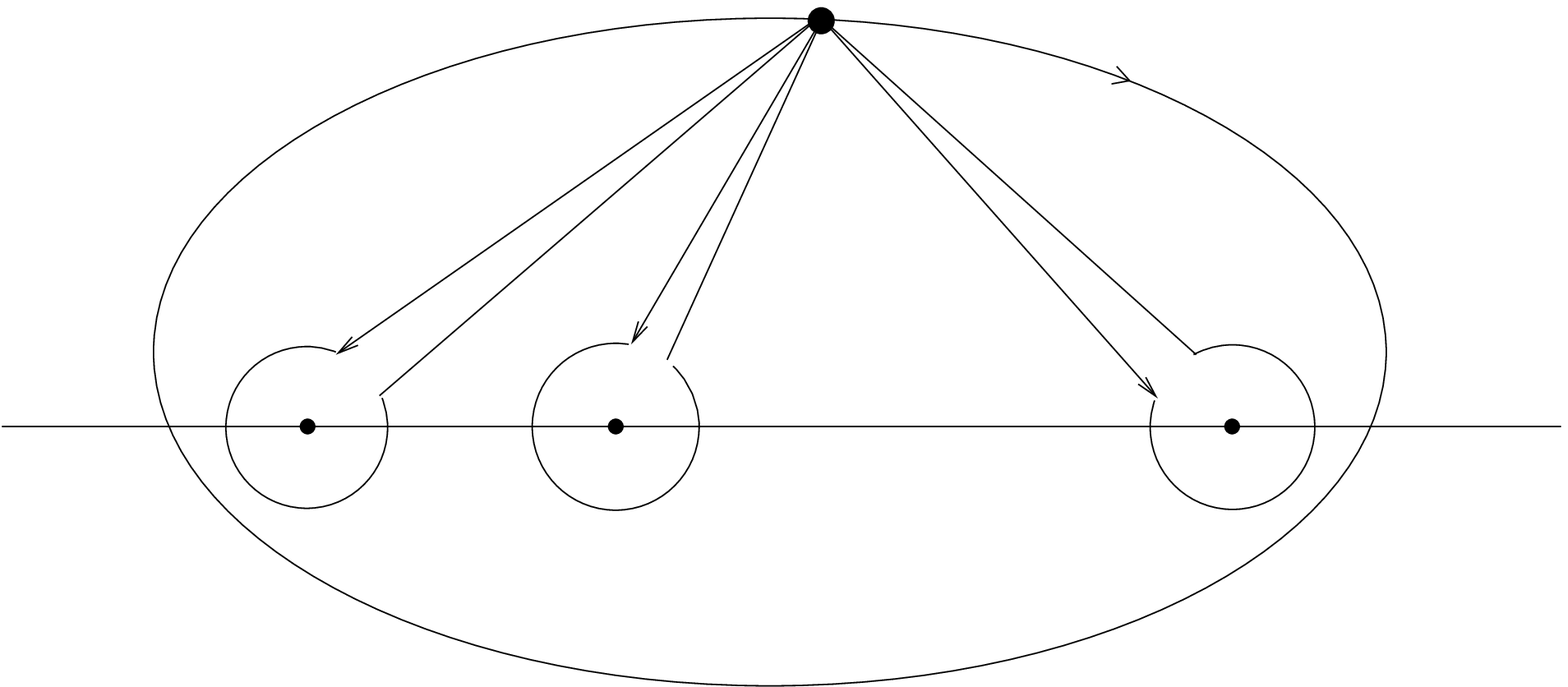}}
\rput(7.4,7.7){$x_0$}
\rput(12.1,6.4){$\ga_{n+3}$}
\rput(2.15,3.5){$t_1$}
\rput(5.3,3.5){$t_2$}
\rput(11.6,3.5){$t_{n+2}$}
\rput(8.5,3.5){$\ldots$}
\rput(3.6,5.2){$\ga_1$}
\rput(5.4,4.8){$\ga_2$}
\end{pspicture}
\caption{Les lacets $\ga_i$}
\label{fig-ga-i}
\end{figure}

L'\'etude pr\'ec\'edente des propri\'et\'es de r\'ealit\'e de l'immersion $X$ et de l'\'equation~\eqref{E*} nous permet de d\'eterminer comment le syst\`eme fondamental $Y_0(x)$ est transform\'e autour de chaque singularit\'e $x=t_i$, c'est-\`a-dire de d\'eterminer un syst\`eme de g\'en\'erateurs de la monodromie de l'\'equation~\eqref{E*}. On fixe un point $x_0$ dans le demi-plan sup\'erieur $\C_+$. Le groupe fondamental $\pi_1\left(\P\ssm S(t),x_0\right)$ est engendr\'e par les classes de lacets $\ga_1,\ldots,\ga_{n+3}$ bas\'es en $x_0$, qui sont repr\'esent\'es \`a la figure~\ref{fig-ga-i}. On note $M_1,\ldots,M_{n+3}$ les matrices de monodromie du syst\`eme fondamental de solutions $Y_0(x)$ le long des classes de lacets $\ga_i$~:
\begin{equation}
 M_i := M_{\ga_i}(Y_0)
 \label{def-Mi}.
\end{equation}
Ces matrices constituent un syst\`eme de g\'en\'erateurs de la monodromie de l'\'equation~\eqref{E*}.

\begin{prop}
Les matrices de monodromie $M_i$ ($i=1,\dotsc,n+3$) du syst\`eme fondamental de solutions $Y_0(x)$ le long des lacets $\ga_i$ s'\'ecrivent
\begin{equation}
    M_i = D_iD_{i-1}^{-1},
\label{def1-Mi}
\end{equation}
o\`u pour tout $i=1,\dotsc,n+3$, la matrice $D_i$ est un relev\'e dans $SU(2)$ du demi-tour vectoriel d'axe $u_i$. \label{prop-mono}
\end{prop}

Par cette proposition, on obtient que la monodromie de l'\'equation~\eqref{E*} est d\'etermin\'ee par les directions des c\^ot\'es du polygone $P$. L'expression des matrices $M_i$ sous la forme de produit de demi-tours successifs n'est donn\'ee ni par Darboux, ni par Garnier. Cette expression sera pourtant essentielle pour \'etablir que les d\'eformations isomonodromiques que l'on va construire d\'efinissent bien des solutions du probl\`eme de Plateau (par la proposition~\ref{prop-sys-reel-2}), fait qui n'est jamais justifi\'e par Garnier.

\begin{proof}
On note $\ga_i\ast Y_0(x)$ le prolongement du syst\`eme fondamental $Y_0(x)$ le long du lacet $\ga_i$. Ce prolongement est \'egalement holomorphe dans $\C_+$, et c'est encore un syst\`eme fondamental de solutions de l'\'equation~\eqref{E*}, \'etant donn\'e que les fonctions $p(x)$ et $q(x)$ sont uniformes dans $\P\ssm S(t)$. La matrice $M_i$ est l'unique matrice inversible qui satisfait
\[
 \ga_i\ast Y_0(x) = Y_0(x)M_i.
\]
Le syst\`eme fondamental $\ga_i\ast Y_0(x)$ constitue les donn\'ees de Weierstrass d'une immersion conforme minimale. Pour d\'eterminer la matrice de monodromie $M_i$, on compare cette immersion \`a l'immersion $X$. Lorsqu'on suit le lacet $\ga_i$, on croise d'abord l'axe r\'eel entre $t_{i-1}$ et $t_i$ et l'immersion $X$ se prolonge donc en d\'efinissant la surface minimale sym\'etrique de $X\left( \C_+\right)$ par rapport au $(i-1)$-i\`eme c\^ot\'e de $P$ ; puis on croise l'axe r\'eel entre $t_i$ et $t_{i+1}$ et on fait un nouveau demi-tour par rapport au $i$-i\`eme c\^ot\'e de la surface obtenue \`a l'\'etape pr\'ec\'edente. L'immersion de donn\'ees de Weierstrass $\ga_i\ast Y_0(x)$ est donc l'image de l'immersion $X$ par le produit de ces deux demi-tours, c'est-\`a-dire par la rotation d'axe $v_i$ et d'angle $2\pi\t_i$. On en d\'eduit que la matrice $M_i$ est un des deux relev\'es de cette rotation. Ceci constitue le r\'esultat qu'obtiennent Darboux et Garnier.

On veut pouvoir comparer les relev\'es des demi-tours intervenant dans des matrices de monodromie successives $M_i$ et $M_{i+1}$, c'est-\`a-dire, en fait, associer un unique relev\'e au demi-tour autour du $i$-i\`eme c\^ot\'e de $P$. On vient de voir que l'immersion $X^i : \C_-\to\R^3$, de donn\'ees de Weierstrass $Y_i:\C_-\to\C^2$ d\'efinies par~\eqref{def-Yi}, repr\'esente la surface minimale sym\'etrique de la surface initiale par rapport au $i$-i\`eme c\^ot\'e de $P$. D'apr\`es le lemme~\ref{lemme-sym-e2}, il existe donc un relev\'e $D_i\in SU(2)$ du demi-tour autour de ce c\^ot\'e tel que pour tout $x\in\C_+$ on ait
\[
  Y_0(x) \cdot D_i = \tau\left( Y_i\right) (x)\cdot J,
\]
ce qui s'\'ecrit
\[
 Y_0(x) = -\tau\left( Y_i\cdot J\cdot\overline D_i\right)(x),
\]
vu que les matrices $A\in SU(2)$ qui sont des relev\'es de demi-tours sont caract\'eris\'ees par l'\'equation $A^2=-\I_2$. En \'ecrivant la relation pr\'ec\'edente pour les syst\`emes $Y_{i-1}(x)$ et $Y_i(x)$, on trouve que pour tout $x\in\C_-$ on a
\[
  Y_{i-1}(x)\cdot J\cdot\overline D_{i-1} = Y_i(x)\cdot J\cdot\overline D_i,
\]
ce qui, par l'identit\'e~\eqref{SU2}, donne
\[
 Y_{i-1}(x) = Y_i(x)\cdot D_i\cdot D_{i-1}^{-1}.
\]
Or la matrice $M_i$ est l'unique matrice qui v\'erifie
\[
 Y_{i-1}(x) = Y_i(x)\cdot M_i,
\]
ce qui donne le r\'esultat annonc\'e.
\end{proof}

\subsection{Exposants en les sommets du polygone}
\label{section-exp}

Pour l'instant, la monodromie de l'\'equation~\eqref{E*} n'est pas enti\`erement d\'etermin\'ee \`a partir du bord polygonal de l'immersion $X$, puisqu'elle d\'epend du choix des relev\'es de chaque demi-tour $D_i$. L'\'etude locale  de l'immersion $X$ au voisinage des singularit\'es $x=t_i$ va nous permettre de lever cette ind\'etermination. Ceci nous permet \'egalement de calculer pr\'ecis\'ement les exposants de l'\'equation~\eqref{E*}, qui ne sont donn\'es par la monodromie qu'\`a un entier pr\`es.

Rappelons que l'immersion $X$ fait au sommet $a_i$ ($i=1,\dotsc,n+2$) un angle de $(1-\e_i\t_i)\pi$, o\`u $ \e_i=\pm 1$, et au sommet $a_{n+3}$, un angle de $(1-\t_{n+3})\pi$, que la surface ait un bout en $a_{n+3}$ ou qu'elle soit born\'ee.

\begin{prop}
Les points $x= t_1,\ldots,t_{n+3}$ sont des singularit\'es fuchsiennes et non logarithmiques de l'\'equation~\eqref{E*}. Pour tout $i=1,\dotsc,n+2$, les exposants en $x=t_i$ sont de la forme
    \[
        -\e_i\frac{\t_i}{2}, \quad r_i+\e_i\frac{\t_i}{2} \quad (r_i\in\N).
    \]
De plus, si $\e_i=-1$, alors $r_i\geq1$. Les exposants au point $x=\infty$ sont de la forme
    \[
        1-\frac{\t_{n+3}}{2}, \quad r_{n+3}-1+\frac{\t_{n+3}}{2} \quad ( r_{n+3}\in\N^*).
    \]
De plus, la surface a un bout h\'elico\"idal en $x=\infty$ si et seulement si $r_{n+3}=1$.
\label{prop-exp-ti}
\end{prop}

\begin{proof}
Montrons tout d'abord que la singularit\'e $x=t_i$ est fuchsienne. Comme on s'int\'eresse \`a pr\'esent \`a des propri\'et\'es locales de l'\'equation~\eqref{E*}, on peut choisir la position du rep\`ere orthonormal de $\R^3$ tel que le vecteur normal $v_i$ co\"incide avec le troisi\`eme vecteur de base $e_3=(0,0,1)$. On note toujours $Y_0=(G,H)$ les donn\'ees de Weierstrass correspondant \`a cette position, et $X$ l'immersion associ\'ee. La matrice de monodromie $M_i$ du syst\`eme $Y_0(x)$  est alors un relev\'e de la rotation d'axe $(O,e_3)$ et d'angle $2\pi\t_i$ et elle s'\'ecrit donc
\begin{equation}
 M_i =\d_i
\begin{pmatrix}
 e^{i\t_i\pi} & 0\\
 0        & e^{-i\t_i\pi}
\end{pmatrix}
, \quad \text{ avec } \d_i=+1 \text{ ou }-1. 
\label{Mi-diago}
\end{equation}
Les fonctions $G(x)$ et $H(x)$ sont donc de la forme
\begin{align*}
    G(x) & = (x-t_i)^{\frac{1-\d_i}{4}+\frac{\t_i}{2}} \;\varphi(x)\\
    H(x) & = (x-t_i)^{-\frac{1-\d_i}{4}-\frac{\t_i}{2}}  \;\psi(x)
\end{align*}
o\`u les fonctions $\varphi(x)$ et $\psi(x)$ sont uniformes au voisinage de $x=t_i$. Comme les primitives
\[
    \int_{x_0}^x G(\xi)^2d\xi,\quad \int_{x_0}^x H(\xi)^2d\xi,\quad \int_{x_0}^x G(\xi)H(\xi)d\xi
\]
qui interviennent dans l'expression de l'immersion $X$ prennent des valeurs finies en $x=t_i$, les fonctions $\varphi(x)$ et $\psi(x)$ n'ont pas de singularit\'e essentielle en $x=t_i$, et sont donc m\'eromorphes en ce point. On en conclut donc que la singularit\'e $x=t_i$ est fuchsienne. Comme la matrice $M_i$ est diagonalisable, cette singularit\'e n'est pas logarithmique. Chacun des deux exposants est d\'etermin\'e \`a un entier pr\`es, et leur somme est un entier relatif.

Pour \^etre plus pr\'ecis sur la valeur des exposants, il faut \'etudier le comportement du syst\`eme fondamental $Y_0(x)$ en utilisant l'expression de l'immersion $X$ aux sommets du polygone. Soient $s^i_1$ et $s^i_2$ les exposants en $x=t_i$, $s^i_1< s^i_2$. Leur somme $r_i:=s_1^i+s_2^i$ est un entier relatif. Supposons tout d'abord $i\neq n+3$. On a vu que le fait de supposer $v_i=e_3$ implique que le syst\`eme fondamental $Y_0(x)$ est canonique en $x=t_i$. Comme la projection st\'er\'eographique nord de $N(x)$ est $-G(x)/H(x)$, on voit que si le vecteur de Gauss $N(t_i)$ est \'egal \`a $e_3$, alors la fonction $H(x)$ est canonique pour l'exposant le plus grand $s^i_2$, et si $N(t_i)$ est \'egal \`a $-e_3$, alors $G(x)$ est canonique pour $s^i_2$. Supposons par exemple que $N(t_i)=e_3$. On a alors en $x=t_i$ les \'equivalents
\[
    G(x) \sim a(x-t_i)^{s^i_1} , \qquad
    H(x) \sim b(x-t_i)^{r_i-s^i_1},
\]
o\`u les constantes $a$ et $b$ sont non nulles. \`A une rotation d'axe $(O,e_3)$ pr\`es, on peut supposer ces constantes r\'eelles. On en d\'eduit, si $r_i\neq-1$,
\begin{equation}
    X(x)-X(t_i) \sim \Re
        \begin{pmatrix}
            \frac{ia^2}{2\alpha_i+1}        (x-t_i)^{2s_1^i+1}\\
            \frac{a^2}{2\alpha_i+1}         (x-t_i)^{2s_1^i+1}\\
            \frac{2iab}{r_i+1}              (x-t_i)^{r_i+1}\\
        \end{pmatrix}.
\label{X-en-ti}
\end{equation}
Mais on ne peut pas avoir $r_i=-1$, car alors l'immersion $X$ serait asymptote \`a une h\'elico\"ide en $x=t_i$ ; on ne peut pas non plus avoir $r_i<-1$, vu que l'immersion $X$ est \`a valeurs finies en $x=t_i$. Lorsque, dans l'\'equivalent pr\'ec\'edent, la quantit\'e $x-t_i$ prend des valeurs r\'eelles infiniment petites, positives puis n\'egatives, on voit que la quantit\'e $(2s_1^i+1)\pi$ est l'angle que fait la surface minimale au sommet $a_i$, c'est-\`a-dire $ 2s_1^i+1 = 1-\e_i\t_i$. On obtient donc
\[
s_1^i = -\e_i\frac{\t_i}{2}, \qquad s_2^i = r_i+\e_i\frac{\t_i}{2},
\]
et lorsque $\e_i=-1$,  l'in\'egalit\'e $s_1^i<s_2^i$ donne la minoration $r_i\geq1$.

Pour d\'eterminer les exposants au point $x=\infty$, on fait le changement de variables $w=1/x$ dans l'immersion $X$
\[
    X\left(\frac{1}{w}\right)=- \Re \int
        \begin{pmatrix}
            i \left(G^2\left(\frac1w\right) - H^2\left(\frac1w\right)\right)\\
            G^2\left(\frac1w\right)         + H^2\left(\frac1w\right)\\
            2i G\left(\frac1w\right) H\left(\frac1w\right)
        \end{pmatrix}
    \frac{dw}{w^2}.
\]
On proc\`ede comme pr\'ec\'edemment en supposant qu'en $w=0$ on a les \'equivalents
\[
    G\left(\tfrac1w\right)  \sim aw^{s^{n+3}_1} , \qquad
    H\left(\tfrac1w\right)  \sim bw^{r_{n+3}-s^{n+3}_1},
\]
avec $a, b\in\R^*$. Si l'entier $r_{n+3}$ est n\'egatif ou nul, alors la surface n'est pas born\'ee au voisinage de $w=0$, et elle n'a pas de bout h\'elico\"idal : ce cas est exclu. Si $r_{n+3}=1$, la surface a un bout h\'elico\"idal. Si $r_{n+3}\geq2$, la situation est la m\^eme que pr\'ec\'edemment. On obtient donc $r_{n+3}\geq1$ et
\[
 2s_1^{n+3}-1 = 1-\t_{n+3}
\]
et on conclut de m\^eme.
\end{proof}

\begin{figure}[!h]
 \begin{minipage}[b]{.47\linewidth}
  \centering\epsfig{figure=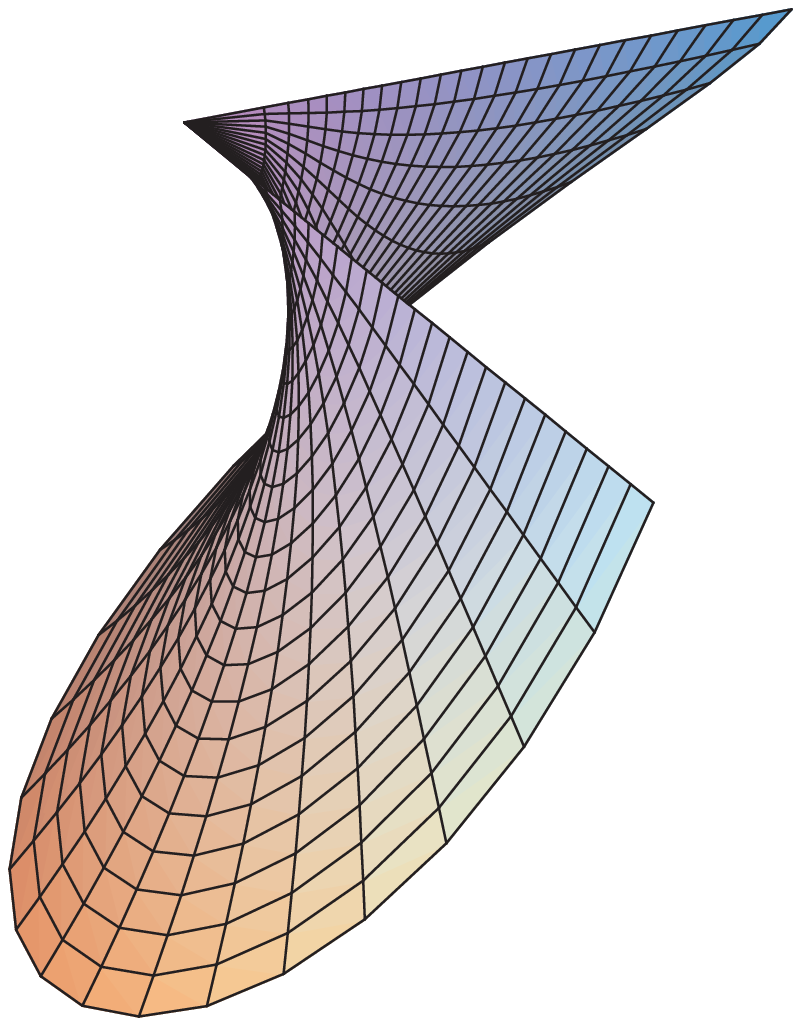,width=\linewidth, height=6cm}
\caption{Situation g\'en\'erique : angle saillant sans ombilic ($\t=2/3,\ \e=1,\ r=0$)}
\label{fig-sommet-aigu-generique}
 \end{minipage} \hfill
 \begin{minipage}[b]{.47\linewidth}
  \centering\epsfig{figure=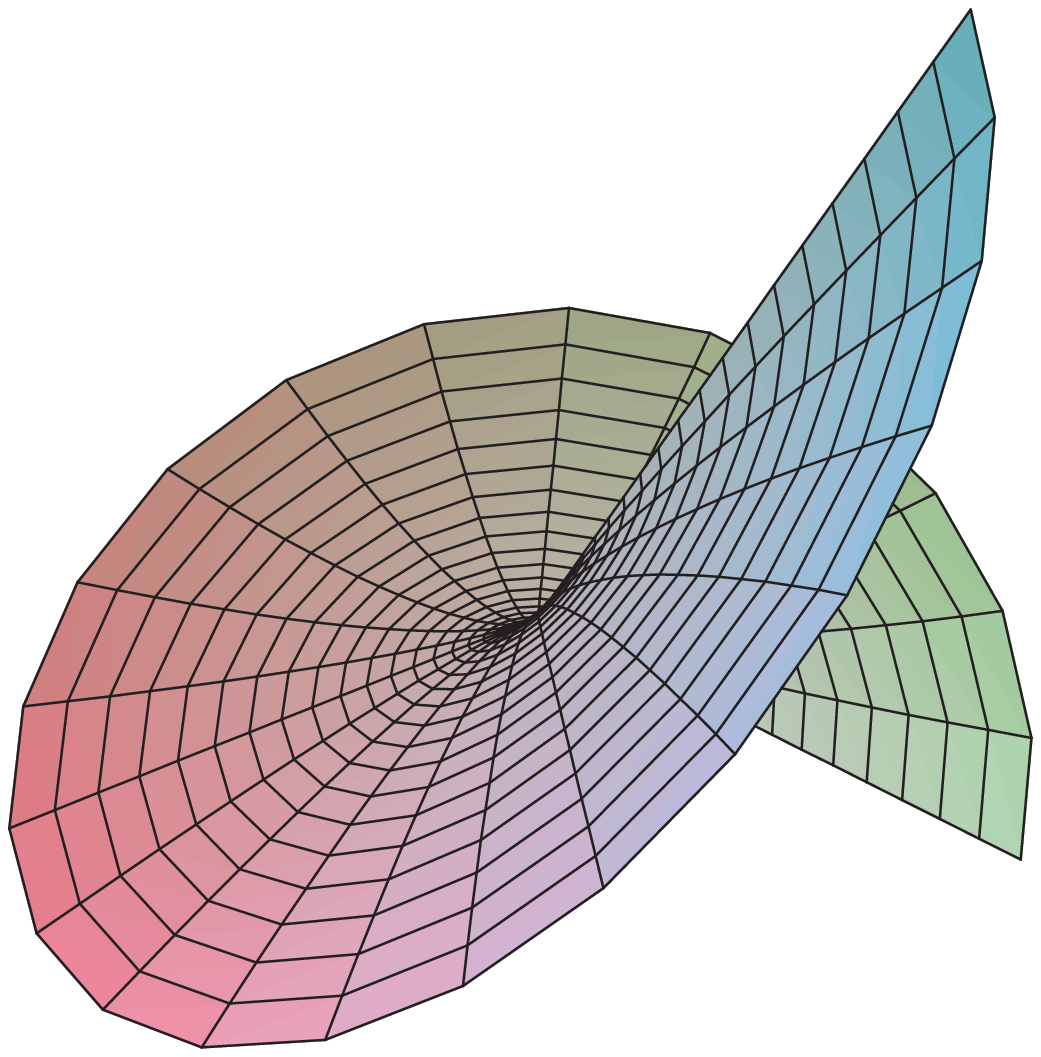,width=\linewidth}
\caption{Angle rentrant sans ombilic ($\t=2/3,\ \e=-1,\ r=1$)}
\label{fig-sommet-obtus-generique}
 \end{minipage}
\end{figure}

\begin{figure}[!h]
 \begin{minipage}[b]{.47\linewidth}
  \centering\epsfig{figure=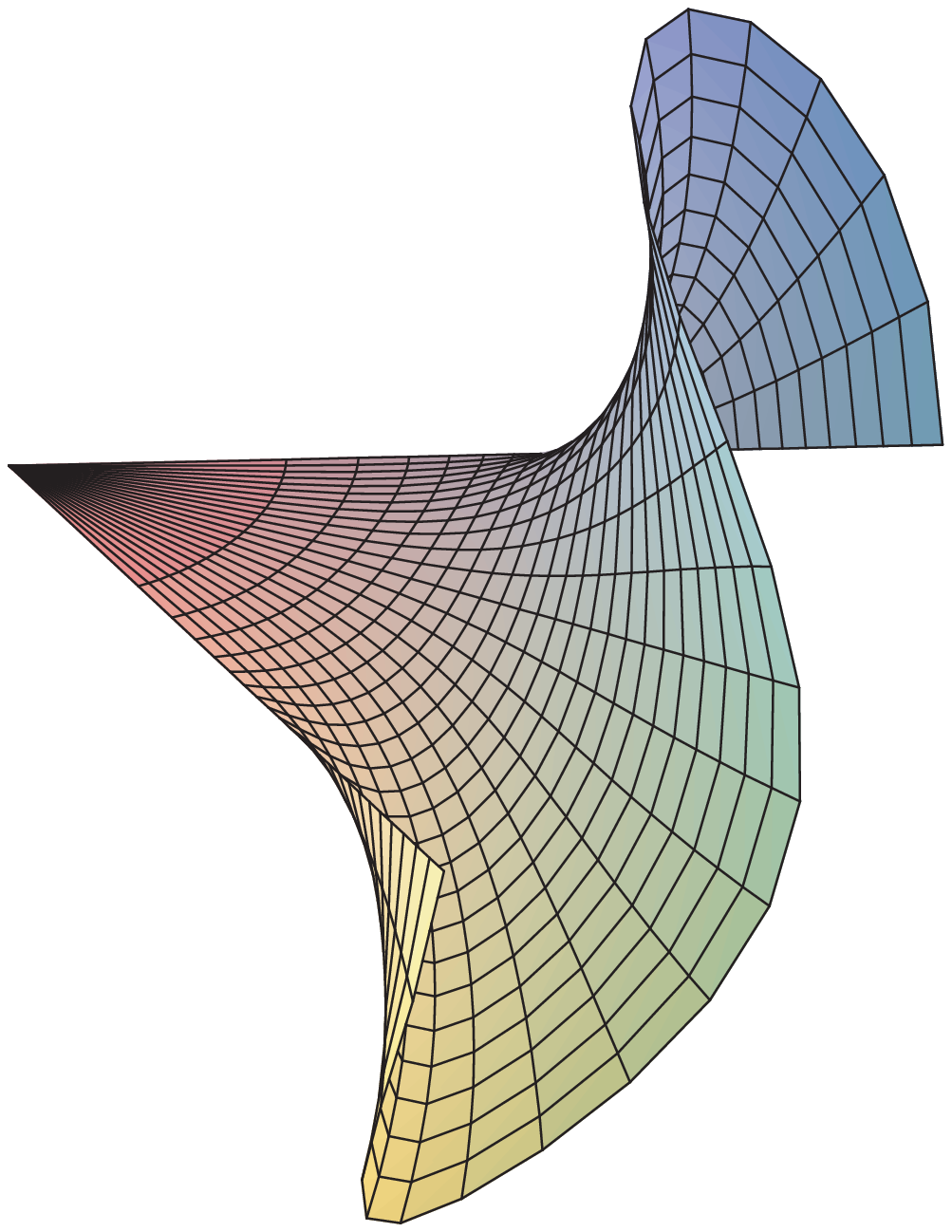,width=\linewidth}
\caption{Angle saillant avec un ombilic ($\t=2/3,\ \e=1,\ r=1$)}
\label{fig-sommet-aigu+ombilic}
 \end{minipage} \hfill
 \begin{minipage}[b]{.47\linewidth}
  \centering\epsfig{figure=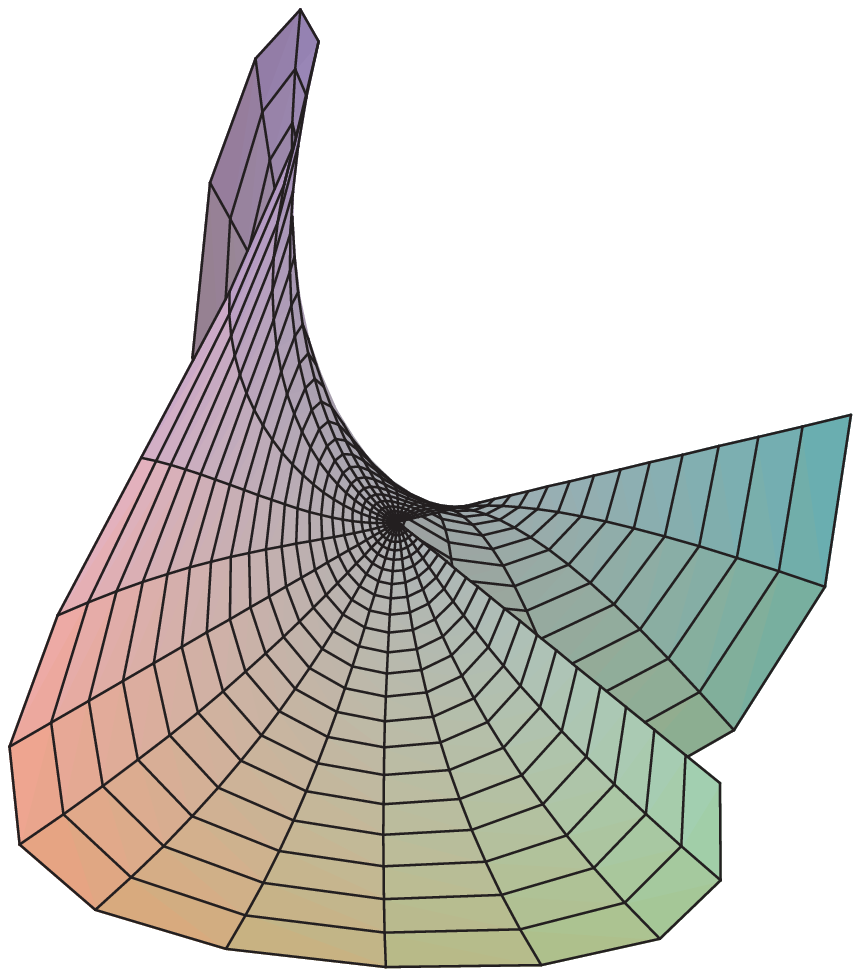,width=\linewidth}
\caption{Angle rentrant avec un ombilic ($\t=2/3,\ \e=-1,\ r=2$)}
\label{fig-sommet-obtus+ombilic}
 \end{minipage}
\end{figure}

On repr\'esente aux figures~\ref{fig-sommet-aigu-generique},~\ref{fig-sommet-obtus-generique},~\ref{fig-sommet-aigu+ombilic} et~\ref{fig-sommet-obtus+ombilic} les diff\'erentes configurations locales possibles pour une surface minimale en un sommet d'un polygone. On a choisi un angle int\'erieur de $\pi/3$ au sommet consid\'er\'e, \ie $\t=2/3$. La surface fait donc un angle saillant de $\pi/3$ ($\e=1$) ou un angle rentrant de $5\pi/3$ ($\e=-1$). Lorsque $\e=1$, les exposants de l'\'equation sont
\[
 -\frac 13, \ r+\frac 13  \qquad (r\geq0).
\]
Lorsque $\e=-1$, les exposants sont
\[
 \frac 13, \ r-\frac 13  \qquad (r\geq1).
\]
Les figures~\ref{fig-sommet-aigu-generique} et~\ref{fig-sommet-obtus-generique} correspondent aux valeurs <<~minimales~>> de l'entier $r$ ($r=0$ lorsque $\e=1$, et $r=1$ lorsque $\e=-1$). Comme on le voit sur les figures~\ref{fig-sommet-aigu+ombilic} et~\ref{fig-sommet-obtus+ombilic}, lorsque l'entier $r$ est sup\'erieur \`a ces valeurs, on peut consid\'erer que le sommet est, en un sens, \'egalement un ombilic.

On dira que la situation en un sommet $a_i\in\R^3$ est \emph{g\'en\'erique} lorsque $\e_i=1$ et $r_i=0$, c'est-\`a-dire lorsque les exposants sont oppos\'es : $-\frac{\t_i}{2}$ et $\frac{\t_i}{2}$  (figure~\ref{fig-sommet-aigu-generique}). On peut alors en d\'eduire les autres configurations possibles en ajoutant un entier naturel \`a l'un des exposants : c'est effectivement ce qui se produira au cours de la d\'eformation isomonodromique. Au sommet $a_{n+3}$, on dira que la situation est g\'en\'erique lorsque $r_{n+3}=1$. En particulier, ceci signifie que g\'en\'eriquement, on a un bout h\'elico\"idal en $a_{n+3}$.


\begin{rema}
Un \emph{point de branchement au bord} de l'immersion $X$ est un point $x_0\in\overline\R$ tel que la norme de $\partial X/\partial x$ tende vers $0$ lorsque $x\in\C_+$ tend vers $x_0$. Par d\'efinition de l'ensemble $\X^n_D$, les seuls points de branchement au bord possibles sont les sommets de $P$. D'apr\`es~\eqref{X-en-ti}, comme le plus petit des exposants en $x=t_i$ est $s^i_1=\e_i\t_i/2$, on a en $x=t_i$
\[
\left\Vert \frac{\partial X}{\partial x}\right\Vert \sim a|x-t_i|^{\e_i\t_i} \qquad (a>0).
\]
Le point  $x=t_i$ est donc un point de branchement si et seulement si $\e_i=-1$. Les uniques points de branchement de l'immersion $X$ sont donc les sommets en lesquels elle a un angle rentrant, et l'ordre de ces points de branchement est $1$ (car $\t_i<1$).
\end{rema}

\begin{rema}
La valeur des exposants donn\'ee \`a la proposition~\ref{prop-exp-ti} implique que les valeurs propres de la matrice $M_i$ sont $\exp(\pm i\pi\t_i)$, c'est-\`a-dire que le signe $\d_i$ intervenant dans sa diagonalis\'ee~\eqref{Mi-diago} est $+1$ (sauf lorsque $\t_i=1/2$, les cas $\d_i=+1$ et $-1$ \'etant alors \'equivalents). On d\'etermine ainsi enti\`erement la matrice $M_i$ \`a partir du polygone $P$, puisqu'on a lev\'e la derni\`ere ind\'etermination, \`a savoir le choix du relev\'e de la rotation d'axe $v_i=(v_i^1,v_i^2,v_i^3)$ et d'angle $2\t_i\pi$ : par~\eqref{def-R_A}, les matrices $M_i$ valent donc
\[
 M_i =\cos(\t_i\pi)\I_2 - i\sin(\t_i\pi)
 \begin{pmatrix}
  -v_i^3     & v_i^1-iv_i^2\\
  v_i^1+iv_i^2 & v_i^3
 \end{pmatrix}.
\]
Cette information suppl\'ementaire provient du fait qu'on a exprim\'e quelles sont les orientations des c\^ot\'es du polygone, et non pas seulement leurs directions. En effet, en \'etudiant le comportement de l'immersion $X$ au voisinage du point $x=t_i$, on a distingu\'e le cas o\`u les c\^ot\'es adjacents au sommet $a_i$ sont dirig\'es par les vecteurs $u_{i-1}$ et $u_i$, du cas o\`u ils sont dirig\'es par les vecteurs $-u_{i-1}$ et $u_i$. Dans le second cas, la normale au sommet $a_i$ est $-v_i$ et l'angle ext\'erieur est $(1-\t_i)\pi$. Ces deux cas d\'efinissent au sommet $a_i$ la m\^eme rotation, mais le choix du relev\'e permet de les distinguer. On en d\'eduit donc \'egalement que les choix des relev\'es $D_i$ des demi-tours sont d\'etermin\'es par les orientations des c\^ot\'es du polygone (\`a une ind\'etermination globale pr\`es, puisque si on remplace toutes les matrices $D_i$ par leurs oppos\'ees, on ne change pas les matrices de monodromie $M_i$). \`A un jeu de directions orient\'ees $D=(D_1,\ldots,D_{n+3})$ correspond donc un $(n+3)$-uplet de relev\'es de demi-tours autour de ces directions, que l'on note \'egalement $D$.
\label{rem-choix-relev\'e}
\end{rema}

Les singularit\'es $t_i$ ($i=1,\ldots,n+3$) sont fuchsiennes. Les autres singularit\'es de l'\'equation~\eqref{E*} sont les ombilics de l'immersion $X$ et leurs conjugu\'es, c'est-\`a-dire des points o\`u le syst\`eme fondamental de solutions $Y_0(x)$ est holomorphe. Ces autres singularit\'es sont donc aussi fuchsiennes. On en d\'eduit donc la proposition suivante.

\begin{prop}
 L'\'equation~\eqref{E*} est une \'equation fuchsienne r\'eelle sur la sph\`ere de Riemann $\P$.
\label{prop-E*-fu}
\end{prop}

On dit que l'\'equation~\eqref{E*} est r\'eelle pour signifier que ses c\oe fficients $p(x)$ and $q(x)$ sont r\'eels sur l'axe r\'eel (proposition~\ref{prop-p-q-reels}).

\section{Singularit\'es apparentes}

Les singularit\'es qui nous reste \`a \'etudier sont les points o\`u les fonctions $G$ et $H$ sont holomorphes, mais o\`u leur Wronskien $GH'-HG'$ s'annule : ce sont les ombilics de l'immersion $X$, et leurs conjugu\'es dans le demi-plan inf\'erieur $\C_-$ (on peut remarquer que, pour une surface minimale, les courbures principales sont nulles en un ombilic). Ces singularit\'es sont fuchsiennes et apparentes (d\'efinition~\ref{def-sg-app}) et leurs exposants sont des entiers naturels. Les deux lemmes suivants pr\'ecisent la valeur de leurs exposants, et le nombre des singularit\'es apparentes.

\begin{lemm}
Les singularit\'es apparentes de l'\'equation~\eqref{E*} sont r\'eelles ou conjugu\'ees deux \`a deux. Deux singularit\'es apparentes qui sont conjugu\'ees ont les m\^emes exposants. Les singularit\'es apparentes de l'\'equation~\eqref{E*} qui sont r\'eelles ou dans $\C_+$ sont les ombilics de l'immersion conforme minimale $X:\C_+\to \R^3$. Les exposants en une de ces singularit\'es $x=\l$ sont $0$ et un entier naturel $m\geq2$, tel que $m-1$ soit l'ordre du z\'ero de la diff\'erentielle de Hopf $Q$ en $x=\l$.
\label{lemme-sg-app}
\end{lemm}

\begin{proof}
Consid\'erons tout d'abord un point r\'egulier quelconque $x=\l$ de l'immersion $X$, $\l\in\overline\C_+\ssm S(t)$. Comme pr\'ec\'edemment, on choisit une position du rep\`ere orthonormal de $\R^3$ tel que le vecteur de Gauss $N(\l)$ de l'immersion $X$ en $x=\l$ co\"incide avec le troisi\`eme vecteur de base $e_3$. Dans cette position, on a
\[
    X(x)-X(\l) = \Re
    \begin{pmatrix}
        (x-\l)   \;\varphi_1(x)\\
        (x-\l)   \;\varphi_2(x)\\
        (x-\l)^{m+1} \;\varphi_3(x)
    \end{pmatrix},
\]
o\`u l'entier $m$ est sup\'erieur ou \'egal \`a $1$, et o\`u les fonctions $\varphi_i(x)$ sont holomorphes au voisinage du point $x=\l$. La fonction $\varphi_3(x)$ ne s'annule pas en $x=\l$, ni l'une ou l'autre des fonctions $\varphi_1(x)$ et $\varphi_2(x)$. Par d\'efinition de la diff\'erentielle de Hopf, l'entier $m-1$ est l'ordre du z\'ero de $Q(x)$ en $x=\l$. Si $m=1$, le point $x=\l$ est un point ordinaire de l'immersion $X$, et si $m\geq2$, c'est un ombilic.

Supposons $m\geq2$. De l'expression de l'immersion $X$ au voisinage de $x=\l$, on d\'eduit que les fonctions $G(x)$ et $H(x)$ satisfont
\begin{itemize}
    \item l'une des primitives $\displaystyle\int_\l^x G(\xi)^2d\xi$ ou $\displaystyle\int_\l^x H(\xi)^2d\xi$ est de la forme $(x-\l)\varphi(x)$,
    \item la primitive $\displaystyle\int_\l^x G(\xi)H(\xi)d\xi$ est de la forme $(x-\l)^{m+1}\varphi(x)$,
\end{itemize}
o\`u $\varphi(x)$ d\'esigne toute fonction holomorphe et non nulle au point $x=\l$. Si on a par exemple
\[
 \int_\l^x G(\xi)^2d\xi = (x-\l)\varphi(x),
\]
alors $G(\l)\neq0$ et $G$ est donc d'exposant $0$. De la deuxi\`eme assertion on d\'eduit alors $H(x) = (x-\l)^m\varphi(x)$.

Dans le demi-plan inf\'erieur $\C_-$, les singularit\'es sont les conjugu\'es des singularit\'es contenues dans $\C_+$ (elles correspondent \`a des points sym\'etriques sur la surface minimale). Comme les exposants en une singularit\'e apparente sont r\'eels, les exposants en deux singularit\'es conjugu\'ees sont les m\^emes.
\end{proof}

L'\'equation~\eqref{E*} a un nombre fini de singularit\'es. Le lemme suivant donne une majoration du nombre $N\in\N$ de singularit\'es apparentes.

\begin{lemm}
 L'\'equation~\eqref{E*} a au plus $n$ singularit\'es apparentes.
\label{lemme-nb-sg-app}
\end{lemm}

\begin{proof}
Il suffit d'appliquer la relation de Fuchs~\eqref{rel-Fu} \`a l'\'equation~\eqref{E*}. On note $\l_1,\ldots,\l_N$ les singularit\'es apparentes, et $m_1,\ldots,m_N$ leurs exposants non nuls respectifs. Par la proposition~\ref{prop-exp-ti} et le lemme~\ref{lemme-sg-app}, et comme l'\'equation~\eqref{E*} a $n+3+N$ singularit\'es, la relation de Fuchs s'\'ecrit
\begin{equation}
 \sum_{i=1}^{n+3}r_i + \sum_{k=1}^N m_k = n+1+N.
\label{rel-Fu-E*}
\end{equation}
Vu les minorations sur les entiers $r_i$ et $m_k$, on obtient $N\leq n$.
\end{proof}

Si le nombre de singularit\'es apparentes est maximal : $N=n$, la valeur des entiers $r_i$ et $m_k$ est d\'etermin\'ee par la relation de Fuchs~\eqref{rel-Fu-E*}, et ils valent alors
\[
  r_i=0 \ (i=1,\ldots,n+2), \qquad r_{n+3}=1, \qquad
    m_k=2 \ (k=1,\ldots,n).
\]
Toutes les singularit\'es de l'\'equation~\eqref{E*} sont donc g\'en\'eriques, et on dit alors que l'immersion $X$ et l'\'equation~\eqref{E*} qui lui est associ\'ee sont elles-m\^emes g\'en\'eriques. En particulier, la surface minimale fait alors en chaque sommet $a_i$ un angle saillant, et elle n'a pas de point de branchement au bord. On peut voir le cas $N<n$ comme provenant de cette situation g\'en\'erique par la fusion de certaines singularit\'es apparentes avec d'autres singularit\'es apparentes ou avec des sommets $t_i$ : c'est effectivement ce qui se produira au cours de la d\'eformation isomonodromique. La fusion d'une singularit\'e apparente d'exposants $0$ et $2$ avec une autre singularit\'e augmente l'un des exposants de cette autre singularit\'e d'une unit\'e. En fait, le sens de ce processus de fusion n'est pas \'evident du point de vue de l'\'equation~\eqref{E*} ; l'utilisation des syst\`emes fuchsiens au chapitre suivant rendra ce processus plus clair et plus simple.

\begin{rema}
Comme on l'a vu \`a la section pr\'ec\'edente, la configuartion g\'en\'erique en le sommet $a_{n+3}$ est d'avoir un bout h\'elico\"idal. Si le bord polygonal de l'immersion $X$ est une courbe ferm\'ee, cela signifie donc qu'une singularit\'e apparente co\"incide avec la singularit\'e $x=\infty$ : ceci transforme les exposants en l'infini de $\left(1-\frac{\t_i}{2} , \frac{\t_i}{2} \right)$ \`a $\left(1-\frac{\t_i}{2} , 1+\frac{\t_i}{2} \right)$. Le nombre maximal de singularit\'es apparentes est alors $n-1$. Ceci explique pourquoi on consid\`ere des disques minimaux ayant pour bord une ligne bris\'ee \emph{pouvant \^etre infinie}. En effet, au chapitre suivant, on sera en particulier amen\'e \`a r\'esoudre le probl\`eme de Riemann--Hilbert pour la monodromie donn\'ee \`a la proposition~\ref{prop-mono}. D'apr\`es le th\'eor\`eme~\ref{thm-boli}, on obtiendra alors des \'equations fuchsiennes ayant au plus $n$ singularit\'es apparantes, et non pas $n-1$. Pour construire des d\'eformations isomonodromiques, que ce soit par le syst\`eme de Garnier ou le syst\`eme de Schlesinger, on a \'egalement besoin g\'en\'eriquement de $n$ singularit\'es apparentes.
\label{rem-N=n}
\end{rema}

\section{Les \'equations fuchsiennes associ\'ees \`a un jeu de directions orient\'ees}
\label{section-cara-equ-fu}

Pour tout jeu de directions orient\'ees $D=(D_1,\ldots,D_{n+3})\in\D^n$, on a montr\'e que pour toute immersion $X\in\X^n_D$, l'unique \'equation diff\'erentielle lin\'eaire du second ordre~\eqref{E*} dont ses donn\'ees de Weierstrass $G$ et $H$ soient solutions satisfait les trois conditions suivantes --- o\`u on note toujours $\t_i\pi$ l'angle ext\'erieur entre les directions $D_{i-1}$ et $D_i$, et o\`u on identifie directions orient\'ees et relev\'es de demi-tours (remarque~\ref{rem-choix-relev\'e}).

\begin{enumerate}[label=(\roman{*}),ref=(\roman{*})]
 \item \label{cond-SdR} L'\'equation~\eqref{E*} est fuchsienne sur la sph\`ere de Riemann $\P$. Elle a $n+3$ singularit\'es non apparentes distinctes $t_1,\ldots,t_n,\ t_{n+1}=0,\ t_{n+2}=1,\ t_{n+3}=\infty$, et au plus $n$ singularit\'es apparentes $\l_1,\dotsc,\l_N$ ($N\leq n$). Son sch\'ema de Riemann est donn\'e par
\begin{equation}
\begin{split}
            &\begin{pmatrix}
                    x=t_i               & x=\infty           & x=\l_k\\
                    -\e_i\frac{\t_i}{2}     & 1-\frac{\t_\infty}{2}          & 0\\
                    r_i+\e_i\frac{\t_i}{2}  & r_\infty-1+\frac{\t_\infty}{2} & m_k
             \end{pmatrix}\\
            &\begin{array}{ccc}
                    i=1,\dotsc,n+2, &   & k=1,\dotsc,N,
             \end{array}
\end{split}
\label{SdR-E*}
\end{equation}
 o\`u $\e_i=\pm1$, les constantes $r_i$ et $m_k$ sont des entiers naturels, qui v\'erifient de plus : $r_\infty\geq1$, $m_k\geq2$ et la relation~\eqref{rel-Fu-E*}.
 \item \label{cond-mono} Un syst\`eme $M_i$ ($i=1,\dotsc,n+3$) de g\'en\'erateurs de la monodromie de l'\'equation~\eqref{E*} le long des lacets $\ga_i$ d\'efinis \`a la figure~\ref{fig-ga-i} s'\'ecrit
\[
 M_i=D_iD_{i-1}^{-1}, \quad \text{ o\`u } D_i \in SU(2), \ D_i^2=-\I_2.
\]
 \item \label{cond-realite} L'\'equation~\eqref{E*} est r\'eelle, et le $n$-uplet de singularit\'es $t=(t_1,\ldots,t_n)$ appartient au simplexe
 \begin{equation}
 \pi^n = \left\{ t\in\R^n \ | \ t_1<\cdots<t_n<0 \right\}.
 \label{def-pi}
 \end{equation}
\end{enumerate}

Remarquons que la condition~\ref{cond-realite}, que l'on appellera \emph{condition de r\'ealit\'e}, assure que les singularit\'es apparentes sont r\'eelles ou conjugu\'ees deux \`a deux. Le fait que les directions $D_i$ ne soient pas toutes coplanaires assure que la monodromie de l'\'equation~\eqref{E*} est irr\'eductible.

\begin{defi}
 Pour tout jeu de directions orient\'ees $D\in\D^n$, on d\'efinit l'ensemble $\E^n_D$ des \'equations fuchsiennes satisfaisant les conditions~\ref{cond-SdR}, \ref{cond-mono} et~\ref{cond-realite} ci-dessus. 
 \label{def-EnD}
\end{defi}

D'apr\`es la proposition~\ref{prop-cara-eq-fu}, et la valeur des exposants, les c\oe fficients $p(x)$ et $q(x)$ d'une \'equation~\eqref{E*} satisfaisant la condition~\ref{cond-SdR} sont de la forme
\begin{align*}
            p(x) &= \sum_{i=1}^{n+2} \dfrac{1-r_i}{x-t_i} + \sum_{k=1}^N \dfrac{1-m_k}{x-\l_k},\\
            q(x) &= -\frac14 \sum_{i=1}^{n+2} \dfrac{\t_i(2\e_ir_i+\t_i)}{(x-t_i)^2} + \dfrac{\kappa}{x(x-1)} - \sum_{i=1}^n\dfrac{t_i(t_i-1)K_i}{x(x-1)(x-t_i)}\\
               & \qquad   + \sum_{k=1}^N\dfrac{\l_k(\l_k-1)\mu_k}{x(x-1)(x-\l_k)},
\end{align*}
o\`u $\kappa=\left(r_\infty-1+\frac{\t_\infty}{2}\right)\left(1- \frac{\t_\infty}{2}\right) +\frac14 \sum_{i=1}^{n+2} \t_i(2\e_ir_i+\t_i) $. Si l'on impose que les $\l_k$ soient des singularit\'es apparentes, alors on obtient que les $K_i$ s'expriment rationnellement en fonction des autres param\`etres $t$, $\l$ et $\mu$ (voir la proposition~\ref{prop-def-Ki}). La strat\'egie que suit Garnier consiste \`a montrer que l'on peut choisir les param\`etres $t$, $\l$ et $\mu$ de telle sorte que l'\'equation ainsi obtenue satisfasse \'egalement la conditions~\ref{cond-mono} et \ref{cond-realite}. On ne d\'etaille pas plus ce point de vue, puisque contrairement \`a Garnier, on va d\`es le chapitre suivant utiliser exclusivement des syst\`emes fuchsiens. 

On peut d\'eduire de l'expression de $p(x)$ que, lorsque $N=n$, c'est-\`a-dire lorsque la surface et l'\'equation sont g\'en\'eriques, la diff\'erentielle de Hopf d'une immersion $X\in\X^n_D$ s'\'ecrit
 \[
 Q = i\frac{\Lambda(x)}{T(x)} \, d x^2,
\]
o\`u les polyn\^omes $\Lambda(x)$ et $T(x)$ sont donn\'es par~\eqref{def-X-T-L}. L'expression g\'en\'erale de la diff\'erentielle de Hopf, lorsque $N<n$, est obtenue en autorisant les $\l_k$ \`a \^etre \'egaux entre eux, et \`a des $t_i$.

D'apr\`es la condition~\ref{cond-mono}, les ensembles $\E^n_D$ sont des ensembles isomonodromiques d'\'equations fuchsiennes. On note
\begin{equation}
  \rho_D:\pi_1\left(\P\ssm S(t),x_0 \right)\to GL(2,\C)
\label{def-rho}
\end{equation}
la repr\'esentation de monodromie engendr\'ee par les matrices $M_i=D_iD_{i-1}^{-1}$. 

La proposition suivante nous dit que les trois conditions ci-dessus caract\'erisent les \'equations diff\'erentielles lin\'eaires du second ordre qui proviennent d'une surface minimale \`a bord polygonal, et donc qu'il est pertinent d'utiliser l'espace $\E^n_D$ pour d\'ecrire les immersions de $\X^n_D$.

\begin{prop}
Soit $D\in\D^n$ un jeu de directions orient\'ees. La correspondance \'etablie par la repr\'esentation de Weierstrass entre les espaces $\X^n_D$ et $\E^n_D$ est bijective. En particulier, toute \'equation de $\E^n_D$ admet un syst\`eme fondamental de solutions $\g$ qui constitue les donn\'ees de Weierstrass d'une immersion conforme minimale $X\in\X^n_D$.
\label{prop-cara-eq}
\end{prop}

Il n'y a aucune traduction naturelle de la longueur des c\^ot\'es du polygone $P$ en terme de propri\'et\'es de l'\'equation fuchsienne~\eqref{E*}. \'Etant donn\'e un jeu de directions orient\'ees $D\in\D^n$, on va donc proc\'eder ainsi pour r\'esoudre le probl\`eme de Plateau : \`a chaque \'equation $(E)$ de l'ensemble $\E^n_D$ est associ\'e par la proposition pr\'ec\'edente un polygone $P_E\in\p^n_D$ pour lequel on sait que le probl\`eme de Plateau admet au moins une solution dans $\X^n_D$. Il s'agit donc de montrer que la famille de polygones $\left(P_E, (E)\in\E^n_D \right)$ ainsi obtenue d\'ecrit enti\`erement l'ensemble $\p^n_D$. On proc\`ede en deux \'etapes : on commence au chapitre~\ref{chapitre-isomono} par d\'ecrire explicitement, en utilisant des d\'eformations isomonodromiques, cette famille de polygones. Puis, on utilise au chapitre~\ref{chapitre-longueur} la description obtenue pour \'etudier leurs rapports de longueur.

\'Enon\c cons d'abord un lemme utile \`a la d\'emonstration de la proposition~\ref{prop-cara-eq} et qui est une cons\'equence imm\'ediate de la m\'ethode de Fr\"obenius en une singularit\'e fuchsienne.

\begin{lemm}
 Soient une \'equation fuchsienne r\'eelle, et $x=x_0$ une singularit\'e r\'eelle et non logarithmique de cette \'equation, d'exposants $\t^-$ et $\t^+$ (qui sont donc r\'eels ou conjugu\'es). Alors, l'\'equation admet en $x=x_0$ un syst\`eme canonique de solutions~:
\[
 g(x)=(x-x_0)^{\t^-}\varphi(x), \qquad
 h(x)=(x-x_0)^{\t^+}\psi(x)
\]
tel que les fonctions $\varphi(x)$ et $\psi(x)$ sont analytiques r\'eelles au voisinage de $x=x_0$.
\label{lemme-sol-cano-r\'eel}
\end{lemm}

\begin{proof}[D\'emonstration de la proposition~\ref{prop-cara-eq}]
Montrons tout d'abord la surjectivit\'e de la correspondance. Soit $(E)$ une \'equation de l'ensemble $\E^n_D$. Remarquons tout d'abord que tout syst\`eme fondamental de solutions $Y_0=\g$ de $(E)$, restreint au demi-plan sup\'erieur $\C_+$, constitue les donn\'ees de Weierstrass d'une immersion conforme minimale $X:\C_+\to\R^3$, d\'efinie \`a translation pr\`es. En effet, les fonctions $G$ et $H$ sont alors holomorphes dans $\C_+$, puisqu'il n'y a pas de singularti\'es non apparentes dans $\C_+$, et elles n'ont pas de z\'ero commun --- sinon, un tel z\'ero serait une singularit\'e apparente de l'\'equation $(E)$ ayant pour exposants deux entiers naturels non nuls, ce qui est exclu par la condition~\ref{cond-SdR}. De plus, cette immersion s'\'etend contin\^ument \`a $\overline\R\ssm S(t)$.

On choisit le syst\`eme fondamental $Y_0(x)$  tel que ses matrices de monodromie le long des lacets $\ga_i$ sont les matrices $M_i$ de la condition~\ref{cond-mono}. Un tel syst\`eme n'est pas unique, l'ensemble des syst\`emes fondamentaux ayant les m\^emes matrices de monodromie sont les syst\`emes $\l\cdot Y_0(x)$ ($\l\in\C^*$). Ceci est une cons\'equence directe de la relation~\eqref{mono-conjuguees} et du fait que les matrices $M_i$ ne sont pas simultan\'ement diagonalisables (car alors les directions $D_i$ seraient toutes coplanaires). Les syst\`emes $\l\cdot Y_0(x)$ d\'efinissent la famille d'immersions $X_\l$. On va montrer que pour un choix convenable $\l_0$ du scalaire $\l$, l'immersion $X_{\l_0}$ est limit\'es par des segments de droite, de directions orient\'ees $D=(D_1,\ldots,D_{n+3})$. On voit qu'un tel scalaire $\l_0$ n'est pas unique, on peut consid\'erer que $\l\in\mathbb S^2$, et que les immersions $X_\l$ sont d\'efinies \`a homoth\'eties de rapport positif pr\`es, \ie sont des \'el\'ements de $\X^n_D$.

Par le lemme~\ref{lemme-realite}, l'immersion $X_\l$ est limit\'ee par des segments de droite si et seulement si, pour tout $i=1,\ldots,n+3$, il existe une matrice $S_i\in SU(2)$ telle que le syst\`eme fondamental $\l\cdot Y_0(x)\cdot S_i$ soit r\'eel ou purement imaginaire sur l'intervalle $]t_i,t_{i+1}[$. On commence par montrer l'existence d'un scalaire $\l$ tel que la condition pr\'ec\'edente soit v\'erifi\'ee pour $i=n+3$. Soit une matrice $S'_\infty\in SU(2)$ telle que
\[
 M_\infty = S'_\infty
    \begin{pmatrix}
     e^{i\t_\infty\pi} & 0\\
     0             & e^{-i\t_\infty\pi}
    \end{pmatrix}
\overline{S'_\infty}^t.
\]
La matrice $S'_\infty$ est un relev\'e d'une rotation envoyant le vecteur normal $v_{n+3}$ sur le vecteur de base $e_3$. Alors le syst\`eme $Y_0(x)\cdot S'_\infty$ est canonique en $x=\infty$, et il s'\'ecrit donc
\[
 Y_0(x)\cdot S'_\infty = \left(a \, g_\infty(x),b \,h_\infty(x)\right)
\]
o\`u $a,\ b\in\C^*$, et o\`u le syst\`eme canonique $(g_\infty(x),h_\infty(x))$ est donn\'e par le lemme~\ref{lemme-sol-cano-r\'eel}. On \'ecrit $a=re^{i(\varphi+\psi)}$ et $b=\rho e^{i(\varphi-\psi)}$, et on choisit
\[
 \l_0:= e^{-i\varphi} \quad \text{et} \quad
 S_\infty := S'_\infty
\begin{pmatrix}
 e^{-i\psi} & 0\\
 0 &  e^{i\psi}
\end{pmatrix}.
\]
Alors la matrice $S_\infty$ est dans $SU(2)$ et on obtient
\[
 \l_0\cdot Y_0(x)\cdot S_\infty = \left(r \, g_\infty(x),\rho \, h_\infty(x)\right) .
\]
Le syst\`eme $\l_0\cdot Y_0(x)\cdot S_\infty$ est donc r\'eel sur l'intervalle $]-\infty,t_1[$.

Montrons \`a pr\'esent qu'il existe une matrice $S_1\in SU(2)$ telle que le syst\`eme $\l_0\cdot Y_0(x)\cdot S_1$ soit r\'eel ou purement imaginaire sur  $]t_1,t_2[$. Par it\'eration, on en d\'eduira le r\'esultat voulu sur chaque intervalle $]t_i,t_{i+1}[$. Le processus d'it\'eration repose sur le fait que d'apr\`es le lemme~\ref{lemme-sol-cano-r\'eel}, il existe pour tout $i=1,\ldots,n+3$ un syst\`eme fondamental canonique au point $x=t_i$
\[
 \mathcal G_i(x) = (g_i(x),h_i(x))
\]
d\'efini et holomorphe dans $\C_+$, qui soit r\'eel sur l'intervalle $]t_i,t_{i+1}[$ et tel que le syst\`eme $(e^{-i\frac{\t_i}{2}\pi}g_i(x),e^{i\frac{\t_i}{2}\pi}h_i(x))$ soit r\'eel sur l'intervalle $]t_{i-1},t_i[$. On en d\'eduit donc que pour tout $i$, il existe une matrice $A_i\in GL(2,\R)$ telle que
\[
 \mathcal G_{i-1}(x) = \mathcal G_i(x)
\begin{pmatrix}
 e^{-i\frac{\t_i}{2}\pi} & 0\\
 0 & e^{i\frac{\t_i}{2}\pi}
\end{pmatrix}
A_i.
\]
Comparons le syst\`eme $\l_0\cdot Y_0(x)$ au syst\`eme $\mathcal G_1(x)$ qui est r\'eel sur $]t_1,t_2[$. Par construction, on a
\[
 \l_0\cdot Y_0(x) = \mathcal G_\infty(x)
\begin{pmatrix}
 r & 0\\
 0 & \rho
\end{pmatrix}
\overline{S_\infty}^t,
\]
et donc
\[
 \l_0\cdot Y_0(x) = \mathcal G_1(x)
\begin{pmatrix}
 e^{-i\frac{\t_1}{2}\pi} & 0\\
 0 & e^{i\frac{\t_1}{2}\pi}
\end{pmatrix}
A_1
\begin{pmatrix}
 r & 0\\
 0 & \rho
\end{pmatrix}
\overline{S_\infty}^t.
\]
Il s'agit de montrer l'existence d'une matrice $B_1\in GL(2,\R)$ et d'une matrice $S_1\in SU(2)$ telles que
\[
 \begin{pmatrix}
 e^{-i\frac{\t_1}{2}\pi} & 0\\
 0 & e^{i\frac{\t_1}{2}\pi}
\end{pmatrix}
A_1
\begin{pmatrix}
 r & 0\\
 0 & \rho
\end{pmatrix}
\overline{S_\infty}^t  =
B_1\overline{S_1}^t.
\]
On obtient ceci en introduisant une matrice $S'_1\in SU(2)$ v\'erifiant
\[
 M_1 = S'_1
    \begin{pmatrix}
     e^{i\t_1\pi} & 0\\
     0            & e^{-i\t_1\pi}
    \end{pmatrix}
 \overline{S'_1}^t,
\]
en comparant, comme \`a l'\'etape pr\'ec\'edente, les syst\`emes $\l_0\cdot Y_0(x)\cdot S'_1$ et $\mathcal G_1(x)$, et enfin en exprimant que le d\'eterminant du produit suivant
\[
 \begin{pmatrix}
 e^{-i\frac{\t_1}{2}\pi} & 0\\
 0 & e^{i\frac{\t_1}{2}\pi}
\end{pmatrix}
A_1
\begin{pmatrix}
 r & 0\\
 0 & \rho
\end{pmatrix}
\overline{S_\infty}^t S'_1
\]
est r\'eel. 

On a donc montr\'e que l'immersion $X_{\l_0}:\C_+\to\R^3$, de donn\'ees de Weierstrass $\l_0\cdot Y_0(x)$, repr\'esente un disque minimal dont le bord est constitu\'e de segments de droite, de longueur \'eventuellement infinie. Vu l'expression des matrices $M_i$ donn\'ee par la condition~\ref{cond-mono}, ces segments de droites sont n\'ecessairement dirig\'es et orient\'es par les $D_i$. Le sch\'ema de Riemann~\eqref{SdR-E*} donne le comportement local de $X_{\l_0}$ au voisinage des points $x=t_i$ : l'immersion $X_{\l_0}$ est born\'ee en les $t_i$ ($i\neq n+3$), et le bord du disque minimal est donc bien un \'el\'ement de $\p^n_D$.

Discutons \`a pr\'esent le comportement de $X_{\l_0}$ au voisinage de $x=\infty$. Si $r_{n+3}\geq2$, alors l'immersion $X_{\l_0}$ se comporte comme en les autres sommets, et le bord du disque minimal repr\'esent\'e par $X_{\l_0}$ est un polygone de $\R^3$. Si $r_{n+3}=1$, l'immersion n'est pas born\'ee au voisinage de $x=\infty$, elle est asymptote \`a une h\'elico\"ide d'axe $v_{n+3}$ contenant les droites passant par les sommets $a_1$ et $a_{n+2}$ et dirig\'ees respectivement par $D_{n+3}$ et $D_{n+2}$. Vu l'\'etude locale r\'ealis\'ee \`a la d\'emonstration de la proposition~\ref{prop-exp-ti}, cette h\'elico\"ide ne peut pas \^etre <<~d\'eg\'en\'er\'ee~>>, \ie elle ne peut pas \^etre plane, et les demi-droites $(a_1,-D_{n+3})$ et $(a_{n+2},D_{n+2})$ ne peuvent pas se couper : le sommet $a_{n+3}$ est en l'infini. Il n'y a pas d'autres comportements possibles en $x=\infty$. L'immersion $X_{\l_0}$ v\'erifie donc bien les conditions de la d\'efinition~\ref{def-XnD} et appartient \`a l'ensemble $\X^n_D$.

Enfin, le caract\`ere injectif de la correspondance entre $\E^n_D$ et $\X^n_D$ provient du fait que d'une part, les immersions de $\X^n_D$ sont d\'efinies \`a translations et homoth\'eties de rapport positif pr\`es, et d'autre part, que dans une famille associ\'ee d'immersions conformes minimales, au plus une immersion repr\'esente une surface minimale \`a bord polygonal.
\end{proof}

\begin{rema}
On observe que dans la d\'emonstration de la proposition~\ref{prop-cara-eq}, pour montrer qu'une \'equation satisfaisant les conditions~\ref{cond-SdR}, \ref{cond-mono} et \ref{cond-realite} d\'efinit une surface minimale limit\'ee par des segments de droite, on a utilis\'e le fait qu'un syst\`eme de g\'en\'erateurs de la monodromie soit dans $SU(2)$, mais nulle part l'\'ecriture en produit de demi-tours successifs de la condition~\ref{cond-mono}. Cette \'ecriture est donc une cons\'equence de la condition de r\'ealit\'e~\ref{cond-realite} et de l'existence d'une repr\'esention unitaire de la monodromie (la condition~\ref{cond-SdR} n'intervient pas dans cette implication). La r\'eciproque de cette assertion nous sera utile dans la r\'esolution du probl\`eme de Plateau et sera discut\'ee au chapitre suivant (proposition~\ref{prop-sys-reel-1}). En fait, on a vu que la condition~\ref{cond-realite} provient du fait que la surface est bord\'ee par un polygone, et que la condition~\ref{cond-mono} est l'expression du principe de r\'eflexion de Schwarz : ces deux conditions ne sont donc pas ind\'ependantes.
\label{rem-iii-implique-ii}
\end{rema}

\chapter{D\'eformations isomonodromiques}
\label{chapitre-isomono}

Le but de ce chapitre est d'obtenir, au moyen de d\'eformations isomonodromiques, une description explicite de l'ensemble $\X^n_D$ des immersions conformes minimales \`a bord polygonal de directions fix\'ees (d\'efinition~\ref{def-XnD}). On va montrer que l'ensemble $\X^n_D$ est param\'etr\'e par le $n$-uplet $t=(t_1,\ldots,t_n)$ des ant\'ec\'edents des sommets, et que la d\'ependance en $t$ des immersions est donn\'ee par le syst\`eme de Schlesinger. Cette description nous sera ensuite utile au chapitre suivant pour r\'esoudre le probl\`eme de Plateau. 

La d\'emarche suivie par Garnier consiste \`a d\'ecrire directement l'ensemble d'\'equations $\E_D^n$ introduit au chapitre pr\'ec\'edent (d\'efinition~\ref{def-EnD}). Les d\'eformations isomonodromiques des \'equations satisfaisant la condition~\ref{cond-SdR} de la section~\ref{section-cara-equ-fu} sont en effet donn\'ees par le syst\`eme de Garnier~\eqref{Gn}. Ce point de vue est tr\`es technique et complexe, principalement parce que le syst\`eme de Garnier n'a pas la propri\'et\'e de Painlev\'e (d\'efinition~\ref{def-prop-P}). On choisit donc plut\^ot de travailler \`a pr\'esent exclusivement avec des syst\`emes fuchsiens, au lieu d'\'equations fuchsiennes.

\bigskip

On va, en se basant sur les r\'esultats de la section~\ref{section-eq->sys}, d\'efinir \`a la section~\ref{section-def-AD} l'ensemble analogue $\A^n_D$ des syst\`emes fuchsiens associ\'es \`a un jeu de directions orient\'ees $D$. La proposition~\ref{prop-eq->sys} nous permet de caract\'eriser les syst\`emes qui appartiennent \`a cet ensemble, en traduisant les conditions~\ref{cond-SdR}, \ref{cond-mono} et \ref{cond-realite}, en des conditions correspondantes~\ref{cond-sys-SdR}, \ref{cond-sys-mono} et \ref{cond-sys-realite} portant sur les syst\`emes.  La condition~\ref{cond-sys-SdR} concerne les singularit\'es et les exposants, la condition~\ref{cond-sys-mono} concerne la monodromie et elle est donc identique \`a la condition~\ref{cond-mono}. La condition~\ref{cond-sys-realite} est toujours une condition de r\'ealit\'e. L'ensemble $\A^n_D$ n'est pas en bijection avec l'ensemble $\X^n_D$, puisque des syst\`emes fuchsiens diff\'erents peuvent d\'efinir la m\^eme \'equation.

Pour d\'ecrire l'ensemble $\A^n_D$, on commence, \`a la section~\ref{section-cond-realite}, par lever une difficult\'e ignor\'ee par Garnier, qu'est la condition de r\'eali\'e~\ref{cond-sys-realite}. On montre que la <<~r\'ealit\'e~>> d'un syst\`eme fuchsien (ou d'une \'equation fuchsienne) peut \^etre caract\'eris\'ee par sa monodromie : on \'etablit une condition n\'ecessaire et suffisante, qu'on appelle \emph{condition \textbf{C1}}, portant sur la monodromie d'un syst\`eme fuchsien pour qu'il v\'erifie la condition~\ref{cond-sys-realite}. En particulier, on montre que la monodromie $\rho_D$ d\'efinie par un jeu $D$ v\'erifie la condition \textbf{C1} : les syst\`emes satisfaisant les conditions~\ref{cond-sys-SdR} et~\ref{cond-sys-mono} v\'erifient donc automatiquement la condition~\ref{cond-sys-realite}.

Enfin, \`a la section~\ref{section-description-sch}, on peut utiliser des d\'eformations isomonodromiques pour d\'ecrire les syst\`emes satisfaisant les conditions~\ref{cond-sys-SdR} et~\ref{cond-sys-mono}. On obtient que l'ensemble $\A^n_D$ contient une famille isomonodromique de syst\`emes fuchsiens $\left(A_D(t), t\in\pi^n\right)$ param\'etr\'ee par les singularit\'es $t=(t_1,\ldots,t_n)$, d\'ecrite par le syst\`eme de Schlesinger~\eqref{schlesinger} et qui est en bijection avec l'ensemble $\X^n_D$. On obtient de plus que la solution $(A_1(t),\ldots,A_{n+2}(t))$ du syst\`eme de Schlesinger correspondant \`a cette famille est holomorphe en tout point du simplexe $\pi^n$ (proposition~\ref{prop-Ai-holo}) : ce r\'esultat, qui simplifiera l'\'etude de la fonction <<~rapports des longueurs~>> au chapitre suivant, est \`a la fois plus fort et plus simple \`a \'etablir que celui obtenu par Garnier pour les \'equations.

Le contenu de ce chapitre est totalement nouveau par rapport \`a l'article de Garnier, et \'egalement beaucoup plus simple que son \'etude des \'equations fuchsiennes de $\E^n_D$.

\section{Les syst\`emes fuchsiens associ\'es \`a un jeu de directions orient\'ees}
\label{section-def-AD}

On souhaite <<~transformer~>> les \'equations de l'ensemble $\E^n_D$ en syst\`emes fuchsiens. On a vu \`a la section~\ref{section-eq->sys} qu'\'etant donn\'e un syst\`eme fuchsien, l'\'equation dont sont solutions les premi\`eres composantes $y_1$ de toute solution $Y=(y_1,y_2)^t$ de ce syst\`eme est fuchsienne (lemme~\ref{lemme-sys->eq}). \`A l'inverse, on a d\'ecrit l'ensemble des syst\`emes fuchsiens normalis\'es en l'infini d\'efinissant, en ce sens, une \'equation fuchsienne donn\'ee (proposition~\ref{prop-eq->sys}). On a donc une correspondance explicite entre \'equations fuchsiennes et syst\`emes fuchsiens normalis\'es en l'infini --- du moins dans le cas \emph{g\'en\'erique}, c'est-\`a-dire lorsque l'\'equation a un nombre maximal $N=n$ de singularit\'es apparentes. Ceci va nous permettre \`a la fois de d\'efinir l'espace analogue $\A^n_D$ des syst\`emes fuchsiens associ\'es \`a un disque minimal \`a bord polygonal, et \'egalement de caract\'eriser les \'el\'ements de cet ensemble par des conditions analogues aux conditions~\ref{cond-SdR}, \ref{cond-mono} et \ref{cond-realite}.

La proposition~\ref{prop-eq->sys} nous dit en particulier qu'un sys\`eme fuchsien non r\'esonnant et normalis\'e en l'infini est d\'etermin\'e par l'\'equation qu'il d\'efinit, par un param\`etre complexe suppl\'ementaire $\xi$, et par le choix d'une normalisation en l'infini parmi deux possibles. Dans la d\'efinition de $\A^n_D$, on impose la normalisation suivante, on verra ensuite pourquoi elle est plus appropri\'ee. Par contre, on a besoin que le param\`etre $\xi$ ne soit pas fix\'e pour pouvoir construire des d\'eformations isomonodromiques.

\begin{defi}
 Pour tout jeu de directions orient\'ees $D\in\D^n$, on d\'efinit l'ensemble $\A^n_D$ des syst\`emes fuchsiens qui d\'efinissent, au sens du lemme~\ref{lemme-sys->eq}, une \'equation qui appartienne \`a l'ensemble $\E^n_D$, et qui soient normalis\'es en l'infini par
  \begin{equation}
 A_\infty = \left(1-\tfrac{\t_\infty}{2}\right)
 \begin{pmatrix}
 1 & 0\\
 0 & -1
 \end{pmatrix}.
 \label{Ainfty}
 \end{equation}
\label{def-AnD}
\end{defi}

Par construction, on obtient donc le r\'esultat suivant.

\begin{prop}
Tout syst\`eme fuchsien appartenant \`a $\A^n_D$ admet une matrice fondamentale de solutions
\[
\Y_0 = 
\begin{pmatrix}
G & H\\
\w G & \w H\\
\end{pmatrix}
\]
dont la premi\`ere ligne $\g$ constitue les donn\'ees de Weierstrass d'une immersion appartenant \`a $\X^n_D$. R\'eciproquement, toute immersion de $\X^n_D$ provient en ce sens d'un syst\`eme de $\A^n_D$.
\label{prop-AnD}
\end{prop}

Remarquons cependant que l'application $\A^n_D \to \X^n_D$ de la proposition pr\'ec\'edente, si elle est toujours bien d\'efinie et surjective, n'est plus injective comme c'\'etait le cas pour l'application analogue $\E^n_D \to \X^n_D$ : comme on n'a pas impos\'e de valeur au param\`etre $\xi\in\C^*$, on a beaucoup plus de syst\`emes que d'\'equations. Pour d\'ecrire l'ensemble $\X^n_D$, il ne sera donc pas n\'ecessaire de d\'ecrire tout l'ensemble $\A^n_D$, mais seulement une partie qui soit en bijection avec $\X^n_D$. En fait, la correspondance entre syst\`emes fuchsiens et disques minimaux est moins naturelle et imm\'ediate que celle entre \'equations fuchsiennes et disques minimaux, puisque il y a beaucoup plus de libert\'e dans le choix d'un syst\`eme associ\'e \`a une immersion. Par exemple, des syst\`emes diff\'erentiels qui ne sont pas fuchsiens d\'efinissent des \'equations qui, elles, sont fuchsiennes, comme le syst\`eme
\[
Y'=
\begin{pmatrix}
0 &1\\
-q(x) & -p(x)
\end{pmatrix}
Y,
\]
o\`u $p(x)$ et $q(x)$ sont les c\oe fficients d'une \'equation fuchsienne.

La proposition~\ref{prop-eq->sys} nous permet d'\'etablir la caract\'erisation suivante.

\begin{theo}
 Pour tout jeu de directions orient\'ees $D\in\D^n$, l'ensemble $\A^n_D$ est l'ensemble des syst\`emes~\eqref{sys-A} qui v\'erifient les trois conditions suivantes.
 \begin{enumerate}[label=(\alph{*}),ref=(\alph{*})]
 \item \label{cond-sys-SdR} Le syst\`eme~\eqref{sys-A} est fuchsien, il a $n+3$ singularit\'es distinctes $t_1,\ldots,t_n$, $t_{n+1}=0$, $t_{n+2}=1$, $t_{n+3}=\infty$, et s'\'ecrit donc~:
\begin{equation}
DY=A(x)Y, \qquad A(x)=\sum_{i=1}^{n+2}\frac{A_i}{x-t_i}.
\tag{$A$}\label{sys-A}
\end{equation}
Pour tout $i=1,\ldots,n+2$, les valeurs propres de la matrice $A_i$ sont $-\frac{\t_i}{2}$ et $\frac{\t_i}{2}$, et \eqref{sys-A} est normalis\'e en l'infini par~\eqref{Ainfty}.
 \item \label{cond-sys-mono} Un syst\`eme $M_i$ ($i=1,\dotsc,n+3$) de g\'en\'erateurs de la monodromie du syst\`eme~\eqref{sys-A} le long des lacets $\ga_i$ d\'efinis \`a la figure~\ref{fig-ga-i} s'\'ecrit
\[
 M_i=D_iD_{i-1}^{-1},
\quad \text{o\`u }D_i\in SU(2),\ {D_i}^2=-\I_2.
\]
 \item \label{cond-sys-realite} Les singularit\'es sont r\'eelles, $t=(t_1,\ldots,t_n)\in\pi^n$, et il existe un nombre r\'eel $\eta$ tel que pour tout $i=1,\ldots,n+2$ la matrice $A_i$ s'\'ecrive
\begin{equation}
  A_i=
\begin{pmatrix}
 a_i&b_ie^{i\eta}\\
 c_ie^{-i\eta}&d_i
\end{pmatrix}
\qquad \text{o\`u }  a_i , d_i\in\R \text{ et }  b_i, c_i \in \left[0,+\infty\right[.
\label{sys-reel}
\end{equation}
\end{enumerate}
\label{thm-AnD}
 \end{theo}
 
Remarquons que la condition~\ref{cond-sys-SdR} est plus simple que la condition analogue~\ref{cond-SdR}. Les syst\`emes v\'erifiant cette condition sont non r\'esonnants car $\t_i\notin\Z$.
 
 \begin{proof}
On \'etablit s\'eparemment chacune des conditions~\ref{cond-sys-SdR}, \ref{cond-sys-mono} et \ref{cond-sys-realite} en traduisant les conditions~\ref{cond-SdR}, \ref{cond-mono} et \ref{cond-realite}.

\vspace{.2cm}

\noindent
\textbf{La condition~\ref{cond-sys-SdR}.}~ Montrons qu'un syst\`eme fuchsien normalis\'e en l'infini par~\eqref{Ainfty} d\'efinit une \'equation satisfaisant la condition~\ref{cond-SdR} si et seulement s'il satisfait la condition~\ref{cond-sys-SdR}. L'\'etude de la section~\ref{section-eq->sys} concerne les \'equations \emph{g\'en\'eriques}, c'est-\`a-dire ayant un nombre maximal $N=n$ de singularit\'es apparentes, et il faut donc la g\'en\'eraliser.

Rappelons que le c\oe fficient $(1,2)$ d'un syst\`eme fuchsien~\eqref{sys-A} normalis\'e en l'infini s'\'ecrit
\[
 A_{12}(x) = \xi \frac{\Lambda(x)}{T(x)},
\]
o\`u
\[
 \Lambda(x)=\prod_{k=1}^n(x-\l_k), \qquad T(x)=\prod_{i=1}^{n+2}(x-t_i),
\]
et $\xi\in\C^*$. En toute g\'en\'eralit\'e, les z\'eros $\l_1,\ldots,\l_n$ de $A_{12}(x)$ peuvent non seulement \^etre \'egaux entre eux, mais \'egalement \^etre confondus avec des p\^oles $x=t_i$. Dans ce cas, la fonction $A_{12}(x)$ est r\'eguli\`ere en $x=t_i$. L'\'equation $(E)$ d\'efinie par un tel syst\`eme est bien d\'efinie, mais ses exposants ne sont pas toujours exactement les valeurs propres des matrices $A_i$. Les z\'eros de $A_{12}(x)$ sont les singularit\'es apparentes de l'\'equation $(E)$. 

On commence par supposer que les valeurs propres des matrices $A_i$ sont $ -\frac{\t_i}{2}$ et $\frac{\t_i}{2}$, c'est-\`a-dire que le syst\`eme~\eqref{sys-A} satisfait la condition~\ref{cond-sys-SdR}, et on montre que l'\'equation~$(E)$ v\'erifie alors la condition~\ref{cond-SdR}. Supposons tout d'abord que les $\l_k$ sont diff\'erents des $t_i$, et \'ecrivons le polyn\^ome $\Lambda(x)$ sous la forme
\[
\Lambda(x)=\prod_{k=1}^N(x-\l_k)^{m_k-1},
\]
o\`u $N\leq n$, $m_k\geq2$, $\sum_{k=1}^N m_k=n-N$, et les $\l_k$ sont \`a pr\'esent suppos\'es distincts. Alors, d'apr\`es le lemme~\ref{lemme-sg-app-sys}, le sch\'ema de Riemann de l'\'equation $(E)$ est
\[
\begin{split}
            &\begin{pmatrix}
                    x=t_i              & x=\infty                      & x=\l_k\\
                    -\frac{\t_i}{2} & 1-\frac{\t_\infty}{2}   & 0\\
                    \frac{\t_i}{2}  & \frac{\t_\infty}{2}       & m_k
             \end{pmatrix}\\
            &\begin{array}{cc}
                    i=1,\dotsc,n+2, &    k=1,\dotsc,N,
             \end{array}
\end{split}
\]
ses singularit\'es sont deux \`a deux distinctes, et les singularit\'es $x=\l_k$ sont apparentes. L'\'equation~$(E)$ v\'erifie alors bien \ref{cond-SdR}.

Consid\'erons \`a pr\'esent que l'un des z\'eros de $A_{12}(x)$ co\"incide avec un p\^ole $x=t_i$, c'est-\`a-dire que $A_{12}(x)$ soit r\'egulier en $x=t_i$ ($i\neq n+3$). V\'erifions que l'un des exposants en $t_i$ de l'\'equation $(E)$ est augment\'e de $1$, \ie que ses exposants sont de la forme $-\e_i\frac{\t_i}{2}$ et $1+\e_i\frac{\t_i}{2}$, o\`u $\e_i=\pm 1$. En effet, on a alors $A^i_{12}=0$ et on choisit $\e_i$ tel que la matrice $A_i$ s'\'ecrive
\[
A_i = 
\begin{pmatrix}
-\e_i\frac{\t_i}{2} & 0\\
A^i_{21} & \e_i\frac{\t_i}{2}
\end{pmatrix}.
\] 
On voit facilement gr\^ace \`a la proposition~\ref{prop-sol-cano} qu'il existe donc une matrice fondamentale de solutions du syst\`eme~\eqref{sys-A} canonique en $x=t_i$
\[
 \Y_i(x) = R_i(x) (x-t_i)^{L_i},
\qquad \text{o\`u } L_i=\e_i\frac{\t_i}{2}
\begin{pmatrix}
  -1 & 0\\
 0 & 1
\end{pmatrix},
\]
telle que la matrice holomorphe $R_i(x)$ soit \'egalement triangulaire inf\'erieure au point $x=t_i$. Ainsi, la premi\`ere ligne de la solution $\Y_i(x)$ a pour exposants $-\e_i\frac{\t_i}{2}$ et $1+\e_i\frac{\t_i}{2}$, bien que $\Y_i(x)$ ait toujours pour exposants les valeurs propres de la matrice $A_i$ : $-\e_i\frac{\t_i}{2}$ et $\e_i\frac{\t_i}{2}$. On peut g\'en\'eraliser cette situation au cas o\`u un nombre arbitraire $r_i\in\N^*$ de singularit\'es apparentes co\"incident avec $t_i$. Alors, la fonction $A_{12}(x)$ a un z\'ero d'ordre $r_i-1$ en $x=t_i$, et il en est de m\^eme pour le c\oe fficient $(1,2)$ de la fonction $R_i(x)$ : l'\'equation~$(E)$ a alors pour exposants $-\e_i\frac{\t_i}{2}$ et $r_i+\e_i\frac{\t_i}{2}$ au point $x=t_i$. Le nombre de ses singularit\'es apparentes a \'et\'e diminu\'e de $r_i$, et son sch\'ema de Riemann est donn\'e par~\eqref{SdR-E*} : l'\'equation~$(E)$ satisfait la condition~\ref{cond-SdR}.

Dans le calcul pr\'ec\'edent, on a vu que seul un des exposants de l'\'equation~$(E)$ en $x=t_i$ peut \^etre augment\'e, et non pas les deux simultan\'ement. En l'infini, on a le m\^eme comportement si des $\l_k$ co\"incident avec le point $x=\infty$, except\'e que seul l'exposant $-1+\frac{\t_\infty}{2}$ peut \^etre augment\'e, puisque la matrice $A_\infty$ est fix\'ee par~\eqref{Ainfty} (alors qu'il y a deux possibilit\'es pour la diagonalis\'ee de $A_i$). 

On remarque que l'on a en fait \'etabli une \'equivalence : pour que l'\'equation~$(E)$ satisfasse la condition~\ref{cond-SdR}, il faut que les valeurs propres des matrices $A_i$ soient $ -\frac{\t_i}{2}$ et $\frac{\t_i}{2}$.
 
\vspace{.2cm}
 
\noindent
\textbf{La condition~\ref{cond-sys-mono}.}~ Il est imm\'ediat que les conditions~\ref{cond-mono} et \ref{cond-sys-mono} soient identiques, puisque un syst\`eme et une \'equation qui sont associ\'es ont la m\^eme monodromie (vu la derni\`ere partie du lemme~\ref{lemme-sys->eq}).

\vspace{.2cm}

\noindent
\textbf{La condition~\ref{cond-sys-realite}.}~ Traduisons \`a pr\'esent la condition de r\'ealit\'e~\ref{cond-realite} pour les syst\`emes. Consid\'erons un syst\`eme fuchsien~\eqref{sys-A} normalis\'e en l'infini dont les singularit\'es $t_i$ sont r\'eelles. Il d\'efinit une \'equation fuchsienne r\'eelle si et seulement s'il d\'efinit la m\^eme \'equation que son syst\`eme conjugu\'e :
\begin{equation}
  DY=\tau(A)(x)Y, \qquad
  \tau(A)(x)=\sum_{i=1}^{n+2}\frac{\overline{A}_i}{x-t_i},
\tag{$\tau A$}\label{sys-conj}
\end{equation}
o\`u l'application $\tau$ est d\'efinie par~\eqref{def-tau}. Le syst\`eme conjugu\'e est \'egalement fuchsien et normalis\'e en l'infini. Son r\'esidu en l'infini est $\overline{A}_\infty$. S'il d\'efinit la m\^eme \'equation que le syst\`eme~\eqref{sys-A}, alors les matrices $A_i$ et $\overline A_i$, $i=1,\ldots,n+2$, ont les m\^emes valeurs propres (puisque elles sont les exposants de l'\'equation), qui sont donc r\'eelles ou conjugu\'ees entre elles. Par contre, les syst\`emes~\eqref{sys-A} et~\eqref{sys-conj} n'ont pas n\'ecessairement la m\^eme normalisation en l'infini, si on suppose seulement que les valeurs propres du r\'esidu $A_\infty$ sont r\'eelles ou conjugu\'ees entre elles. Si on suppose que les valeurs propres de $A_\infty$ sont r\'eelles (et c'est bien le cas ici), alors les syst\`emes~\eqref{sys-A} et~\eqref{sys-conj} ont la m\^eme normalisation en l'infini. Alors, par la proposition~\ref{prop-eq->sys}, pour qu'ils d\'efinissent la m\^eme \'equation, il faut et il suffit qu'il existe un nombre complexe non nul $\xi $ tel que pour tout $i=1,\ldots,n+2$, on ait
\[
 \overline{A}_i=
    \begin{pmatrix}
     A_{11}^i               & \xi  A_{12}^i\\
     \frac1\xi  A_{21}^i & A_{22}^i
    \end{pmatrix}.
\]
Alors $\left| \xi\right| =1$, et la condition pr\'ec\'edente est \'equivalente \`a~\eqref{sys-reel}.
\end{proof}

La d\'emonstration de la condition~\ref{cond-sys-SdR} nous permet de justifier le choix de la normalisation en l'infini~\eqref{Ainfty} : pour que le bord polygonal de l'immersion associ\'ee \`a une \'equation de $\E^n_D$ soit une courbe ferm\'ee, il faut que les exposants de l'\'equation soient de la forme : $1-\frac{\t_\infty}{2}$ et $r_\infty-1+\frac{\t_\infty}{2}$, avec $r_\infty\geq2$. C'est donc ce deuxi\`eme exposant que l'on veut pouvoir augmenter.

\bigskip

L'introduction de l'ensemble des syst\`emes fuchsiens associ\'es \`a un jeu de directions orient\'ees apporte un point de vue nouveau \`a la m\'ethode de Garnier. L'approche suivie par Garnier est la suivante : il d\'ecrit l'ensemble des \'equations satisfaisant les conditions~\ref{cond-SdR} et \ref{cond-mono} au moyen du syst\`eme de Garnier (en oubliant la condition de r\'ealit\'e~\ref{cond-realite}), et il obtient ainsi une famille d'\'equations $\left(E_D(t), t\in\pi^n \right) $ param\'etr\'ee par $t$. Cependant, le syst\`eme de Garnier n'ayant pas la propri\'et\'e de Painlev\'e, il est oblig\'e \`a plusieurs reprises d'utiliser le syst\`eme de Schlesinger pour \'etudier cette famille d'\'equations : Garnier ne voit les syst\`emes fuchsiens que comme un outil ponctuel permettant de lever certaines difficult\'es rencontr\'ees avec les \'equations fuchsiennes, principalement pour \'etudier la r\'egularit\'e de la fonction <<~rapports des longueurs~>>. \`A chaque fois, Garnier fait une sorte d'aller-retour entre \'equations et syst\`emes. Ce travail est long et complexe, il repose sur l'\'etude de la transformation du syst\`eme de Garnier en le syst\`eme de Schlesinger, qui a depuis \'et\'e expos\'e en d\'etail dans~\cite{IKSY} (chapitre 3, section 6).

On a choisi au contraire d'adapter les r\'esultats du chapitre~\ref{chapitre-equ-fu} de mani\`ere \`a obtenir directement une correspondance entre les disques minimaux \`a bord polygonal et les syst\`emes fuchsiens, puis de travailler exclusivement avec ces derniers. Cette utilisation syst\'ematique des syst\`emes fuchsiens pr\'esente de nombreux avantages : d'une part, comme on l'a dit, elle permet d'\'eviter d'\'etudier la transformation du syst\`eme de Garnier en le syst\`eme de Schlesinger. D'autre part, comme les syst\`emes ont une structure plus canonique que les \'equations, cette approche permet de multiples simplifications : notamment gr\^ace \`a la propri\'et\'e de Painlev\'e, mais pas uniquement, comme la proposition~\ref{prop-Ai-holo}.

Un autre point que l'on va d\'evelopper dans ce chapitre et qui est compl\`etement absent de l'article de Garnier est l'\'etude de la condition de r\'ealit\'e~\ref{cond-sys-realite}. Il semble que Garnier consid\`ere que la famille isomonodromique $\left(E_D(t), t\in\pi^n \right) $ qu'il a construite v\'erifie automatiquement la condition~\ref{cond-realite}, et il lui donne un sens g\'eom\'etrique en terme de surfaces minimales --- bien qu'il n'ait pas non plus \'etabli de r\'esultat analogue \`a la proposition~\ref{prop-cara-eq-fu}. Cette interpr\'etation est malgr\'e tout exacte, puisque on va montrer \`a la section suivante que la condition de r\'ealit\'e~\ref{cond-sys-realite} est une cons\'equence des conditions~\ref{cond-sys-SdR} et \ref{cond-sys-mono}.

\section{La condition de r\'ealit\'e}
\label{section-cond-realite}

Cette section ne concerne pas uniquement l'ensemble $\A^n_D$, on va \'etablir des r\'esultats g\'en\'eraux sur les syst\`emes fuchsiens non r\'esonnants et normalis\'es en l'infini. On a vu \`a la remarque~\ref{rem-iii-implique-ii} que la condition~\ref{cond-mono} et la condition r\'ealit\'e~\ref{cond-realite} ne sont pas ind\'ependantes. On va montrer que pour les syt\`emes fuchsiens, la condition de r\'ealit\'e~\ref{cond-sys-realite} est \'equivalente \`a une condition, que l'on appellera \emph{condition \textbf{C1}}, portant uniquement sur la monodromie, et que cette condition est v\'erifi\'ee en particulier par une monodromie satisfaisant la condition~\ref{cond-sys-mono}. Pour cela, on \'etablit d'abord un r\'esultat d'unicit\'e classique pour les syst\`emes fuchsiens non r\'esonnants.

\subsection{Un r\'esultat d'unicit\'e}

\begin{lemm}
Soient deux syst\`emes fuchsiens non r\'esonnants
\begin{align}
 DY&=A(x)Y \tag{$A$}\label{A}\\
 DZ&=B(x)Z \tag{$B$}\label{B}
\end{align}
Les syst\`emes~\eqref{A} et~\eqref{B} ont les m\^emes singularit\'es, les m\^emes exposants et la m\^eme monodromie si et seulement s'il existe une matrice inversible $C$ telle que
\[
 B(x)=CA(x)C^{-1}.
\]

Si, de plus, les syst\`emes~\eqref{A} et~\eqref{B} sont normalis\'es
en l'infini, alors il existe un nombre complexe non nul $\xi $ tel
que la matrice $C$ soit \'egale \`a
\[
 C=
\begin{pmatrix}
 1&0\\
 0&\xi 
\end{pmatrix}
\text{ ou }
\begin{pmatrix}
 0&1\\
 \xi &0
\end{pmatrix}.
\]
\label{lemme-unicite}
\end{lemm}

\begin{proof}
La condition suffisante est \'evidente. Supposons que les syst\`emes~\eqref{A} et~\eqref{B} ont les m\^emes singularit\'es $t_1,\ldots, t_n$, $t_{n+1}=0$, $t_{n+2}=1$, $t_{n+3}=\infty$, les m\^emes exposants et la m\^eme monodromie. Alors ils s'\'ecrivent
\[
 A(x)=\sum_{i=1}^{n+2}\frac{A_i}{x-t_i}, \qquad B(x)=\sum_{i=1}^{n+2}\frac{B_i}{x-t_i},
\]
et les matrices $A_i$ et $B_i$ ont les m\^emes valeurs propres. On note
\[
 L_i = \begin{pmatrix}
        \t_i^+ & 0\\
    0      & \t_i^-
       \end{pmatrix}
\]
la diagonalis\'ee des matrices $A_i$ et $B_i$. Il existe deux matrices fondamentales $\Y(x)$ et $\mathbf Z(x)$ de solutions, respectivement, du syst\`eme~\eqref{A} et du syst\`eme~\eqref{B}, qui ont les m\^emes matrices de monodromie. On pose alors pour tout $x$ dans le rev\^etement universel de l'ensemble $\P\ssm S(t)$
\[
 C(x):=\mathbf Z(x) \cdot \Y(x)^{-1}.
\]
La matrice $C(x)$ est donc m\'eromorphe dans le rev\^etement universel de $\P\ssm S(t)$ ; on va montrer qu'elle est holomorphe dans $\P$, c'est-\`a-dire constante. Remarquons tout d'abord que $C(x)$ est uniforme dans $\P\ssm S(t)$ : en effet, pour tout $\gamma \in\pi_1(\P\ssm S(t),x_0)$, vu que $M_\gamma(\Y) =M_\gamma(\mathbf Z)$, on
a
\begin{align*}
 \gamma *C(x) =\left( \mathbf Z(x) \cdot M_\gamma(\mathbf Z)\right)  \cdot \left(\Y(x) \cdot M_\gamma(\Y) \right) ^{-1} = C(x).
\end{align*}
De plus, la matrice $C(x)$ n'est singuli\`ere qu'aux points o\`u $\det\Y(x)$ s'annule. Ceci est impossible, car la fonction $\det\Y(x)$ v\'erifie
\begin{align*}
 D\left(\det \Y(x)\right) &=\det \Y(x) \tr \left( D\Y(x) \cdot \Y(x)^{-1} \right) \\
    &= \det \Y(x)\tr A(x) \\
    &=\det \Y(x)\sum_{i=1}^{n+2} \frac{\tr L_i}{x-t_i},
\end{align*}
ce qui donne
\[
 \det \Y(x) = K\prod_{i=1}^{n+2}(x-t_i)^{\tr L_i}
\]
($K\in\C^*$). La matrice $C(x)$ est donc holomorphe dans $\P\ssm S(t)$.

\'Etudions \`a pr\'esent le comportement de $C(x)$ au voisinage d'une singularit\'e $x=t_i$ ($i=1,\ldots,n+2$). Soit $M_i$ la matrice de monodromie des matrices fondamentales $\Y(x)$ et $\mathbf Z(x)$ autour de la singularit\'e $x=t_i$~:
\[
 M_i = C_i
\begin{pmatrix}
 e^{2i\pi\t_i^+} & 0\\
 0  &      e^{2i\pi\t_i^-}
\end{pmatrix}
C_i^{-1}
\]
o\`u $C_i\in GL(2,\C)$. Alors les matrices fondamentales de solutions $\Y(x)\cdot C_i$ et $\mathbf Z(x)\cdot C_i$ sont canoniques au point $x=t_i$~:
\begin{align*}
 \mathbf Y(x)\cdot C_i &= R_i(x)(x-t_i)^{L_i}\\
 \mathbf Z(x)\cdot C_i &= S_i(x)(x-t_i)^{L_i}
\end{align*}
o\`u les matrices $R_i(x)$ et $S_i(x)$ sont holomorphes et inversibles au point $x=t_i$. On en d\'eduit
\[
 C(x)=S_i(x)R_i(x)^{-1},
\]
et $C(x)$ est holomorphe en $x=t_i$. On montrerait de m\^eme que la matrice $C(x)$ est holomorphe en $x=\infty$. Elle est donc holomorphe sur la sph\`ere de Riemann $\P$ : elle est ind\'ependante de $x$.

Si on suppose de plus que les matrices $A_\infty$ et $B_\infty$ sont diagonales, alors on note
\[
A_\infty =
\begin{pmatrix}
    \t_\infty^+ & 0\\
    0 & \t_\infty^-
\end{pmatrix}
\]
et donc
\[
 B_\infty = A_\infty \text{ ou }
\begin{pmatrix}
    \t_\infty^- & 0\\
    0 & \t_\infty^+
\end{pmatrix},
\]
\ie $B_\infty=A_\infty$ ou $B_\infty=JA_\infty J^{-1}$ avec
\[
 J=\begin{pmatrix}
    0&-1\\
    1&0
\end{pmatrix}.
\]
Comme par ailleurs $B_\infty=CA_\infty C^{-1}$ et comme $\t_\infty^+\neq \t_\infty^-$, on en d\'eduit dans le premier cas que la matrice $C$ est diagonale, et dans le second, qu'elle est anti-diagonale.
\end{proof}

\subsection{Syst\`emes fuchsiens <<~r\'eels~>>}

On vient de voir qu'un syst\`eme fuchsien non
r\'esonnant et normalis\'e en l'infini est enti\`erement d\'etermin\'e par
ses singularit\'es $t_1,\ldots,t_n$, par les valeurs propres des
matrices $A_i$, par sa monodromie et par un param\`etre
suppl\'ementaire $\xi \in\C^*$. On va d\'eterminer \`a pr\'esent \`a quelle
condition sur ces donn\'ees le syst\`eme~\eqref{A00} v\'erifie la
condition de r\'ealit\'e~\ref{cond-sys-realite}. Pour les singularit\'es
et les valeurs propres, la r\'eponse est imm\'ediate : les
singularit\'es doivent \^etre r\'eelles ou conjugu\'ees deux \`a deux ; les
valeurs propres en une singularit\'e r\'eelle doivent \^etre r\'eelles ou
conjugu\'ees entre elles, et les valeurs propres en deux
singularit\'es conjugu\'ees doivent \^etre conjugu\'ees. On ne s'int\'eresse
ici qu'au cas o\`u les singularit\'es $t_1,\ldots,t_n$ sont r\'eelles
(on obtiendrait le m\^eme r\'esultat dans le cas o\`u elles sont
seulement r\'eelles ou conjugu\'ees deux \`a deux, mais la d\'emonstration
est un peu plus technique). Par souci de simplicit\'e, on suppose que $t=(t_1,\ldots,t_n)$ est dans le simplexe $\pi^n$ d\'efini par~\eqref{def-pi}.

On reprend les notations de la section~\ref{section-sys-fu}. On consid\`ere un syst\`eme fuchsien non r\'esonnant et normalis\'e en l'infini
\begin{equation}
DY=A(x)Y, \qquad A(x)=\sum_{i=1}^{n+2}\frac{A_i}{x-t_i},
\tag{$A_0$}\label{A00}
\end{equation}
et on note $\t_i^+$ et $\t_i^-$ les valeurs propres des matrices de r\'esidu $A_i$ ($i=1,\ldots,n+3$).

\begin{prop}
On suppose que le $n$-uplet de singularit\'es $(t_1,\ldots,t_n)$ du syst\`eme fuchsien~\eqref{A00} est dans le simplexe $\pi^n$, que les valeurs propres $\t_i^+$ et $\t_i^-$ ($i=1,\ldots,n+2$) sont r\'eelles ou conjugu\'ees entre elles, et que les valeurs propres $\t_\infty^+$ et $\t_\infty^-$ sont r\'eelles. Alors les trois assertions suivantes sont \'equivalentes~:
\begin{itemize}
\item l'\'equation fuchsienne associ\'ee au sens du lemme~\ref{lemme-sys->eq} au syst\`eme~\eqref{A00} est r\'eelle ;
\item les matrices $A_i$ sont de la forme~\eqref{sys-reel} ;
\item pour tout syst\`eme de g\'en\'erateurs $\left(M_1,\ldots,M_{n+3}\right)$ de la monodromie le long des lacets $\ga_1,\ldots,\ga_{n+3}$, il existe une matrice $C\in GL_2(\C)$ telle que pour tout $j=1,\ldots,n+3$ on ait
\begin{equation}
 C^{-1}\overline{M_j}C = (M_j\ldots M_1)^{-1} M_j^{-1} (M_j\ldots M_1).
\label{C1}
\end{equation}
On appelle la derni\`ere de ces assertions la \emph{condition \textbf{C1}}.
\end{itemize}
\label{prop-sys-reel-1}
\end{prop}

On remarque qu'il existe un syst\`eme de g\'en\'erateurs $\left(M_1,\ldots,M_{n+3}\right)$ pour lequel la matrice $C$ est la matrice identit\'e $\I_2$.

\begin{proof}
On a d\'ej\`a vu que les deux premi\`eres assertions sont \'equivalentes. On consid\`ere le syst\`eme conjugu\'e~\eqref{tau-A0} au syst\`eme~\eqref{A00} :
\begin{equation}
  DY=\tau(A)(x)Y, \qquad
  \tau(A)(x)=\sum_{i=1}^{n+2}\frac{\overline{A}_i}{x-t_i}.
\tag{$\tau A_0$}\label{tau-A0}
\end{equation}
Le syst\`eme~\eqref{A00} d\'efinit une \'equation fuchsienne r\'eelle s'il d\'efinit la m\^eme \'equation que le syst\`eme conjugu\'e~\eqref{tau-A0}. On a vu \'egalement que ceci \'equivaut \`a l'existence d'un nombre $\xi \in\C^*$ tel que pour tout $i=1,\ldots,n+2$, on ait
\[
\overline A_i =
\begin{pmatrix}
 1 & 0\\
 0 & \xi 
\end{pmatrix}
A_i
\begin{pmatrix}
 1 & 0\\
 0 & \xi 
\end{pmatrix}
^{-1}.
\]
Par hypoth\`ese, les syst\`emes~\eqref{A00} et~\eqref{tau-A0} ont les m\^emes singularit\'es, les m\^emes exposants et la m\^eme normalisation en l'infini. D'apr\`es le lemme~\ref{lemme-unicite}, ils d\'efinissent donc la m\^eme \'equation si et seulement s'ils ont la m\^eme monodromie.

Soit $\Y(x)$ une matrice fondamentale de solutions du syst\`eme~\eqref{A00} d\'efinie et holomorphe dans le demi-plan sup\'erieur $\C_+$. On note ses matrices de monodromie $M_i$ :
\begin{equation*}
 M_{\ga_i}(\Y)=M_i.
\end{equation*}
On note $\Y_i(x)$ le prolongement \`a $\C_-$ de la matrice fondamentale $\Y(x)$ \`a travers l'intervalle $]t_i,t_{i+1}[$ (c'est-\`a-dire le long de tout chemin joignant un point de $\C_+$ \`a un point de $\C_-$ et croisant l'axe r\'eel une seule fois entre $t_i$ et $t_{i+1}$) ; la matrice fondamentale $\Y_i(x)$ est d\'efinie et holomorphe sur l'ouvert simplement connexe
\[
 U_i :=\C_+ \cup \C_- \cup \, ]t_i,t_{i+1}[.
\]
La matrice $\tau(\Y_i)(x)$, elle aussi holomorphe et inversible sur $U_i$, est une matrice fondamentale de solutions du syst\`eme~\eqref{tau-A0}. Pour que les syst\`emes~\eqref{A00} et~\eqref{tau-A0} aient la m\^eme
monodromie, il faut et il suffit que pour une valeur de $i$, les matrices de monodromie des solutions fondamentales $\Y_i(x)$ et $\tau(\Y_i)(x)$ soient conjugu\'ees, c'est-\`a-dire qu'il existe une matrice inversible $C$ telle que pour tout $j=1,\ldots,n+3$ on ait
\[
 M_{\ga_j}\left(\tau(\Y_i)\right)=CM_jC^{-1}.
\]
On choisit le prolongement $\Y_{n+3}(x)$ (\ie $i=n+3$).

\begin{figure}[!h]
\centering
\begin{pspicture}(0,0)(16,8)
\rput(8,4){\includegraphics[height=7cm]{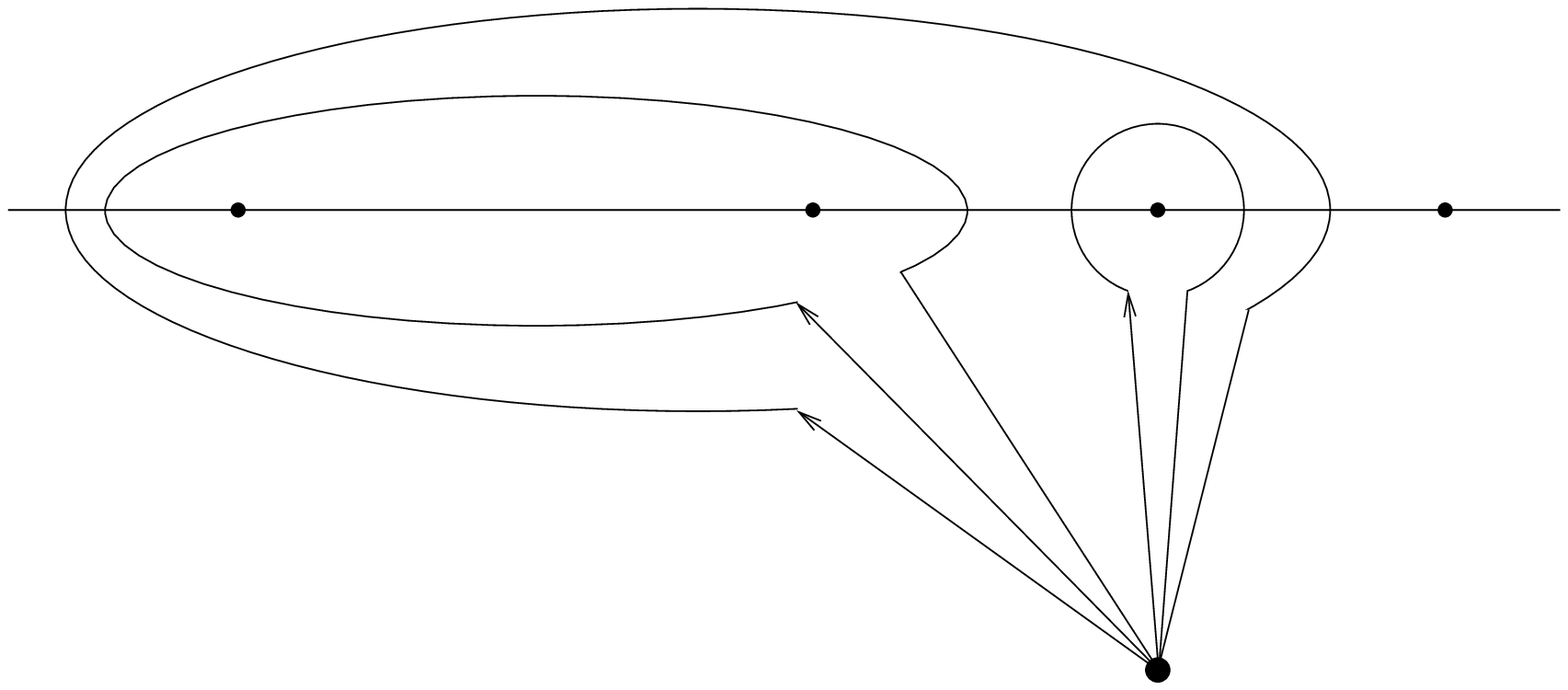}}
\rput(11.85,0.2){$\bar x_0$}
\rput(11.2,4.3){$\overline\ga_j$}
\rput(2.35,5.7){$t_1$}
\rput(8.3,5.7){$t_{j-1}$}
\rput(11.9,5.7){$t_j$}
\rput(14.8,5.7){$t_{j+1}$}
\rput(5.3,5.7){$\ldots$}
\rput(5.4,3.1){$\a$}
\rput(6.6,4.45){$\b$}
\end{pspicture}
\caption{Les lacets $\overline{\ga}_j$, $\a$ et $\b$}
\label{fig-lacet2}
\end{figure}

Il faut exprimer en fonction des matrices $M_j$ les matrices :
\begin{equation*}
 M_{\ga_j}(\tau(\Y_{n+3}))=\overline{M_{\overline{\ga}_j}(\Y_{n+3})}.
\end{equation*}
Le lacet $\overline{\ga}_j$ a pour point de base $\bar{x}_0$ et tourne en sens inverse du sens trigonom\'etrique autour de $t_j$. Pour calculer $M_{\overline{\ga}_j}(\Y_{n+3})$, la difficult\'e vient de ce qu'on sait comment est transform\'ee, en g\'en\'eral, la matrice fondamentale $\Y_i(x)$ le long des lacets $\overline{\ga}_i$ et $\overline{\ga}_{i+1}$ (comme on le voit sur la figure~\ref{fig-lacet3}), mais pas le long d'un lacet $\overline{\ga}_j$ quelconque. On va donc proc\'eder par it\'eration. On d\'ecompose $\overline{\ga}_j$ en le produit de deux lacets. Soient $\a, \b \in \pi_1(\P\ssm S(t), \bar{x}_0)$ les deux classes de lacets orient\'es n\'egativement et qui encerclent respectivement les singularit\'es $t_1,\ldots,t_j$ et $t_1,\ldots,t_{j-1}$ (l\`a encore, les indices s'entendent modulo $n+3$ : si $j=1$, le lacet $\b$ est homotope \`a un point). Les lacets $\a$ et $\b$ sont repr\'esent\'es \`a la figure~\ref{fig-lacet2}. Alors
\[
 \overline{\ga}_j = \a \b^{-1},
\]
et donc
\[
 M_{\overline{\ga}_j}(\Y_{n+3}) = M_{\a}(\Y_{n+3})M_{\b}(\Y_{n+3})^{-1}.
\]
Montrons que
\[
 M_{\a}(\Y_{n+3}) =M_1^{-1}\ldots M_j^{-1}.
\]
On remarque que, vu la d\'efinition des matrices fondamentales $\Y_i(x)$, on a pour tout $i=1,\ldots,n+3$
\[
 \overline{\ga}_i\ast \Y_{i-1}(x) = \Y_i(x).
\]
Donc, comme $\a=\overline{\ga}_j \cdots \overline{\ga}_1$, on obtient par it\'eration
\[
 \a\ast \Y_{n+3}(x) = \Y_j(x).
\]

\begin{figure}[!h]
\centering
\begin{pspicture}(0,0)(6,5)
\rput(3,2.5){\includegraphics{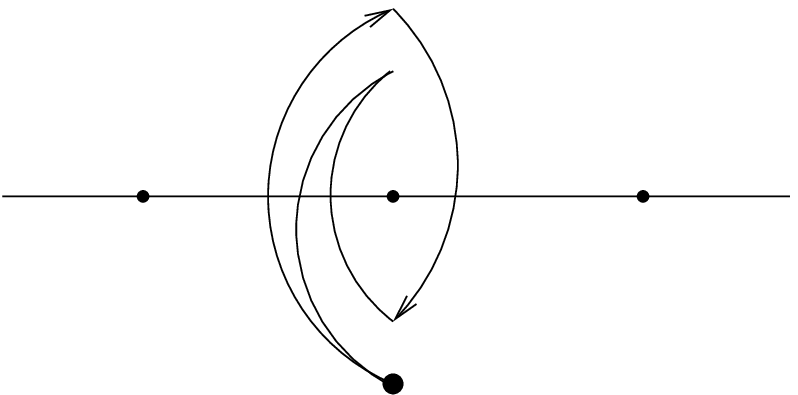}}
\rput(2.95,0.2){$\bar x_0$}
\rput(3.5,4.3){$\overline\ga_i$}
\rput(0.4,2.8){$t_{i-1}$}
\rput(2.95,2.8){$t_i$}
\rput(5.45,2.8){$t_{i+1}$}
\end{pspicture}
\caption{On d\'ecompose le lacet $\overline{\ga}_j$}
\label{fig-lacet3}
\end{figure}

\noindent
Par ailleurs, comme on a aussi $\overline{\ga}_i\ast \Y_{i-1}(x) = \Y_{i-1}(x)M_i^{-1}$ (voir la d\'ecomposition du lacet $\overline\ga_j$ \`a figure~\ref{fig-lacet3}), alors pour tout $i=1,\ldots,n+3$
\[
 \Y_i(x)=\Y_{i-1}(x)M_i^{-1},
\]
donc $\Y_j(x)=\Y_{j-1}(x)M_j^{-1}=\cdots = \Y_{n+3}(x)(M_1^{-1} \ldots M_j^{-1}) $, ce qui donne la formule annonc\'ee pour $M_{\a}(\Y_{n+3}) $. De m\^eme, on a
\[
 M_{\b}(\Y_{n+3}) =M_1^{-1}\ldots M_{j-1}^{-1}.
\]
Finalement, on obtient
\[
 M_{\ga_j}(\tau(\Y_{n+3}))=\overline{(M_j\ldots M_1)^{-1} M_j^{-1} (M_j\ldots M_1)},
\]
ce qui donne bien l'\'equivalence annonc\'ee.
\end{proof}

\begin{rema}
On a un r\'esultat analogue \`a la proposition~\ref{prop-sys-reel-1} pour les \'equations fuchsiennes : \emph{une \'equation fuchsienne sans singularit\'e logarithmique, dont les singularit\'es et les exposants sont r\'eels, est r\'eelle si et seulement si sa monodromie v\'erifie la condition \textbf{C1}}. Dans~\cite{Desideri-these}, on en a d\'eduit que la condition \textbf{C1} est \'egalement une condition n\'ecessaire et suffisante pour qu'une solution $\left( \l(t),\mu(t) \right) = \left( \l_1(t),\dots,\l_n(t),\mu_1(t),\ldots,\mu_n(t) \right)$ du syst\`eme de Garnier~\eqref{Gn} (voir l'appendice~\ref{annexe-garnier}) soit r\'eelle ou conjugu\'ee deux \`a deux, c'est-\`a-dire que $\overline \l(\bar t)$ et $\overline \mu(\bar t)$ soient obtenus \`a partir respectivement de $\l(t)$ et $\mu(t)$ par une m\^eme permutation de leurs indices (corollaire 3.17. de~\cite{Desideri-these}).
\end{rema}

\subsection{Cas o\`u la monodromie est unitarisable}

Dans le cas o\`u il existe un syst\`eme de g\'en\'erateurs $\left(M_1,\ldots,M_{n+3}\right)$ de la monodromie du syst\`eme~\eqref{A00} qui soit contenu dans le groupe des matrices unitaires $U(2)$, ou dans le groupe $U(1,1)$, on peut simplifier l'\'ecriture de la condition \textbf{C1}.

\begin{prop}
Sous les m\^emes hypoth\`eses qu'\`a la proposition~\ref{prop-sys-reel-1}, si un syst\`eme de g\'en\'erateurs $\left(M_1,\ldots,M_{n+3}\right)$ de la monodromie du syst\`eme~\eqref{A00} est contenu dans $U(2)$ ou dans $U(1,1)$, alors le syst\`eme~\eqref{A00} v\'erifie l'une des trois assertions \'equivalentes de la proposition~\ref{prop-sys-reel-1} si et seulement s'il existe $n+3$ matrices inversibles $D_1,\ldots,D_{n+3}$ telles que
\[
\begin{cases}
     M_j=D_jD_{j-1}^{-1} \qquad (j=1,\ldots,n+3)\\
     \frac{1}{\d_1}{D_1}^2 =\cdots= \frac{1}{\d_{n+3}}{D_{n+3}}^2
\end{cases}
\]
o\`u on a not\'e $\d_j = \det D_j$ pour tout $j=1,\ldots,n+3$. On appelle cette condition la \emph{condition \textbf{C2}}.
\label{prop-sys-reel-2}
\end{prop}

Rappelons que le groupe $U(1,1)$ est le groupe des matrices $M\in M(2,\C)$ telles que
\[
 M
 \begin{pmatrix}
  1 & 0\\
  0 & -1
 \end{pmatrix}
 \overline M^t =
 \begin{pmatrix}
  1 & 0\\
  0 & -1
 \end{pmatrix} .
\]

\begin{proof}
Pour toute matrice $M\in U(2)$, on a
\[
 J^{-1} M J = \det(M) \overline M
\]
(ce qui redonne la relation~\eqref{SU2} lorsque $M\in SU(2)$). Si les matrices $M_1,\ldots,M_{n+3}$ sont dans le groupe unitaire $U(2)$, alors la condition \textbf{C1} est \'equivalente \`a l'existence d'une matrice inversible $C$ telle que pour tout $j=1,\ldots,n+3$, on ait
\begin{equation}
 (JC)^{-1}M_j(JC) = \det(M_j)(M_j\ldots M_1)^{-1} M_j^{-1} (M_j\ldots M_1).
\label{C1'}
\end{equation}
(\emph{condition \textbf{C1'}}). On a la m\^eme expression lorsque les matrices $M_1,\ldots,M_{n+3}$ sont dans le groupe $ U(1,1)$, en remplaçant la matrice $J$ par la matrice
\[
\begin{pmatrix}
 0 & i\\
 i & 0
\end{pmatrix}.
\]
Pour la d\'emonstration, on se limitera donc au cas o\`u le syst\`eme de g\'en\'erateurs est dans le groupe $U(2)$.

Montrons que les deux conditions \textbf{C1'} et \textbf{C2} sont \'equivalentes. Pour tout choix de la matrice inversible $D_{n+3}$, par la relation $M_{n+3}\cdots M_1 = \I_2$, il existe des matrices inversibles $D_1,\ldots,D_{n+2}$, d\'etermin\'ees de mani\`ere unique, telles que pour tout $j=1,\ldots,n+3$, on ait
\[
 M_j=D_jD_{j-1}^{-1}
\]
(o\`u les indices sont consid\'er\'es modulo $n+3$). Alors on a
\[
 \det(M_j) = \frac{\d_j}{\d_{j-1}}
\]
o\`u $\d_j=\det D_j$. La relation~\eqref{C1'}$_j$ se r\'ecrit alors de la façon suivante
\[
 (JC)^{-1}D_jD_{j-1}^{-1}(JC) = \frac{\d_j}{\d_{j-1}}D_{n+3}D_j^{-1}D_{j-1}D_{n+3}^{-1}.
\]
Si les matrices $D_{n+3}$ et $C$ v\'erifient
\[
 D_{n+3}^{-1}=JC,
\]
alors la relation~\eqref{C1'}$_j$ est \'equivalente \`a
\[
 \frac{1}{\d_{j-1}}{D_{j-1}}^2 = \frac{1}{\d_j}{D_j}^2,
\]
et on obtient ainsi l'\'equivalence annonc\'ee.
\end{proof}

On en d\'eduit donc que pour tout jeu de directions orient\'ees $D\in\D^n$, les syst\`emes fuchsiens dont la monodromie soit la classe de la repr\'esentation $\rho_D:\pi_1\left(\P\ssm S(t),x_0\right)\to GL(2,\C)$ d\'efinie par $D$, et dont les singularit\'es et les exposants sont r\'eels v\'erifient automatiquement la condition de r\'ealit\'e~\ref{cond-sys-realite}. L'ensemble $\A^n_D$ est ainsi simplement l'ensemble des syst\`emes v\'erifiant les conditions~\ref{cond-sys-SdR} et~\ref{cond-sys-mono} et dont les singularit\'es sont r\'eelles et ordonn\'ees.

\section{Description par le syst\`eme de Schlesinger}
\label{section-description-sch}

On va maintenant utiliser des d\'eformations isomonodromiques par le syst\`eme de Schlesinger~\eqref{schlesinger} pour d\'ecrire une partie de l'ensemble $\A^n_D$. On v\'erifiera ensuite que cette partie convient, \ie qu'elle est en bijection avec l'ensemble $\X^n_D$ des disques minimaux \`a bord polygonal. Enfin, on montrera un r\'esultat de r\'egularit\'e pour cette description.

\subsection{Le choix d'une famille isomonodromique}
\label{section-def-ADt}

Soit un jeu de directions orient\'ees $D\in\D^n$. On fixe arbitrairement un point $t^0\in\pi^n$, et on consid\`ere un syst\`eme fuchsien ($A_0$) dont la monodromie est la classe de $\rho_D$ et dont la position des singularit\'es est donn\'ee par $t^0$. Un tel syst\`eme existe toujours, puisque pour les syt\`emes de taille $2\times2$, le probl\`eme de Riemann--Hilbert admet toujours une solution (on peut se reporter au livre d'Anosov et Bolibruch~\cite{AB}, ou \`a l'article de Beauville~\cite{Beauville} pour une pr\'esentation synth\'etique des r\'esultats connus sur le probl\`eme de Riemann--Hilbert). On peut toujours supposer que le syst\`eme est normalis\'e en l'infini, et qu'il v\'erifie la condition~\ref{cond-sys-SdR}. Soit $U\subset\B^n$ un voisinage simplement connexe du simplexe $\pi^n$, o\`u l'ensemble $\B^n$ est d\'efini par \eqref{def-B}. Les r\'esidus $\left( A_1^0,\ldots,A_{n+2}^0\right)$ du syst\`eme ($A_0$) sont une condition initiale du Schlesinger~\eqref{schlesinger}, qui est compl\`etement int\'egrable (th\'eor\`eme~\ref{thm-schlesinger}). On obtient donc ainsi une famille isomonodromique de syst\`emes fuchsiens $\left( A_D(t) , t\in U \right)$ d\'ecrite par le syst\`eme de Schlesinger, telle que $(A_D(t^0))=(A_0)$. Les conditions~\ref{cond-sys-SdR} et~\ref{cond-sys-mono} sont satisfaites par le syst\`eme ($A_0$), et sont conserv\'ees au cours de la d\'eformation. D'apr\`es la proposition~\ref{prop-sys-reel-2}, on en d\'eduit
\begin{equation}
\left( A_D(t) , t\in\pi^n \right)  \subset  \A^n_D.
\label{ADt}
\end{equation}
Tous les choix possibles pour la solution ($A_0$) du probl\`eme de Riemann--Hilbert induit de cette mani\`ere une famille isomonodromique de syst\`emes fuchsiens contenue dans l'ensemble  $\A^n_D$ (lorsque $t\in\pi^n$), et bien s\^ur tous les \'el\'ements de $\A^n_D$ appartiennent \`a une telle famille.

Consid\'erons \`a pr\'esent deux de ces familles $\left( A^1_D(t) , t\in U \right)$ et $\left( A^2_D(t) , t\in U \right)$. Pour chaque valeur $t\in U$, les syst\`emes fuchsiens $\left( A^1_D(t)\right)$ et $\left( A^2_D(t)\right)$ ont les m\^emes singularit\'es, les m\^emes exposants et la m\^eme monodromie, et leurs normalisations en l'infini sont identiques (donn\'ees par~\eqref{Ainfty}). D'apr\'es le lemme~\ref{lemme-unicite} et la proposition~\ref{prop-eq->sys}, ces deux syst\`emes d\'efinissent la m\^eme \'equation, que l'on note $\left( E_D(t)\right)$ : ils correspondent \`a des valeurs diff\'erentes du param\`etre $\xi$. \'Etant donn\'e que toute \'equation de $\E^n_D$ provient d'un syst\`eme de $ \A^n_D$, qui appartient lui-m\^eme \`a une famille isomonodromique~\eqref{ADt}, la famille isomonodromique d'\'equations fuchsiennes $\left( E_D(t) , t\in\pi^n \right)$ d\'ecrit ainsi enti\`erement l'ensemble $\E^n_D$, qui est donc param\'etr\'e par $t$~:
\[
 \E^n_D = \left( E_D(t) , t\in\pi^n \right).
\]
En fait, cette d\'ependance en $t$ est \'egalement donn\'ee par le syst\`eme de Garnier~\eqref{Gn}, mais on n'utilisera pas ce point de vue.

Finalement, on choisit arbitrairement une famille isomonodromique $\left( A_D(t) , t\in\pi^n \right)$, donn\'ee par une sous-vari\'et\'e d'une vari\'et\'e int\'egrale du syst\`eme de Schlesinger, et on note
\begin{equation}
 DY = A_D(x,t) Y,\qquad A_D(x,t)= \sum_{i=1}^{n+2} \frac{A_{D,i}(t)}{x-t_i}.
\tag{$A_D(t)$}\label{A_D(t)}
\end{equation}
Cette famille est en bijection avec l'ensemble $\X^n_D$, et permet de le d\'ecrire ainsi : d'apr\`es la proposition~\ref{prop-AnD}, pour tout $t\in\pi^n$, il existe une solution fondamentale $\Y_0(x,t)$ du syst\`eme $\left( A_D(t)\right)$ dont la premi\`ere ligne $\gxt$ constitue les donn\'ees de Weierstrass d'une immersion de $\X^n_D$, que l'on note $X_D(t) $, et on a
\[
 \X^n_D = \left( X_D(t) , t\in\pi^n \right).
\]
On note \'egalement $P_D(t)\in\p^n_D$ le bord polygonal du disque repr\'esent\'e par $X_D(t)$. La famille $\left( P_D(t) , t\in\pi^n \right)$ est exactement la famille des polygones de direction $D$ qui sont le bord d'au moins un disque minimal. L'objet du chapitre suivant est d'utiliser cette description par le syst\`eme de Schlesinger pour montrer qu'elle d\'ecrit enti\`erement l'ensemble $\p^n_D$. Remarquons que la solution fondamentale $\Y_0(x,t)$ est $M$-invariante, \ie que sa repr\'esentation de monodromie est ind\'ependante de $t$, puisqu'il s'agit de la repr\'esentation $\rho_D$.

\begin{rema}
On a vu que deux solutions diff\'erentes du probl\`eme de Riemann--Hilbert appartenant \`a l'ensemble $\A^n_D$ se distinguent par leur valeur du param\`etre $\xi$
\[
\xi =\sum_{i=1}^{n+2} t_i A^i_{12}.
\]
On peut montrer que si les matrices $\left( A_1(t),\ldots A_{n+2}(t)\right)$ sont solutions du syst\`eme de Schlesinger, alors le param\`etre $\xi (t)$ satisfait le syst\`eme de Pfaff
\[
 \frac{\partial \xi }{\partial t_i} = (\t_\infty-1) A^i_{12}(t), \qquad (i=1,\ldots,n)
\]
qui permet de d\'ecrire les relations entre le syst\`eme de Schlesinger et le syst\`eme de Garnier (voir~\cite{IKSY}).
\end{rema}

\begin{rema}
La preuve au chapitre suivant que la famille de polygones $\left( P_D(t) , t\in\pi^n \right)$ d\'ecrit l'ensemble $\p^n_D$ tout entier re-montrera \emph{a posteriori} que l'on a bien choisi la famille $\left( A_D(t) , t\in\pi^n \right)$, ainsi que les l'ensembles $\E^n_D$ et $\A^n_D$. On peut remarquer que pour r\'esoudre le probl\`eme de Plateau, on aurait p\^u ne pas utiliser d'\'equations fuchsiennes, et introduire directement l'ensemble $\A^n_D$ comme l'ensemble des syst\`emes fuchsiens satisfaisant les conditions~\ref{cond-sys-SdR}, \ref{cond-sys-mono} et \ref{cond-sys-realite}. On aurait p\^u alors seulement montrer qu'un tel syst\`eme d\'efinit bien une immersion $X$ qui appartient \`a $\X^n_D$ (\ie un r\'esultat analogue \`a la proposition~\ref{prop-cara-eq}), sans v\'erifier qu'on les obtient toutes ainsi, puisque c'est une cons\'equence du th\'eor\`eme~\ref{thm-F-surj}. Il y a plusieurs raisons pour lesquelles on n'a pas proc\'ed\'e ainsi. Tout d'abord, il n'est pas clair comment on peut d\'efinir directement un syst\`eme diff\'erentiel \`a partir seulement des donn\'ees de Weierstrass $\g$ : il y a beaucoup trop de choix possibles. On a choisi d'utiliser des syst\`emes fuchsiens, ce qui r\'eduit consid\'erablement le nombre de syst\`emes diff\'erentiels concern\'es, uniquement parce qu'on a d\'emontr\'e que l'unique \'equation de solution fondamentale $\g$ est fuchsienne. Par ailleurs, comment obtenir les conditions~\ref{cond-sys-SdR} et \ref{cond-sys-realite} sans utiliser d'\'equations ? En particulier, le comportement locale d'une immersion $X\in\X^n_D$ est donn\'ee par les exposants de l'\'equation associ\'ee, et non pas par les valeurs propres des matrices $A_i$ (ceci est l'objet de la premi\`ere partie de la d\'emonstration du th\'eor\`eme~\ref{thm-AnD}). L'utilisation d'\'equations fuchsiennes semble \^etre un d\'etour n\'ecessaire.
\end{rema}

\subsection{Singularit\'es mobiles des solutions r\'eelles du syst\`eme de Schlesinger}

Par le th\'eor\`eme~\ref{thm-sch-prop-P}, toute solution du syst\`eme de Schlesinger est m\'eromorphe dans le rev\^etement universel de l'ensemble $\B^n$. On \'etablit \`a pr\'esent un r\'esultat plus fort de r\'egularit\'e pour les solutions du syst\`eme de Schlesinger provenant d'une monodromie satisfaisant la condition \textbf{C1}, ou de mani\`ere \'equivalente, d\'efinissant une famille isomonodromique de syst\`emes fuchsiens v\'erifiant la condition de r\'ealit\'e~\ref{cond-sys-realite}. Ce r\'esultat s'appliquera donc \`a la famille $\left( A_D(t) , t\in\pi^n \right)$. En se restreignant aux syst\`emes v\'erifiant~\ref{cond-sys-realite}, on obtient un r\'esultat plus fort que celui de Garnier (puisqu'il ne parvient pas \`a exclure l'existence de p\^oles doubles en $t\in\pi^n$), et beaucoup plus simple \`a \'etablir.

\begin{prop}
Soit une solution $\left(A_1(t),\ldots,A_{n+2}(t) \right) $ du syst\`eme de Schlesinger~\eqref{schlesinger} d\'efinie dans un ouvert simplement connexe $U\subset \B^n$ contenant le simplexe $\pi^n$, et soit $(A_t)$ le syst\`eme fuchsien associ\'e. On suppose que les valeurs propres $\t_i^{\pm}$ ($i=1,\ldots,n+2$) sont r\'eelles ou conjugu\'ees, et que les valeurs propres $\t_\infty^{\pm}$ sont r\'eelles. S'il existe une valeur $t^0\in \pi^n$ telle que la monodromie du syst\`eme fuchsien~$(A_{t^0})$ v\'erifie la condition \textbf{C1}, alors pour tout $t\in \pi^n$ les matrices $A_1(t),\ldots,A_{n+2}(t)$ s'\'ecrivent sous la forme~\eqref{sys-reel}, et elles sont holomorphes en tout point de $\pi^n$.
\label{prop-Ai-holo}
\end{prop}

\begin{proof}
La premi\`ere partie de la proposition est \'evidente. On peut supposer que les matrices $A_i(t)$ sont \`a trace nulle. En effet, pour toutes constantes $k_1,\ldots,k_{n+2}\in\R$ les matrices
\[
 B_i(t) := A_i(t) + k_i\I_2 \qquad (i=1,\ldots,n+2)
\]
constituent \'egalement une solution du syst\`eme de Schlesinger, et sont encore sous la forme~\eqref{sys-reel}. Quitte \`a transformer ainsi les matrices $A_i(t)$, on peut donc supposer que pour les valeurs r\'eelles de $t$, elles s'\'ecrivent
\[
A_i(t) =
 \begin{pmatrix}
 a_i(t) &  b_i(t)e^{i\eta(t)}\\
 c_i(t)e^{-i\eta(t)} &  -a_i(t)
\end{pmatrix},
\]
o\`u la fonction $a_i(t)$ est \`a valeurs r\'eelles et o\`u les fonctions $b_i(t)$ et $c_i(t)$ sont \`a valeurs positives dans $\pi^n$. La matrice $A_i(t)$ est m\'eromorphe dans le rev\^etement universel de l'ensemble $\B^n$. On note $\frac{\t_i}{2}$ et $-\frac{\t_i}{2}$ ses valeurs propres ; elles sont ind\'ependantes de $t$ et pour tout $t$ r\'eel, on a
\[
 \frac{\t_i^2}{4}=a_i(t)^2+b_i(t)c_i(t).
\]
On en d\'eduit que $a_i(t)$ et le produit $b_i(t)c_i(t)$ sont born\'es dans $\pi^n$. Les fonctions $A^i_{11}(t)=a_i(t)$ et $A^i_{12}(t)A^i_{21}(t)=b_i(t)c_i(t)$ sont donc holomorphes en tout point de $\pi^n$.

Montrons que les fonctions $A^i_{12}(t)$ sont holomorphes dans $\pi^n$. Soit un point $t^0\in\pi^n$. On \'etudie le comportement en la variable $t_j$ au point $t_j^0$, les autres $t_k$, $k\neq j$, \'etant fix\'es en $t_k^0$. On raisonne par l'absurde, et on choisit $i\neq j$ tel que $A^i_{12}(t)$ ait un p\^ole d'ordre $p_i\geq1$ maximal en  $t_j=t_j^0$ (par rapport aux autres $A^l_{12}(t)$, $l\neq j$). D'apr\`es le syst\`eme de Schlesinger, on a
\[
 \frac{\partial A^i_{12}}{\partial t_j} = \frac{2A^j_{11}}{t_j-t_i}A^i_{12}-\frac{2A^i_{11}}{t_j-t_i}A^j_{12}.
\]
Comme les fonctions $\frac{A^j_{11}(t)}{t_j-t_i}$ et $\frac{A^i_{11}(t)}{t_j-t_i}$ sont holomorphes en tout point de $\pi^n$, on voit que $A^j_{12}(t)$ a en $t_j=t_j^0$ un p\^ole d'ordre au minimum $p_i+1$. Or ceci contredit l'\'equation
\[
 \frac{\partial A^j_{12}}{\partial t_j} = -\sum_{\substack{l=1\\ l\neq j}}^{n+2}\frac{\partial A^l_{12}}{\partial
 t_j}.
\]
Les fonctions $A^l_{12}(t)$, $l\neq j$, sont donc holomorphes en $t_j=t_j^0$, et vu l'\'equation pr\'ec\'edente, $A^j_{12}(t)$ l'est alors \'egalement. On proc\'ederait de m\^eme, et on aboutirait au m\^eme r\'esultat pour les fonctions $A^i_{21}(t)$.
\end{proof}

On peut donc d\'eduire de cette proposition que la solution $\left(A_{D,1}(t),\ldots,A_{D,n+2}(t) \right)$ du syst\`eme de Schlesinger associ\'ee \`a un jeu de directions orient\'ees $D\in\D^n$ est holomorphe dans un voisinage simplement connexe $U\subset\B^n$ du simplexe $\pi^n$. Ceci va simplifier l'\'etude de la r\'egularit\'e \`a l'int\'erieur de $\pi^n$ de la fonction <<~rapports des longueurs~>> $F_D(t)$, comme on le verra au chapitre~\ref{chapitre-longueur}.

\chapter{Rapports de longueurs des c\^ot\'es}
\label{chapitre-longueur}
\selectlanguage{francais}

On suppose toujours fix\'e un jeu de directions orient\'ees $D\in\D^n$. On a obtenu au chapitre pr\'ec\'edent que l'ensemble $\X^n_D$ des immersions conformes repr\'esentant des disques minimaux \`a bord polygonal de direction $D$ est une famille $\left(X_D(t), t\in\pi^n \right) $, param\'etr\'ee par le $n$-uplet de singularit\'es $t=(t_1,\ldots,t_n)$, qui sont \'egalement les ant\'ec\'edents par l'immersion $X_D(t):\C_+\to\R^3$ des sommets de leur bord polygonal. La d\'ependance en $t$ des immersions $X_D(t)$ est donn\'ee par le syst\`eme de Schlesinger~\eqref{schlesinger}. Pour chaque valeur de $t\in\pi^n$, les donn\'ees de Weierstrass $\gxt$ de l'immersion $X_D(t)$ constituent la premi\`ere ligne d'une matrice fondamentale de solutions $\Y_0(x,t)$, qui est $M$-invariante, du syst\`eme fuchsien $\left(A_D(t)\right)$. Comme cette solution fondamentale est d\'efinie \`a multiplication scalaire r\'eelle pr\`es, les immersions de $\X^n_D$ sont d\'efinies non seulement \`a translation pr\`es, mais \'egalement \`a homoth\'eties de rapport positif pr\`es. On a  not\'e $P_D(t)\in\p^n_D$ le bord polygonal du disque repr\'esent\'e par $X_D(t)$, et le but de ce chapitre est de montrer l'\'egalit\'e suivante
\[
\p^n_D = \left(P_D(t), t\in\pi^n \right).
\]
Un syst\`eme de coordonn\'ees sur $\p_D^n$ est donn\'e par $n$ rapports de longueurs de c\^ot\'es. Vu l'expression~\eqref{metrique-ff2} de la m\'etrique induite des immersions $X_D(t)$, les rapports de longueurs des c\^ot\'es de tout repr\'esentant du polygone $P_D(t)$ s'\'ecrivent
\[
 r_i(t)  = \frac{\displaystyle\int_{t_i}^{t_{i+1}} \left(|G(x,t)|^2 + |H(x,t)|^2\right) dx}{\displaystyle\int_0^1 \left(|G(x,t)|^2 + |H(x,t)|^2 \right) dx}
\]
($i=1,\ldots,n$). On d\'efinit la fonction <<~rapports des longueurs~>> $F_D(t)$ associ\'ee au jeu de directions orient\'ees $D$ ainsi
\[
 F_D:\pi^n \to \,]0,+\infty[\,^n, \qquad F_D(t)=(r_1(t),\ldots,r_n(t)).
\]
Le but de ce chapitre est donc d'\'etablir le th\'eor\`eme suivant, qui conclut la d\'emonstration du th\'eor\`eme~\ref{thm-Plateau}, et qui en est la partie la plus difficile.

\begin{theo}
 \'Etant donn\'e un jeu de directions orient\'ees $D\in\D^n$, la fonction <<~rapports des longueurs~>> $F_D:\pi^n \to \,]0,+\infty[\,^n$ est surjective.
\label{thm-F-surj}
\end{theo}

En quelque sorte, on veut montrer que les directions orient\'ees des c\^ot\'es sont param\'etr\'ees par la monodromie des syst\`emes fuchsiens, tandis que la position $t$ de leurs singularit\'es code les longueurs des c\^ot\'es. Mais la d\'etermination des longueurs par $t$ est moins explicite.

\bigskip

\`A la section~\ref{section-fct-rapport}, on commence par pr\'eciser le choix de la solution fondamentale $M$-invariante $\Y_0(x,t)$ de mani\`ere \`a pouvoir \'etudier sa d\'ependance en $t$ --- de nouveau, l'utilisation de syst\`emes au lieu d'\'equations simplifiera cette d\'etermination. On en d\'eduit ensuite, gr\^ace \`a la proposition~\ref{prop-Ai-holo}, que la fonction $F_D(t)$ est analytique r\'eelle dans le simplexe $\pi^n$ (proposition~\ref{prop-F-holo-dans-pi}).

On expose \`a la section~\ref{section-rec} la m\'ethode que l'on va suivre pour d\'emontrer le th\'eor\`eme~\ref{thm-F-surj}. Elle repose sur l'\'etude de la fonction $F_D(t)$ au bord du simplexe $\pi^n$ et sur une r\'ecurrence portant sur le nombre $n+3$ de c\^ot\'es des polygones. En identifiant les simplexes $\pi^n$ et $]0,+\infty[\,^n$, on d\'efinit une fonction
\[
 \w F_D:\,]0,+\infty[\,^n \to\,]0,+\infty[\,^n.
\]
Pour montrer que la fonction $F_D$ est surjective, on va montrer que la fonction $\w F_D$ est de degr\'e $1$, c'est-\`a-dire homotope \`a l'identit\'e. On \'etablit un r\'esultat de topologie (proposition~\ref{prop-F-surj}) qui nous permet de nous ramener \`a montrer que la fonction $\w F_D$ est continue et de degr\'e $1$ au bord de $]0,+\infty[\,^n$. Pour obtenir cela, il faut interpr\'eter la fonction $F_D\big|_{\partial \pi^n}$ en terme de nouvelles fonctions <<~rapports des longueurs~>> de dimension inf\'erieure. Le bord du simplexe $\pi^n$ est constitu\'e de simplexes de dimension inf\'erieure. Regardons par exemple ce qui se passe lorsque la singularit\'e $t_n$ tend vers $0$, \ie en un point de la face $\mathcal F$ du bord de $\pi^n$
\[
\mathcal F = \left\{ (t_1,\ldots,t_n)\in\R^n\ |\ t_1<\cdots<t_{n-1}<t_n=0\right\}\simeq\pi^{n-1}.
\]
Il para\^it naturel de s'attendre \`a ce que le $n$-i\`eme c\^ot\'e $a_n(t)a_{n+1}(t)$ du polygone $P_D(t)$ <<disparaisse>>, c'est-\`a-dire que le rapport de longueur $r_n(t)$ tende vers $0$. On montre de plus que lorsque $t_n=0$ et que $t'=(t_1,\ldots,t_{n-1})$ d\'ecrit le simplexe $\pi^{n-1}$, on obtient la famille de polygones $P_{D'}(t')$ d\'efinie par les directions orient\'ees
\[
 D'=\left(D_1,\ldots,D_{n-1},D_{n+1},D_{n+2},D_{n+3}\right) \in \D^{n-1}.
\]
Ceci signifie que la fonction $F_D(t)$ s'\'etend contin\^ument \`a la face $\mathcal F$ du bord de $\pi^n$ et que pour tout $t'\in\pi^{n-1}$, on a
\[
 F_D(t',0) = \left(F_{D'}(t'), 0\right).
\]
On g\'en\'eralise cette assertion \`a toutes les faces du simplexe $\pi^n$ : c'est la proposition~\ref{prop-F-cont-au-bord}, dont la d\'emonstration constitue la majeure partie de ce chapitre. On proc\`ede ensuite par r\'ecurrence, en faisant l'hypoth\`ese qu'au rang $n-1$, pour tout entier $k\leq n-1$ et tout jeu de directions orient\'ees $D'\in\D^k$, la fonction $\w F_{D'}:\,]0,+\infty[\,^k \to \,]0,+\infty[\,^k$ est de degr\'e $1$. Les propositions~\ref{prop-F-surj} et~\ref{prop-F-cont-au-bord} assurent l'h\'er\'edit\'e de l'hypoth\`ese de r\'ecurrence. L'initialisation au rang $n=1$ (cas d'un bord quadrilat\'eral) est imm\'ediate une fois que l'on a obtenu la proposition~\ref{prop-F-cont-au-bord}.

Les sections~\ref{section-pseudo-choc} et~\ref{section-cas-reel} sont consacr\'ees \`a la d\'emonstration de la proposition~\ref{prop-F-cont-au-bord}. La partie la plus difficile est d'obtenir la continuit\'e de la fonction $F_D(t)$ au bord, et non pas son interpr\'etation g\'eom\'etrique. \`A la section~\ref{section-pseudo-choc}, on reprend des r\'esultats g\'en\'eraux sur les singularit\'es fixes du syst\`eme de Schlesinger, que Garnier appelle les pseudo-chocs, c'est-\`a-dire en les points tels que $t_i=t_j$, $i\neq j$. Ces r\'esultats sont une partie plus connue du travail de Garnier~\cite{Garnier26}, et ont \'et\'e d\'evelopp\'es et g\'en\'eralis\'es par Sato, Miwa et Jimbo~\cite{SMJ}. On reprend ces r\'esultats en en approfondissant des aspects qui nous seront utiles pour \'etudier l'holomorphie de la fonction $F_D(t)$ en les pseudo-chocs. \`A la section~\ref{section-cas-reel}, on applique cette \'etude g\'en\'erale aux solutions particuli\`eres du syst\`eme de Schlesinger qui nous int\'eresse, c'est-\`a-dire au cas r\'eel. En rassemblant et en adaptant les r\'esultats de la section pr\'ec\'edente, on \'etablit la proposition~\ref{prop-F-cont-au-bord}.

\bigskip

La d\'emonstration propos\'ee dans ce chapitre est tr\`es diff\'erente de celle de Garnier, m\^eme s'il utilise aussi le comportement de la famille de syst\`emes $\left(A_D(t), t\in\pi^n\right)$ au bord du simplexe $\pi^n$ et une r\'ecurrence sur le nombre de c\^ot\'es des polygones. Mais son hypoth\`ese de r\'ecurrence n'est pas la m\^eme, car il ne s'appuie pas sur un r\'esultat de topologie global tel que la proposition~\ref{prop-F-surj}. C'est pourquoi son \'etude est plus complexe. De plus, Garnier n'introduit pas la fonction <<~rapports des longueurs~>>, il cherche d'abord \`a refermer les polygones $P_D(t)$ en faisant dispara\^itre une singularit\'e apparente <<~de trop~>> (remarque~\ref{rem-N=n}), puis \`a ajuster $n-1$ rapports de longueurs. Il \'ecrit ces conditions sous la forme d'un syst\`eme $\mathcal S^n$ \`a $n$ \'equations. Il montre que le syst\`eme $\mathcal S^n$ tend vers un syst\`eme analogue de dimension inf\'erieure $\mathcal S^{n-1}$ lorsque $t_n\to0$, et ceci passe en particulier par l'\'etude compliqu\'ee du syst\`eme de Garnier~\eqref{Gn} lorsque $t_n\to0$. Il proc\`ede ensuite par r\'ecurrence : il prolonge une solution du syst\`eme $\mathcal S^{n-1}$ en une solution du syst\`eme $\mathcal S^n$. Il utilise pour cela le th\'eor\`eme d'inversion locale, et doit montrer que le jacobien d'une fonction (qui est quasiment $F_D(t)$) n'est pas nul au bord et \`a l'int\'erieur du simplexe $\pi^n$. La d\'emonstration de ce dernier point est obscure, voire peu convaincante. De plus, l'initialisation de la r\'ecurrence pour le cas du quadrilat\`ere est tr\`es elliptique, comme l'attestent les propres travaux ult\'erieurs de Garnier : il \'etudie dans les ann\'ees 1950 et 1960 le cas du quadrilat\`ere dans les articles~\cite{Garnier51},~\cite{Garnier62a} et~\cite{Garnier62b}, et y soul\`eve plusieurs difficult\'es qu'il ne mentionne pas dans~\cite{Garnier28}.

\section{La fonction << rapports des longueurs >> \texorpdfstring{$F_D(t)$}{}}
\label{section-fct-rapport}

\subsection{D\'efinition}

Consid\'erons la famille isomonodromique de syst\`emes fuchsiens $\left(A_D(t), t\in\pi^n\right)$ associ\'ee \`a un jeu de directions orient\'ees $D\in\D^n$, que l'on a introduite au chapitre pr\'ec\'edent (section~\ref{section-def-ADt}). Pour \'etudier la d\'ependance en $t$ de la solution fondamentale $\Y_0(x,t)$, on va la comparer \`a une solution que l'on conna\^it mieux, la solution fondamentale canonique en l'infini $\Y_\infty(x,t)$. 

Pour tout $t\in\pi^n$, comme le syst\`eme $\left(A_D(t)\right)$ est normalis\'e en l'infini, il admet une unique matrice fondamentale de solutions canonique en l'infini de la forme
\[
 \mathbf Y_\infty(x,t)= R_\infty\left(\frac1x,t\right) x^{-L_\infty}, \qquad \text{avec } L_\infty=A_\infty=
  \left(1-\tfrac{\t_\infty}{2}\right)
 \begin{pmatrix}
 1 & 0\\
 0 & -1
 \end{pmatrix},
\]
o\`u la matrice $R_\infty(w,t)$ est holomorphe en $w=0$ et v\'erifie $R_\infty(0,t)=\I_2$. D'apr\`es le th\'eor\`eme~\ref{thm-schlesinger}, cette solution est $M$-invariante. De plus, comme la partie principale en $x=\infty$ est ind\'ependante de $t$~:
\[
 \mathbf Y_\infty(x,t) \sim x^{-L_\infty},
\]
la d\'ependance en $t$ de la solution fondamentale $\Y_\infty(x,t)$ est enti\`erement d\'etermin\'ee par la d\'ependance en $t$ de la matrice
\[
 A_D(x,t) =\sum_{i=1}^{n+2}\frac{A_{D,i}(t)}{x-t_i}.
\]
On a vu au chapitre pr\'ec\'edent que les matrices $A_{D,i}(t)$ sont holomorphes en tout point $t\in \pi^n$ (proposition~\ref{prop-Ai-holo}). Il existe donc un ouvert simplement connexe $U$ de l'ensemble $\B^n$ qui contient le simplexe $\pi^n$ tel que les matrices $A_{D,i}(t)$ sont holomorphes dans $U$. On obtient donc le lemme suivant.

\begin{lemm}
La solution fondamentale $\Y_\infty(x,t)$ est holomorphe dans tout ouvert simplement connexe de l'ensemble $\left( \P\times U\right) \ssm S$.

Pour tout $i=1,\ldots,n+2$, il existe une matrice $C_i\in GL(2,\C)$ ind\'ependante de $t$ telle que
\[
 \Y_\infty(x,t) = R_i(x,t)(x-t_i)^{L_i} \cdot C_i,
\]
o\`u la matrice $R_i(x,t)$ est holomorphe et inversible dans un voisinage de l'hypersurface $x=t_i$ de $\P\times U$. La matrice $R_i(x,t)$ se prolonge analytiquement le long de toute courbe de $\P\times U$ ne coupant aucune des hypersurfaces $x=t_j$ ($j\neq i$).
\label{lemme-Yinfty-en-ti}
\end{lemm}

On rappelle que l'ensemble $S\subset \P\times U$ est l'ensemble des singularit\'es de la famille de syst\`emes $(A_D(t), t\in U)$
\[
 S = \bigcup_{t\in U} S(t)\times \{t\}
\]
avec
\[
 S(t) = \{t_1,\ldots,t_{n+3}\}.
\]

\begin{proof}
La premi\`ere partie du lemme en \'evidente. Pour tout $i=1,\ldots,n+2$, au voisinage de la singularit\'e $x=t_i$, il existe par la proposition~\ref{prop-sol-cano} des matrices fondamentales de solutions de la forme
\[
 R_i(x,t)(x-t_i)^{L_i},
\]
o\`u la matrice $R_i(x,t)$ est holomorphe en $x$ au point $x=t_i$ et $R_i^0(t):=R_i(x,t)\big|_{x=t_i}$ est inversible et v\'erifie
\[
 A_{D,i}(t) = R_i^0(t)L_iR_i^0(t)^{-1}.
\]
Comme la matrice $A_{D,i}(t)$ est holomorphe dans $U$, il existe des matrices $R_i^0(t)\in GL(2,\C)$ qui diagonalisent $A_{D,i}(t)$ et qui soient holomorphes dans $U$. On en d\'eduit que la matrice $R_i(x,t)$, d\'efinie par une condition initiale $R_i^0(t)$ holomorphe, est holomorphe au voisinage de l'hypersurface $x=t_i$ de $\P\times U$.

A priori, la matrice de connexion entre les matrices fondamentales $R_i(x,t)(x-t_i)^{L_i}$ et $\Y_\infty(x,t)$ d\'epend de $t$. Comme la matrice de monodromie $M_i(\Y_\infty)$ de la solution fondamentale $\Y_\infty(x,t)$ est ind\'ependante de $t$, il existe une matrice $C_i\in GL(2,\C)$ telle que
\[
 M_i(\Y_\infty)=C_i^{-1} e^{2i\pi L_i}C_i.
\]
Alors, les solutions fondamentales $R_i(x,t)(x-t_i)^{L_i}$ et $\Y_\infty(x,t)\cdot C_i^{-1}$ ont la m\^eme matrice de monodromie $e^{2i\pi L_i}$ au point $x=t_i$, qui est diagonale et non scalaire. On montre facilement que ceci implique qu'il existe une matrice diagonale $\Delta_i(t)$ inversible et holomorphe dans $U$ telle que
\begin{align*}
 \Y_\infty(x,t) &= R_i(x,t)(x-t_i)^{L_i}\cdot \Delta_i(t)\cdot C_i\\
&= R_i(x,t)\Delta_i(t)(x-t_i)^{L_i}\cdot C_i,
\end{align*}
et la matrice $R_i(x,t)\Delta_i(t)$ convient.
\end{proof}

La solution fondamentale $\Y_0(x,t)$ dont la premi\`ere ligne constitue les donn\'ees de Weierstrass d'un disque minimal \`a bord polygonal est $M$-invariante, puisque sa repr\'esentation de monodromie est engendr\'ee par les matrices $M_i$ de la condition~\ref{cond-sys-mono}. La proposition~\ref{prop-sol-M-inv} permet de caract\'eriser l'ensemble des matrices fondamentales de solutions qui sont $M$-invariantes, sous r\'eserve que la monodromie des syst\`emes $\left(A_D(t)\right)$ soit irr\'eductible --- et c'est bien le cas, car les directions $D_i$ ne sont pas toutes coplanaires. On a donc
\begin{equation}
 \Y_0(x,t) = \mu(t)\Y_\infty(x,t)\cdot C,
 \label{(Y0-Yinfty)}
\end{equation}
o\`u la matrice inversible $C$, qui est ind\'ependante de $t$, est une matrice de conjugaison entre les matrices de monodromie de la solution $\Y_\infty(x,t)$ et les matrices $M_i$, et o\`u la fonction $\mu :U\to \C^*$ est holomorphe. Comme la solution $\Y_0(x,t)$ est d\'efinie \`a multiplication pr\`es par une fonction r\'eelle de $t$ jamais nulle, on peut la choisir comme suit.

\begin{lemm}
Il existe une matrice $C_0\in GL(2,\R)$ ind\'ependante de $t$ telle que la premi\`ere ligne de la solution fondamentale $\Y_0(x,t)=\Y_\infty(x,t)\cdot C_0$ constitue les donn\'ees de Weierstrass d'une immersion appartenant \`a $\X^n_D$.
\label{lemme-realite-Yinfty}
\end{lemm}

\begin{proof}
Par les propri\'et\'es de r\'ealit\'e du syst\`eme $\left(A_D(t)\right)$, pour tout $t\in\pi^n$, la solution fondamentale $\Y_\infty(x,t)$ est \`a valeurs r\'eelles d\`es que $x\in]-\infty,t_1[$. Si, par souci de simplicit\'e, on choisit d\'efinitivement une position du rep\`ere de $\R^3$ telle que la direction $D_{n+3}$ est dirig\'ee par le second vecteur de base $e_2$, alors on a vu au chapitre~\ref{chapitre-equ-fu} que la premi\`ere ligne de la solution fondamentale $\Y_0(x,t)$ est r\'eelle ou purement imaginaire d\`es que $x\in]-\infty,t_1[$. Quitte \`a inverser l'orientation de $D_{n+3}$, on peut supposer qu'elle est r\'eelle, et on en conclut donc que
\[
 \forall t\in\pi^n \quad \mu(t)C\in GL(2,\R).
\]
En particulier, les \'el\'ements de la matrice $C$ ont tous le m\^eme argument ; il existe donc un nombre r\'eel $\varphi$ tel que la matrice $C_0 :=e^{i\varphi}C$ soit dans $GL(2,\R)$. Alors $\mu_0(t):= e^{-i\varphi}\mu(t)$ est \`a valeurs r\'eelles dans $\pi^n$, et les solutions fondamentales $\Y_0(x,t)=\mu_0(t)\Y_\infty(x,t)\cdot C_0$ et $\Y_\infty(x,t)\cdot C_0$ d\'efinissent des surfaces minimales homoth\'etiques. On peut donc supposer que $\Y_0(x,t)=\Y_\infty(x,t)\cdot C_0$.
\end{proof}

Remarquons que l'expression obtenue dans le lemme pr\'ec\'edent pour les donn\'ees de Weierstrass est beaucoup plus simple que celle obtenue par Garnier \`a partir de solutions d'\'equations fuchsiennes.

Notons, pour une matrice $\Y \in M(2,\C)$,
\[
 L_1(\Y) := \sqrt{|y_1|^2 + |z_1|^2}
 \qquad \text{ o\`u }
 \Y =\begin{pmatrix}
         y_1 & z_1\\
         y_2 & z_2\\
     \end{pmatrix}
\]
(bien que ce ne soit pas une norme). Alors les longueurs des c\^ot\'es du disque minimal d\'efini par la solution fondamentale $\Y_0(x,t)$ sont donn\'ees, pour tout $t\in\pi^n$, par
\[
\ell_i(t) = \int_{t_i}^{t_{i+1}} L_1\left(\Y_\infty(x,t)\cdot C_0\right)^2
 dx \ \in \,]0,+\infty[
\]
$(i=1,\ldots,n+1)$. Elles sont bien d\'efinies de par leur interpr\'etation g\'eom\'etrique, mais aussi parce que les exposants en $x=t_i$ du syst\`eme fuchsien $\left(A_D(t)\right)$ sont strictement sup\'erieurs \`a $-1/2$. De plus, les fonctions $\ell_i(t)$ ne peuvent s'annuler dans $\pi^n$, car alors la premi\`ere ligne de la solution $\Y_0(x,t)$ serait nulle pour tout $x$ dans l'intervalle $]t_i,t_{i+1}[$, ce qui est impossible. Les rapports de longueurs s'\'ecrivent donc, pour tout $i=1,\ldots,n$,
\begin{equation}
 r_i(t)  = \frac{\displaystyle\int_{t_i}^{t_{i+1}} L_1\left(\Y_0(x,t)\right)^2 dx}{\displaystyle\int_0^1 L_1\left(\Y_0(x,t)\right)^2 dx} =  \frac{\displaystyle\int_{t_i}^{t_{i+1}} L_1\left(\mathbf Y_\infty(x,t)\cdot C_0\right)^2 dx}{\displaystyle\int_0^1 L_1\left(\mathbf Y_\infty(x,t)\cdot C_0\right)^2 dx}.
\label{rapport-Yinfty}
\end{equation}

\subsection{Holomorphie}

On veut \'etendre la fonction $F_D(t)$ en une fonction holomorphe dans un voisinage du simplexe $\pi^n$. Pour cela, il faut obtenir les rapports $r_i(t)$ par l'int\'egration de fonctions holomorphes en $t$, c'est-\`a-dire, en particulier, ne comportant pas de module. La solution fondamentale $\mathbf Y_0(x,t)$ a des propri\'et\'es de r\'ealit\'e qui permettent de se d\'ebarrasser des modules dans l'expression~\eqref{rapport-Yinfty}. Gr\^ace au lemme~\ref{lemme-realite-Yinfty}, on a choisi des donn\'ees de Weierstrass qui h\'eritent \`a la fois des propri\'et\'es de r\'ealit\'e de la solution $\Y_0(x,t)$, et de la r\'egularit\'e de la solution $\Y_\infty(x,t)$, puisqu'on a p\^u \'eliminer la d\'ependance en $t$ due \`a la fonction $\mu(t)$ dans l'expression~\eqref{(Y0-Yinfty)} de $\Y_0(x,t)$. On obtient ainsi une expression de la fonction $F_D(t)$ qui sera aussi utile pour l'\'etude en les pseudo-chocs.

\begin{prop}
Soit un jeu de directions orient\'ees $D\in\D^n$. Il existe un ouvert simplement connexe de l'ensemble $\B^n$ contenant $\pi^n$ et contenu dans $U$, que l'on note encore $U$, et une fonction $\underline F_D:U\to\C^n$ holomorphe dans $U$ qui prolonge la fonction <<~rapports des longueurs~>> $F_D:\pi^n\to\,]0,+\infty[\, ^n$ :
\[
 \underline F_D \big|_{\pi^n} = F_D.
\]
\label{prop-F-holo-dans-pi}
\end{prop}

\begin{proof}
On fixe $i\in\{1,\ldots,n+1\}$. Consid\'erons une matrice $S_i\in SU(2)$ qui soit un relev\'e d'une rotation envoyant la direction $D_i$ sur le second axe de coordonn\'ees. Alors la premi\`ere ligne de la solution fondamentale $\mathbf Y_0(x,t)\cdot S_i$ constitue les donn\'ees de Weierstrass d'une surface minimale bord\'ee par un polygone dont le $i$-\`eme c\^ot\'e est parall\`ele au second axe de coordonn\'ees. On a vu qu'alors cette premi\`ere ligne est r\'eelle ou purement imaginaire lorsque $x\in]t_i,t_{i+1}[$. On peut choisir $S_i$ telle qu'elle soit r\'eelle. Consid\'erons la solution fondamentale
\begin{equation}
 \Y_i(x,t): = \Y_0(x,t)\cdot S_i=\Y_\infty(x,t)\cdot C_0\cdot S_i.
\label{def-Y_i}
\end{equation}
La premi\`ere ligne $( g_i(x,t),h_i(x,t))$ de la solution $\Y_i(x,t)$ est donc \'egalement r\'eelle lorsque $x\in]t_i,t_{i+1}[$. Comme $S_i\in SU(2)$, on a pour tout $t\in\pi^n$ et tout $x\in]t_i,t_{i+1}[$
\[
  L_1\left( \Y_0(x,t)\right)^2
    = L_1\left( \Y_0(x,t)\cdot S_i\right)^2
    =g_i(x,t)^2 + h_i(x,t)^2.
\]
On a donc
\begin{equation}
 \ell_i(t) =\displaystyle\int_{t_i}^{t_{i+1}}\left(g_i(x,t)^2 + h_i(x,t)^2\right)dx,
 \label{r_i-par-Y_i}
\end{equation}
et $r_i(t) = \ell_i(t) / \ell_{n+1}(t)$.

On peut donc \'etendre la fonction $F_D(t)$ \`a l'ouvert $U$. En effet, quitte \`a diminuer l'ouvert $U$, on supposer que pour tout $t$ dans $U$, pour tout $i=1,\ldots,n$, les points $t_j$ ($j\neq i,i+1$) n'appartiennent pas au segment de droite limit\'e par $t_i$ et $t_{i+1}$. On peut donc toujours calculer les int\'egrales pr\'ec\'edentes le long des segments joignant $t_i$ et $t_{i+1}$. Alors pour tout $t\in U$ on a
\[
 \ell_i(t) =(t_{i+1}-t_i)\displaystyle\int_0^1 \left(g_i(t_i+\xi(t_{i+1}-t_i),t)^2 + h_i(t_i+\xi(t_{i+1}-t_i),t)^2\right)d\xi.
\]
Pour tout $t\in U$, la fonction $\ell_{n+1}(t)$ n'est jamais nulle, vu que ceci forcerait les fonctions $g_{n+1}(\cdot,t)$ et $h_{n+1}(\cdot,t)$ \`a \^etre identiquement nulles sur l'intervalle $]0,1[$.

Montrons que les fonctions $\ell_i(t)$ ($i=1,\ldots,n+1$) sont holomorphes en un point $t^0\in\pi^n$. D'apr\`es le lemme~\ref{lemme-Yinfty-en-ti}, comme la matrice $ C_iC_0S_i$ est ind\'ependante de $t$, la fonction
\[
 \mathcal G_i(\xi,t) :=g_i\left(t_i+\xi(t_{i+1}-t_i),t\right)^2 + h_i\left( t_i+\xi(t_{i+1}-t_i),t\right) ^2
\]
est holomorphe en $t$ au point $t=t^0$ pour tout $\xi$ fix\'e, $0<\xi<1$, et donc il suffit de la dominer par une fonction int\'egrable ind\'ependante de $t$, pour tout $t$ dans un voisinage de $t^0$. Soit $\e>0$ tel que la boule
\[
 B_\e\left(t^0\right) = \left\lbrace t\in\C^n\ |\ \forall i=1,\ldots,n \quad |t_i-t_i^0|<\e\right\rbrace
\]
soit contenue dans l'ouvert $U$. On scinde l'intervalle d'int\'egration
\[
 \ell_i(t) =\ell_i^-(t) + \ell_i^+(t),
\]
avec
\[
 \ell_i^-(t)= (t_{i+1}-t_i)\displaystyle\int_0^\frac12 \mathcal G_i(\xi,t)d\xi \quad \text{ et } \quad \ell_i^+(t)= (t_{i+1}-t_i)\displaystyle\int_\frac12^1 \mathcal G_i(\xi,t)d\xi.
\]
Consid\'erons la fonction $\ell_i^-(t)$. Il faut choisir $\e$ tel que
pour tout $\xi\in\, [0,\frac12]$ et pour tout $t\in
B_\e\left(t^0\right)$, la quantit\'e $\xi(t_{i+1}-t_i)$ soit contenue
dans un disque centr\'e en $0$ de rayon $\eta_i$ ind\'ependant de
$\xi$ et de $t$ et qui ne contienne aucune des valeurs singuli\`eres
$t_j-t_i$, $j\neq i$. On n'entre pas dans les d\'etails de calculs ;
si on suppose que $\e<(t_{i+1}^0-t_i^0)/6$ ($i=1,\ldots,n$), alors
\[
 \eta_i= \frac23\left|t_{i+1}^0-t_i^0\right|
\]
convient. Toujours par le lemme~\ref{lemme-Yinfty-en-ti} et parce que la matrice $C_0S_i C_i^{-1}$ est ind\'ependante de $t$, les fonctions $g_i(x,t)$ et $h_i(x,t)$ sont au voisinage de $x=t_i$ des combinaisons lin\'eaires \`a c\oe fficients ind\'ependants de $t$ de fonctions de la forme
\[
 (x-t_i)^{-\frac{\t_i}{2}}\varphi_i(x-t_i,t) \ \text{ et } \ (x-t_i)^{\frac{\t_i}{2}}\psi_i(x-t_i,t)
\]
o\`u les fonctions $\varphi_i(y,t)$ et $\psi_i(y,t)$ sont holomorphes en $t\in U$ et en $y$ tant que $y\neq t_j-t_i$ ($j\neq i$). Ces fonctions $\varphi_i(y,t)$ et $\psi_i(y,t)$ sont donc born\'ees pour tout $y$ tel que $|y|<\eta_i$ et pour tout $t\in B_\e\left(t^0\right)$. Il existe donc des constantes $K_0,\ K_1,\ K_{-1}>0$ telles que pour tout $\xi\in\, [0,\frac12]$ et tout $t\in B_\e(t^0)$, on ait
\begin{align*}
 |\mathcal G_i(\xi,t)| &\leq K_0 + K_{-1}|t_{i+1}-t_i|^{-\t_i}\xi^{-\t_i} + K_1|t_{i+1}-t_i|^{\t_i}\xi^{\t_i}\\
        &\leq K_0 + K_{-1}(2\eta_i)^{-\t_i}\xi^{-\t_i} + K_1(2\eta_i)^{\t_i}\xi^{\t_i}.
\end{align*}
On obtient donc que la fonction $\ell_i^-(t)$ est holomorphe au point $t^0$. On proc\`ederait de m\^eme pour $\ell_i^+(t)$. La fonction $\underline F_D(t)$ est donc holomorphe en tout point du simplexe $\pi^n$. Elle est donc holomorphe dans un ouvert simplement connexe $\underline U$ de $\B^n$ contenant $\pi^n$, et on appelle toujours $U$ l'intersection $U\cap \underline U$.
\end{proof}

\section{La d\'emonstration par r\'ecurrence}
\label{section-rec}

\subsection{La proposition fondamentale}

D'apr\`es la proposition~\ref{prop-F-holo-dans-pi}, la fonction $F_D:\pi^n \to\,]0,+\infty[\,^n$ est continue dans $\pi^n$. Par identification naturelle des simplexes $\pi^n$ et $]0,+\infty[\,^n$ (identification que l'on va pr\'eciser dans la suite), on obtient une fonction continue 
\[
 \w F_D:\,]0,+\infty[\,^n \to\,]0,+\infty[\,^n.
\]
Pour montrer que la fonction $F_D$ est surjective, on va montrer que la fonction $\w F_D$ est homotope \`a l'identit\'e, c'est-\`a-dire de degr\'e $1$. Le point essentiel pour \'etablir ce r\'esultat est l'\'etude du comportement de $F_D$ au bord du simplexe $\pi^n$. On commence par \'etablir la proposition suivante, qui, une fois obtenu ce comportement au bord, nous permettra de conclure gr\^ace \`a un raisonnement par r\'ecurrence.

\begin{prop}
 Soient un ensemble convexe et compact $K$ de $\R^n$, et une fonction continue $f :K \to K$ telle que $f(\partial K) \subset \partial K$. Si la fonction $f\big|_{\partial K} : \partial K\to \partial K$ est de degr\'e $1$, alors la fonction $f:K \to K$ est de degr\'e $1$ dans $K$.
\label{prop-F-surj}
\end{prop}

La notion standard de degr\'e concerne les applications diff\'erentiables (ou seulement continues) entre vari\'et\'es sans bord. On peut n\'eanmoins l'\'etendre aux vari\'et\'es ayant un bord, \`a la condition que les applications pr\'eservent le bord. Cependant, la notion importante ici est le fait que pour une application continue, \^etre de degr\'e $1$ est \'equivalent \`a \^etre homotope \`a l'identit\'e : on veut montrer que la fonction $F_D$ pr\'eserve la structure simpliciale du bord des domaines $\pi^n$ et $]0,+\infty[\,^n$ (apr\`es identification naturelle). 

\begin{proof}
On commence par montrer cette proposition lorsque le convexe compact $K$ co\"incide avec la boule unit\'e ferm\'ee $\bar B:=\bar B_1(0)$ de $\R^n$ pour la norme euclidienne $\Vert\cdot\Vert$. On proc\`ede par d\'eformations homotopiques. Par hypoth\`ese, il existe une fonction continue $h:[0,1]\times \partial B \to \partial B$ telle que 
\begin{align*}
 h\left(0,\cdot\right) &= f\big|_{\partial B}\\
 h\left(1,\cdot\right) &= \id_{\partial B}.
\end{align*}
On va construire une fonction continue $H:[0,1]\times\bar B \to\bar B$ telle que 
\begin{equation}
\begin{split}
  H\left(0,\cdot\right) &= \id_B\\
 H\left(1,\cdot\right) &= f.
\end{split}
\label{homotopie}
\end{equation}
On proc\`ede en deux \'etapes, suivant la valeur de $t$ (voir figure~\ref{fig-degre}). On d\'efinit tout d'abord la fonction $H(t,\cdot):\bar B \to\bar B$ pour $t$ fix\'e, $0<t\leq\frac12$, en faisant une r\'etractation de $f$ de $\bar B_t(0)$ dans $\bar B_t(0)$, puis en la transformant au bord par la fonction $h$ pour obtenir l'identit\'e.

\begin{figure}[!ht]
\centering
\begin{pspicture}(0,0)(16,8)
\rput(8,4){\includegraphics[height=6cm]{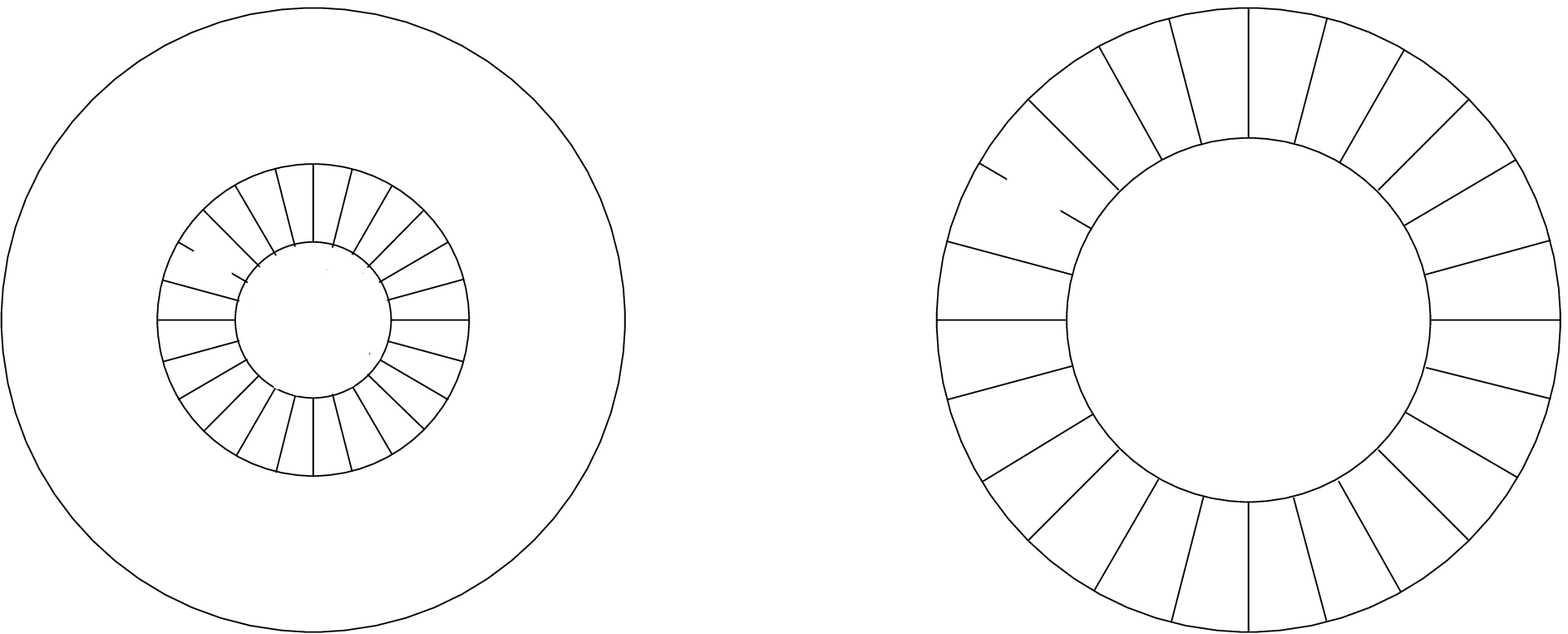}}
\rput(8,7){$B$}
\psline{->}(7.4,6.9)(5.6,6.2)
\psline{->}(8.6,6.9)(10.4,6.2)
\rput(3.5,0.4){$0<t\leq\frac 12$}
\rput(3.5,4){$f$}
\rput(1.7,5.5){$\id$}
\rput(2.55,4.5){$h$}
\rput(7.3,5.4){$B_t(0)$}
\psline{->}(6.8,5.4)(4.2,4.3)
\rput(7.4,2.4){$B_{2t}(0)$}
\psline{->}(6.8,2.5)(4.9,3.5)
\rput(12.5,0.4){$\frac 12\leq t\leq1$}
\rput(12.5,4){$f$}
\rput(10.4,5.2){$h$}
\rput(9.3,1.4){$B_t(0)$}
\psline{->}(9.8,1.6)(11,3)
\end{pspicture}
\caption{La fonction $H(t,\cdot)$ suivant la valeur de $t$.}
\label{fig-degre}
\end{figure}

\noindent
Plus pr\'ecis\'ement, on pose
\begin{align*}
 \forall x\in B_t(0)                \quad H(t,x) &= tf\left(\frac xt\right)\\
 \forall x\in B_{2t}(0) \ssm B_t(0) \quad H(t,x) &= \Vert x\Vert \, h\left(\frac{\Vert x\Vert}{t}-1, \frac{x}{\Vert x\Vert}\right)\\
 \forall x\in \bar B \ssm B_{2t}(0) \quad H(t,x) &= x.
\end{align*}
De m\^eme, pour $\frac12\leq t\leq1$, on se contente de grossir et de tronquer le cas pr\'ec\'edent~:
\begin{align*}
 \forall x\in B_t(0)             \quad H(t,x) &= tf\left(\frac xt\right)\\
 \forall x\in \bar B \ssm B_t(0) \quad H(t,x) &= \Vert x\Vert \, h\left(\frac{\Vert x\Vert}{t}-1, \frac{x}{\Vert x\Vert}\right).
\end{align*}
La fonction $H$ v\'erifie la condition~\eqref{homotopie} et est continue en tout point de $[0,1]\times \bar B\ssm\{(0,0)\}$. Pour v\'erifier qu'elle est continue au point $(0,0)$ et que $H(0,0)=0$, il suffit de remarquer que pour tout $0< t\leq\frac12$, on a
\begin{align*}
 \forall x\in B_t(0)                \quad \Vert H(t,x)\Vert &\leq t\\
 \forall x\in B_{2t}(0) \ssm B_t(0) \quad \Vert H(t,x)\Vert &\leq \Vert x\Vert.
\end{align*}

Dans le cas g\'en\'eral, lorsque le convexe compact $K$ est quelconque, il existe un hom\'eomorphisme $\varphi: K\to \bar B$ qui envoie le bord de $K$ sur la sph\`ere $\partial B$. Alors en appliquant le r\'esultat qu'on vient d'\'etablir \`a la fonction
\[
 g :=\varphi\circ f\circ\varphi^{-1}:\bar B \to\bar B,
\]
on obtient que la fonction $f$ est de degr\'e $1$.
\end{proof}

On va montrer que la fonction $\w F_D$ s'\'etend contin\^ument au bord de $]0,+\infty[\,^n$ et que $\w F_D\left(\partial (\,]0,+\infty[\,^n)\right) \subset\partial (\,]0,+\infty[\,^n)$. On proc\`edera par r\'ecurrence pour obtenir que la fonction
\[
 \w F_D\big|_{\partial (\,]0,+\infty[\,^n)} :\partial (\,]0,+\infty[\,^n) \to \partial (\,]0,+\infty[\,^n)
\]
est de degr\'e $1$. Commençons par compactifier les simplexes $\pi^n$ et $]0,+\infty[\,^n$ dans $\overline\R^n=\left(\R\cup\{-\infty,+\infty\}\right)^n$ et par expliciter leur bord et la mani\`ere de les identifier. On \'ecrit
\begin{align*}
 ]0,+\infty[\,^n &= \left\{ (x_1,\ldots,x_n)\in\R^n\ |\ 0<x_n<x_{n-1}+x_n<\cdots<x_1+\cdots+x_n<+\infty\right\}\\
    & =\left\{ (x_1,\ldots,x_n)\in\R^n\ |\ 0<s_n(x)<\cdots<s_1(x)<+\infty\right\}
\end{align*}
o\`u on a pos\'e pour $i=1,\ldots,n$
\[
 s_i(x)=x_i+\cdots+x_n
\]
et $s_0(x)=+\infty$, $s_{n+1}(x)=0$. Les adh\'erences sont donn\'ees par
\begin{align*}
 &\overline{\pi^n} = \left\{ t\in\overline\R^n \ |\  -\infty\leq t_1 \leq\cdots\leq t_n\leq 0\right\}  \\
 &\overline{]0,+\infty[\,^n} = [0,+\infty]\,^n =\left\{ x\in\overline\R^n\ |\ 0\leq s_n(x)\leq\cdots\leq s_1(x)\leq +\infty\right\} .
\end{align*}
Les bords $\partial\pi^n$ et $\partial (\,]0,+\infty[\,^n)$ sont constitu\'es de simplexes de dimensions $0$ \`a $n-1$. On param\`etre ces simplexes de la façon suivante. Soit $\Delta := \{0,1\}^{n+1} \ssm \{(0,\ldots,0),(1,\ldots,1)\}$. On note $\delta= (\delta_0,\ldots,\delta_n)$ les \'el\'ements de $\Delta$. On a la r\'eunion disjointe de simplexes
\[
\partial\pi^n = \bigsqcup_{\delta\in\Delta}P^\delta,
\]
avec
\[
P^\delta = \left\lbrace (t_1,\ldots,t_n)\in \overline{\pi^n} \ |\ \forall i=0,\ldots,n \ \ t_i=t_{i+1}
\Leftrightarrow \delta_i=0\right\rbrace
\]
o\`u on note $t_{n+1}=0$ et $t_0=t_{n+3}=-\infty$. Pour tout $\delta \in \Delta$, on a un isomorphisme naturel
\[
\varphi_\delta : P^\delta \to \pi^{|\delta|}
\]
o\`u la dimension du simplexe est donn\'ee par
\[
 |\delta|=\displaystyle\sum_{i=0}^n\delta_i -1.
\]
Cet isomorphisme est obtenu en <<enlevant>> les composantes $t_i$ telles que $\delta_i=0$ ($i=1,\ldots,n$) et celles qui valent $-\infty$. De m\^eme
\[
\partial(\,]0,+\infty[\,^n) = \bigsqcup_{\delta\in\Delta}R^\delta,
\]
avec
\[
R^\delta = \left\{ (x_1,\ldots,x_n)\in [0,+\infty] \, ^n \ |\ \forall i=0,\ldots,n \ \ s_{i+1}(x)=s_i(x)
\Leftrightarrow \delta_i=0\right\}.
\]
On a \'egalement les isomorphismes
\[
\psi_\delta : R^\delta \to \,]0,+\infty[\,^{|\delta|}.
\]
De m\^eme, on note $D^\delta \in D^{|\delta|}$ le jeu de directions orient\'ees obtenu \`a partir de $D\in \D^n $ en <<enlevant>> les directions orient\'ees $D_i$ telles que $\delta_i=0$ ($i=0,\ldots,n$). Les deux directions orient\'ees $D_{n+1}$ et $D_{n+2}$ ne peuvent donc jamais disparaître. Grâce \`a la d\'efinition~\ref{def-Dn} de l'ensemble $\D^n$, on voit que le jeu de directions orient\'ees $D^\delta$ appartient bien \`a $\D^{|\delta|}$. Alors
\[
 F_{D^\delta} : \pi^{|\delta|} \to \,]0,+\infty[\,^{|\delta|}.
\]
Le but des sections suivantes va \^etre d'\'etablir la proposition fondamentale~:

\begin{prop}
 Pour tout $\delta \in \Delta$, la fonction <<~rapports des longueurs~>> $F_D(t)$ associ\'ee \`a un jeu de direction $D\in\D^n$ s'\'etend contin\^ument \`a la face $P^\delta$ de $\pi^n$ et
\begin{equation}
 F_D\big|_{P^\delta}=\psi_\delta^{-1} \circ F_{D^\delta} \circ \varphi_\delta.
\end{equation}
\label{prop-F-cont-au-bord}
\end{prop}

Pour tout $n\in\N^*$, on consid\`ere un hom\'eomorphisme
\[
 \Phi_n : \,]0,+\infty[\,^n \to \pi^n
\]
tel que pour tout $\d\in\Delta$ on ait
\[
 \Phi_n\left(R^\d\right) = P^\d.
\]
On pose alors
\[
 \w F_D:= F_D \circ\Phi_n, \qquad \w F_D:\,]0,+\infty[\,^n \to\,]0,+\infty[\,^n.
\]
\'Etant donn\'e les propositions~\ref{prop-F-surj} et~\ref{prop-F-cont-au-bord}, pour montrer que la fonction $\w F_D$ est de degr\'e $1$, on va faire une r\'ecurrence forte, et la bonne hypoth\`ese est~:

\vspace{0.3cm}
\noindent
\textbf{Hypoth\`ese de r\'ecurrence au rang $n$ :} pour tout $k=1,\ldots,n$, pour tout jeu de directions orient\'ees $D\in\mathcal D^k$ la fonction
\[
 \w F_D= F_D \circ\Phi_k, \qquad \w F_D:\,]0,+\infty[\,^k \to\,]0,+\infty[\,^k
\]
est de degr\'e $1$.

\vspace{0.3cm}
\noindent
Pour tout $\delta \in \Delta$, comme $|\delta|<n$, on obtient ainsi, gr\^ace \`a l'hypoth\`ese de r\'ecurrence au rang $n-1$, que $\w F_D\big|_{R^\delta}:R^\delta\to R^\delta$ est de degr\'e $1$, et on a donc
\[
 \w F_D\big|_{\partial (\,]0,+\infty[\,^n)} :\partial (\,]0,+\infty[\,^n) \to \partial (\,]0,+\infty[\,^n)
\]
est de degr\'e $1$. Par la proposition~\ref{prop-F-surj}, on peut alors en conclure que la fonction $\w F_D:\,]0,+\infty[\,^n\to\,]0,+\infty[\,^n$ est de degr\'e $1$, et l'h\'er\'edit\'e de la r\'ecurrence est \'etablie.

\subsection{Le cas du quadrilat\`ere (\texorpdfstring{$n=1$}{n=1})}

L'initialisation de la r\'ecurrence au rang $n=1$ est imm\'ediate \`a partir de la proposition~\ref{prop-F-cont-au-bord}. Dans ce cas, pour tout $D=(D_t,D_0,D_1,D_\infty)\in\D^1$, la fonction <<~rapports des longueurs~>>
\[
 F_D :\, ]-\infty,0[\,\to\,]0,+\infty[
\]
est le rapport de la longueur du premier c\^ot\'e (de direction $D_t$) sur la longueur du deuxi\`eme (de direction $D_0$). Ici, $\Delta=\{\delta_1,\delta_2\}$ avec $\delta_1=\{0,1\}$ et $\delta_2=\{1,0\}$, et $P^{\delta_1}=\{-\infty\}$, $R^{\delta_1}=\{+\infty\}$, $P^{\delta_2}=\{0\}$ et $R^{\delta_2}=\{0\}$. La proposition~\ref{prop-F-cont-au-bord} nous donne donc ce \`a quoi on pouvait raisonnablement s'attendre~:
\[
 \lim_{t\to0} F_D(t)=0 \quad \text{ et } \quad \lim_{t\to-\infty} F_D(t)=+\infty.
\]
On peut choisir
\[
 \Phi_1:\,]0,+\infty[\,\to\,]-\infty,0[, \quad \Phi_1(t)=-t
\]
c'est-\`a-dire
\[
 \w F_D :\,]0,+\infty[\,\to\,]0,+\infty[, \quad \w F_D(t) = F_D(-t).
\]
On en d\'eduit donc que la fonction $\w F_D$ est de degr\'e $1$ (cas particulier \'evident de la dimension $1$ de la proposition~\ref{prop-F-surj}). 

\begin{figure}[!ht]
\centering
\begin{pspicture}(0.5,0)(16.5,11)
\rput(8,6){\includegraphics[height=10cm]{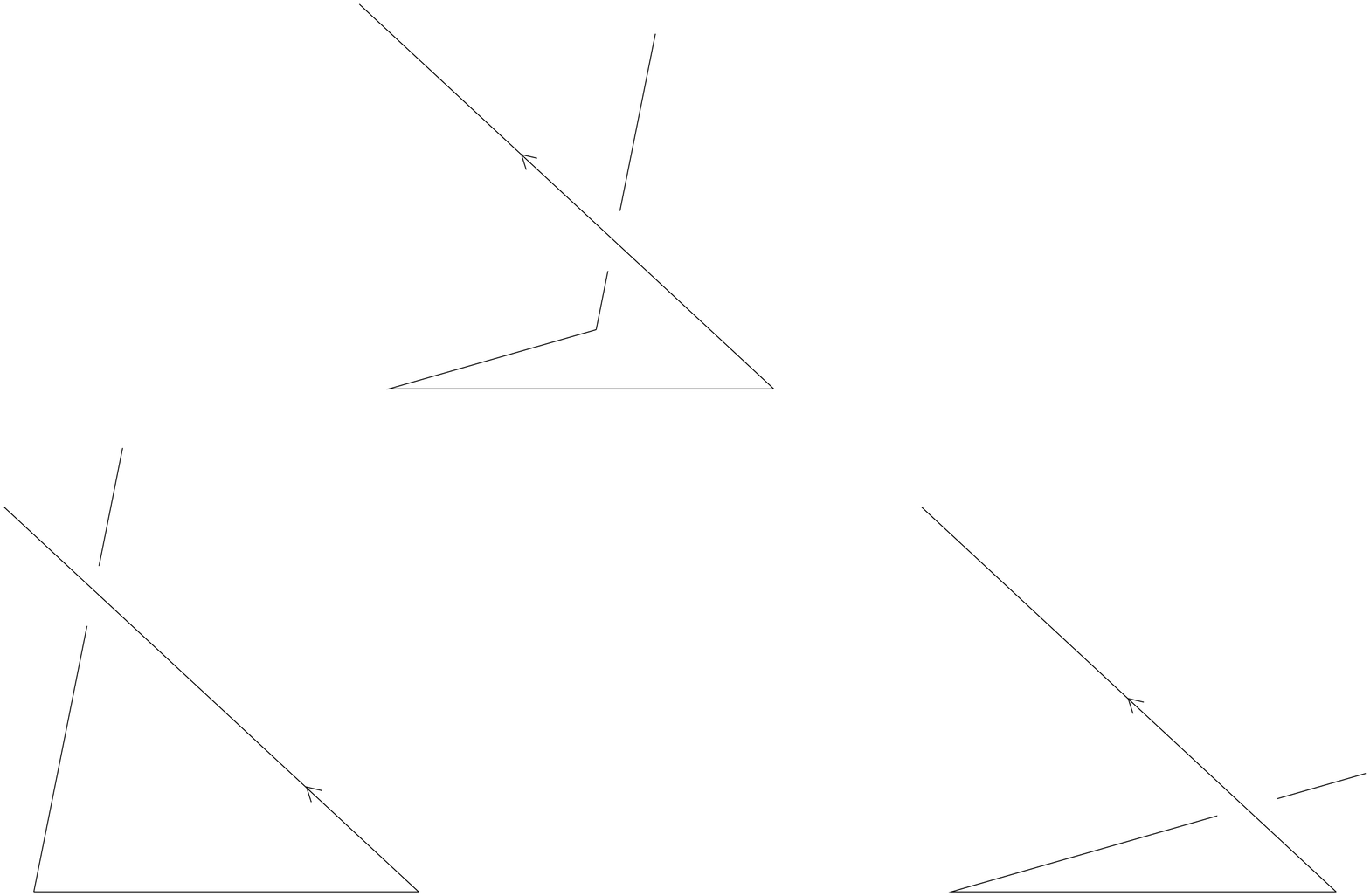}}
\rput(6.8,5.8){$-\infty<t<0$}
\rput(6.7,7.5){$a_t$}
\rput(4.4,6.5){$a_0$}
\rput(9.4,6.5){$a_1$}
\rput(6.4,10.8){$a_\infty$}
\rput(5.8,7.2){$D_t$}
\rput(7.6,6.85){$D_0$}
\rput(4.7,10.2){$D_1$}
\rput(7.8,9.5){$D_\infty$}
\rput(2.8,0.2){$t=0$}
\rput(0.4,0.9){$a_0$}
\rput(5.4,0.9){$a_1$}
\rput(1,5.6){$a_\infty$}
\rput(3.5,1.25){$D_0$}
\rput(2.6,3.7){$D_1$}
\rput(0.8,3.2){$D_\infty$}
\rput(13.2,0.2){$t=-\infty$}
\rput(10.8,0.9){$a_0$}
\rput(15.7,0.9){$a_1$}
\rput(14,4){$a_\infty$}
\rput(15.2,2.5){$D_t$}
\rput(14,1.25){$D_0$}
\rput(11.4,4.2){$D_1$}
\end{pspicture}
\caption{D\'eformation du quadrilat\`ere d\'efini par un jeu de directions orient\'ees}
\label{fig-quadri}
\end{figure}

On repr\'esente \`a la figure~\ref{fig-quadri} les variations lorsque $t\to-\infty$ et $t\to0$ du quadrilat\`ere $P_D(t)$ d\'efini par le jeu de directions orient\'ees $D$, et pour lequel le probl\`eme de Plateau admet une solution. On note $a_t=X(t)$, $a_0=X(0)$, $a_1=X(1)$ et $a_\infty=X(\infty)$ les sommets de ce quadrilat\`ere. Les sommets $a_0$ et $a_1$ ne peuvent pas disparaître au cours de la d\'eformation. Comme les quadrilat\`eres $\left(P_D(t), t\in\, ]-\infty,0[\right)$ sont d\'efinis \`a translation et homoth\'etie de rapport positif pr\`es, et comme la direction $D_0$ est fix\'ee, on peut supposer que la position des sommets $a_0$ et $a_1$ est fixe. 

Aux cas limites, lorsque $t=-\infty$ ou $t=0$, les donn\'ees de Weierstrass d'une surface minimale limit\'ee par un triangle ayant un sommet en l'infini sont des solutions d'une \'equation hyperg\'eom\'etrique.

\begin{rema}
 Si les directions orient\'ees $D_\infty$, $D_t$ et $D_0$ sont dans un m\^eme plan, et si la direction $D_1$ n'appartient pas ce plan, alors ces directions ne sont pas les directions d'un quadrilat\`ere de $\R^3$, et il n'existe aucune valeur de $t$ pour laquelle le quadrilat\`ere $P_D(t)$ <<~se referme~>>. Par contre, suivant l'orientation des directions $D_\infty$, $D_t$ et $D_0$, il peut exister une valeur de $t$ telle que la demi-droite $\left(a_t(t),-D_\infty\right)$ passe par le sommet $a_1(t)$ (qui devient donc aussi le sommet $a_\infty(t)$) : on obtient un triangle de $\R^3$.
\end{rema}

\subsection{Le changement de variables}

On va d\'etailler uniquement le cas des faces $P^\delta$ o\`u $\delta =(1,\ldots,1,0,\ldots,0)$ : on fixe un entier $p$, $1\leq p\leq n$, et on \'etudie la fonction $F_D(t)$ lorsque $t_p,t_{p+1},\ldots,t_n$ tendent vers $t_{n+1}=0$, les autres variables $t_1,\ldots,t_{p-1}$ demeurant \`a distance mutuelle sup\'erieure \`a un r\'eel strictement positif. De mani\`ere g\'en\'erale, on va noter par $\a$ les indices prenant les valeurs $1,\ldots,p-1,n+2$, et par $\b$ ceux variant entre $p$ et $n+1$. Pour tout $t\in\B^n$, on fait le changement de variables suivant
\begin{equation}
 \tau:=t_p, \qquad \nu_\b:=\frac{t_\b}{\tau} \quad (p\leq \b\leq n+1),
\label{chang-var-t}
\end{equation}
et on note $t'=(t_1,\ldots,t_{p-1})$ et $\nu=(1,\nu_{p+1},\ldots,\nu_n)$. Par abus de notation, on identifiera $\nu$ et $(\nu_{p+1},\ldots,\nu_n)$. En particulier, on dira que $\nu\in\B^{n-p}$ pour signifier que $(\nu_{p+1},\ldots,\nu_n)\in\B^{n-p}$. Alors
\[
 t=\left( t',\tau\cdot\nu\right) .
\]
\`A $t'\in\B^{p-1}$ et $\nu\in\B^{n-p}$ fix\'es, le $n$-uplet
$(t',\tau\cdot\nu)$ est dans $\B^n$ d\`es que $|\tau|$ est
suffisamment petit. On d\'efinit l'image $V$ de l'ouvert $U$ par le
changement de variables
\begin{equation}
 V:=\left\{ \left(t',\nu,\tau\right)\in\C^n \ |\ \left(t',\tau\cdot\nu\right)\in U\right\} .
\label{def-V}
\end{equation}
Dans le cas r\'eel, c'est-\`a-dire lorsque la variable $t$ est dans le simplexe $\pi^n$, la variable $t'$ est dans $\pi^{p-1}$ et la variable $\nu$ est dans le simplexe $\w\pi^{n-p}$ d\'efini par
\[
 \w\pi^k:=\left\{ (\nu_1,\ldots,\nu_k)\in\R^k\ |\ 0<\nu_k<\cdots<\nu_1<1 \right\} .
\]

On consid\`ere un voisinage simplement connexe $U'$ du simplexe $\pi^{p-1}$ contenu dans $\B^{p-1}$, et un voisinage simplement connexe $\w U$ du simplexe $\w\pi^{n-p}$ contenu dans $\B^{n-p}$ tels que pour tout $\left(t',\nu\right)\in U'\times\w U$, il existe $\tau\in\C^*$ tel que le $n$-uplet $\left(t',\nu,\tau\right)$ soit dans l'ensemble $V$. On suppose de plus que l'ouvert $\w U$ est born\'e : ceci est possible puisque le simplexe $\w\pi^{n-p}$ l'est.

La proposition suivante rassemble les r\'esultats que l'on va \'etablir dans les deux sections suivantes. Elle donne le comportement de la fonction $F_D$ en les variables $\left(t',\nu,\tau\right)$ aux points $\left(t'^0,\nu^0,0\right)$, avec $\left(t'^0,\nu^0\right)\in U'\times\w U$. Comme on va \'etudier le comportement de la fonction $F_D$ en chacune des variables $t'$, $\nu$ et $\tau$ s\'eparemment, on utilisera pour conclure le th\'eor\`eme de l'analyticit\'e s\'epar\'ee d'Hartogs. C'est pourquoi on a eu besoin d'\'etendre la fonction $F_D(t)$ \`a l'ouvert $U$. On verra ensuite que cette proposition nous permet de d\'eduire la continuit\'e de $F_D$ en la variable $t$ en les points du bord du simplexe $\pi^n$. 

\begin{prop}
Soient un jeu de directions orient\'ees $D\in\D^n$ et un entier $p$, $1\leq p\leq n$. On d\'efinit le jeu de directions orient\'ees $D'\in\D^{p-1}$ par
\[
  D'=(D_1,\ldots,D_{p-1},D_{n+1},D_{n+2},D_{n+3}),
\]
et on note $\s\pi$ la mesure de l'angle ext\'erieur entre les directions orient\'ees $D_{p-1}$ et $D_{n+1}$ telle que $0<\s<1$. Soit un ouvert $\Omega'$ de $U'$ tel que pour tout $\a=1,\ldots,p-1$, sa projection $\Omega'_\a$ sur la $\a$-i\`eme coordonn\'ee v\'erifie
\[
 \text{dist}(\Omega'_\a,0)>0.
\]
Alors il existe  $\e>0$ tel que pour tout secteur $S_{\e,\varphi}=\left\{ \tau \in\C \ | \ 0<|\tau|<\e , \ |\arg\tau|<\varphi\right\}$, le produit cart\'esien
\[
 \Omega'\times \w U \times S_{\e,\varphi}
\]
soit contenu dans $V$ et que dans ce produit la fonction $\underline F_D\left(t',\tau\cdot\nu\right)$ v\'erifie
\[
 \underline F_D\left(t',\tau\cdot\nu\right)=\mathcal H\left(t',\nu,\tau^\sigma,\tau^{1-\sigma}\right),
\]
o\`u $\mathcal H(t',\nu,u,v)$ est une fonction holomorphe en $(t',\nu,u,v)$ au voisinage de chacun des points $\left(t'^0,\nu^0,0,0\right)$, avec $t'^0\in\Omega'$ et $\nu^0\in\w U$.

De plus, pour tout $\left(t',\nu\right)\in U'\times\w U$, on a
\[
 \lim_{\tau\to0} \underline F_D\left(t',\tau\cdot\nu\right) = \left(\underline F_{D'}\left(t'\right),0\ldots,0\right).
\]
\label{prop-F-holo-en-t'-nu-tau}
\end{prop}

\begin{rema}
 On proc\`ederait de m\^eme pour les autres faces du simplexe $\pi^n$, et on obtiendrait des r\'esultats analogues, en faisant des changements de variables adapt\'es, par exemple~:
\begin{align*}
 t &= \left(t_1,\ldots,t_{p-1}, \tau +t_q,\ldots,\tau\nu_{q-1} +t_q,t_q,\ldots,t_n\right), \qquad \tau\to 0\\
 t &= \left(t_1,\ldots,t_{p-1},\tau+t_{p+1},t_{p+1},\ldots,t_{q-1},\tau\nu+t_{q+1},t_{q+1},\ldots,t_n\right), \qquad \tau\to 0\\
 t &= \left(\frac{\nu_1}{\tau},\ldots,\frac1\tau,t_p,\ldots,t_n\right), \qquad \tau\to 0.
\end{align*}
\label{rem-chang-var}
\end{rema}

La proposition~\ref{prop-F-holo-en-t'-nu-tau} permet d'\'etablir la proposition fondamentale~\ref{prop-F-cont-au-bord}.

\begin{proof}[D\'emonstration de la proposition~\ref{prop-F-cont-au-bord}]
Pour \'etendre de mani\`ere continue la fonction $F_D(t)$ en chacune des faces $P^\delta$ du bord de $\pi^n$, on va proc\'eder par r\'ecurrence sur la codimension $n-|\delta|$ de $P^\delta$.

Soit une <<~hyper-face~>> $P^\d$ de $\pi^n$, c'est-\`a-dire telle qu'il existe un entier $p\in\{0,\ldots,n\}$ v\'erifiant $\d_p=0$ et $\d_i=1$ pour tout $i\neq p$. Soit $t^0$ un point de $P^\d$. Alors $t_p^0=t_{p+1}^0$. Dans ce cas, le changement de variables adapt\'e est
\[
 \tau:=t_p-t_{p+1}^0, \quad t':=(t_i)_{1\leq i\leq n,\ i\neq p}.
\]
Alors une variante adapt\'ee au point $t^0$ de la proposition~\ref{prop-F-holo-en-t'-nu-tau} nous assure que la fonction $F_D(t',\tau+t_{p+1}^0)$ est holomorphe en $\left( t',\tau^\sigma,\tau^{1-\sigma}\right) $ au point $t'= t'^0$, $\tau=0$, et on obtient donc que la fonction $F_D(t)$ est continue en $t^0$.

Supposons que la fonction $F_D(t)$ se prolonge contin\^ument \`a
toutes les faces de codimension inf\'erieure ou \'egale \`a $q-1$. Soit
$t^0$ un point d'une face $P^\delta$ de codimension
$n-|\delta|=q$. Pour simplifier l'\'ecriture de la d\'emonstration, on
va supposer encore $\delta=(1,\ldots,1,0,\ldots,0)$, c'est-\`a-dire
$t^0=(t'^0,0,\ldots,0)$, avec $t'^0\in\pi^{p-1}$ et
$p=|\delta|+1=n-q+1$. Soit $K'$ un compact de $\pi^{p-1}$ tel que
$t'^0$ soit \`a l'int\'erieur de $K'$. Alors, par la
proposition~\ref{prop-F-holo-en-t'-nu-tau}, on sait qu'il existe
$\e>0$ tel que pour tous $t'\in K'$, $\nu\in\w\pi^{n-p}$,
$-\e<\tau<0$, on ait
\begin{align*}
 F_D\left(t',\tau\cdot\nu\right) &=\mathcal H\left(t',\nu,\tau^\s,\tau^{1-\s}\right)\\
    &= \left(F_{D'}\left(t'\right),0\ldots,0\right) + \tau^\sigma \mathcal H_1\left(t',\nu,\tau^\sigma,\tau^{1-\sigma}\right) + \tau^{1-\sigma} \mathcal H_2\left(t',\nu,\tau^\sigma,\tau^{1-\sigma}\right)
\end{align*}
o\`u les fonctions $\mathcal H_i(t',\nu,u,v)$ ont les m\^emes propri\'et\'es que la fonction $\mathcal H$. Par l'hypoth\`ese de r\'ecurrence, \'etant donn\'e que la codimension des faces de $\w\pi^{n-p}$ est inf\'erieure
ou \'egale \`a $n-p=q-1$, la fonction $F_D(t',\tau\cdot\nu)$ se prolonge contin\^ument en tous les points $t=(t',\tau\cdot\nu)$ tels que
\[
 t'\in K',\ \nu\in\partial\w\pi^{n-p},\ -\e<\tau<0.
\]
La fonction $\mathcal H(t',\nu,\tau^\s,\tau^{1-\s})$ est donc continue dans le compact
\[
 K'\times \overline{\widetilde\pi^{n-p}}\times [-\e,0].
\]
On en conclut donc qu'il existe deux constantes $C_1,\ C_2>0$ telles que pour tout $(t',\nu,\tau)$ dans ce compact, on ait
\[
 \left\Vert \mathcal H_i\left(t',\nu,\tau^\sigma,\tau^{1-\sigma}\right)\right\Vert\leq C_i
\]
($i=1,2$). Et donc
\[
 \left\Vert F_D\left(t',\tau\cdot\nu\right)-\left(F_{D'}\left(t'\right),0\ldots,0\right)\right\Vert\leq C_1|\tau|^\s + C_2|\tau|^{1-\s}.
\]
L'ensemble $\left\{ t=(t',t_p,\ldots,t_n)\in\pi^n \ |\ t'\in K', -\e<t_p<0 \right\} $ est bien un voisinage de $t^0$ dans $\pi^n$, et pour tout $t$ dans cet ensemble, on a
\begin{align*}
 \left\Vert F_D(t)-\left(F_{D'}\left(t'^0\right),0\ldots,0\right)\right\Vert &\leq \left\Vert F_D(t)-\left(F_{D'}\left(t'\right),0\ldots,0\right)\right\Vert \\
    &\qquad+ \left\Vert \left(F_{D'}\left(t'\right),0\ldots,0\right) - \left(F_{D'}\left(t'^0\right),0\ldots,0\right)\right\Vert\\
    &\leq C_1|t_p|^\sigma + C_2|t_p|^{1-\sigma} + C_0\left\Vert t'-t'^0\right\Vert,
\end{align*}
o\`u la derni\`ere in\'egalit\'e provient du fait que la fonction $F_{D'}(t')$ est lipschitzienne dans le compact $K'$. La fonction $F_D(t)$ est donc bien continue au point $t^0$.
\end{proof}

\section{Les pseudo-chocs}
\label{section-pseudo-choc}

Dans cette section, on rappelle des r\'esultats connus sur le comportement du syst\`eme de Schlesinger au voisinage des singularit\'es que Garnier appelle <<pseudo-chocs~>>, c'est-\`a-dire lorsque plusieurs $t_i$ viennent se confondre. On ne se limite pas ici au cas r\'eel, ni aux syst\`emes fuchsiens dont la monodromie v\'erifie une condition du type~\ref{cond-sys-mono}. Ces r\'esultats sont une partie connue du travail de Garnier. Ils ont \'et\'e modernis\'es et approfondis par M. Sato, T. Miwa et M. Jimbo dans~\cite{SMJ}. On les adapte \`a la situation qui nous int\'eresse : le but de cette section est d'obtenir la d\'ependance en $\tau$ de la fonction $F_D(t',\tau\cdot\nu)$ au point $\tau=0$. On donne \`a l'appendice~\ref{annexe} les d\'emonstrations des principaux r\'esultats de~\cite{SMJ} et~\cite{Jimbo} que l'on va utiliser, et on \'etablit dans ce chapitre uniquement les propri\'et\'es nouvelles dont on a besoin.

On consid\`ere une famille isomonodromique de syst\`emes fuchsiens non r\'esonnants et normalis\'es en l'infini
\[
 \frac{dY}{dx}=A(x,t)Y, \qquad \text{ o\`u } A(x,t)=\sum_{i=1}^{n+2}\frac{A_i(t)}{x-t_i}
\]
o\`u les matrices $(A_1(t),\ldots,A_{n+2}(t))$ sont solutions du syst\`eme de Schlesinger~\eqref{schlesinger}. On suppose que les matrices $A_i(t)$ ($i=1,\ldots,n+2$) sont \`a trace nulle. On note
\[
 -\frac{\t_i}{2}, \quad \frac{\t_i}{2}
\]
les valeurs propres de la matrice $A_i(t)$, qui sont constantes, ainsi que
\[
 A_\infty = -\sum_{i=1}^{n+2}A_i(t) =
\left(1-\tfrac{\t_\infty}{2}\right)
\begin{pmatrix}
 1&0\\
 0&-1
\end{pmatrix}.
\]
On fixe un entier $p$, $1\leq p\leq n$, et on \'etudie le comportement des matrices $A_i(t)$ lorsque $t_p,\ldots,t_n$ tendent vers $0$, les autres variables $t_1,\ldots,t_{p-1}$ demeurant \`a distance mutuelle sup\'erieure \`a un nombre strictement positif. On fait le changement de variables~\eqref{chang-var-t}. Le syst\`eme pr\'ec\'edent s'\'ecrit alors
\begin{equation}
 \frac{dY}{dx} = \left( \sum_\a \frac{A_\a(t',\tau\cdot\nu)}{x-t_\a} + \sum_\b \frac{A_\b(t',\tau\cdot\nu)}{x-\tau\nu_\b}\right) Y.
\label{sys-restreint}
\end{equation}
Dans cette section (\`a l'exception de la
proposition~\ref{prop-sch-lim}), on va supposer les variables
$( t',\nu) \in\B^{p-1}\times\B^{n-p}$ fix\'ees. On pose
\begin{equation}
 \begin{split}
   r &=\min\left\lbrace |t_\a|,\ \a=1,\ldots,p-1,n+2\right\rbrace  >0,\\
   R &=\max\left\lbrace |\nu_\b|,\ \b=p,\ldots,n\right\rbrace  \geq1.
 \end{split}
\label{def-r-R}
\end{equation}
D\`es que $|\tau|<r/R$, le $n$-uplet $(t',\tau\cdot\nu)$ est dans $\B^n$. En fixant $( t',\nu)$, on va donc pour chaque valeur $t'^0$ de $t'$, limiter l'\'etude le long de toute droite passant par le point $(t'^0,0,\ldots,0)$ et contenue dans le sous-espace $t'=t'^0$. Ces droites sont param\'etr\'ees par la variable $\nu$. Quand il n'y a pas d'ambigu\"it\'e, on ne note plus la d\'ependance en $t'$ et en $\nu$. Les transformations isomonodromiques de param\`etre $\tau$ du syst\`eme~\eqref{sys-restreint} sont donn\'ees par le syst\`eme de Schlesinger restreint~:
\begin{equation}
 \begin{split}
  \frac{dA_\a}{d\tau} &= \sum_\b \frac{\nu_\b}{\tau\nu_\b-t_\a} \left[A_\b(\tau),A_\a(\tau)\right]\\
  \frac{dA_\b}{d\tau} &= \sum_\a \frac{\nu_\b}{\tau\nu_\b-t_\a} [A_\a(\tau),A_\b(\tau)] + \frac1\tau \sum_{\b'(\neq \b)}[A_{\b'}(\tau),A_\b(\tau)].
 \end{split}
\label{sch-restreint}
\end{equation}

\subsection{Les solutions du syst\`eme de Schlesinger}

On \'etudie le comportement des solutions du syst\`eme de Schlesinger restreint~\eqref{sch-restreint} lorsque $\tau$ tend vers $0$. Ceci nous permettra ensuite d'en d\'eduire celui des solutions du syst\`eme fuchsien~\eqref{sys-restreint}. Le th\'eor\`eme suivant est \'etabli par Garnier dans~\cite{Garnier26} quand $p=n$, et dans~\cite{Garnier28} dans le cas r\'eel pour $p$ quelconque. Il est repris et g\'en\'eralis\'e dans~\cite{SMJ}, en particulier aux autres changements de variables de la remarque~\ref{rem-chang-var} et aux syst\`emes de dimension quelconque.

\begin{theo}[\cite{SMJ}]
Soient $A_\a^0$ ($\a=1,\ldots,p-1,n+2$) et $A_\b^0$ ($\b=p,\ldots,n+1$) des matrices constantes dont les valeurs propres sont respectivement $\left(-\tfrac{\t_\a}{2}, \tfrac{\t_\a}{2}\right)$ et $\left(-\tfrac{\t_\b}{2}, \tfrac{\t_\b}{2}\right)$. On suppose de plus que
\[
 \sum_\a A_\a^0 + \sum_\b A_\b^0 = -A_\infty
\]
et que les valeurs propres $\mu$ et $-\mu$ de la matrice
\[
 \Lambda :=\sum_\b A_\b^0
\]
v\'erifient : $0<2\Re (\mu)<1$. On note $\s=2\Re (\mu)$. Soient $\s_1$ et $K$ deux constantes telles que
\[
 \s<\s_1<1 \qquad \text{ et } \qquad|A_\a^0|<K, \ |A_\b^0|<K.
\]
Alors il existe $\e>0$ tel que dans tout secteur $S_{\e,\varphi}=\left\{ \tau \in\C \ | \ 0<|\tau|<\e , \ |\arg\tau|<\varphi\right\}$, il existe une unique solution $A_\a(\tau)$ ($\a=1,\ldots,p-1,n+2$), $A_\b(\tau)$ ($\b=p,\ldots,n+1$) du syst\`eme~\eqref{sch-restreint} v\'erifiant~:
\begin{equation}
\begin{split}
 &|A_\a(\tau)-A_\a^0|  \leq K|\tau|^{1-\sigma_1},\\
 &|\tau^{-\Lambda} A_\b(\tau)\tau^\Lambda-A_\b^0|  \leq K|\tau|^{1-\sigma_1}.
\end{split}
\label{A_i-in\'egalit\'es-tau=0}
\end{equation}
\label{thm-SMJ}
\end{theo}

On donne la d\'emonstration du th\'eor\`eme~\eqref{thm-SMJ} \`a l'appendice~\ref{annexe}. La proposition suivante, qui n'est pas dans~\cite{SMJ}, se d\'eduit ais\'ement de cette d\'emonstration. On pose
\begin{align*}
 \w A_\a(\tau) &= \tau^{-\Lambda}A_\a(\tau)\tau^\Lambda\\
 \w A_\b(\tau) &= \tau^{-\Lambda}A_\b(\tau)\tau^\Lambda.
\end{align*}

\begin{prop}
 Les matrices $A_\a(\tau)$ et $A_\b(\tau)$ du th\'eor\`eme~\ref{thm-SMJ}, ainsi que les matrices $\w A_\a(\tau)$ et $\w A_\b(\tau)$ v\'erifient dans tout secteur $S_{\e,\varphi}$, o\`u $\e>0$ est donn\'e au th\'eor\`eme~\ref{thm-SMJ}, les propri\'et\'es suivantes
\begin{align}
 A_\a(\tau)-A_\a^0 &= \tau^{1-\sigma} \mathcal H(\tau^\sigma,\tau^{1-\sigma})\label{Aa-holo}\\
 \tau^{-\Lambda}\left( A_\a(\tau)-A_\a^0\right)\tau^\Lambda &= \tau^{1-\sigma} \mathcal H(\tau^\sigma,\tau^{1-\sigma})\label{Aa-tilde-holo}\\
 \tau^{-\Lambda}A_\b(\tau)\tau^\Lambda -A_\b^0 &= \tau^{1-\sigma} \mathcal H(\tau^\sigma,\tau^{1-\sigma})\label{Ab-tilde-holo}\\
 A_\b(\tau) &= \tau^{-\sigma} \mathcal H(\tau^\sigma,\tau^{1-\sigma})\label{Ab-holo}
\end{align}
o\`u $\mathcal H(u,v)$ d\'esigne toute fonction holomorphe en $(u,v)$ dans un voisinage du point $(0,0)$ contenu dans $\C^2$.
\label{prop-A-holo-en-tau}
\end{prop}

\begin{proof}
Remarquons tout d'abord que la propri\'et\'e~\eqref{Ab-holo} est une cons\'equence imm\'ediate de~\eqref{Ab-tilde-holo}, puisque si une matrice $A(\tau)$ est holomorphe en $\tau^\s,\ \tau^{1-\s}$, alors on a
\[
 \tau^{-\Lambda}A(\tau)\tau^\Lambda = \tau^{-\s}\mathcal H(\tau^\sigma,\tau^{1-\sigma}).
\]

\`A la d\'emonstration du th\'eor\`eme~\ref{thm-SMJ}, qui se trouve \`a l'appendice~\ref{annexe}, on construit la solution $A_\a(\tau)$ et $A_\b(\tau)$ du syst\`eme~\eqref{sch-restreint} par int\'egrations successives. On rappelle cette construction. Il faut r\'ecrire le syst\`eme de Schlesinger restreint~\eqref{sch-restreint} avec les matrices $A_\a(\tau)$ ($\a=1,\ldots,p-1,n+2$), $\w A_\b(\tau)$ ($\b=p,\ldots,n+1$) comme inconnues :
\begin{equation*}
\begin{split}
  \frac{dA_\a}{d\tau} = \sum_\b \frac{\nu_\b}{\tau\nu_\b-t_\a} \left[\tau^\Lambda\w A_\b(\tau)\tau^{-\Lambda},A_\a(\tau)\right]\\
  \frac{d\w A_\b}{d\tau} = \sum_\a \frac{\nu_\b}{\tau\nu_\b-t_\a} \left[\tau^{-\Lambda}A_\a(\tau)\tau^\Lambda,\w A_\b(\tau)\right]\\
    + \frac1\tau \sum_{\b'}\left[\left(\w A_{\b'}(\tau)-A_{\b'}^0\right),\w A_\b(\tau)\right].
\end{split}
\end{equation*}
On construit la solution recherch\'ee en proc\'edant par it\'eration. On pose
\[
 A_\a^{(0)}(\tau) = A_\a^0, \qquad  \w A_\b^{(0)}(\tau) = A_\b^0,
\]
et pour tout entier naturel $k$, on d\'efinit les matrices $A_\a^{(k)}(\tau)$ et $\w A_\b^{(k)}(\tau)$ \`a partir de $A_\a^{(k-1)}(\tau)$ et $\w A_\b^{(k-1)}(\tau)$ par :
\begin{align*}
  A_\a^{(k)}(\tau) = A_\a^0 + \sum_\b \int_0^\tau \frac{\nu_\b}{s\nu_\b-t_\a} \left[s^\Lambda\w A_\b^{(k-1)}(s)s^{-\Lambda},A_\a^{(k-1)}(s)\right]ds\\
  \w A_\b^{(k)}(\tau) = A_\b^0 + \sum_\a \int_0^\tau \frac{\nu_\b}{s\nu_\b-t_\a} \left[s^{-\Lambda}A_\a^{(k-1)}(s)s^\Lambda,\w A_\b^{(k-1)}(s)\right]ds\\
    + \sum_{\b'} \int_0^\tau \frac1s \left[\left(\w A_{\b'}^{(k-1)}(s)-A_{\b'}^0\right),\w A_\b^{(k-1)}(s)\right]ds.
\end{align*}
Les int\'egrales sont calcul\'ees le long du segment joignant $0$ et $\tau$~:
\[
 \{s=re^{i\psi}\ | \ 0<r<|\tau|,\ \psi=\arg \tau\}.
\]
On a montr\'e ensuite par r\'ecurrence que les matrices $A_\a^{(k)}(\tau)$ et $\w A_\b^{(k)}(\tau)$ sont bien d\'efinies et qu'elles convergent uniform\'ement dans tout secteur $S_{\e,\varphi}$, o\`u $\e$ est bien choisi. Leurs limites constituent la solution recherch\'ee. Pour montrer la proposition~\ref{prop-A-holo-en-tau}, il suffit donc de montrer que les matrices $A_\a^{(k)}(\tau)$ et $A_\b^{(k)}(\tau)$ v\'erifient pour tout $k$ les propri\'et\'es~\eqref{Aa-holo},~\eqref{Aa-tilde-holo} et~\eqref{Ab-tilde-holo}.

On proc\`ede \'egalement par r\'ecurrence. L'initialisation est imm\'ediate. Si les matrices $A_\a^{(k-1)}(\tau)$ et $A_\b^{(k-1)}(\tau)$ v\'erifient les propri\'et\'es~\eqref{Aa-holo} et~\eqref{Ab-tilde-holo}, alors on voit que les matrices $A_\a^{(k)}(\tau)$ et $\w A_\b^{(k)}(\tau)$ sont obtenues par l'int\'egration de fonctions de la forme
\[
 \tau^{-\sigma} \mathcal H(\tau^\sigma,\tau^{1-\sigma}).
\]
Elles sont donc elles-m\^emes de la forme
\[
 \tau^{1-\sigma} \mathcal H(\tau^\sigma,\tau^{1-\sigma}),
\]
\ie les matrices $A_\a^{(k)}(\tau)$ et $A_\b^{(k)}(\tau)$ v\'erifient les propri\'et\'es~\eqref{Aa-holo} et~\eqref{Ab-tilde-holo}. Elles v\'erifient \'egalement la propri\'et\'e~\eqref{Aa-tilde-holo}, \'etant donn\'e qu'on a
\begin{multline*}
 \tau^{-\Lambda}\left( A_\a(\tau)-A_\a^0\right)\tau^\Lambda  = \\
\sum_\b \int_0^\tau \frac{\nu_\b}{s\nu_\b-t_\a} \left[\left(\frac s\tau\right)^\Lambda\w A_\b^{(k-1)}(s)\left(\frac s\tau\right)^{-\Lambda},\tau^{-\Lambda}A_\a^{(k-1)}(s)\tau^\Lambda\right]ds.
\end{multline*}
La propri\'et\'e~\eqref{Aa-tilde-holo} est donc une cons\'equence de~\eqref{Aa-holo} et~\eqref{Ab-tilde-holo}.
\end{proof}

Garnier~\cite{Garnier26} \'etablit le r\'esultat suivant, qui ne figure pas sous une forme aussi g\'en\'erale dans~\cite{SMJ}.

\begin{prop}
 Toute solution $A_\a(\tau)$ ($\a=1,\ldots,p-1,n+2$), $A_\b(\tau)$ ($\b=p,\ldots,n,n+1$) du syst\`eme de Schlesinger restreint~\eqref{sch-restreint} admet une limite quand $\tau\to0$ au sens de~\eqref{A_i-in\'egalit\'es-tau=0}.
\label{prop-G-representation}
\end{prop}

Je ne donne pas la d\'emonstration (compliqu\'ee) de Garnier. Comme on se limite au cas des syst\`emes de taille $2\times2$, cas o\`u le probl\`eme de Riemann--Hilbert admet toujours une solution, on d\'eduira ais\'ement cette proposition de la proposition~\ref{prop-mono-a-la-lim}, c'est-\`a-dire de la monodromie des syst\`emes fuchsiens associ\'es \`a chaque solution du syst\`eme de Schlesinger restreint~\eqref{sch-restreint}. On n'utilisera la proposition~\ref{prop-G-representation} qu'\`a la section suivante.

On donne \`a pr\'esent la d\'ependance en $t'$ et en $\nu$ au voisinage de $\tau=0$ des matrices $A_\a(t',\tau\cdot\nu)$ et $A_\b(t',\tau\cdot\nu)$. On sait d\'ej\`a que lorsque $\tau\neq0$, ces matrices sont m\'eromorphes en $t'$ et en $\nu$ tant que la variable $t=(t',\tau\cdot\nu)$ reste dans $\B^n$ (par la propri\'et\'e de Painlev\'e). La proposition suivante permet d'\'etendre ce r\'esultat aux matrices
\begin{align*}
 A_\a^0 &=A_\a^0(t',\nu)\\
 A_\b^0 &=A_\b^0(t',\nu)\\
 \Lambda &=\Lambda(t',\nu).
\end{align*}
Sa d\'emonstration est donn\'ee \`a l'appendice~\ref{annexe}.


\begin{prop}
Les matrices $A_\a^0(t',\nu)$ ($\a=1,\ldots,p-1,n+2$) et $\Lambda(t',\nu)$ sont solutions du syst\`eme de Schlesinger suivant
\begin{align*}
\begin{cases}
 d'A_\a' &= \displaystyle\sum_{\a'\neq\a} \left[A'_{\a'},A'_\a\right]d'\log(t_\a-t_{\a'})\\
 d_\nu A'_\a &= 0
\end{cases}
\end{align*}
o\`u on a pos\'e $A_{n+1}^0(t',\nu) :=\Lambda(t',\nu)$, et o\`u $d'$ d\'esigne la diff\'erentiation par rapport \`a $t'=(t_1,\ldots,t_{p-1})$ et $d_\nu$ la diff\'erentiation par rapport \`a $\nu=(\nu_{p+1},\ldots,\nu_n)$.

Les matrices $A^0_\a(t',\nu)$ ($\a=1,\ldots,p-1,n+2$) et $A_\b^0(t',\nu)$ ($\b=p,\ldots,n+1$) sont solutions du syst\`eme
\begin{align*}
\begin{cases}
 d'A_\b^0 &= -\displaystyle\sum_{\a=1}^{p-1} \left[ A_\b^0,A_\a^0\right] d'\log(t_\a)\\
 d_\nu A_\b^0 &= \displaystyle\sum_{\b'\neq\b} \left[A_{\b'}^0,A_\b^0\right] d_\nu\log(\nu_\b-\nu_{\b'})
\end{cases}.
\end{align*}
\label{prop-sch-lim}
\end{prop}

En particulier, les matrices $A_\a^0(t',\nu)$ et $\Lambda(t',\nu)$ sont ind\'ependantes de $\nu$ et sont solutions du syst\`eme de Schlesinger~\eqref{schlesinger} de dimension $p-1$.

\subsection{Les solutions du syst\`eme fuchsien}

Pour toute matice fondamentale de solutions $\Y(x,\tau)$ du syst\`eme~\eqref{sys-restreint}, la matrice $\wY(y,\tau)=\tau^{-\Lambda} \Y(\tau y,\tau)$ est une matrice fondamentale de solutions du syst\`eme fuchsien non r\'esonnant
\begin{equation}
 \frac{dY}{dy} = \widetilde A(y,\tau) Y,
\label{sys-fuchsien-tilde}
\end{equation}
o\`u la matrice $\w A(y,\tau)$ est d\'efinie par
\[
 \w A(y,\tau) = \sum_\a \frac{\w A_\a(\tau)}{y-\dfrac{t_\a}{\tau}} + \sum_\b \frac{\w A_\b(\tau)}{y-\nu_\b} = \tau \left(\tau^{-\Lambda}A(\tau y,\tau)\tau^\Lambda\right).
\]
Le syst\`eme~\eqref{sys-fuchsien-tilde} n'est pas normalis\'e en l'infini.

\begin{prop}
\begin{enumerate}
 \item La solution fondamentale $\Y_\infty(x,\tau)$ normalis\'ee en l'infini du syst\`eme~\eqref{sys-restreint} est holomorphe en $\tau^\sigma$, $\tau^{1-\sigma}$ au point $\tau=0$ pour tout $x\neq0$ fix\'e. Sa limite $\displaystyle\lim_{\tau\to0} \Y_\infty(x,\tau)$ existe donc et est solution du syst\`eme fuchsien
\begin{equation}
 \frac{dY}{dx} = \left( \sum_\a \frac{A_\a^0}{x-t_\a} + \frac{\Lambda}{x} \right) Y
\label{sys-lim1-tau=0}.
\end{equation}
 \item La solution fondamentale $\wY(y,\tau):=\tau^{-\Lambda}\Y_\infty(\tau y,\tau)$ du syst\`eme~\eqref{sys-fuchsien-tilde} est holomorphe en $\tau^\sigma$, $\tau^{1-\sigma}$ au point $\tau=0$ pour tout $y\in\C$ fix\'e. Sa limite $\displaystyle\lim_{\tau\to0} \wY(y,\tau)$ existe donc et est solution du syst\`eme fuchsien
\begin{equation}
 \frac{dY}{dy} = \sum_\b \frac{A_\b^0}{y-\nu_\b} \ Y
\label{sys-lim2-tau=0}.
\end{equation}
\end{enumerate}
\label{prop-Yinfty-holo-en-tau}
\end{prop}

\begin{proof}
On ne montre que l'assertion (2) ; l'assertion (1) se montre de la m\^eme mani\`ere. On calcule la d\'eriv\'ee de la matrice $\wY(y,\tau)$ par rapport \`a $\tau$, pour $y$ fix\'e. On suppose que $|\tau|<r/|y|$. Sachant que
\[
 \frac{\partial}{\partial t_i}\Y_\infty(x,t) = -\frac{A_i(t)}{x-t_i} \ \Y_\infty (x,t),
\]
(lemme~\ref{lemme-omega-sch}), on trouve
\begin{align*}
 \frac{d}{d\tau}\Y_\infty(\tau y,\tau) &= \left( yA(\tau y,\tau) -\sum_\b\frac{\nu_\b A_\b(\tau)}{\tau y-\tau\nu_\b} \right) \Y_\infty (\tau y,\tau)\\
    &=\left( y\sum_\a\frac{A_\a(\tau)}{\tau y-t_\a} + \frac1\tau \sum_\b A_\b(\tau) \right) \Y_\infty (\tau y,\tau).
\end{align*}
Et comme $A_\infty=-\displaystyle\sum_\a A_\a(\tau) -\sum_\b A_\b(\tau) = -\sum_\a A_\a^0 - \Lambda$, on obtient
\[
 \frac{d}{d\tau}\wY(y,\tau) = \left( -\frac1\tau \sum_\a\tau^{-\Lambda}\left( A_\a(\tau)-A_\a^0\right) \tau^\Lambda +y\sum_\a\frac{\tau^{-\Lambda}A_\a(\tau)\tau^\Lambda}{\tau y-t_\a} \right) \wY(y,\tau),
\]
c'est-\`a-dire, vu~\eqref{Aa-tilde-holo}
\[
 \frac{d}{d\tau}\wY(y,\tau) = \tau^{-\sigma}\mathcal H\left( y,\tau^\sigma,\tau^{1-\sigma}\right) \wY(y,\tau),
\]
o\`u $\mathcal H(y,u,v)$ d\'esigne une fonction holomorphe au voisinage de $(y_0,0,0)$ pour tout $y_0\in\C$. On en conclut donc qu'il existe une matrice 
\[
 Q(y,\tau)=\I_2+\tau^{1-\sigma}\mathcal H_1\left( y,\tau^\sigma,\tau^{1-\sigma}\right),
\]
o\`u la fonction $\mathcal H_1$ a les m\^emes propri\'et\'es que $\mathcal H$, et une matrice $\wY^0(y)$ ind\'ependante de $\tau$ telles que
\[
 \wY(y,\tau)=Q(y,\tau)\wY^0(y).
\]

Il ne reste donc qu'\`a prouver que la matrice $\wY^0(y)$ est solution du syst\`eme fuchsien~\eqref{sys-lim2-tau=0}. Pour cela, il suffit de v\'erifier que la matrice
\[
 \widetilde A(y,\tau) = \sum_\a \frac{\widetilde A_\a(\tau)}{y-\dfrac{t_\a}{\tau}} + \sum_\b \frac{\widetilde A_\b(\tau)}{y-\nu_\b}  
\]
tend en $\tau=0$ vers la matrice
\[
 \sum_\b \frac{A_\b^0}{y-\nu_\b}.
\]
Ceci est \'evident, \'etant donn\'e que la matrice
\[
 \sum_\a \frac{\widetilde A_\a(\tau)}{y-\dfrac{t_\a}{\tau}} = \tau^{1-\sigma} \sum_\a \frac{\tau^\sigma\widetilde A_\a(\tau)}{\tau y-t_\a}
\]
tend vers la matrice nulle par l'assertion~\eqref{Aa-tilde-holo} de la proposition~\ref{prop-A-holo-en-tau}. Pour la deuxi\`eme partie de l'assertion (1), on aurait montr\'e de m\^eme que la matrice
\[
 A(x,\tau) = \sum_\a \frac{A_\a(\tau)}{x-t_\a} + \sum_\b \frac{A_\b(\tau)}{x-\tau\nu_\b} 
\]
tend en $\tau=0$ vers la matrice
\[
 \sum_\a \frac{A_\a^0}{x-t_\a} + \frac{\Lambda}{x}
\]
en remarquant que
\begin{align*}
 \sum_\b \frac{A_\b(\tau)}{x-\tau\nu_\b} &= \frac1x \sum_\b A_\b(\tau) + \tau\sum_\b \frac{\nu_\b A_\b(\tau)}{x(x-\tau\nu_\b)}\\
&= -\frac1x \left( A_\infty + \sum_\a A_\a(\tau)\right) + \tau^{1-\sigma}\sum_\b \frac{\nu_\b\tau^\sigma A_\b(\tau)}{x(x-\tau\nu_\b)}.
\end{align*}
\end{proof}

Le syst\`eme fuchsien~\eqref{sys-lim1-tau=0} est non r\'esonnant et normalis\'e en l'infini, \'etant donn\'e que
\[
 -\sum_\a A_\a^0 - \Lambda = -\sum_\a A_\a^0 - \sum_\b A_\b^0 = A_\infty.
\]
Soit $\Y^0_\infty(x)$ sa matrice fondamentale de solutions normalis\'ee en l'infini. Comme $0<\s<1$, le comportement local de $\Y^0_\infty(x)$ au voisinage des singularit\'es du syst\`eme~\eqref{sys-lim1-tau=0} est donn\'e par
\begin{align}
\Y^0_\infty(x)
  & =\left(S_\a^0 +\mathcal O(x-t_\a)\right)(x-t_\a)^{L_\a}\cdot C_\a^0  & x\to t_\a\nonumber\\
  & =\left(\I_2 +\mathcal O(x)\right)x^\Lambda\cdot C^0                   & x\to 0 \ \label{comp-local-Y0infty}\\
  & =\left(\I_2 +\mathcal O\left(x^{-1}\right)\right)x^{-L_\infty}                 & x\to \infty\nonumber
\end{align}
o\`u les matrices $S_\a^0$, $C_\a^0$ et $C^0$ sont inversibles, les matrices $L_\a$ d\'esignent comme pr\'ec\'edemment les diagonalis\'ees de $A_\a(\tau)$ (et donc aussi de $A_\a^0$) et $L_\infty=A_\infty$. Le syst\`eme~\eqref{sys-lim2-tau=0} n'est pas normalis\'e en l'infini, puisque la matrice $\Lambda$ n'est pas diagonale, mais il existe de m\^eme une unique matrice fondamentale de solutions $\wY^0_\infty(y)$ dont le comportement local est donn\'e par
\begin{align}
\wY^0_\infty(y)
  & =\left(\w S_\b^0 +\mathcal O(y-\nu_\b)\right)(y-\nu_\b)^{L_\b}\cdot \w C_\b^0 & y\to \nu_\b \label{comp-local-Ytilde-infty}\\
  & =\left(\I_2 +\mathcal O\left(y^{-1}\right)\right)y^{-\Lambda}                           & y\to\infty\nonumber
\end{align}
avec $\w S_\b^0$, $\w C_\b^0 \in GL(2,\C)$. La proposition suivante se trouve dans~\cite{Jimbo}. On ne donne pas sa d\'emonstration, qui proc\`ede des m\^eme m\'ethodes que celle du th\'eor\`eme~\ref{thm-SMJ}.

\begin{prop}[\cite{Jimbo}]
 On a
\[
 \lim_{\tau\to0} \Y_\infty(x,\tau) = \Y^0_\infty(x),
 \ \ \
 \lim_{\tau\to0} \tau^{-\Lambda} \Y_\infty(\tau y,\tau) = \mathbf{\widetilde Y}^0_\infty(y)\cdot C^0.
\]
 De plus, pour $\tau \neq 0$, le comportement local de la matrice fondamentale $\Y_\infty(x,\tau)$ est donn\'e par
\begin{align}
 \Y_\infty(x,\tau)
& =\left(S_\a(\tau) +\mathcal O(x-t_\a)\right)(x-t_\a)^{L_\a}\cdot C_\a^0  & x\to t_\a\nonumber\\
& =\left(\w S_\b(\tau)+\mathcal O(x-\tau\nu_\b)\right)(x-\tau\nu_\b)^{L_\b}\cdot \w C_\b^0\cdot C^0 & x\to \tau\nu_\b\label{comp-local-Yinfty}\\
& =\left(\I_2 +\mathcal O\left(x^{-1}\right)\right)x^{-L_\infty}  & x\to\infty\nonumber
\end{align}
o\`u les matrices $S_\a(\tau)$ et $\w S_\b(\tau)$ sont inversibles, et les matrices $C^0$, $C_\a^0$ et $\w C_\b^0$ sont d\'efinies ci-dessus.
\label{prop-mono-a-la-lim}
\end{prop}

La proposition~\ref{prop-mono-a-la-lim} nous permet d'\'etablir simplement la proposition~\ref{prop-G-representation}.

\begin{proof}[D\'emonstration de la proposition~\ref{prop-G-representation}]
Soit une solution quelconque $A_\a(\tau)$ ($\a=1,\ldots,p-1,n+2$), $A_\b(\tau)$ ($\b=p,\ldots,n,n+1$) du syst\`eme de Schlesinger restreint~\eqref{sch-restreint} telle que la somme $\sum_\a A_\a(\tau) + \sum_\b A_\b(\tau)$ soit constante et diagonale. Soit $\Y_\infty(x,\tau)$ l'unique solution fondamentale normalis\'ee en l'infini du syst\`eme fuchsien~\eqref{sys-restreint} d\'efini par les matrices $A_\a(\tau)$ et $A_\b(\tau)$. Cette solution est $M$-invariante. Il existe donc des matrices inversibles $C^0$, $C_\a^0$ et $\widetilde C_\b^0$ ind\'ependantes de $\tau$ telles que le comportement local de la solution $\Y_\infty(x,\tau)$ soit donn\'e par~\eqref{comp-local-Yinfty}. On d\'efinit la matrice $\Lambda$ de mani\`ere \`a ce que les probl\`emes de Riemann--Hilbert~\eqref{comp-local-Y0infty} et~\eqref{comp-local-Ytilde-infty} v\'erifient bien
\[
 M_\infty\left( \Y_\infty^0\right)\cdot (C^0)^{-1}e^{2i\pi \Lambda}C^0\cdot M_{p-1}\left( \Y_\infty^0\right) \cdots M_1\left( \Y_\infty^0\right) =\I_2,
\]
et
\[
 e^{2i\pi \Lambda}\cdot \widetilde M_{n+1}\left( \wY_\infty^0\right)\cdots \widetilde M_p\left( \wY_\infty^0\right) =\I_2,
\]
et que les valeurs propres de $\Lambda$ soit oppos\'ees : $\mu$ et $-\mu$ et v\'erifient $0<2\Re(\mu)<1$. Les deux conditions pr\'ec\'edentes sont \'equivalentes par la relation
\[
 M_{n+3}\left( \Y_\infty\right) \cdots M_1\left( \Y_\infty\right) =\I_2.
\]
Alors, comme on consid\`ere des syst\`emes de taille $2\times2$, on sait que les probl\`emes de Riemann--Hilbert~\eqref{comp-local-Y0infty} et~\eqref{comp-local-Ytilde-infty} admettent respectivement d'uniques
solutions $\Y_\infty^0(x)$ et $\wY_\infty^0(y)$. On d\'efinit les matrices constantes $A_\a^0$ ($\a=1,\ldots,p-1,n+2$), $A_\b^0$ ($\b=p,\ldots,n,n+1$) respectivement associ\'ees aux solutions $\Y_\infty^0(x)$ et $\wY_\infty^0(y)$. Par le th\'eor\`eme~\ref{thm-SMJ}, ces matrices $A_\a^0$, $A_\b^0$ sont les conditions initiales au sens de~\eqref{A_i-in\'egalit\'es-tau=0} d'une unique solution du syst\`eme de Schlesinger restreint~\eqref{sch-restreint}. Cette solution est n\'ecessairement la solution $A_\a(\tau)$, $A_\b(\tau)$ par unicit\'e de la matrice $\Y_\infty(x,\tau)$ satisfaisant le probl\`eme de Riemann--Hilbert~\eqref{comp-local-Yinfty}.
\end{proof}

\section{Le cas r\'eel}
\label{section-cas-reel}

On consid\`ere \`a pr\'esent la limite d'une famille isomodromique de syst\`emes fuchsiens $\left(A_D(t), t\in U\right)$, associ\'ee \`a un jeu de directions orient\'ees $D\in\D^n$ et d\'ecrite par le syst\`eme de Schlesinger, que l'on a introduite \`a la section~\ref{section-def-ADt}. L'ouvert simplement connexe $U$ est un voisinage contenu dans $\B^n$ du simplexe $\pi^n$
\[
 \pi^n=\left\{ (t_1,\ldots,t_n)\in\R^n\ |\ t_1<\cdots<t_n<0\right\}  ,
\]
tel que la solution du syst\`eme de Schlesinger $\left(A_{D,1}(t),\ldots,A_{D,n+2}(t)\right)$ correspondant \`a cette famille est holomorphe dans $U$. D'apr\`es la proposition~\ref{prop-G-representation}, cette solution admet une limite
\[
 A_{D,\a}^0(t') \ (\a=1,\ldots,p-1,n+2), \qquad A_{D,\b}^0(t',\nu) \ (\b=p,\ldots,n+1)
\]
au sens de~\eqref{A_i-in\'egalit\'es-tau=0} lorsque $\tau$ tend vers $0$. D'apr\`es la proposition~\ref{prop-sch-lim}, les matrices $A_{D,\a}^0(t')$ ($\a=1,\ldots,p-1,n+2$) et $\Lambda(t')$ $(t'\in U')$ sont solutions du syst\`eme de Schlesinger de dimension $p-1$. 

D'apr\`es la proposition~\ref{prop-Yinfty-holo-en-tau}, pour chaque valeur de $(t',\nu)\in U'\times\w U$, le syst\`eme fuchsien $\left(A_D(t',\tau\cdot\nu)\right)$ tend lorsque $\tau$ tend vers $0$ vers le syst\`eme fuchsien limite ind\'ependant de $\nu$ suivant
\begin{equation}
 \frac{dY}{dx} = \left( \sum_\a \frac{A_\a^0(t')}{x-t_\a} + \frac{\Lambda(t')}{x} \right) Y.
\tag{$A^0_D(t')$}\label{A0(t')}
\end{equation}
La famille de syst\`emes fuchsiens limites $\left(A^0_D(t'), t'\in U'\right)$ est donc isomonodromique et d\'ecrite par le syst\`eme de Schlesinger. Les syst\`emes~$A^0_D(t')$ sont non r\'esonnants. Pour tout $\a=1,\ldots,p-1,n+2$, les valeurs propres de la matrice $A_{D,\a}^0(t')$ sont ind\'ependantes de $t'$ et valent
\[
 -\frac{\t_\a}{2}, \ \frac{\t_\a}{2}
\]
et les valeurs propres de la matrice $\Lambda(t')$ sont $-\mu$ et $\mu$, avec $\s=2\Re(\mu)$. Les syst\`emes~$\left(A^0_D(t')\right)$ sont normalis\'es en l'infini et ils ont la m\^eme normalisation que les syst\`emes $\left(A_D(t)\right)$.

\begin{lemm}
Soient un jeu de directions orient\'ees $D\in\D^n$ et un entier $p$, $1\leq p\leq n$. La famille isomonodromique
\[
 \left(A^0_D(t'), t'\in\pi^{p-1}\right)
\]
est contenue dans l'ensemble $\A_{D'}^{p-1}$ des syst\`emes fuchsiens associ\'es au jeu de directions orient\'ees $D'\in\D^{p-1}$ d\'efini par
\begin{equation}
 D'=(D_1,\ldots,D_{p-1},D_{n+1},D_{n+2},D_{n+3}).
\label{def-D'}
\end{equation}
On note donc le syst\`eme $\left(A^0_D(t')\right)$ par $\left(A_{D'}(t')\right)$.

De plus, la fonction <<~rapports des longueurs~>> $F_{D'}(t')=(r'_1(t'),\ldots,r'_{p-1}(t'))$ associ\'ee au jeu de directions orient\'ees $D'\in\D^{p-1}$ est donn\'ee par
\begin{equation}
 r'_\a(t') =\frac{\displaystyle\int_{t_\a}^{t_{\a+1}} L_1\left(\Y_\infty^0(x,t')\cdot C_0\right)^2 dx}{\displaystyle\int_0^1 L_1\left(\Y_\infty^0(x,t')\cdot C_0\right)^2 dx}
\label{rapport-r'-Yinfty0}
\end{equation}
($\a=1,\ldots,p-1$), o\`u la solution fondamentale $\Y_\infty^0(x,t')$ est la solution normalis\'ee en l'infini du syst\`eme~$A_{D'}(t')$ et la matrice $C_0$ est d\'efinie au lemme~\ref{lemme-realite-Yinfty}.
\label{lemme-sysFu-lim-tau-0}
\end{lemm}

\begin{proof}
Pour la premi\`ere partie du lemme, il suffit de v\'erifier que la monodromie du syst\`eme~$\left(A^0_D(t')\right)$ est engendr\'ee par les matrices $M_\a^0$ d\'efinies par
\[
 M_\a^0:=M_\a=D_\a D_{\a-1}^{-1} \qquad (\a=1,\ldots,p-1,n+2,n+3)
\]
et 
\[
 M_{n+1}^0:=D_{n+1}D_{p-1}^{-1}.
\]
Par la proposition~\ref{prop-mono-a-la-lim}, pour tout $\a=1,\ldots,p-1,n+2$, les monodromies des solutions fondamentales $\Y_\infty(x,\tau)$ et $\Y^0_\infty(x)$ autour de la singularit\'e $t_\a$ sont les m\^emes~:
\[
 M_\a\left( \Y^0_\infty\right)={C_\a^0}^{-1}e^{2i\pi L_\a}C_\a^0=M_\a\left( \Y_\infty\right)
\]
et donc, vu la condition~\ref{cond-sys-mono} et le lemme~\ref{lemme-realite-Yinfty},
\[
 M_\a\left( \Y^0_\infty\right)=C_0 M_\a C_0^{-1}=C_0 D_\a D_{\a-1}^{-1}C_0^{-1}.
\]
De m\^eme, en $t_{n+3}=\infty$~:
\[
 M_\infty\left( \Y^0_\infty\right)=e^{2i\pi L_\infty}=C_0D_{n+3}D_{n+2}^{-1}C_0^{-1}.
\]
Il ne reste plus qu'\`a d\'eterminer la monodromie autour de la singularit\'e $t_{n+1}=0$~:
\begin{align*}
 M_{n+1}\left( \Y^0_\infty\right)
    &=\left( M_{p-1}\left( \Y^0_\infty\right) \cdots M_1\left( \Y^0_\infty\right)M_{n+3}\left( \Y^0_\infty\right)M_{n+2}\left( \Y^0_\infty\right)\right)^{-1}\\
    &=\left(C_0D_{p-1}D_{p-2}^{-1}D_{p-2}\cdots D_{n+1}^{-1}C_0^{-1}\right)^{-1}\\
    &=C_0D_{n+1}D_{p-1}^{-1}C_0^{-1}.
\end{align*}
On a donc montr\'e que pour tout $\a=1,\ldots,p-1,n+1,n+2,n+3$, on a
\[
 M_\a\left( \Y^0_\infty\right)=C_0 M_\a^0C_0^{-1}
\]
o\`u la matrice de conjugaison $C_0$ est la m\^eme qu'entre les matrices $M_i\left( \Y_\infty\right)$ et les matrices $M_i$. Gr\^ace \`a cela, en proc\'edant exactement comme \`a la d\'emonstration du lemme~\ref{lemme-realite-Yinfty}, on obtient l'expression~\eqref{rapport-r'-Yinfty0} des rapports $r'_\a(t')$.
\end{proof}

On d\'eduit en particulier de ce lemme que les valeurs propres de la matrice $\Lambda(t')$ sont r\'eelles et valent
\[
 -\frac\s2, \ \frac\s2
\]
o\`u $\s\pi$ est la mesure de l'angle ext\'erieur entre les directions orient\'ees $D_{p-1}$ et $D_{n+1}$ telle que $0<\s<1$.

Quitte \`a diminuer l'ouvert simplement connexe $U'$, on peut supposer gr\^ace \`a la proposition~\ref{prop-Ai-holo} que les matrices $A_{D,\a}^0(t')$ et $\Lambda(t')$ sont holomorphes dans $U'$.

\begin{lemm}
Soient un jeu de directions orient\'ees $D\in\D^n$ et un entier $p$, $1\leq p\leq n$. Pour tout $\left(t',\nu\right)\in U'\times\w U$ fix\'e, il existe $\e>0$ tel que le prolongement de la fonction <<~rapports des longueurs~>> $\underline F_D(t',\tau\cdot\nu)$ soit holomorphe en $\tau^\sigma$, $\tau^{1-\sigma}$ au point $\tau=0$ dans tout secteur $S_{\e,\varphi}$. De plus, on a
\[
 \lim_{\tau\to0} \underline F_D\left(t',\tau,\tau\nu_{p+1},\ldots,\tau\nu_n\right) = \left(\underline F_{D'}\left(t'\right),0,\ldots,0\right)
\]
o\`u le jeu de directions orient\'ees $D'\in\D^{p-1}$ est donn\'e par~\eqref{def-D'}.
\label{lemme-F-holo-en-tau}
\end{lemm}

\begin{proof}
On choisit $\e>0$ tel que pour tout $\a$ on ait $|t_\a|>\e$. Consid\'erons l'expression~\eqref{r_i-par-Y_i} de la fonction $\underline F_D(t)$ \`a partir des solutions fondamentales $\Y_i(x,t',\tau\cdot\nu)$ d\'efinies par~\eqref{def-Y_i} : pour tout $i=1,\ldots,n$
\[
 r_i(t',\tau\cdot\nu)= \frac{\ell_i(t',\tau\cdot\nu)}{\ell_{n+1}(t',\tau\cdot\nu)}
\]
o\`u pour $\a=1,\ldots,p-2,\ n+1$
\[
 \ell_\a(t',\tau\cdot\nu)=\displaystyle\int_{t_\a}^{t_{\a+1}}\left(g_\a(x,t',\tau\cdot\nu)^2 + h_\a(x,t',\tau\cdot\nu)^2\right)dx,
\]
et
\[
 \ell_{p-1}(t',\tau\cdot\nu)=\displaystyle\int_{t_{p-1}}^{\tau}\left(g_{p-1}(x,t',\tau\cdot\nu)^2 + h_{p-1}(x,t',\tau\cdot\nu)^2\right)dx,
\]
et pour $\b=p,\ldots,n$
\[
 \ell_\b(t',\tau\cdot\nu)=\displaystyle\int_{\tau\nu_\b}^{\tau\nu_{\b+1}}\left(g_\b(x,t',\tau\cdot\nu)^2 + h_\b(x,t',\tau\cdot\nu)^2\right)dx
\]
o\`u les fonctions $\left(g_i(x,t',\tau\cdot\nu),h_i(x,t',\tau\cdot\nu)\right)$ constituent la premi\`ere ligne de la solution fondamentale $\Y_i(x,t',\tau\cdot\nu)$. Les int\'egrales sont calcul\'ees le long des segments joignant respectivement $t_i$ et $t_{i+1}$. On ne d\'etaille pas le cas de la fonction $\ell_{p-1}(t',\tau\cdot\nu)$ ; il faudrait, comme \`a la d\'emonstration de la proposition~\ref{prop-F-holo-dans-pi}, la d\'ecomposer en $\ell_{p-1}=\ell_{p-1}^-+\ell_{p-1}^+$ avec
\[
 \ell_{p-1}^- = \int_{t_{p-1}}^{\frac r\e\tau} \quad \text{ et } \quad \ell_{p-1}^+ = \int_{\frac r\e\tau}^{\tau}\ ,
\]
puis \'etudier la fonction $\ell_{p-1}^-$ comme les fonctions $\ell_\a$ et la fonction $\ell_{p-1}^+$ comme les fonctions $\ell_\b$.

Pour tout $\a=1,\ldots,p-2,\ n+1$, d'apr\`es l'assertion (1) de la proposition~\ref{prop-Yinfty-holo-en-tau}, les solutions fondamentales $\Y_\a(x,t',\tau\cdot\nu)$ sont holomorphes en $\tau^\sigma$, $\tau^{1-\sigma}$ au point $\tau=0$ d\`es que $x\neq0$, et on en d\'eduit donc que les fonctions $\ell_\a(t',\tau\cdot\nu)$ sont \'egalement  holomorphes en $\tau^\sigma$, $\tau^{1-\sigma}$ (la situation est plus simple ici qu'\`a la d\'emonstration de la proposition~\ref{prop-F-holo-dans-pi}, \'etant donn\'e que les bornes d'int\'egration et le facteur $(x-t_\a)^{L_\a}$
sont ind\'ependants de $\tau$). On obtient de m\^eme que la fonction $\ell_{n+1}(t',\tau\cdot\nu)$ ne s'annule jamais pour $|\tau|<\e$. De plus, par les propositions~\ref{prop-F-holo-dans-pi} et~\ref{prop-mono-a-la-lim}, les solutions $\Y_\a(x,t',\tau\cdot\nu)$ ont une limite ind\'ependante de $\nu$ quand $\tau\to0$ qui est solution du syst\`eme $\left(A_{D'}(t')\right)$ et qui v\'erifie :
\begin{align*}
 \Y^0_\a(x,t')
    &:=\lim_{\tau\to0} \Y_\a(x,t',\tau\cdot\nu)\\
    &= \lim_{\tau\to0} \left(\Y_\infty(x,t',\tau\cdot\nu)\cdot C_0\cdot S_\a\right)\\
    &=\Y^0_\infty(x,t')\cdot C_0\cdot S_\a.
\end{align*}
On note par $(g_\a^0(x,t'),h_\a^0(x,t'))$ la premi\`ere ligne de la solution fondamentale $\Y^0_\a(x,t')$, et on obtient donc
\[
 \lim_{\tau\to0} \ell_\a(t',\tau\cdot\nu) =\displaystyle\int_{t_\a}^{t_{\a+1}}\left(g_\a^0(x,t')^2 + h_\a^0(x,t')^2\right)dx.
\]
D'apr\`es l'expression~\eqref{rapport-r'-Yinfty0} des rapports $r'_\a(t')$, comme la matrice $S_\a$ est dans $SU(2)$, on a
\[
 r'_\a(t') = \frac{\displaystyle\int_{t_\a}^{t_{\a+1}}\left(g_\a^0(x,t')^2 + h_\a^0(x,t')^2\right)dx}{\displaystyle\int_0^1\left(g_{n+1}^0(x,t')^2 + h_{n+1}^0(x,t')^2\right)dx},
\]
ce qui donne
\[
 \lim_{\tau\to0} r_\a(t',\tau\cdot\nu) = r'_\a(t').
\]

Pour tout $\b=p,\ldots,n$, on exprime les fonctions $\ell_\b(t',\tau\cdot\nu)$ \`a partir des solutions fondamentales $\wY_\b(y,t',\tau\cdot\nu):=\tau^{-\Lambda} \Y_\b(\tau y,t',\tau\cdot\nu)$ du syst\`eme fuchsien $\left(\w A_D(t)\right)$, qui est le syst\`eme~\eqref{sys-fuchsien-tilde} associ\'e au syst\`eme $\left(A_D(t)\right)$. Pour all\'eger les notations, on ne note plus la d\'ependance en $t'$ et en $\nu$. En faisant le changement de variables
\[
 y=\frac{x}{\tau},
\]
on obtient
\[
 \ell_\b(\tau) = \tau\displaystyle\int_{\nu_\b}^{\nu_{\b+1}}\left(g_\b(\tau y,\tau)^2 + h_\b(\tau y,\tau)^2\right)dy.
\]
On note
\[
 \tau^\Lambda =
\begin{pmatrix}
 a(\tau)&b(\tau)\\
 c(\tau)&d(\tau)
\end{pmatrix}
\]
et
\[
 \wY_\b(y,\tau) =
\begin{pmatrix}
 \tilde y_1(y,\tau) & \tilde z_1(y,\tau)\\
 \tilde y_2(y,\tau) & \tilde z_2(y,\tau)
\end{pmatrix}.
\]
Alors
\begin{multline*}
 g_\b(\tau y,\tau)^2  + h_\b(\tau y,\tau)^2 = a(\tau)^2\left(\tilde y_1(y,\tau)^2 + \tilde z_1(y,\tau)^2\right)
 + b(\tau)^2\left(\tilde y_2(y,\tau)^2 + \tilde z_2(y,\tau)^2\right) \\
 + 2a(\tau)b(\tau)\left(\tilde y_1(y,\tau) \tilde y_2(y,\tau)+ \tilde z_1(y,\tau)\tilde z_2(y,\tau)\right).
\end{multline*}
Comme les \'el\'ements de la matrice $\tau^\Lambda$ sont de la forme $c_1\tau^\frac\s2 +c_{-1}\tau^{-\frac\s2}$ ($c_h \in\C$), les quantit\'es suivantes
\[
 \tau a(\tau)^2, \quad \tau b(\tau)^2, \quad \tau a(\tau)b(\tau)
\]
sont polynomiales en $\tau^\sigma$ et $\tau^{1-\sigma}$ et s'annulent en $\tau=0$. Par l'assertion (2) de la proposition~\ref{prop-Yinfty-holo-en-tau}, la solution fondamentale $\wY_\b(y,\tau)=\wY(y,\tau)\cdot C_0\cdot S_\b$ est holomorphe en $\tau^\sigma$, $\tau^{1-\sigma}$ lorsque $|y|<r/\e$, et donc en particulier quand $y$ appartient \`a l'intervalle $\left]\nu_\b,\nu_{\b+1}\right[$. Les int\'egrales
\[
 \displaystyle\int_{\nu_\b}^{\nu_{\b+1}}\left(\tilde y_k(y,\tau)^2 + \tilde z_k(y,\tau)^2\right)dy
\]
($k=1,2$) et
\[
 \displaystyle\int_{\nu_\b}^{\nu_{\b+1}}\left(\tilde y_1(y,\tau) \tilde y_2(y,\tau)+ \tilde z_1(y,\tau)\tilde z_2(y,\tau)\right)dy
\]
sont donc holomorphes en $\tau^\sigma$, $\tau^{1-\sigma}$ (l\`a encore, par les m\^emes arguments qu'\`a la d\'emonstration de la proposition~\ref{prop-F-holo-dans-pi}). On peut donc en conclure que les fonctions $\ell_\b(t',\tau\cdot\nu)$ sont holomorphes en $\tau^\sigma$, $\tau^{1-\sigma}$ et qu'elles v\'erifient~:
\[
 \lim_{\tau\to0} \ell_\b(t',\tau\cdot\nu) = 0.
\]
\end{proof}

On peut enfin \'etablir la proposition~\ref{prop-F-holo-en-t'-nu-tau}

\begin{proof}[D\'emonstration de la proposition~\ref{prop-F-holo-en-t'-nu-tau}]
Au vu des r\'esultats pr\'ec\'edents, il s'agit simplement d'appliquer le th\'eor\`eme de l'analyticit\'e s\'epar\'ee d'Hartogs. Le lemme~\ref{lemme-F-holo-en-tau} nous donne le comportement en $\tau$ de la fonction $\underline F_D(t',\tau\cdot\nu)$ \`a $\left(t',\nu\right)\in U'\times\w U$ fix\'e. Il ne reste plus qu'\`a v\'erifier qu'en $\tau=0$, cette fonction est holomorphe en $\left(t',\nu\right)$. Comme en $\tau=0$, la fonction $\underline F_D(t',\tau\cdot\nu)$ vaut
\[
 \left(\underline F_{D'}\left(t'\right),0\ldots,0\right),
\]
elle est donc ind\'ependante de $\nu$ et holomophe en $t'\in U'$ par la proposition~\ref{prop-F-holo-dans-pi} appliqu\'ee \`a la dimension $p-1$, et par le choix de l'ouvert simplement connexe $U'$ tel que la solution du syst\`eme de Schlesinger $\left(A_{D',1},\ldots,A_{D',p+1}\right)$ soit holomorphe dans $U'$.
\end{proof}

\appendix
\chapter{Le syst\`eme de Garnier}
\label{annexe-garnier}

Pour \^etre complet, on introduit le syst\`eme de Garnier, qui d\'ecrit les d\'eformations isomonodromiques des \'equations fuchsiennes qui n'ont pas de singularit\'e logarithmique. M\^eme si la r\'esolution du probl\`eme de Plateau propos\'ee dans ce m\'emoire n'utlise pas, contrairement \`a celle de Garnier, le syst\`eme de Garnier, on est malgr\'e tout amen\'e \`a le mentionner \`a plusieurs reprises, ne serait-ce que pour comparer les deux points de vue. 

\bigskip

On consid\`ere une \'equation fuchsienne sur la sph\`ere de Riemann $\P$
\begin{equation}
 D^2y + p(x)Dy+ q(x)y=0
\label{E-G}
\end{equation}
de singularit\'es deux \`a deux distinctes $t_1,\dotsc,t_n$, $t_{n+1}=0$, $t_{n+2}=1$, $t_{n+3}=\infty$ et $\l_1,\dotsc,\l_n$ et de sch\'ema de Riemann
\begin{eqnarray*}
 &\begin{pmatrix}
             x=t_i & x=\infty      & x=\lambda_k\\
             0     & \a            & 0\\
             \t_i  & \a + \t_\infty & 2
  \end{pmatrix}\\
 &\begin{array}{ccc}
                    i=1,\dotsc,n+2, & \,  & k=1,\dotsc,n .
  \end{array}
\end{eqnarray*}
On suppose que les singularit\'es $x=\l_k$ sont apparentes (d\'efinition~\ref{def-sg-app}) et que les exposants v\'erifient $\t_i \notin \Z$, $i=1,\dotsc,n+3$ (on note parfois $\t_{n+3}$ pour $\t_\infty$). L'\'equation~\eqref{E-G} n'a donc aucune singularit\'e logarithmique. La relation de Fuchs~\eqref{rel-Fu} impose
\[
 \a = \frac12 \left(1- \sum_{i=1}^{n+3} \t_i \right) .
\]
Le th\'eor\`eme~\ref{thm-boli} nous assure que pour toute monodromie irr\'eductible, il existe une \'equation de ce type ayant cette monodromie. Le but de cette section est de d\'ecrire les transformations isomonodromiques de l'\'equation~\eqref{E-G}. On commence par pr\'eciser l'expression de ses c\oe ficients $p(x)$ et $q(x)$.

D'apr\`es la proposition~\ref{prop-cara-eq-fu}, les c\oe fficients $p(x)$ et $q(x)$ de l'\'equation~\eqref{E-G} s'\'ecrivent
\begin{align*}
        p(x) & =       \sum_{i=1}^{n+2} \dfrac{a_i}{x-t_i}     + \sum_{k=1}^n \dfrac{c_k}{x-\l_k},\\
        q(x) & =       \sum_{i=1}^{n+2} \dfrac{b_i}{(x-t_i)^2} + \sum_{k=1}^n \dfrac{d_k}{(x-\l_k)^2}
              - \sum_{i=1}^{n+2} \dfrac{K_i}{x-t_i}     + \sum_{k=1}^n \dfrac{\mu_k}{x-\l_k},
\end{align*}
avec
\begin{equation}
 -\sum_{i=1}^{n+2} K_i + \sum_{k=1}^n \mu_k = 0.
\label{(eq-G-infty)}
\end{equation}
Les exposants de l'\'equation nous permettent de calculer certaines des constantes intervenant dans l'expression de $p(x)$ et $q(x)$. L'\'equation caract\'eristique en $x=t_i$ est
\[
    s^2+(a_i-1)s+b_i=0,
\]
et ses racines sont $0$ and $\t_i$. On en d\'eduit que  $ a_i=1-\t_i$ et $ b_i=0$. De m\^eme, on obtient $ c_k=-1$ et $ d_k=0$. De l'\'equation caract\'eristique en l'infini, on d\'eduit
\[
 \a(\a+\t_\infty) = -\sum_{i=1}^{n+2} t_iK_i + \sum_{k=1}^n \l_k\mu_k.
\]
De cette relation et de~\eqref{(eq-G-infty)}, on d\'eduit $K_{n+1}$ et $K_{n+2}$ en fonction des autres constantes et on obtient l'expression suivante des c\oe fficients $p(x)$ et $q(x)$~:
\begin{equation}
        \begin{cases}
            p(x) =& \displaystyle\sum_{i=1}^{n+2} \dfrac{1-\t_i}{x-t_i} - \sum_{k=1}^n \dfrac{1}{x-\l_k}\\
            q(x) =& \dfrac{\a(\a+\t_\infty)}{x(x-1)} - \displaystyle\sum_{i=1}^n\dfrac{t_i(t_i-1)K_i}{x(x-1)(x-t_i)}\\
                  & + \displaystyle\sum_{k=1}^n\dfrac{\l_k(\l_k-1)\mu_k}{x(x-1)(x-\l_k)}
        \end{cases}\label{pq}
\end{equation}
o\`u les $K_i$, $\mu_k$ sont des constantes inconnues
\begin{align*}
 K_i   &=-\res \left(q(x),x=t_i  \right)\\
 \mu_k &= \res \left(q(x),x=\l_k \right).
\end{align*}
Pour chaque valeur fix\'ee de $\t=(\t_1,\ldots,\t_{n+3})$, l'\'equation~\eqref{E-G} de c\oe fficients~\eqref{pq} d\'epend donc d'au plus $4n$ param\`etres
\[
    t_1,\dotsc,t_n,\; \l_1,\dotsc,\l_n,\; \mu_1,\dotsc,\mu_n,\;
    K_1,\dotsc,K_n.
\]
Cependant, toutes les valeurs de ces param\`etres ne d\'efinissent pas
n\'ecessairement une \'equation ayant des singularit\'es apparentes en
les $\l_k$ (vu les exposants en ces singularit\'es, elles peuvent
\^etre logarithmiques). La proposition suivante, obtenue en
appliquant la m\'ethode de Fr\"obenius aux points $x=\l_k$, donne une
condition n\'ecessaire et suffisante pour que l'\'equation~\eqref{E-G}
n'ait aucune singularit\'e logarithmique. Sa d\'emonstration se trouve dans~\cite{IKSY}.

\begin{prop}
Les points $\l_1,\dotsc,\l_n$ sont des singularit\'es non logarithmiques de l'\'equation~\eqref{E-G} de c\oe fficients $p(x)$ et $q(x)$ d\'efinis par~\eqref{pq} si et seulement si les r\'esidus $K_i$ sont donn\'es par
    \[
        K_i = M_i \sum_{k=1}^n M^{k,i}\left(\mu_k^2 - \sum_{j=1}^{n+2} \frac{\t_j-\delta_{ij}}{\l_k-t_j}\mu_k +
        \frac{\a(\a+\t_\infty)}{\l_k(\l_k-1)}\right),
    \]
o\`u $M_i$ et $M^{k,i}$ sont d\'efinis par
    \[
        M_i     =- \frac{\Lambda(t_i)}{T'(t_i)} \quad \text{et} \quad
        M^{k,i} = \frac{T(\l_k)}{(\l_k-t_i)\Lambda'(\l_k)},
    \]
o\`u les polyn\^omes $\Lambda(x)$ et $T(x)$ sont donn\'es par~\eqref{def-X-T-L} .
\label{prop-def-Ki}
\end{prop}

Les r\'esidus $K_i$ sont donc des fractions rationnelles de $(\t,\l,\mu,t)$. Les \'equations~\eqref{E-G} v\'erifiant les hypoth\`eses souhait\'ees d\'ependent uniquement des param\`etres $(\t,\l,\mu,t)$, on les note donc $\eq$. On cherche \`a quelle condition des variations de ces param\`etres pr\'eservent la monodromie d'une telle \'equation. Les exposants $\t=(\t_1,\ldots,\t_{n+3})$ sont n\'ecessairement constants pendant une d\'eformation isomonodromique continue. On pose
\[
 \B^n= \left\{ (t_1,\ldots,t_n) \in (\C^* \ssm \{1\})^n \quad |\quad \forall i\neq j \quad t_i \neq t_j \right\},
\]
et on cherche \`a caract\'eriser les sous-vari\'et\'es $M$ de $\C^n \times
\C^n \times \B^n$ telles que la famille d'\'equations $\eq$, $\left(
(\l,\mu,t)\in M\right) $ soit isomonodromique.
Dans~\cite{Garnier12}, Garnier donne le syst\`eme d'\'equations aux
d\'eriv\'ees partielles qui d\'ecrit les d\'eformations isomonodromiques
des \'equations $\eq$. Le param\`etre de la d\'eformation est le
param\`etre $t$, et le syst\`eme d\'ecrit les variations des param\`etres
$\l_k(t)$ en fonction de $t$, tandis que les r\'esidus $\mu_k(t)$,
vus \'egalement comme des fonctions de $t$, s'expriment
rationnellement \`a partir des $\l_k(t)$ et de leurs d\'eriv\'ees
premi\`eres. Okamoto~\cite{Okamoto} a mis en \'evidence la structure
hamiltonienne de ce syst\`eme, et lui a donn\'e le nom de
\emph{syst\`eme de Garnier}. C'est sous cette forme qu'il est connu aujourd'hui.

\begin{defi}
Le \emph{syst\`eme de Garnier} $(\mathcal{G}_n)$ de dimension $n$ est le syst\`eme hamiltonien
    \begin{equation}
        \begin{cases}
            \dfrac{\partial \l_i}{\partial t_j}  =
            \dfrac{\partial K_j}{\partial \mu_i} \\
            \dfrac{\partial \mu_i}{\partial t_j}  = -
            \dfrac{\partial K_j}{\partial \l_i}   \\
        \end{cases}
    \label{Gn}
    \end{equation}
($i,j=1,\dotsc,n$), o\`u les Hamiltoniens $K_i=K_i(\t,\l,\mu,t) $ sont donn\'es \`a la proposition~\ref{prop-def-Ki}.
\end{defi}

On a alors

\begin{theo}
Soit $\t=(\t_1,\ldots,\t_{n+3})\in\left(\C\ssm\Z\right)^{n+3}$.

\begin{enumerate}
 \item Le syst\`eme~$(\mathcal{G}_n)$ est compl\`etement int\'egrable.
 \item Soit $M$ une sous-vari\'et\'e de $\C^n \times \C^n \times \B^n$. Alors la famille d'\'equations  $\eq$, $\left( (\l,\mu,t)\in M\right) $ est isomonodromique si et seulement si $M$ est une sous-vari\'et\'e d'une vari\'et\'e int\'egrale du syst\`eme de Garnier~$(\mathcal{G}_n)$.
\end{enumerate}
\label{thm-sys-G}
\end{theo}

Une solution $(\l(t),\mu(t))$ du syst\`eme de Garnier~$(\mathcal{G}_n)$ est d\'etermin\'ee par la donn\'ee d'une monodromie pour l'\'equation $\eq$.

\begin{rema}
Dans le cas o\`u $n=1$, en notant $\left(\l,\mu,t,K\right) $ les quantit\'es $\left(\l_1,\mu_1,t_1,K_1\right) $, on obtient que l'Hamiltonien $K(\l,\mu,t)$ est donn\'e par
\begin{align*}
 K(\l,\mu,t) =& \frac{1}{t(t-1)} \Big[ \l(\l-1)(\l-t)\mu^2 - \Big(\t_2(\l-1)(\l-t)\\
& + \t_3\l(\l-t) + (\t_1-1)\l(\l-1) \Big)\mu + \kappa (\l-t) \Big]
\end{align*}
o\`u
\[
 \kappa = \frac14\left(\left(\t_1+\t_2+\t_3-1\right)^2-\t_4^2\right).
\]
En \'eliminant la variable conjugu\'ee $\mu$, on trouve que le syst\`eme de Garnier $(\mathcal{G}_1)$ est \'equivalent \`a la sixi\`eme \'equation de Painlev\'e~\eqref{PVI}~:
\begin{equation}
\begin{split}
 \frac{d^2\l}{dt^2} = & \frac12 \left(\frac1\l +\frac{1}{\l-1} +\frac{1}{\l-t}\right) \left(\frac{d\l}{dt}\right)^2 - \left(\frac1t +\frac{1}{t-1} +\frac{1}{\l-t}\right)\frac{d\l}{dt}\\
& + \frac{\l(\l-1)(\l-t)}{t^2(t-1)^2} \left(\a + \b\frac{t}{\l^2} + \gamma\frac{t-1}{(\l-1)^2} + \d\frac{t(t-1)}{(\l-t)^2} \right)
\end{split}
 \tag{$\PVI$}\label{PVI}
\end{equation}
avec
\[
 \a=\frac12\t_4^2, \quad \b=-\frac12\t_2^2, \quad \d=\frac12\t_3^2, \quad \gamma=\frac12\left(1-\t_1^2\right).
\]
En ce sens, le syst\`eme~$(\mathcal{G}_n)$ constitue une g\'en\'eralisation de l'\'equation~\eqref{PVI} en un syst\`eme aux d\'eriv\'ees partielles compl\`etement int\'egrable.
\end{rema}

\chapter{D\'emonstrations de r\'esultats utilis\'es au chapitre~5}
\label{annexe}

On va donner les d\'emonstrations des r\'esultats dus \`a Sato, Miwa et
Jimbo~\cite{SMJ}, ainsi qu'\`a Jimbo~\cite{Jimbo} dont on a eu
besoin au chapitre~\ref{chapitre-longueur} pour \'etudier la
fonction <<~rapports des longueurs~>> $F_D(t)$. On ne d\'emontre que
le th\'eor\`eme~\ref{thm-SMJ} et la proposition~\ref{prop-sch-lim}. La
d\'emonstration de la proposition~\ref{prop-mono-a-la-lim} proc\`ede
des m\^emes m\'ethodes que celle du th\'eor\`eme~\ref{thm-SMJ}.

\subsection*{D\'emonstration du th\'eor\`eme~\ref{thm-SMJ}}

On rappelle l'\'enonc\'e du th\'eor\`eme~\ref{thm-SMJ}.

\begin{theo*}
 Soient $A_\a^0$ ($\a=1,\ldots,p-1,n+2$) et $A_\b^0$ ($\b=p,\ldots,n+1$) des matrices constantes dont les valeurs propres sont respectivement $\left(-\t_\a/2, \t_\a/2\right)$ et $\left(-\t_\b/2, \t_\b/2\right)$. On suppose de plus que
\[
 \sum_\a A_\a^0 + \sum_\b A_\b^0 = -A_\infty
\]
et que les valeurs propres $\mu$ et $-\mu$ de la matrice
\[
 \Lambda :=\sum_\b A_\b^0
\]
v\'erifient : $0<2\Re (\mu)<1$. On note $\s=2\Re (\mu)$. Soient $\s_1$ et $K$ deux constantes telles que
\[
 \s<\s_1<1 \qquad \text{ et } \qquad|A_\a^0|<K, \ |A_\b^0|<K.
\]
Alors il existe $\e>0$ tel que dans tout secteur $S_{\e,\varphi}=\left\{ \tau \in\C \ | \ 0<|\tau|<\e , \ |\arg\tau|<\varphi\right\}$, il existe une unique solution $A_\a(\tau)$ ($\a=1,\ldots,p-1,n+2$), $A_\b(\tau)$ ($\b=p,\ldots,n+1$) du syst\`eme~\eqref{sch-restreint} v\'erifiant~:
\begin{equation}
\begin{split}
 &|A_\a(\tau)-A_\a^0|  \leq K|\tau|^{1-\s_1},\\
 &|\tau^{-\Lambda} A_\b(\tau)\tau^\Lambda-A_\b^0|  \leq K|\tau|^{1-\s_1}.
\end{split}
\end{equation}
\end{theo*}

\begin{proof}
On pose, pour tout $\b=p,\ldots,n+1$,
\[
 \w A_\b(\tau) = \tau^{-\Lambda}A_\b(\tau)\tau^\Lambda .
\]
On r\'ecrit le syst\`eme de Schlesinger restreint~\eqref{sch-restreint} avec les matrices $A_\a(\tau)$ ($\a=1,\ldots,p-1,n+2$), $\w A_\b(\tau)$ ($\b=p,\ldots,n+1$) comme inconnues :
\begin{equation}
\begin{split}
  \frac{dA_\a}{d\tau} = \sum_\b \frac{\nu_\b}{\tau\nu_\b-t_\a} \left[\tau^\Lambda\w A_\b(\tau)\tau^{-\Lambda},A_\a(\tau)\right]\\
  \frac{d\w A_\b}{d\tau} = \sum_\a \frac{\nu_\b}{\tau\nu_\b-t_\a} \left[\tau^{-\Lambda}A_\a(\tau)\tau^\Lambda,\w A_\b(\tau)\right]\\
    + \frac1\tau \sum_{\b'}\left[\left(\w A_{\b'}(\tau)-A_{\b'}^0\right),\w A_\b(\tau)\right].
\end{split}
\label{nv-sch-restreint}
\end{equation}
On construit la solution recherch\'ee en proc\'edant par it\'eration. On pose
\[
 A_\a^{(0)}(\tau) = A_\a^0, \qquad  \w A_\b^{(0)}(\tau) = A_\b^0,
\]
et pour tout entier naturel $k$, on d\'efinit les matrices $A_\a^{(k)}(\tau)$ et $\w A_\b^{(k)}(\tau)$ \`a partir de $A_\a^{(k-1)}(\tau)$ et $\w A_\b^{(k-1)}(\tau)$ par :
\begin{align*}
  A_\a^{(k)}(\tau) = A_\a^0 + \sum_\b \int_0^\tau \frac{\nu_\b}{s\nu_\b-t_\a} \left[s^\Lambda\w A_\b^{(k-1)}(s)s^{-\Lambda},A_\a^{(k-1)}(s)\right]ds\\
  \w A_\b^{(k)}(\tau) = A_\b^0 + \sum_\a \int_0^\tau \frac{\nu_\b}{s\nu_\b-t_\a} \left[s^{-\Lambda}A_\a^{(k-1)}(s)s^\Lambda,\w A_\b^{(k-1)}(s)\right]ds\\
    + \sum_{\b'} \int_0^\tau \frac1s \left[\left(\w A_{\b'}^{(k-1)}(s)-A_{\b'}^0\right),\w A_\b^{(k-1)}(s)\right]ds.
\end{align*}
Les int\'egrales sont calcul\'ees le long du segment joignant $0$ et $\tau$~:
\[
 \{s=re^{i\psi}\ | \ 0<r<|\tau|,\ \psi=\arg \tau\}.
\]
Soit une constante $\d$ telle que $0<\d<1$. 
On va montrer par r\'ecurrence que les matrices les $A_\a^{(k)}(\tau)$ et $\w A_\b^{(k)}(\tau)$ sont bien d\'efinies et qu'elles convergent uniform\'ement dans tout voisinage compact de $\tau=0$. Pour cela, on va montrer qu'il existe un nombre $\e>0$ ne d\'ependant que de $\s$, $\s_1$, $\d$, $K$, $r$ et $R$ tel que les matrices $A_\a^{(k)}(\tau)$ et $\w A_\b^{(k)}(\tau)$ v\'erifient pour tout $\tau$ dans le secteur $S_{\e,\varphi}$ les conditions asymptotiques suivantes~:
\begin{align}
 \left| A_\a^{(k)}(\tau)-A_\a^0 \right| &\leq K|\tau|^{1-\s_1}
\label{I1}\\
 \left| \w A_\b^{(k)}(\tau) -A_\b^0 \right| &\leq K|\tau|^{1-\s_1}
\label{I3}
\end{align}
et
\begin{align}
 \left| A_\a^{(k)}(\tau)-A_\a^{(k-1)}(\tau) \right| &\leq K\d^{k-1}|\tau|^{1-\s_1}
\label{I4}\\
 \left| \w A_\b^{(k)}(\tau) -A_\b^{(k-1)}(\tau) \right| &\leq K\d^{k-1}|\tau|^{1-\s_1}
\label{I6}
\end{align}

L'initialisation est \'evidente. Supposons que les matrices $A_\a^{(k)}(\tau)$ et $\w A_\b^{(k)}(\tau)$ sont bien d\'efinies et qu'elles v\'erifient les majorations~\eqref{I1}$_k$,\ldots, \eqref{I6}$_k$. On doit avoir
\[
 \e<\frac rR
\]
o\`u les constantes $r$ et $R$ sont d\'efinies par~\eqref{def-r-R}. On choisit de plus
\[
 \e<1.
\]
On a alors par les majorations~\eqref{I1}$_k$ et~\eqref{I3}$_k$ et par d\'efinition de la constante $K$
\begin{equation}
 \left|A_\a^{(k)}(\tau) \right|< 2K, \qquad \left|\w A_\b^{(k)}(\tau) \right|< 2K.
\label{I7}
\end{equation}
Or pour toute matrice $C\in M(2,\C)$, les \'el\'ements des matrices $\tau^{-\Lambda}C\tau^\Lambda$ et $\tau^\Lambda C \tau^{-\Lambda}$ sont des polynômes du premier degr\'e en $\tau^\s$ et $\tau^{-\s}$, et donc
\[
 \left|\tau^{-\Lambda}C\tau^\Lambda\right|\leq |C| |\tau|^{-\s}, \qquad
 \left|\tau^\Lambda C \tau^{-\Lambda}\right|\leq |C| |\tau|^{-\s}.
\]
On peut donc d\'eduire de~\eqref{I7}$_k$
\[
 \left|\tau^{-\Lambda}A_\a^{(k)}(\tau)\tau^\Lambda\right|\leq 2K |\tau|^{-\s}, \qquad
 \left|\tau^\Lambda \w A_\b^{(k)}(\tau) \tau^{-\Lambda}\right|\leq 2K |\tau|^{-\s}.
\]
On peut d\'eduire des ces majorations et des majorations~\eqref{I1}$_k$ et~\eqref{I3}$_k$ que les matrices $A_\a^{(k+1)}(\tau)$ et $\w A_\b^{(k+1)}(\tau)$ sont bien d\'efinies.

\'Etablissons les majorations~\eqref{I1}$_{k+1}$ et~\eqref{I3}$_{k+1}$. On remarque tout d'abord que l'on a pour tout $\tau$ dans le secteur $S_{\e,\varphi}$
\[
 \left| \frac{\nu_\b}{\tau\nu_\b-t_\a}\right| \leq \left(\frac rR -\e\right)^{-1}<\left(\frac rR -1\right)^{-1}
\]
Alors
\begin{align*}
 \left| A_\a^{(k+1)}(\tau)-A_\a^0 \right| &\leq 2\sum_\b \int_0^{|\tau|} \left|\frac{\nu_\b}{s\nu_\b-t_\a}\right| \left|s^\Lambda\w A_\b^{(k)}(s)s^{-\Lambda}\right|\left|A_\a^{(k)}(s)\right|ds\\
&\leq 8(n-p+2)\left(\tfrac rR -1\right)^{-1}K^2 \int_0^{|\tau|} \frac{ds}{s^\s}\\
&\leq K|\tau|^{1-\s_1} \left[ \frac{8K(n-p+2)}{(1-\s)\left(\tfrac rR -1\right)}  \right]\e^{\s_1-\s}.
\end{align*}
De m\^eme
\begin{align*}
 \left|\w A_\b^{(k+1)}(\tau) - A_\b^0 \right| &\leq 2\sum_\a \int_0^{|\tau|} \left|\frac{\nu_\b}{s\nu_\b-t_\a} \right| \left| s^{-\Lambda}A_\a^{(k)}(s)s^\Lambda\right|\left|\w A_\b^{(k)}(s)\right|ds\\
    & \quad + 2\sum_{\b'} \int_0^{|\tau|} \frac1s \left|\w A_{\b'}^{(k)}(s)-A_{\b'}^0\right|\left|\w A_\b^{(k)}(s)\right|ds\\
&\leq K|\tau|^{1-\s_1} \left[  \frac{4K}{(1-\s)} \left( \frac{2p}{\left(\tfrac rR -1\right)} +(n-p+2) \right) \right]\e^{\s_1-\s}.
\end{align*}
Il suffit donc de choisir $\e$ tel que $\e^{\s_1-\s}$ soit inf\'erieur \`a la plus grande des deux quantit\'es suivantes
\[
 \frac{(1-\s)\left(\tfrac rR -1\right)}{8K(n-p+2)}, \qquad \qquad
 \frac{(1-\s)}{4K} \left[ \frac{2p}{\left(\tfrac rR -1\right)} +(n-p+2) \right]^{-1}.
\]

On obtient de m\^eme les majorations~\eqref{I4}$_{k+1}$ et~\eqref{I6}$_{k+1}$. On en d\'eduit donc que les suites $A_\a^{(k)}(\tau)$ et $\w A_\b^{(k)}(\tau)$ convergent uniform\'ement dans tout voisinage compact de $\tau=0$. On note
\begin{align*}
 A_\a(\tau):= \lim_{k\to+\infty} A_\a^{(k)}(\tau)\\
 \w A_\b(\tau):= \lim_{k\to+\infty} \w A_\b^{(k)}(\tau).
\end{align*}
Alors les matrices $A_\a(\tau)$ et $A_\b(\tau) = \tau^\Lambda \w A_\b(\tau)\tau^{-\Lambda}$ constitue une solution du syst\`eme de Schlesinger restreint~\eqref{sch-restreint}. Cette solution v\'erifie les conditions asymptotiques~\eqref{A_i-in\'egalit\'es-tau=0}. L'unicit\'e de cette solution se montrerait de m\^eme par r\'ecurrence.
\end{proof}

\subsection*{D\'emonstration de la proposition~\ref{prop-sch-lim}}

On rappelle l'\'enonc\'e de la proposition~\ref{prop-sch-lim}.

\begin{prop*}
Les matrices $A_\a^0(t',\nu)$ ($\a=1,\ldots,p-1,n+2$) et $\Lambda(t',\nu)$ sont solutions du syst\`eme de Schlesinger suivant
\begin{equation}
\begin{cases}
 d'A_\a' &= \displaystyle\sum_{\a'\neq\a} \left[A'_{\a'},A'_\a\right]d'\log(t_\a-t_{\a'})\\
 d_\nu A'_\a &= 0
\end{cases}
\label{A-sch-lim1}
\end{equation}
o\`u on a pos\'e $A_{n+1}^0(t',\nu) :=\Lambda(t',\nu)$, et o\`u $d'$ d\'esigne la diff\'erentiation par rapport \`a $t'=(t_1,\ldots,t_{p-1})$ et $d_\nu$ la diff\'erentiation par rapport \`a $\nu=(\nu_{p+1},\ldots,\nu_n)$.

Les matrices $A^0_\a(t',\nu)$ ($\a=1,\ldots,p-1,n+2$) et $A_\b^0(t',\nu)$ ($\b=p,\ldots,n+1$) sont solutions du syst\`eme
\begin{equation}
\begin{cases}
 d'A_\b^0 &= -\displaystyle\sum_{\a=1}^{p-1} \left[ A_\b^0,A_\a^0\right] d'\log(t_\a)\\
 d_\nu A_\b^0 &= \displaystyle\sum_{\b'\neq\b} \left[A_{\b'}^0,A_\b^0\right] d_\nu\log(\nu_\b-\nu_{\b'})
\end{cases}.
\label{A-sch-lim2}
\end{equation}
\end{prop*}

\begin{proof}
On \'etablit uniquement le syst\`eme~\eqref{A-sch-lim2}. Le syst\`eme~\eqref{A-sch-lim1} se montre de la m\^eme mani\`ere, et il est plus simple \`a \'etablir. D'apr\`es le syst\`eme de Schlesinger~\eqref{schlesinger}, on a
\[
 d'A_\b = \sum_\a [A_\a,A_\b] d'\log(t_\a-\tau\nu_\b),
\]
et
\[
 d_\nu A_\b = -\tau\sum_\a [A_\a,A_\b] d_\nu\log(t_\a-\tau\nu_\b) + \sum_{\b'\neq\b} [A_{\b'},A_\b] d_\nu\log(\nu_{\b'}-\nu_\b).
\]
On va en d\'eduire les \'equations v\'erifi\'ees par les matrices $\w A_\b(\tau) = \tau^{-\Lambda} A_\b(\tau) \tau^\Lambda$. Pour cela, il faut v\'erifier que pour tout $\a=1,\ldots,p-1$
\begin{equation}
 d'\Lambda = -\sum_\a\left[ \Lambda, A_\a^0\right]d'\log t_\a.
\label{eq-Lambda}
\end{equation}
Comme
\[
 A_\infty = -\sum_\a A_\a^0 -\sum_\b A_\b^0 = -\sum_\a A_\a(\tau)-\sum_\b A_\b(\tau),
\]
on a
\[
 \Lambda = \lim_{\tau\to0} \sum_\b A_\b(\tau).
\]
Or
\begin{align*}
 \left[ \sum_\b A_\b , A_\a\right]
&= \sum_\b (\tau\nu_\b-t_\a)\frac{\partial A_\b}{\partial t_\a}\\
&= \tau^{1-\sigma}\sum_\b \nu_\b\frac{\partial\left( \tau^\sigma A_\b\right) }{\partial t_\a} - t_\a\frac{\partial}{\partial t_\a}\sum_\b A_\b.
\end{align*}
Grâce \`a l'assertion~\eqref{Ab-holo} de la proposition~\ref{prop-A-holo-en-tau}, lorsque $\tau$ tend vers $0$, on obtient~\eqref{eq-Lambda}. On en d\'eduit
\begin{equation}
 d'\tau^\Lambda = -\sum_\a\left[ \tau^\Lambda, A_\a^0\right] d'\log t_\a , \qquad
 d'\tau^{-\Lambda} = -\sum_\a\left[ \tau^{-\Lambda}, A_\a^0\right] d'\log t_\a,
\label{eq-tau^Lambda}
\end{equation}
vu que
\begin{align*}
 \frac{\partial\tau^\Lambda}{\partial t_\a}
&= \log(\tau)\int_0^1\tau^{(1-u)\Lambda} \frac{\partial\Lambda}{\partial t_\a} \tau^{u\Lambda} du\\
&= -\frac{1}{t_\a}\int_0^1\tau^{(1-u)\Lambda} \left[ \log(\tau)\Lambda, A_\a^0\right] \tau^{u\Lambda} du\\
&= -\frac{1}{t_\a} \left[\tau^{(1-u)\Lambda} A_\a^0 \tau^{u\Lambda} \right]_{u=0}^{u=1}.
\end{align*}
On obtient donc, d'une part,
\[
 d_\nu \widetilde A_\b = -\tau^{1-\sigma}\sum_\a \left[\tau^\sigma\widetilde A_\a,\widetilde A_\b\right] d_\nu\log(t_\a-\tau\nu_\b) + \sum_{\b'\neq\b} \left[\widetilde A_{\b'},\widetilde A_\b\right] d_\nu\log(\nu_{\b'}-\nu_\b),
\]
et donc, vu l'assertion~\eqref{Aa-tilde-holo} de la proposition~\ref{prop-A-holo-en-tau}, quand $\tau$ tend vers $0$, on obtient
\[
 d_\nu A_\b^0 = \sum_{\b'\neq\b} \left[A_{\b'}^0,A_\b^0\right] d_\nu\log(\nu_{\b'}-\nu_\b).
\]
D'autre part,
\begin{align*}
 d'\widetilde A_\b
&= \sum_\a\left( \frac{\left[\widetilde A_\a,\widetilde A_\b\right]}{1-\tau\frac{\nu_\b}{t_\a}} - \left(\left[ \tau^{-\Lambda}, A_\a^0\right]A_\b\tau^\Lambda + \tau^{-\Lambda} A_\b\left[ \tau^\Lambda A_\a^0\right]\right) \right) d'\log t_\a ,\\
&= \sum_\a\left(\frac{\left[\widetilde A_\a,\widetilde A_\b\right]}{1-\tau\frac{\nu_\b}{t_\a}} - \left[ \widetilde A_\b,A_\a^0-\tau^{-\Lambda}A_\a^0\tau^\Lambda\right] \right) d'\log t_\a \\
&= -\sum_\a\left(\left[\widetilde A_\b,A_\a^0\right] + \left[\widetilde A_\b,\tau^{-\Lambda}\left( A_\a-A_\a^0\right) \tau^\Lambda\right] + \e(\tau) \left[\widetilde A_\b,\tau^\sigma\widetilde A_\a\right] \right) d'\log t_\a,
\end{align*}
o\`u $\e(\tau)$ est une fonction qui tend vers $0$ avec $\tau$. On a finalement \`a la limite, de nouveau par la proposition~\ref{prop-A-holo-en-tau},
\[
 d'A_\b^0 = -\displaystyle\sum_{\a=1}^{p-1} \left[ A_\b^0,A_\a^0\right] d'\log(t_\a).
\]
\end{proof}


\nocite{Schwarz} \nocite{Weierstrass2} 

\bibliographystyle{alpha-fr} 
\bibliography{biblio}

\end{document}